%% file: main.tex
\documentclass[10pt]{simplearticle}

\makeatletter

\input{preamble.tex}

\input{metainfo.tex}

\makeatother

\begin{document}

\maketitlepage

\tableofcontents{}
\input{sec0.tex}
\input{sec1.tex}
\input{sec2.tex}
\input{sec3.tex}
\input{sec4.tex}
\input{sec5.tex}
\input{sec6.tex}

\input{biblio.tex}
\end{document}

%% file: preamble.tex
\usepackage[inline]{enumitem}

\DeclareMathAlphabet{\mathbbold}{U}{bbold}{m}{n} % blackboard bold font:  \mathbb 1, \mathbb 2, etc
\usepackage{mathrsfs}			% calligraphic fonts

\usepackage{mathtools}			% amsmath extension
\usepackage{amsthm}				% theorems
\usepackage{amssymb}			% symbols

\usepackage{mquiver}
\usepackage{tikz-cd}			% commutative diagrams
\usetikzlibrary{decorations.pathreplacing,shapes.misc}
\usetikzlibrary{nfold}

\usepackage{xstring}%!DEL	

\usepackage{tikz}
\newcommand{\tobjbase}{%
  \mathord{%
    \tikz[baseline=(zero.base)]{
      \node[inner sep=0pt, outer sep=0pt] (zero) {$0$};
      \node[inner sep=0pt, outer sep=0pt] at (zero.center) {\tiny$\ast$};
    }%
  }%
}

\newcommand{\tcbase}[1]{%
  \mathord{%
    \tikz[baseline=(char.base)]{
      \node[draw,circle,inner sep=0pt,outer sep=0pt,minimum size=1em,font=\small,line width=0.2pt] (char) {#1};
    }%
  }%
}

\newcommand{\norm}[1]{\left\Vert #1\right\Vert }
\newcommand{\ev}{\operatorname{ev}}
\newcommand{\as}{\operatorname{as}}
\newcommand{\Const}{\operatorname{Const}}
\newcommand{\const}{\operatorname{const}}
\newcommand{\diag}{\operatorname*{diag}}
\newcommand{\pr}{\operatorname{pr}}
\newcommand{\colim}{\operatorname*{colim}}
\newcommand{\dist}{\operatorname{dist}}
\newcommand{\Ad}{\operatorname*{Ad}}
\newcommand{\id}{\operatorname{id}}
\newcommand{\Id}{\operatorname{Id}}
\newcommand{\GEFC}{\mathbb{\mathtt{GEFC}}}
\newcommand{\hGEFC}{\mathtt{hGEFC}}
\newcommand{\DEFC}{\mathbb{\mathtt{DEFC}}}
\newcommand{\LEFC}{\mathbb{\mathtt{LEFC}}}
\newcommand{\Cstar}{\mathtt{C^{*}}}
\newcommand{\hCstar}{\mathtt{hC^{*}}}
\newcommand{\hAsy}{\mathtt{hAsy}}
\newcommand{\EFC}{\Cstar^{\Cstar}}
\newcommand{\Ch}{\operatorname{Ch}}
\newcommand{\chex}{\Ch_{\operatorname{ex}}}
\newcommand{\mathsb}[1]{\boldsymbol{\mathsf{#1}}}
\newcommand{\C}{\mathfrak{C}^{\,}}

\newcommand{\asadj}{\dashv_{as}}
\newcommand{\mmph}{\text{{␣}}}
\newcommand{\tc}[1]{\tcbase{#1}}
\newcommand{\spec}{\operatorname{sp}}
\newcommand{\uu}{\dot{\oplus}}
\newcommand{\au}{\alpha^{\oplus}}
\newcommand{\lu}{\lambda^{\oplus}}
\newcommand{\ru}{\rho^{\oplus}}
\newcommand{\bzero}{\boldsymbol{0}}
\newcommand{\tensorunit}{\mathbbold1}
\newcommand{\bbzero}{\mathbbold0}
\newcommand{\bpp}[1]{(#1^{#1})_{\textup{cart}}}
\newcommand{\terminal}{\operatorname{!}}
\newcommand{\Terminal}{\boldsymbol{!}}
\newcommand{\tobj}{\tobjbase}
\newcommand{\rc}{\operatorname{rc}}
\newcommand{\lc}{\operatorname{lc}}
\newcommand{\inv}{\operatorname{inv}}
\newcommand{\fr}{\textrm{fr}}
\newcommand{\Br}{\operatorname{\mathsb{Br}}}
\renewcommand{\i}{\operatorname{\mathsb i}}
\renewcommand{\u}{\operatorname{\mathsb u}}
\newcommand{\EV}{\mathsb{EV}}
\newcommand{\oK}{\mathbb{K}}
\newcommand{\fK}{\mathbf{K}}
\newcommand{\fM}{\mathbf{M}}
\newcommand{\oM}{\mathbb{M}}
\newcommand{\ttt}{\operatorname{tt}}
\newcommand{\kapparev}{\overleftarrow{\kappa}}
\newcommand{\mmdd}{.}
\newcommand{\mmdc}{,}
\newcommand{\hspc}[1]{\hspace{#1ex}}
\newcommand{\myllbracket}{[[}
\newcommand{\myrrbracket}{]]}
\newcommand{\bbrackets}[1]{\myllbracket#1\myrrbracket}
\newcommand{\Vect}{\operatorname{Vect}}
\newcommand{\modrel}[1]{/\hspace{-0.6ex}#1\hspace{0.1ex}}
\newcommand{\noloc}{\nobreak
\mspace{4mu plus 1mu}
{:}
\nonscript\mkern-\thinmuskip
\mathpunct{}
\mspace{2mu}
}
\newcommand{\singledot}{\,\cdot\,}

%% file: metainfo.tex
\title{On generalized morphisms associated with \\
endofunctors of $C^*$-algebras
}

\authors{
		\ainfo{Georgii~S.~Makeev}
		{School of Mathematical Sciences, Key Laboratory of MEA (Ministry of
Education) and Shanghai Key Laboratory of PMMP, East China Normal
University, Shanghai 200241, People’s Republic of China.
}
		{makeev.gs@gmail.com}
		{}
}

\keywords{bimonoidal categories; endofunctors; homotopy; $E$-theory; $C^*$-algebras; asymptotic algebra}

\date{September 2025}

\begin{abstract}
We introduce a class of good endofunctors of $C^{*}$-algebras, endow
it with a structure of a bimonoidal category, and define homotopies
of natural transformations between such endofunctors. For every pair
of $C^{*}$-algebras and a good endofunctor, we construct a commutative
monoid of generalized morphisms, and endow these monoids with a bilinear
composition. This construction generalizes the homotopy category of
asymptotic homomorphisms used in the definition of the Connes-Higson
$E$-theory. We also introduce the notion of asymptotically adjoint
good endofunctors, which has interesting applications to $E$-theory
and $K$-homology.
\end{abstract}

%% file: sec0.tex
\section*{Introduction\addcontentsline{toc}{section}{Introduction}}

This paper establishes categorical foundations for the future series
of works on generalized homotopies of $*$-homomorphisms between $C^{*}$-algebras,
which
we first introduced in~\cite{makeev_unsuspended}.
To clarify our motivation, denote by $F$ a map (not necessarily functorial)
that sends a $C^{*}$-algebra $B$ to the $C^{*}$-algebra $F(B)$,
and
a surjective $*$-homomorphism $A\xrightarrow{\psi}B$
to a $*$-homomorphism $F(A)\xrightarrow{F(\psi)}F(B)$. We say that
two $*$-homomorphisms $\varphi_{0},\varphi_{1}\colon A\to F(B)$
are \emph{$F$-homotopic} (written $\varphi_{0}\simeq_{F}\varphi_{1}$)
if there is a $*$-homomorphism $\Phi\colon A\to F(IB)$ making the
following diagram commutative for~$j=0,1$:
\[
\input{diagrams/d346.tikz}
\]
where $IB\coloneqq C([0,1],B)$ denotes the $C^{*}$-algebra of continuous
functions on $[0,1]$ taking values in $B$, and $e_{t}\colon IB\to B$
stands for the evaluation at $t\in[0,1]$. Generally speaking, $F$-homotopy
need not be an equivalence relation, but in some nice cases considered
below, it is. For such $F$, we introduce the set of \emph{generalized
morphisms} $\left[A,F,B\right]\coloneqq\hom(A,F(B))\modrel{\simeq_{F}}$.
It turns out that many Kasparov-type theories can be described within
this framework:
\begin{itemize}
\item \textbf{$K_{1}$-theory} (see, e.g.,~\cite{GHT}): $K_{1}(B)=\left[C_{0}(\mathbb{R}),\Id,B\right]$,
where $\Id$ denotes the identity functor.
\item \textbf{The Kasparov $KK_{0}$-theory} (see~\cite{blackadar1998}):
$KK_{0}(A,B)=[A,\mathcal{D},B]$, where $\mathcal{D}(B)$ is the \emph{double},
defined as the pullback
\[
\input{diagrams/d347.tikz}
\]
in which $\mathcal{M}(B)$ denotes the multipliers of $B$, and $q_{B}\colon\mathcal{M}(B)\to\mathcal{Q}(B)$
is the quotient $*$-homomorphism onto the \emph{corona} $\mathcal{Q}(B)\coloneqq\mathcal{M}(B)/B$.
Since
$\mathcal{M}$ is an endofunctor in the category
of $C^{*}$-algebras and surjective $*$-homomorphisms~\cite{lance1995},
there is a well-defined $*$-homomorphism $\mathcal{M}(e_{j})$, which
allows us to define $\mathcal{Q}(e_{j})$ and $\mathcal{D}(e_{j})$
by a diagram chase.
\item \textbf{Weak homotopy classes of extensions} (see~\cite{houghton1999universal}):
$\textrm{ext}(A,B)=\left[A,\mathcal{Q},B\right]$, where $\mathcal{Q}$
is the corona defined above.
\item \textbf{The Connes--Higson $E$-theory} (see~\cite{GHT}): $E_{0}(A,B)=\left[SA,\mathfrak{A},SB\right]$,
where $S\coloneqq C_{0}(\mathbb{R},\mmph)$ is the suspension functor,
and $\mathfrak{A}$ is the asymptotic algebra functor.
\end{itemize}
For simplicity, the $C^{*}$-algebras $A$ and $B$ in the examples
above are assumed to
be separable and stable.

The last example, which in fact inspired the idea of $F$-homotopy
(cf.~\cite[Definition~2.2]{GHT}), serves as the starting point of
this research and is worth a more detailed discussion.

For a $C^{*}$-algebra $B$, denote by $\mathfrak{A}B$ the \emph{asymptotic
algebra} of $B$ defined as $\mathfrak{A}B\coloneqq C_{b}([0,\infty),B)/C_{0}([0,\infty),B)$,
where $C_{b}([0,\infty),B)$ stands for the $C^{*}$-algebra of continuous
$B$-valued functions on $[0,\infty)$, and $C_{0}([0,\infty),B)$
denotes the ideal in $C_{b}([0,\infty),B)$ consisting of the functions
vanishing at infinity. One can show that $\mathfrak{A}$ is an endofunctor
of the category of $C^{*}$-algebras and $*$-homomorphisms.

For every $n\geq0$, let $\mathfrak{A}^{n}$ denote the $n$-fold
asymptotic algebra functor $\mathfrak{A}^{n}\coloneqq\mathfrak{A}\dots\mathfrak{A}$
($n$~times). An \emph
{asymptotic homomorphism} from $A$ to $B$
is a $*$-homomorphism $A\to\mathfrak{A}^{n}B$ for some $n\geq0$.
The $\mathfrak{A}^{n}$-homotopy classes of such $*$-homomorphisms
are denoted by $\bbrackets{A,B}_{n}\coloneqq\hom(A,\mathfrak{A}^{n}B)\modrel{\simeq_{\mathfrak{A}^{n}}}$;
and the class of $\varphi\colon A\to\mathfrak{A}^{n}B$ in $\bbrackets{A,B}_{n}$
is written as $\bbrackets{\varphi}_{n}$.

Let $\alpha\colon\Id\Rightarrow\mathfrak{A}$ be the natural transformation
whose component at the object $B$ sends $b$ to the class of the
constant function $t\mapsto b$ in $\mathfrak{A}B$. For every $n\geq0$,
there is a map of sets
\[
\bbrackets{A,B}_{n}\to\bbrackets{A,B}_{n+1}\colon\bbrackets{\varphi}_{n}\mapsto\bbrackets{\alpha\mathfrak{A}^{n}B\circ\varphi}_{n+1},
\]
where $\alpha\mathfrak{A}^{n}B$ stands for the component of $\alpha$
at $\mathfrak{A}^{n}B$. Introduce the set
$\bbrackets{A,B}\coloneqq\colim_{n}\bbrackets{A,B}_{n}$,
and denote by $\bbrackets{\varphi}$ the class of $\varphi\colon A\to\mathfrak{A}^{n}B$
in $\bbrackets{A,B}$. The \emph
{homotopy category of asymptotic
homomorphisms} (written $\hAsy$) is the category with objects $C^{*}$-algebras,
arrows from $A$ to $B$ elements of the set $\bbrackets{A,B}$, and
composition given by
\begin{equation}
\bbrackets{B\xrightarrow{\psi}\mathfrak{A}^{m}C}\circ\bbrackets{A\xrightarrow{\varphi}\mathfrak{A}^{n}B}\coloneqq\bbrackets{A\xrightarrow{\mathfrak{A}^{n}\psi\circ\varphi}\mathfrak{A}^{n+m}C}.\label{eq:ashom-cps}
\end{equation}
If we restrict to the category of separable $C^{*}$-algebras, the
colimit construction in the definition of $\bbrackets{A,B}$ can be
discarded (see~\cite[Theorem~2.16]{GHT}), i.e., if $A$ is separable,
then
\[
\bbrackets{A,B}\cong\bbrackets{A,B}_{1}=[A,\mathfrak{A},B].
\]

In case of stable\footnote{A $C^{*}$-algebra $B$ is called \emph{stable} if $B$ is isomorphic
to $B\otimes\oK$, where $\oK$ stands for the $C^{*}$-algebra of
compact operators in a separable infinite dimensional Hilbert space.} $B$, the set $\bbrackets{A,B}$ can be endowed with the structure
of a commutative monoid with the sum $\bbrackets{\varphi\colon A\to\mathfrak{A}^{n}B}+\bbrackets{\psi\colon A\to\mathfrak{A}^{n}B}$
defined as the class of the asymptotic homomorphism
\begin{equation}
A\to\mathfrak{A}^{n}B\colon a\mapsto\mathfrak{A}^{n}\theta_{2}\left(\begin{array}{cc}
\varphi(a) & 0\\
0 & \psi(a)
\end{array}\right),\label{eq:ashom-compos}
\end{equation}
where $\theta_{2}\colon M_{2}B\to B$ is a stability isomorphism,
with $M_{2}B$ standing for the $C^{*}$-algebra of $2\textrm{-by-}2$
matrices with entries in~$B$.

The \emph
{Connes--Higson $E$-theory category} is a category with
objects separable $C^{*}$-algebras, and arrows from $A$ to $B$
elements of the set $E_{0}(A,B)\coloneqq\bbrackets{S\oK A,S\oK B}$,
where $S\coloneqq C_{0}(\mathbb{R},\mmph)$ stands for the \emph
{suspension}
functor, which determines an abelian group structure on $E_{0}(A,B)$:
the inverse element is induced by the pre-composition (or post-composition)
with the homomorphism $f\mapsto f(-\singledot)$. Similarly to operator
$K$-theory, one can inductively define the sequence of bifunctors
$E_{n}(A,B)\coloneqq E_{n-1}(SA,B)$ from $C^{*}$-algebras to abelian
groups, which forms a generalized (co)homology theory satisfying the
Bott periodicity $E_{n+2}(A,B)\cong E_{n}(A,B)$. In fact, for separable
nuclear $C^{*}$-algebras, $E$-theory is isomorphic to the Kasparov
$KK$-theory~\cite[Theorem~25.6.3]{blackadar1998}.

The idea of replacing asymptotic homomorphisms with even more general
morphisms was introduced
in our first work~\cite{makeev_2019}
on this subject, where we proved that for separable $C^{*}$-algebras
$A$ and $B$, there is an isomorphism
\[
\left[SA,\mathfrak{A},B\right]\cong\left[A,\mathfrak{M}^{u}_{\mathbb{Z}}\mathfrak{A},B\right],
\]
where $\mathfrak{M}^{u}_{\mathbb{Z}}$ is the endofunctor associated
with the Roe algebra of integers equipped with its natural metric.
This result suggested that there should be a class of sufficiently
nice endofunctors which, similarly to $\mathfrak{A}$, give rise to
bivariant theories for $C^{*}$-algebras, and that the Roe functor
$\mathfrak{M}^{u}_{\mathbb{Z}}$ behaves, roughly speaking, as a right
adjoint to the suspension $S$. These ideas were formalized in~\cite{makeev_unsuspended},
where we introduced the notion of \emph{good endofunctor} --- an
endofunctor of $C^{*}$-algebras satisfying the two properties, both
inspired by the asymptotic algebra functor~$\mathfrak{A}$. Namely,
we called $F$ a \emph{good endofunctor} if:%--!--
\begin{enumerate}[label=\textup{\textbf{(G\arabic*)}}, ref=\textbf{(G\arabic*)}]

\item\label{enu:o1}$F$ preserves pullbacks along a pair of $*$-homomorphisms,
one of which is surjective (cf.~\cite[Lemma~2.5]{GHT});
\item\label{enu:o2}there exists a natural transformation $IF\Rightarrow FI$
making the following diagram commutative for~$j=0,1$:
\begin{equation}
\input{diagrams/d579.tikz}\label{eq:proto-kappa}
\end{equation}
(cf.~\cite[proof of Lemma~2.9]{GHT}).
\end{enumerate}
The first property ensured that $F$-homotopy is an equivalence relation,
and the second one yielded a well-defined composition of generalized
morphisms
\begin{equation}
[B\xrightarrow{\psi}GC]\circ[A\xrightarrow{\varphi}FB]\coloneqq[A\xrightarrow{F\psi\circ\varphi}FGC],\label{eq:composition}
\end{equation}
which is analogous to the composition of asymptotic homomorphisms,
cf.~(\ref{eq:ashom-cps}). In that article we also introduced \emph{homotopies
of natural transformations} and\ 
\emph{asymptotic adjunction},
and showed that for a pair of asymptotically adjoint good endofunctors
$L\asadj R$, there is an isomorphism
\begin{equation}
\colim_{n}[LA,\mathfrak{A}^{n},B]\cong\colim_{n}[A,R\mathfrak{A}^{n},B]\label{eq:isomorphism}
\end{equation}
natural in $A$ and $B$.

Unfortunately, axioms~\ref{enu:o1} and~\ref{enu:o2} did not provide
any tools for working with the additive structure on $[A,F,\oK\otimes B]$
defined as in~(\ref{eq:ashom-compos}). For example, we could not
show within that framework that for stable $C^{*}$-algebras $A$,
$B$, and $C$, the composition~(\ref{eq:composition}) is bilinear,
which is a highly desirable property. In this paper, we address this
problem by proposing an upgraded version of these axioms and redefining
the notion of ``good endofunctors'' from~\cite{makeev_unsuspended}.
This generalization, which requires more sophisticated categorical
techniques, lays the foundations for future publications, where we
construct a right asymptotic adjoint for the functor of tensoring
with continuous functions on a locally compact scalable metric space,
and propose a new model for the relative $K$-homology of compact
metric spaces.

This article is organized as follows.

Section~1 covers the categorical machinery we rely upon, such as
(bi)monoidal categories, monoids and rigs, quotient categories, and
some strictification results for unitary bimonoidal functors.

Section~2 discusses certain functor categories endowed with natural
bimonoidal structure, and defines the left tensoring functor which
plays a crucial role in the definition of good endofunctors.

Section~3 deals with the category of $C^{*}$-algebras and $*$-homomorphisms
and yields the bimonoidal category of \emph{decent endofunctors},
which generalizes the class of endofunctors satisfying Axiom~\ref{enu:o1}.
This section also discusses \emph{tensor-type endofunctors}, provides
a strictification theorem for the category of $C^{*}$-algebras, and
introduces the rig of compact operators which is essential for proving
bilinearity of composition.

In Section~4, we introduce \emph{labeled endofunctors} of $C^{*}$-algebras
whose axioms are designed to replace Axiom~\ref{enu:o2}. Labeled
endofunctors which are also decent yield \emph{good endofunctors }allowing
us to define homotopies of natural transformations and componentwise
quotients. We also show that good endofunctors and their homotopy
version naturally carry bimonoidal structures. Using the strictification
result from Section~1, we promote the rig structure of compact operators
from the previous section to the level of endofunctors.

In Section~5, for every good endofunctor $F$ and every pair of $C^{*}$-algebras
$A$ and $B$, we obtain the commutative monoid $[A,F\fK,B]$ of \emph{generalized
morphisms} from $A$ to $B$, with $\fK$ standing for the stabilization
functor. We also prove that the composition~(\ref{eq:composition})
is bilinear and that~(\ref{eq:isomorphism}) is an isomorphism of
monoids. We finish this section with a discussion of \emph{homotopy
symmetric good endofunctors} which yield invertible generalized morphisms.

\begin{acknowledgement}
The author would like to thank the anonymous referee whose detailed
comments greatly helped improve the manuscript.
\end{acknowledgement}

%% file: diagrams/d346.tikz
\begin{tikzcd}[ampersand replacement=\&]
{A} \& {F(IB)} \\
 \& {F(B),}
\arrow["\Phi", from=1-1, to=1-2]
\arrow["\varphi_{j}"', from=1-1, to=2-2]
\arrow["F(e_{j})", from=1-2, to=2-2]
\end{tikzcd}

%% file: diagrams/d347.tikz
\begin{tikzcd}[ampersand replacement=\&]
{\mathcal{D}(B)} \& {\mathcal{M}(B)} \\
{\mathcal{M}(B)} \& {\mathcal{Q}(B),}
\arrow[from=1-1, to=1-2]
\arrow[from=1-1, to=2-1]
\arrow["q_{B}", from=1-2, to=2-2]
\arrow["q_{B}", from=2-1, to=2-2]
\end{tikzcd}

%% file: diagrams/d579.tikz
\begin{tikzcd}[ampersand replacement=\&]
{IF} \&  \& {FI} \\
 \& {F,}
\arrow[Rightarrow, from=1-1, to=1-3]
\arrow["\ev_{j}F"', Rightarrow, from=1-1, to=2-2]
\arrow["F\ev_{j}", Rightarrow, from=1-3, to=2-2]
\end{tikzcd}

%% file: sec1.tex
\section{\label{sec:moncats}Monoidal and bimonoidal categories}

In this section, we provide the necessary definitions and facts concerning
monoidal and bimonoidal categories. The main references for this section
are~\cite{JY_bimon24,JY2021,mclane-categories}.

Let us first fix some basic notation. Consider the categories, functors,
and natural transformations as in the diagram below:
\[
\input{diagrams/d582.tikz}
\]
We denote by
\begin{itemize}
\item $FB$ and $F\varphi$ the values of the functor $F$ respectively
at the object $B$ and the morphism $\varphi$ in the category~$\mathcal{C}_{1}$;
\item $GF$ the composition of the functors $F$ and $G$;
\item $\alpha B\colon FB\to F'B$ the component of the natural transformation
$\alpha$ at the object~$B$;
\item $\beta F$, $G\alpha$ the whiskering of $\beta$ with $F$, and of
$G$ with $\alpha$, respectively;
\item $\beta\alpha\colon GF\Rightarrow G'F'$ the horizontal composition
of $\beta$ and~$\alpha$ (see~\cite[II.5 and~XII.3]{mclane-categories});
\item $\alpha'\circ\alpha\colon F\Rightarrow F''$ the vertical composition
of $\alpha$ and~$\alpha'$ (see~\cite[II.5 and~XII.3]{mclane-categories}).
\end{itemize}
We write $B\in\in\mathcal{C}$ and $f\in\mathcal{C}$ to express the
fact that $B$ is an object and $f$ is an arrow in the category~$\mathcal{C}$.
\begin{rem}
The binary operation of composing two functors will be referred to
as $\singledot$. We also use $GF$ and $G\cdot F$ interchangeably.
\end{rem}

\subsection{Monoidal categories}
\begin{defn}
A \emph{monoidal category} is a tuple $(\mathcal{C},\otimes,\tensorunit,\alpha,\lambda,\rho)$
consisting of
\begin{itemize}
\item a category $\mathcal{C}$;
\item a functor $\otimes\colon\mathcal{C}\times\mathcal{C}\to\mathcal{C}$
called the \emph{monoidal product};
\item an object $\tensorunit\in\in\mathcal{C}$ called the \emph{monoidal
unit};
\item a natural isomorphism $(X\otimes Y)\otimes Z\xrightarrow{\alpha_{X,Y,Z}}X\otimes(Y\otimes Z)$
for all $X,Y,Z\in\in\mathcal{C}$ called the \emph{associativity isomorphism};
\item natural isomorphisms $\tensorunit\otimes X\xrightarrow{\lambda_{X}}X\xleftarrow{\rho_{X}}X\otimes\tensorunit$
for all $X\in\in\mathcal{C}$ called respectively the \emph{left unit
isomorphism} and the \emph{right unit isomorphism}.
\end{itemize}
These data are required to make the following diagrams commutative
for all $W,X,Y,Z\in\in\mathcal{C}$:
\begin{alignat*}{3}
 & \textrm{(Pentagon axiom)} & \qquad &\ 
& \input{diagrams/d325.tikz}\\
 & \textrm{(Middle unity axiom)} & \qquad &\ 
& \input{diagrams/d324.tikz}
\end{alignat*}
A monoidal category is called \emph{strict} if the associativity and
unit isomorphisms are identities. In this case, $(\mathcal{C},\otimes,\tensorunit,\alpha,\lambda,\rho)$
will be abbreviated to $(\mathcal{C},\otimes,\tensorunit)$.
\end{defn}

\begin{example}\label{exa:str-moncat}
Given
a category $\mathcal{C}$,
denote by $\mathcal{C}^{\mathcal{C}}$ the category of endofunctors
of $\mathcal{C}$ and natural transformations between them. Composition
in this category is vertical composition of natural transformations.
If we denote by $\Id$ the identity endofunctor of $\mathcal{C}$,
then $(\mathcal{C}^{\mathcal{C}},\singledot,\Id)$ is a strict monoidal
category where the monoidal product $\singledot$ is given on objects
by composition of endofunctors, and on morphisms by horizontal composition
of natural transformations.
\end{example}

\begin{defn}\label{def:symmoncat}
A \emph{symmetric monoidal category} is a pair
$\left(\mathcal{C},\xi\right)$ in which:
\begin{itemize}
\item $\mathcal{C}=(\mathcal{C},\otimes,\tensorunit,\alpha,\lambda,\rho)$
is a monoidal category;
\item $\xi$ is a natural isomorphism $X\otimes Y\xrightarrow{\xi_{X,Y}}Y\otimes X$
for $X,Y\in\in\mathcal{C}$, called the \emph{symmetry isomorphism},
\end{itemize}
such that the following diagrams commute for all $X,Y,Z\in\in\mathcal{C}$:
\[
\input{diagrams/d431.tikz}
\]
\begin{alignat*}{2}
\input{diagrams/d432.tikz} & \qquad\qquad & \input{diagrams/d433.tikz}
\end{alignat*}

\end{defn}

\begin{example}\label{exa:Vect-moncat}
The category of real vector spaces $\Vect_{\mathbb{R}}$
can be endowed with two different symmetric monoidal structures: $\Vect^{\otimes}_{\mathbb{R}}\coloneqq(\Vect_{\mathbb{R}},\otimes,\mathbb{R},\alpha^{\otimes},\lambda^{\otimes},\rho^{\otimes},\xi^{\otimes})$
and $\Vect^{\oplus}_{\mathbb{R}}\coloneqq(\Vect_{\mathbb{R}},\oplus,\{0\},\alpha^{\oplus},\lambda^{\oplus},\rho^{\oplus},\xi^{\oplus})$,
where $\otimes$ and $\oplus$ denote respectively the tensor product
and direct sum, and $\alpha^{\otimes}$, $\dots$, $\xi^{\oplus}$
are defined in the obvious way. For example, $\xi^{\oplus}\colon U\oplus V\to V\oplus U\colon(u,v)\mapsto(v,u)$.
\end{example}

\begin{thm}
[Coherence Theorem,~{\cite[VII.2]{mclane-categories}}]\label{thm:coherence}
Let
$\mathcal{D}$ be a diagram in a monoidal category,
involving only the associativity isomorphisms, the unit isomorphisms,
their inverses, identity morphisms, the monoidal products, and composites.
Then $\mathcal{D}$ commutes.
\end{thm}

\begin{defn}\label{def:monoidal-functor}
 For monoidal categories $\mathcal{C}$
and $\mathcal{C}'$, a \emph{monoidal functor}
\[
(F,F^{2},F^{0})\colon\mathcal{C}\to\mathcal{C}'
\]
consists of:
\begin{itemize}
\item a functor $F\colon\mathcal{C}\to\mathcal{C}'$;
\item a natural transformation $F^{2}_{X,Y}\colon FX\otimes FY\to F(X\otimes Y)$
where $X,Y\in\in\mathcal{C}$;
\item a morphism $F^{0}\colon\tensorunit_{\mathcal{C}'}\to F\tensorunit_{\mathcal{C}}$.
\end{itemize}
These data are required to make the following diagrams commutative
for all $X,Y,Z\in\in\mathcal{C}$:
\[
\input{diagrams/d327.tikz}
\]
\begin{alignat*}{2}
\input{diagrams/d328.tikz} & \qquad\qquad & \input{diagrams/d328.tikz}
\end{alignat*}
A monoidal functor $(F,F^{2},F^{0})$ is called \emph{strong} (resp.
\emph{strict}) if both $F^{0}$ and $F^{2}$ are isomorphisms (resp.
identities). A monoidal functor between symmetric monoidal categories
is called \emph{symmetric} if the following diagram commutes:
\[
\input{diagrams/d345.tikz}
\]
For brevity, $(F,F^{2},F^{0})$ will sometimes be abbreviated to $F$.
\end{defn}

\begin{rem}
We shall often abuse notation and omit subscripts for $F^{2}$, $F^{0}$,
as well as for $\alpha$, $\lambda$, $\rho$, $\xi$, etc.
\end{rem}

\begin{defn}\label{def:monfun-compos}
Let $\mathcal{C}$, $\mathcal{C}'$ and
$\mathcal{C}''$ be monoidal categories, and let $(F,F^{2},F^{0})\colon\mathcal{C}\to\mathcal{C}'$
and $(G,G^{2},G^{0})\colon\mathcal{C}'\to\mathcal{C}''$ be two monoidal
functors. Their composition is defined as $(GF,(GF)^{2},(GF)^{0})$
where
\begin{alignat*}{1}
(GF)^{2}\colon & GFX\otimes GFY\xrightarrow{G^{2}}G(FX\otimes FY)\xrightarrow{GF^{2}}GF(X\otimes Y);\\
(GF)^{0}\colon & 1_{\mathcal{C}''}\xrightarrow{G^{0}}G1_{\mathcal{C}'}\xrightarrow{GF^{0}}GF1_{\mathcal{C}}.
\end{alignat*}
One can check that the composition of monoidal functors is again a
monoidal functor.
\end{defn}

\begin{defn}\label{def:mon-nt}
 For monoidal functors $F,G\colon\mathcal{C}\to\mathcal{C}'$,
a \emph{monoidal natural transformation} $\theta\colon F\to G$ is
a natural transformation between the underlying functors such that
the diagrams
\begin{alignat}{1}
\input{diagrams/d330.tikz}\qquad\qquad & \input{diagrams/d331.tikz}\label{eq:mon-nt}
\end{alignat}
commute for all $X,Y\in\in\mathcal{C}$.
\end{defn}

\begin{defn}
A \emph{monoid} in a monoidal category $\mathcal{C}$ is a triple
$(X,\mu,1)$ where
\begin{itemize}
\item $X$ an object in $\mathcal{C}$;
\item $\mu\colon X\otimes X\to X$ a morphism, called the \emph{multiplication};
\item $1\colon\tensorunit\to X$ a morphism, called the \emph{unit}.
\end{itemize}
These data are required to make the following diagrams commutative:
\begin{alignat*}{2}
\input{diagrams/d332.tikz} & \qquad\qquad & \input{diagrams/d333.tikz}
\end{alignat*}

\end{defn}

\begin{defn}
A \emph{commutative monoid} in a symmetric monoidal category $(\mathcal{C},\xi)$
is a monoid $(X,\mu,1)$ in $\mathcal{C}$ such that the following
diagram commutes:
\[
\input{diagrams/d334.tikz}
\]
\end{defn}

\begin{example}
Any real vector space $V$ can be regarded as a commutative monoid
$(V,+,0)$ in $\Vect^{\oplus}_{\mathbb{R}}$, where $+\colon V\oplus V\to V$
denotes the addition, and $0\colon\{0\}\to V$ stands for the unique
linear map from the zero space.
\end{example}

\subsection{Bimonoidal categories}
\begin{defn}
A \emph{bimonoidal category} is a tuple
\[
(\mathcal{C},(\oplus,\text{\ensuremath{\bbzero}},\alpha^{\oplus},\lambda^{\oplus},\rho^{\oplus},\xi^{\oplus}),(\otimes,\tensorunit,\alpha^{\otimes},\lambda^{\otimes},\rho^{\otimes}),(\lambda^{\bullet},\rho^{\bullet}),(\delta^{l},\delta^{r}))
\]
consisting of the following data:
\begin{itemize}
\item $(\mathcal{C},\oplus,\text{\ensuremath{\bbzero}},\alpha^{\oplus},\lambda^{\oplus},\rho^{\oplus},\xi^{\oplus})$
is a symmetric monoidal category (called the \emph{additive structure});
\item $(\mathcal{C},\otimes,\tensorunit,\alpha^{\otimes},\lambda^{\otimes},\rho^{\otimes})$
is a monoidal category (called the \emph{multiplicative structure});
\item $\lambda^{\bullet}$ and $\rho^{\bullet}$ are natural isomorphisms
$\bbzero\otimes A\xrightarrow{\lambda^{\bullet}_{A}}\bbzero\xleftarrow{\rho^{\bullet}_{A}}A\otimes\bbzero$
(called, respectively, the \emph{left} and \emph{right multiplicative
zero});
\item $\delta^{l}$ and $\delta^{r}$ are natural monomorphisms
\begin{alignat*}{1}
A\otimes(B\oplus B')\xrightarrow{\delta^{l}_{A,B,B'}}(A\otimes B)\oplus(A\otimes B') & ,\\
(B\oplus B')\otimes A\xrightarrow{\delta^{r}_{B,B',A}}(B\otimes A)\oplus(B'\otimes A)
\end{alignat*}
(called, respectively, the \emph{left} and \emph{right distributivity
morphism}).
\end{itemize}
The above data are required to satisfy the 22 axioms (known as \emph{Laplaza
axioms}) which the reader can find in
the Appendix.
A bimonoidal category is called \emph{tight} if $\delta^{l}$ and
$\delta^{r}$ are isomorphisms.

A \emph{symmetric bimonoidal category} is defined similarly, but with
an extra isomorphism $\xi^{\otimes}$ which makes the multiplicative
structure a symmetric monoidal category $(\mathcal{C},\otimes,\tensorunit,\alpha^{\otimes},\lambda^{\otimes},\rho^{\otimes},\xi^{\otimes})$;
furthermore, these data are required to satisfy two extra axioms (24
axioms in total; see the Appendix).
\end{defn}

\begin{example}\label{exa:Vect-bimoncat}
The category of real vector spaces $\Vect_{\mathbb{R}}$
carries the structure of symmetric bimonoidal category:
\[
(\Vect_{\mathbb{R}},(\oplus,\{0\},\alpha^{\oplus},\lambda^{\oplus},\rho^{\oplus},\xi^{\oplus}),(\otimes,\mathbb{R},\alpha^{\otimes},\lambda^{\otimes},\rho^{\otimes}),(\lambda^{\bullet},\rho^{\bullet}),(\delta^{l},\delta^{r})).
\]
with the multiplicative and additive structure as in the Example~\ref{exa:Vect-moncat},
and $\lambda^{\bullet}$, $\rho^{\bullet}$, $\delta^{l}$, and $\delta^{r}$
defined by the obvious formulas.
\end{example}

\begin{defn}\label{def:bimon-fun}
 Let $\mathcal{C}$ and $\mathcal{C}'$ be
bimonoidal categories. A \emph{bimonoidal functor} from $\mathcal{C}$
to $\mathcal{C}'$ is a tuple
\[
(F,F^{2}_{\oplus},F^{0}_{\oplus},F^{2}_{\otimes},F^{0}_{\otimes})\colon\mathcal{C}\to\mathcal{C}'
\]
consisting of the following data:
\begin{itemize}
\item $(F,F^{2}_{\oplus},F^{0}_{\oplus})\colon\mathcal{C}\to\mathcal{C}'$
is a symmetric monoidal functor from the additive structure of $\mathcal{C}$
to the additive structure of $\mathcal{C}'$;
\item $(F,F^{2}_{\otimes},F^{0}_{\otimes})\colon\mathcal{C}\to\mathcal{C}'$
is a monoidal functor from the multiplicative structure of $\mathcal{C}$
to the multiplicative structure of $\mathcal{C}'$.
\end{itemize}
These data are required to make the following diagrams in $\mathcal{C}'$
commutative for all $A,B,C\in\in\mathcal{C}$:
\begin{alignat}{2}
 & \input{diagrams/d352.tikz} & \qquad & \input{diagrams/d350.tikz}\label{eq:bmf1}
\end{alignat}
\begin{alignat}{2}
 & \input{diagrams/d351.tikz} & \qquad & \input{diagrams/d349.tikz}\label{eq:bmf2}
\end{alignat}
A bimonoidal functor $F$ is called
\begin{itemize}
\item \emph{strong} if $F^{2}_{\oplus}$, $F^{0}_{\oplus}$, $F^{2}_{\otimes}$,
$F^{0}_{\otimes}$ are isomorphisms;
\item \emph{unitary} if $F^{2}_{\oplus}$, $F^{2}_{\otimes}$ are isomorphisms,
and $F^{0}_{\oplus}$, $F^{0}_{\otimes}$ are identities;
\item \emph{strict} if $F^{2}_{\oplus}$, $F^{0}_{\oplus}$, $F^{2}_{\otimes}$,
$F^{0}_{\otimes}$ are identities.
\end{itemize}
\end{defn}

\begin{defn}\label{def:bimon-nt}
 A \emph{bimonoidal} natural transformation between
two bimonoidal functors $F,G\colon\mathcal{C}\to\mathcal{C}'$ is
a natural transformation $\theta\colon F\to G$ between the underlying
functors such that both
\begin{alignat*}{2}
 & \theta\colon(F,F^{2}_{\oplus},F^{0}_{\oplus})\to(G,G^{2}_{\oplus},G^{0}_{\oplus}) & \qquad\textrm{and}\qquad & \theta\colon(F,F^{2}_{\otimes},F^{0}_{\otimes})\to(G,G^{2}_{\otimes},G^{0}_{\otimes})
\end{alignat*}
are monoidal natural transformations.
\end{defn}

\begin{defn}\label{def:bimon-equiv}
Two
(bi)monoidal categories
$\mathcal{C}$ and $\mathcal{C}'$ are \emph{(bi)monoidally equivalent},
if there exist (bi)monoidal functors $F\colon\mathcal{C}\to\mathcal{C}'$
and $G\colon\mathcal{C}'\to\mathcal{C}$ and (bi)monoidal natural
isomorphisms $FG\cong1_{\mathcal{C}'}$ and $GF\cong1_{\mathcal{C}}$.
\end{defn}

\begin{defn}
[\cite{nlab-ring-obj}]\label{def:rig} A
\emph{rig}
in a bimonoidal category $\mathcal{C}$ is a tuple $(R,\mu^{\otimes},\mu^{\oplus},\eta^{\otimes},\eta^{\oplus})$
where:
\begin{itemize}
\item $(R,\mu^{\oplus},\eta^{\oplus})$ is a commutative monoid in $(\mathcal{C},\oplus,\text{\ensuremath{\bbzero}},\alpha^{\oplus},\lambda^{\oplus},\rho^{\oplus},\xi^{\oplus})$;
\item $(R,\mu^{\otimes},\eta^{\otimes})$ is a monoid in $(\mathcal{C},\otimes,\tensorunit,\alpha^{\otimes},\lambda^{\otimes},\rho^{\otimes})$,
\end{itemize}
such that the following diagrams commute:
\begin{alignat}{2}
\input{diagrams/d335.tikz} & \qquad\qquad & \input{diagrams/d337.tikz}\label{eq:rig-1}
\end{alignat}
\begin{alignat}{2}
\input{diagrams/d336.tikz} & \qquad\qquad & \input{diagrams/d338.tikz}\label{eq:rig-2}
\end{alignat}
A
rig $(R,\mu^{\otimes},\mu^{\oplus},\eta^{\otimes},\eta^{\oplus})$
in a symmetric bimonoidal category $\mathcal{C}$ is called \emph{commutative}
if its multiplicative part $(R,\mu^{\otimes},\eta^{\otimes})$ is
a commutative monoid in the multiplicative part of~$\mathcal{C}$.

\end{defn}

\begin{example}
Any
real unital algebra $A$ can be regarded as a
rig $(A,\singledot,1,+,0)$ in the bimonoidal category $\Vect_{\mathbb{R}}$
from Example~\ref{exa:Vect-bimoncat}, with $\singledot$ and $+$
denoting respectively the multiplication and addition in $A$, and
$1$ and $0$ being the following linear maps:
\begin{alignat*}{2}
 & 1\colon\mathbb{R}\to A\colon x\mapsto x\cdot1_{A}, & \qquad\qquad & 0\colon\{0\}\to A\colon0\mapsto0_{A}.
\end{alignat*}
\end{example}

\subsection{Cartesian additive structure}

In
this paper, we shall mostly need bimonoidal categories
with a specific additive structure defined as follows.
\begin{defn}
Let $\mathcal{C}$ be a category with binary products\footnote{This means that any two objects admit a product.}
and a terminal object, denoted here respectively by $\oplus$ and
$\tobj$ (such categories are called \emph
{cartesian}). Then $\mathcal{C}$
can be endowed with the \emph{cartesian symmetric monoidal structure}
$(\mathcal{C},\oplus,\tobj,\alpha^{\oplus},\lambda^{\oplus},\rho^{\oplus},\xi^{\oplus})$
with $\alpha^{\oplus}$, $\lambda^{\oplus}$, $\rho^{\oplus}$, $\xi^{\oplus}$
the canonical isomorphisms induced by the universal property of products
and terminal objects.
\end{defn}

\begin{defn}\label{def:codistributive}
We call a monoidal category $(\mathcal{C},\otimes,\tensorunit,\alpha^{\otimes},\lambda^{\otimes},\rho^{\otimes})$
\emph{codistributive} if:
\begin{itemize}
\item $\mathcal{C}$ has binary products (denoted here by $\oplus$);
\item $\mathcal{C}$ has a terminal object (which we fix and denote by $\tobj$);
\item The natural morphisms
\begin{alignat}{2}
A\otimes(B_{1}\oplus B_{2}) & \xrightarrow{\delta^{l}}(A\otimes B_{1})\oplus(A\otimes B_{2}),\qquad\qquad &\ 
& \tobj\otimes A\xrightarrow{\terminal^{l}}\tobj,\label{eq:22}\\
(B_{1}\oplus B_{2})\otimes A & \xrightarrow{\delta^{r}}(B_{1}\otimes A)\oplus(B_{2}\otimes A), &\ 
& A\otimes\tobj\xrightarrow{\terminal^{r}}\tobj\label{eq:23}
\end{alignat}
induced respectively by the universal properties of $\oplus$ and
$\tobj$, are isomorphisms.
\end{itemize}
\end{defn}

\begin{rem}
The notion from the previous definition is a dual non-symmetric version
of a \emph{distributive symmetric monoidal category}; see~\cite[Definition~2.3.1]{JY_bimon24}.
\end{rem}

\begin{prop}\label{prop:cart-bicat}
Let $(\mathcal{C},\otimes,\tensorunit,\alpha^{\otimes},\lambda^{\otimes},\rho^{\otimes})$
be a codistributive monoidal category. Then there is a tight bimonoidal
category
\begin{equation}
(\mathcal{C},(\oplus,\tobj,\alpha^{\oplus},\lambda^{\oplus},\rho^{\oplus},\xi^{\oplus}),(\otimes,\tensorunit,\alpha^{\otimes},\lambda^{\otimes},\rho^{\otimes}),(\terminal^{l},\terminal^{r}),(\delta^{l},\delta^{r})),\label{eq:induced-bimon}
\end{equation}
where:
\begin{itemize}
\item $(\oplus,\tobj,\alpha^{\oplus},\lambda^{\oplus},\rho^{\oplus},\xi^{\oplus})$
is the cartesian symmetric monoidal structure;
\item $\terminal^{l},\terminal^{r},\delta^{l},\delta^{r}$ are the isomorphisms
as in Definition~\ref{def:codistributive}.
\end{itemize}
\end{prop}

\begin{proof}
The 22 axioms follow from the universal property of products and the
terminal object, and the naturality of the multiplicative structure.
Tightness follows from codistributivity.
\end{proof}

\begin{example}
The monoidal category $\Vect^{\otimes}_{\mathbb{R}}$ defined in Example~\ref{exa:Vect-moncat}
is codistributive, and the bimonoidal category $\Vect_{\mathbb{R}}$
defined in Example~\ref{exa:Vect-bimoncat} has cartesian additive
structure.
\end{example}

\begin{rem}
Proposition~\ref{prop:cart-bicat} has an obvious symmetric version.
\end{rem}

\begin{rem}
Proposition~\ref{prop:cart-bicat} also has a dual symmetric version;
cf.~\cite[Proposition~2.3.2]{JY_bimon24}.
\end{rem}

\subsection{Quotients}
\begin{defn}
An equivalence relation $\sim$ on parallel arrows of a bimonoidal
category is called a \emph{bimonoidal congruence} if $f\sim f'$ and
$g\sim g'$ imply $f\circ g\sim f'\circ g'$ (whenever $f$ and $g$
are composable), $f\otimes g\sim f'\otimes g'$, and $f\oplus g\sim f'\oplus g'$.
\end{defn}

Given a bimonoidal congruence and a morphism $f$, denote by $[f]$
its congruence class.
\begin{defn}\label{def:aj}
Let $\mathcal{C}$ be a bimonoidal category, and $\sim$
a bimonoidal congruence. Introduce the \emph{bimonoidal quotient category}
\begin{equation}
(\mathcal{C}\modrel{\sim},(\oplus,\bbzero,[\alpha^{\oplus}],[\lambda^{\oplus}],[\rho^{\oplus}],[\xi^{\oplus}]),(\otimes,\tensorunit,[\alpha^{\otimes}],[\lambda^{\otimes}],[\rho^{\otimes}]),([\lambda^{\bullet}],[\rho^{\bullet}]),([\delta^{l}],[\delta^{r}]))\label{eq:ad}
\end{equation}
where $\mathcal{C}\modrel{\sim}$ is the quotient category, and $\oplus$
and $\otimes$ are defined on morphisms as $[f]\otimes[g]\coloneqq[f\otimes g]$
and $[f]\oplus[g]\coloneqq[f\oplus g]$.
\end{defn}

\begin{rem}\label{rem:quot-tight-strict}
It is straightforward to check that
(\ref{eq:ad}) is indeed a bimonoidal category. It is clear that $\mathcal{C}\modrel{\sim}$
is tight whenever $\mathcal{C}$ is, and that $\mathcal{C}\modrel{\sim}$
has a strict multiplicative/additive structure whenever $\mathcal{C}$
does.
\end{rem}

\begin{lem}
The following correspondence is a strict bimonoidal functor:
\[
\mathcal{C}\to\mathcal{C}\modrel{\sim}\colon(A\xrightarrow{f}B)\mapsto(A\xrightarrow{[f]}B).
\]
\end{lem}

\begin{proof}
Straightforward.
\end{proof}

\begin{lem}\label{lem:ak}
Let $\mathcal{C}$, $\mathcal{C}'$ be bimonoidal categories,
$\sim$ a bimonoidal congruence on $\mathcal{C}$, and $F\colon\mathcal{C}\to\mathcal{C}'$
a bimonoidal functor, such that $Ff=Fg$ whenever $f\sim g$. Then
$F$ factors through the quotient:
\[
\input{diagrams/d418.tikz}
\]
with $\overline{F}$ a bimonoidal functor. If $F$ is full, then so
is $\overline{F}$. If $F$ is strict, then so is $\overline{F}$.
\end{lem}

\begin{proof}
Define the functor
$\overline{F}\colon\mathcal{C}\modrel{\sim}\to\mathcal{C}'\colon(A\xrightarrow{[f]}B)\mapsto(FA\xrightarrow{Ff}FB)$.
Then
$(\overline{F},F^{2}_{\otimes},F^{0}_{\otimes},F^{2}_{\oplus},F^{0}_{\oplus})$
is the required bimonoidal functor.
\end{proof}

\subsection{Bracketing}

In this section, we describe a method for constructing strict bimonoidal
functors from unitary bimonoidal ones. This is useful for transferring
diagrams to another category while preserving all coherence relations,
which
will allow us to easily establish the rig structure
of the stabilization endofunctor in Proposition~\ref{prop:rig-K}.
The approach is inspired by~\cite{becerra2024strictification}, from
which we also adopt some notation.
\begin{defn}
Let $X$ be a set. Denote by $X^{\fr}$ the \emph{free $\{\oplus,\otimes\}$-algebra}
generated by $X$, which is defined inductively as follows:
\begin{itemize}
\item $X\subset X^{\fr}$;
\item if $a,b\in X^{\fr}$, then $a\otimes b\in X^{\fr}$ and $a\oplus b\in X^{\fr}$.
\end{itemize}
We also define the \emph{length} of an element in $X^{\fr}$ inductively
as
\[
|x|\coloneqq\begin{cases}
1, & x\in X;\\
|y|+|z|, & x=y\otimes z~\textrm{or}~x=y\oplus z~~\textrm{for some}~y,z\in X^{\fr}.
\end{cases}
\]
Denote by $\mathscr{W}\coloneqq\{\mmph\}^{\fr}$ the free $\{\oplus,\otimes\}$-algebra
generated by the singleton $\{\mmph\}$.
\end{defn}

\begin{defn}
A \emph{bracketed object} in a bimonoidal category $\mathcal{C}$
is a tuple $(S,w)$, where $S\coloneqq(A_{1},A_{2},\dots,A_{n})$
is a non-empty finite sequence of objects of $\mathcal{C}$ and $w\in\mathscr{W}$
such that $|w|=|S|$. For every bracketed object $(S,w)$, we define
the \emph{underlying object} $\u(S,w)\in\in\mathcal{C}$ by replacing
the placeholders $\mmph$ in $w$ with the items from $S$, for example,
\[
\u((A,B,C,D),(\mmph\oplus\mmph)\otimes(\mmph\oplus\mmph))=(A\oplus B)\otimes(C\oplus D).
\]
Abusing notation, we write $(A,\mmph)$ instead of $((A),\mmph)$.
\end{defn}

\begin{defn}
For a bimonoidal category $\mathcal{C}$, define the bimonoidal category
$\Br(\mathcal{C})$ as follows:
\begin{itemize}
\item Objects of $\Br(\mathcal{C})$ are the bracketed objects in $\mathcal{C}$;
\item Morphisms in $\Br(\mathcal{C})$ are the morphisms between the underlying
objects, i.e.,
\[
\Br(\mathcal{C})((S,w),(S',w'))\coloneqq\mathcal{C}(\u(S,w),\u(S',w'));
\]
\item The unit and zero are defined as $(\tensorunit,\mmph)$ and $(\bbzero,\mmph)$,
respectively;
\item For $\square\in\{\otimes,\oplus\}$, we define
\[
(S_{1},w_{1})\square(S_{2},w_{2})\coloneqq(S_{1}*S_{2},w_{1}\square w_{2})
\]
where $S_{1}*S_{2}$ stands for the concatenation of $S_{1}$ and
$S_{2}$;
\item The morphisms $\alpha^{\oplus},~\lambda^{\oplus},~\rho^{\oplus},~\alpha^{\otimes},~\lambda^{\otimes},~\rho^{\otimes},~\lambda^{\bullet},~\rho^{\bullet},~\delta^{l},~\delta^{r}$
are defined as those for the underlying objects; for example, $\alpha^{\otimes}$
is given as
\begin{multline*}
\u\left(((S_{1},w_{1})\otimes(S_{2},w_{2}))\otimes(S_{3},w_{3})\right)=(\u(S_{1},w_{1})\otimes\u(S_{2},w_{2}))\otimes\u(S_{3},w_{3})\\
\xrightarrow{\alpha^{\otimes}}\u(S_{1},w_{1})\otimes(\u(S_{2},w_{2})\otimes\u(S_{3},w_{3}))=\u\left((S_{1},w_{1})\otimes((S_{2},w_{2})\otimes(S_{3},w_{3}))\right).
\end{multline*}
\end{itemize}
It is straightforward to verify that $\Br(\mathcal{C})$ is well-defined.
Furthermore, we can define the fully faithful functor $\i\colon\mathcal{C}\to\Br(\mathcal{C})$
sending $A$ to $(A,\mmph)$, which allows us to regard $\mathcal{C}$
as a full subcategory of $\Br(\mathcal{C})$ and call $\i$ the inclusion
functor.

\end{defn}

\begin{prop}\label{prop:reflector}
The
correspondences $\i$
and $\u$ are\ 
bimonoidal functors, with $\i$ unitary and $\u$
strict. Furthermore, $\i\colon\mathcal{C}\rightleftarrows\Br(\mathcal{C})\noloc\u$
is a bimonoidal equivalence.
\end{prop}

\begin{proof}
Straightforward.
\end{proof}

Let $\mathcal{C}$ and $\mathcal{C}'$ be bimonoidal categories, and
$F\colon\mathcal{C}\to\mathcal{C}'$ a unitary bimonoidal functor.
For every sequence of objects $S=(A_{1},A_{2},\dots,A_{n})$, let
us define the sequence
\[
F[S]\coloneqq(FA_{1},FA_{2},\dots,FA_{n}).
\]
For every $(S,w)\in\in\Br(\mathcal{C})$, introduce the isomorphism
$\Omega_{F}\colon\u(F[S],w)\xrightarrow{\cong}F\u(S,w)$ defined inductively
as:
\begin{alignat*}{1}
\Omega_{F}\colon & \u(FA,\mmph)\xrightarrow{=}F\u(A,\mmph),\\
\Omega_{F}\colon & \u(F[S_{1}*S_{2}],w_{1}\square w_{2})=\u(F[S_{1}],w_{1})\square\u(F[S_{2}],w_{2})\xrightarrow{\Omega_{F}\square\Omega_{F}}F\u(S_{1},w_{1})\square F\u(S_{2},w_{2})\\
 & \qquad\qquad\xrightarrow{F^{2}_{\square}}F(\u(S_{1},w_{1})\square\u(S_{2},w_{2}))=F\u(S_{1}*S_{2},w_{1}\square w_{2}),~~~~~\textrm{for }\square\in\{\otimes,\oplus\}.
\end{alignat*}

\begin{thm}\label{thm:defFbar}
Every
unitary bimonoidal functor
$F\colon\mathcal{C}\to\mathcal{C}'$ gives rise to the strict bimonoidal
functor
\[
\overline{F}\colon\Br(\mathcal{C})\to\Br(\mathcal{C}')
\]
defined at object $(S,w)$ as $(F[S],w)$, and at arrow $(S,w)\xrightarrow{\varphi}(S',w')$
as the composite
\[
\u(F[S],w)\xrightarrow{\Omega_{F}}F\u(S,w)\xrightarrow{F\varphi}F\u(S',w')\xrightarrow{\Omega^{-1}_{F}}\u(F[S'],w').
\]
Furthermore, if $F$ is full or faithful, then so is $\overline{F}$.
\end{thm}

\begin{proof}
It is clear that $\overline{F}$ is a functor. Furthermore, we have
\begin{alignat*}{2}
\overline{F}(S_{1},w_{1})\otimes\overline{F}(S_{2},w_{2})=\overline{F}((S_{1},w_{1})\otimes(S_{2},w_{2})), & \qquad & \overline{F}(\tensorunit,\mmph)=(F\tensorunit,\mmph)=(\tensorunit,\mmph),\\
\overline{F}(S_{1},w_{1})\oplus\overline{F}(S_{2},w_{2})=\overline{F}((S_{1},w_{1})\oplus(S_{2},w_{2})), &\ 
& \overline{F}(\bbzero,\mmph)=(F\bbzero,\mmph)=(\bbzero,\mmph).
\end{alignat*}
That $\overline{F}$ preserves the associativity isomorphisms follows
from the diagram
\[
\input{diagrams/d489.tikz}
\]
in which the left and right subdiagrams commute by the definition
of $\Omega_{F}$; the upper rectangle commutes by the naturality of
$\alpha^{\otimes}$; and the lower rectangle does since $F$ is a
monoidal functor. Hence, $\overline{F}\alpha^{\otimes}=\Omega^{-1}_{F}\circ F\alpha^{\otimes}\circ\Omega_{F}=\alpha^{\otimes}$.
The rest of the proof is similar.

The last statement follows from the equality $\overline{F}\varphi=\Omega^{-1}_{F}\circ F\varphi\circ\Omega_{F}$.
\end{proof}

\begin{prop}\label{prop:natOmega}
The morphisms $\Omega_{F}$ defined above give
rise to the natural isomorphism $\Omega_{F}\colon\u\overline{F}\xRightarrow{\cong}F\u$.
\end{prop}

\begin{proof}
This follows from the definition of $\overline{F}$.
\end{proof}

\begin{prop}\label{prop:brquot}
Every bimonoidal congruence $\sim$ on a bimonoidal
category $\mathcal{C}$ gives rise to a bimonoidal congruence on $\Br(\mathcal{C})$
defined as follows: $\varphi\sim\psi$ in $\Br(\mathcal{C})$ if and
only if $\u\varphi\sim\u\psi$ in $\mathcal{C}$. Furthermore, $\Br(\mathcal{C}\modrel{\sim})=\Br(\mathcal{C})\modrel{\sim}$.
\end{prop}

\begin{proof}
It is clear that $\sim$ is indeed a bimonoidal congruence on $\Br(\mathcal{C})$.
To show that $\Br(\mathcal{C}\modrel{\sim})$ and $\Br(\mathcal{C})\modrel{\sim}$
coincide, just note that objects of both categories are precisely
the bracketed
objects of $\mathcal{C}$, and morphisms
of both categories are precisely the congruence classes of the morphisms
between the underlying objects in $\mathcal{C}$.
\end{proof}

%% file: diagrams/d582.tikz
\begin{tikzcd}[ampersand replacement=\&]
{\mathcal{C}_{3}} \&  \& {\mathcal{C}_{2}} \&  \& {\mathcal{C}_{1}.}
\arrow[""{name=0, anchor=center, inner sep=0}, "G'"{description}, from=1-3, to=1-1]
\arrow[""{name=1, anchor=center, inner sep=0}, "G''"{description}, shift left, curve={height=-30pt}, from=1-3, to=1-1]
\arrow[""{name=2, anchor=center, inner sep=0}, " G"{description}, shift right, curve={height=30pt}, from=1-3, to=1-1]
\arrow[""{name=3, anchor=center, inner sep=0}, "F'"{description}, from=1-5, to=1-3]
\arrow[""{name=4, anchor=center, inner sep=0}, "F''"{description}, shift left, curve={height=-30pt}, from=1-5, to=1-3]
\arrow[""{name=5, anchor=center, inner sep=0}, " F"{description}, shift right, curve={height=30pt}, from=1-5, to=1-3]
\arrow["\beta", Rightarrow, fixedlength,  from=2, to=0]
\arrow["\beta'", Rightarrow, fixedlength,  from=0, to=1]
\arrow["\alpha", Rightarrow, fixedlength,  from=5, to=3]
\arrow["\alpha'", Rightarrow, fixedlength,  from=3, to=4]
\end{tikzcd}

%% file: diagrams/d325.tikz
\begin{tikzcd}[column sep=tiny, ampersand replacement=\&]  %d325
 \& {{(W\otimes(X\otimes Y))\otimes Z}} \&  \\
{{((W\otimes X)\otimes Y)\otimes Z}} \&  \& {{W\otimes((X\otimes Y)\otimes Z)}} \\
{{(W\otimes X)\otimes(Y\otimes Z)}} \&  \& {{W\otimes(X\otimes(Y\otimes Z)),}}
\arrow["{\alpha_{W,X\otimes Y,Z}}", from=1-2, to=2-3]
\arrow["{\alpha_{W,X,Y}\otimes Z}", from=2-1, to=1-2]
\arrow["{\alpha_{W\otimes X,Y,Z}}"', from=2-1, to=3-1]
\arrow["{W\otimes\alpha_{X,Y,Z}}", from=2-3, to=3-3]
\arrow["{\alpha_{W,X,Y\otimes Z}}", from=3-1, to=3-3]
\end{tikzcd}

%% file: diagrams/d324.tikz
\begin{tikzcd}[column sep=tiny, ampersand replacement=\&]  %d324
{{(X\otimes\tensorunit)\otimes Y}} \& {{\phantom{(W\otimes X)\otimes(Y\otimes Z)}}} \& {{X\otimes(\tensorunit\otimes Y)}} \\
{{\phantom{(W\otimes X)\otimes(Y\otimes Z)}}} \& {{X\otimes Y.}} \& {{\phantom{(W\otimes X)\otimes(Y\otimes Z)}}}
\arrow["{\alpha_{X,\tensorunit,Y}}", from=1-1, to=1-3]
\arrow["{\rho_{X}\otimes Y}"', from=1-1, to=2-2]
\arrow["{X\otimes\lambda_{Y}}", from=1-3, to=2-2]
\end{tikzcd}

%% file: diagrams/d431.tikz
\begin{tikzcd}[ampersand replacement=\&]
{X\otimes(Z\otimes Y)} \& {X\otimes(Y\otimes Z)} \\
{(X\otimes Z)\otimes Y} \& {(X\otimes Y)\otimes Z} \\
{Y\otimes(X\otimes Z)} \& {(Y\otimes X)\otimes Z\mmdc}
\arrow["X\otimes\xi_{Z,Y}", from=1-1, to=1-2]
\arrow["\alpha_{X,Y,Z}^{-1}", from=1-2, to=2-2]
\arrow["\alpha_{X,Z,Y}", from=2-1, to=1-1]
\arrow["\xi_{X\otimes Z,Y}"', from=2-1, to=3-1]
\arrow["\alpha_{Y,X,Z}^{-1}"', from=3-1, to=3-2]
\arrow["\xi_{Y,X\otimes Z}"', from=3-2, to=2-2]
\end{tikzcd}

%% file: diagrams/d432.tikz
\begin{tikzcd}[column sep=tiny, ampersand replacement=\&]  %d432
 \& {Y\otimes X} \\
{X\otimes Y} \&  \& {X\otimes Y\mmdc}
\arrow["\xi_{Y,X}", from=1-2, to=2-3]
\arrow["\xi_{X,Y}", from=2-1, to=1-2]
\arrow[equals, from=2-3, to=2-1]
\end{tikzcd}

%% file: diagrams/d433.tikz
\begin{tikzcd}[column sep=small, ampersand replacement=\&]  %d433
{X\otimes\tensorunit} \&  \& {\tensorunit\otimes X} \\
 \& {X\mmdd}
\arrow["\xi_{X,\tensorunit}", from=1-1, to=1-3]
\arrow["\rho_{X}"', from=1-1, to=2-2]
\arrow["\lambda_{X}", from=1-3, to=2-2]
\end{tikzcd}

%% file: diagrams/d327.tikz
\begin{tikzcd}[ampersand replacement=\&]
{(FX\otimes FY)\otimes FZ} \&  \& {FX\otimes(FY\otimes FZ)} \\
{F(X\otimes Y)\otimes FZ} \&  \& {FX\otimes F(Y\otimes Z)} \\
{F((X\otimes Y)\otimes Z)} \&  \& {F(X\otimes(Y\otimes Z))\mmdc}
\arrow["\alpha_{FX,FY,FZ}", from=1-1, to=1-3]
\arrow["F_{X,Y}^{2}\otimes FZ"', from=1-1, to=2-1]
\arrow["FX\otimes F_{Y,Z}^{2}", from=1-3, to=2-3]
\arrow["F_{X\otimes Y,Z}^{2}"', from=2-1, to=3-1]
\arrow["F_{X,Y\otimes Z}^{2}", from=2-3, to=3-3]
\arrow["F\alpha_{X,Y,Z}"', from=3-1, to=3-3]
\end{tikzcd}

%% file: diagrams/d328.tikz
\begin{tikzcd}[ampersand replacement=\&]
{FX\otimes\tensorunit_{\mathcal{C}}} \& {FX} \\
{FX\otimes F\tensorunit_{\mathcal{C}'}} \& {F(X\otimes\tensorunit_{\mathcal{C}'})\mmdd}
\arrow["\rho_{FX}", from=1-1, to=1-2]
\arrow["FX\otimes F^{0}"', from=1-1, to=2-1]
\arrow["F_{X,\tensorunit}^{2}"', from=2-1, to=2-2]
\arrow["F\rho_{X}"', from=2-2, to=1-2]
\end{tikzcd}

%% file: diagrams/d345.tikz
\begin{tikzcd}[column sep=large, ampersand replacement=\&]  %d345
{FX\otimes FY} \& {FY\otimes FX} \\
{F(X\otimes Y)} \& {F(Y\otimes X)\mmdd}
\arrow["\xi_{FX,FY}", from=1-1, to=1-2]
\arrow["F_{X,Y}^{2}"', from=1-1, to=2-1]
\arrow["F_{Y,X}^{2}", from=1-2, to=2-2]
\arrow["F\xi_{X,Y}", from=2-1, to=2-2]
\end{tikzcd}

%% file: diagrams/d330.tikz
\begin{tikzcd}[column sep=large, ampersand replacement=\&]  %d330
{FX\otimes FY} \& {GX\otimes GY} \\
{F(X\otimes Y)} \& {G(X\otimes Y)}
\arrow["\theta X\otimes\theta Y", from=1-1, to=1-2]
\arrow["F^{2}"', from=1-1, to=2-1]
\arrow["F^{2}", from=1-2, to=2-2]
\arrow["\theta(X\otimes Y)"', from=2-1, to=2-2]
\end{tikzcd}

%% file: diagrams/d331.tikz
\begin{tikzcd}[ampersand replacement=\&]
 \& {\tensorunit_{\mathcal{C}'}} \\
{F\tensorunit_{\mathcal{C}}} \&  \& {G\tensorunit_{\mathcal{C}}}
\arrow["F^{0}"', from=1-2, to=2-1]
\arrow["G^{0}", from=1-2, to=2-3]
\arrow["\theta\tensorunit_{\mathcal{C}}"', from=2-1, to=2-3]
\end{tikzcd}

%% file: diagrams/d332.tikz
\begin{tikzcd}[column sep=small, ampersand replacement=\&]  %d332
{X\otimes(X\otimes X)} \&  \& {(X\otimes X)\otimes X} \\
{X\otimes X} \& {X} \& {X\otimes X\mmdc}
\arrow["\alpha", from=1-1, to=1-3]
\arrow["X\otimes\mu"', from=1-1, to=2-1]
\arrow["\mu\otimes X", from=1-3, to=2-3]
\arrow["\mu", from=2-1, to=2-2]
\arrow["\mu"', from=2-3, to=2-2]
\end{tikzcd}

%% file: diagrams/d333.tikz
\begin{tikzcd}[ampersand replacement=\&]
{\tensorunit\otimes X} \& {X\otimes X} \& {X\otimes\tensorunit} \\
 \& {X\mmdd}
\arrow["1\otimes X", from=1-1, to=1-2]
\arrow["\lambda"', from=1-1, to=2-2]
\arrow["\mu", from=1-2, to=2-2]
\arrow["X\otimes1"', from=1-3, to=1-2]
\arrow["\rho", from=1-3, to=2-2]
\end{tikzcd}

%% file: diagrams/d334.tikz
\begin{tikzcd}[ampersand replacement=\&]
{X\otimes X} \&  \& {X\otimes X} \\
 \& {X\mmdd}
\arrow["\xi", from=1-1, to=1-3]
\arrow["\mu"', from=1-1, to=2-2]
\arrow["\mu", from=1-3, to=2-2]
\end{tikzcd}

%% file: diagrams/d352.tikz
\begin{tikzcd}[ampersand replacement=\&]
{FA\otimes(FB\oplus FC)} \& {(FA\otimes FB)\oplus(FA\otimes FC)} \\
{FA\otimes F(B\oplus C)} \& {F(A\otimes B)\oplus F(A\otimes C)} \\
{F(A\otimes(B\oplus C))} \& {F((A\otimes B)\oplus(A\otimes C))\mmdc}
\arrow["\delta^{l}", from=1-1, to=1-2]
\arrow["FA\otimes F_{\oplus}^{2}"', from=1-1, to=2-1]
\arrow["F_{\otimes}^{2}\oplus F_{\otimes}^{2}", from=1-2, to=2-2]
\arrow["F_{\otimes}^{2}"', from=2-1, to=3-1]
\arrow["F_{\oplus}^{2}", from=2-2, to=3-2]
\arrow["F\delta^{l}"', from=3-1, to=3-2]
\end{tikzcd}

%% file: diagrams/d350.tikz
\begin{tikzcd}[column sep=tiny, ampersand replacement=\&]  %d350
{\bbzero\otimes FA} \&  \& {F\bbzero\otimes FA} \\
{\bbzero} \&  \& {F(\bbzero\otimes A)} \\
 \& {F\bbzero\mmdc}
\arrow["F_{\oplus}^{0}\otimes FA", from=1-1, to=1-3]
\arrow["\lambda^{\bullet}"', from=1-1, to=2-1]
\arrow["F_{\otimes}^{2}", from=1-3, to=2-3]
\arrow["F_{\oplus}^{0}"', from=2-1, to=3-2]
\arrow["F\lambda^{\bullet}", from=2-3, to=3-2]
\end{tikzcd}

%% file: diagrams/d351.tikz
\begin{tikzcd}[ampersand replacement=\&]
{(FA\oplus FB)\otimes FC} \& {(FA\otimes FC)\oplus(FB\otimes FC)} \\
{F(A\oplus B)\otimes(FC)} \& {F(A\otimes C)\oplus F(B\otimes C)} \\
{F((A\oplus B)\otimes C)} \& {F((A\otimes C)\oplus(B\otimes C))\mmdc}
\arrow["\delta^{r}", from=1-1, to=1-2]
\arrow["F_{\oplus}^{2}\otimes FC"', from=1-1, to=2-1]
\arrow["F_{\otimes}^{2}\oplus F_{\otimes}^{2}", from=1-2, to=2-2]
\arrow["F_{\otimes}^{2}"', from=2-1, to=3-1]
\arrow["F_{\oplus}^{2}", from=2-2, to=3-2]
\arrow["F\delta^{r}"', from=3-1, to=3-2]
\end{tikzcd}

%% file: diagrams/d349.tikz
\begin{tikzcd}[column sep=tiny, ampersand replacement=\&]  %d349
{FA\otimes\bbzero} \&  \& {FA\otimes F\bbzero} \\
{\bbzero} \&  \& {F(A\otimes\bbzero)} \\
 \& {F\bbzero\mmdd}
\arrow["FA\otimes F_{\oplus}^{0}", from=1-1, to=1-3]
\arrow["\rho^{\bullet}"', from=1-1, to=2-1]
\arrow["F_{\otimes}^{2}", from=1-3, to=2-3]
\arrow["F_{\oplus}^{0}"', from=2-1, to=3-2]
\arrow["F\rho^{\bullet}", from=2-3, to=3-2]
\end{tikzcd}

%% file: diagrams/d335.tikz
\begin{tikzcd}[column sep=tiny, ampersand replacement=\&]  %d335
{R\otimes(R\oplus R)} \&  \& {(R\otimes R)\oplus(R\otimes R)} \\
{R\otimes R} \& {R} \& {R\oplus R\mmdc}
\arrow["\delta^{l}", from=1-1, to=1-3]
\arrow["R\otimes\mu^{\oplus}"', from=1-1, to=2-1]
\arrow["\mu^{\otimes}\oplus\mu^{\otimes}", from=1-3, to=2-3]
\arrow["\mu^{\otimes}"', from=2-1, to=2-2]
\arrow["\mu^{\oplus}", from=2-3, to=2-2]
\end{tikzcd}

%% file: diagrams/d337.tikz
\begin{tikzcd}[ampersand replacement=\&]
{\bbzero\otimes R} \& {\bbzero} \\
{R\otimes R} \& {R\mmdc}
\arrow["\lambda^{\bullet}", from=1-1, to=1-2]
\arrow["\eta^{\oplus}\otimes R"', from=1-1, to=2-1]
\arrow["\eta^{\oplus}", from=1-2, to=2-2]
\arrow["\mu^{\otimes}"', from=2-1, to=2-2]
\end{tikzcd}

%% file: diagrams/d336.tikz
\begin{tikzcd}[column sep=tiny, ampersand replacement=\&]  %d336
{(R\oplus R)\otimes R} \&  \& {(R\otimes R)\oplus(R\otimes R)} \\
{R\otimes R} \& {R} \& {R\oplus R\mmdc}
\arrow["\delta^{r}", from=1-1, to=1-3]
\arrow["\mu^{\oplus}\otimes R"', from=1-1, to=2-1]
\arrow["\mu^{\otimes}\oplus\mu^{\otimes}", from=1-3, to=2-3]
\arrow["\mu^{\otimes}"', from=2-1, to=2-2]
\arrow["\mu^{\oplus}", from=2-3, to=2-2]
\end{tikzcd}

%% file: diagrams/d338.tikz
\begin{tikzcd}[ampersand replacement=\&]
{R\otimes\bbzero} \& {\bbzero} \\
{R\otimes R} \& {R\mmdd}
\arrow["\rho^{\bullet}", from=1-1, to=1-2]
\arrow["R\otimes\eta^{\oplus}"', from=1-1, to=2-1]
\arrow["\eta^{\oplus}", from=1-2, to=2-2]
\arrow["\mu^{\otimes}"', from=2-1, to=2-2]
\end{tikzcd}

%% file: diagrams/d418.tikz
\begin{tikzcd}[ampersand replacement=\&]
{\mathcal{C}} \&  \& {\mathcal{C}'} \\
{\mathcal{C}\modrel{\sim}}
\arrow[" F", from=1-1, to=1-3]
\arrow[from=1-1, to=2-1]
\arrow["\overline{F}"', curve={height=12pt}, from=2-1, to=1-3]
\end{tikzcd}

%% file: diagrams/d489.tikz
\begin{tikzcd}[column sep=small, ampersand replacement=\&]  %d489
{\big(\u(F[S_{1}],w_{1})\otimes\u(F[S_{2}],w_{2})\big)\otimes\u(F[S_{3}],w_{3})} \& {\u(F[S_{1}],w_{1})\otimes\big(\u(F[S_{2}],w_{2})\otimes\u(F[S_{3}],w_{3})\big)} \\
{\big(F\u([S_{1}],w_{1})\otimes F\u([S_{2}],w_{2})\big)\otimes F\u([S_{3}],w_{3})} \& {F\u([S_{1}],w_{1})\otimes\big(F\u([S_{2}],w_{2})\otimes F\u([S_{3}],w_{3})\big)} \\
{F\big(\big(\u([S_{1}],w_{1})\otimes\u([S_{2}],w_{2})\big)\otimes\u([S_{3}],w_{3})\big)} \& {F\big(\u([S_{1}],w_{1})\otimes\big(\u([S_{2}],w_{2})\otimes\u([S_{3}],w_{3})\big)\big)}
\arrow["\alpha^{\otimes}", from=1-1, to=1-2]
\arrow["(\Omega_{F}\otimes\Omega_{F})\otimes\Omega_{F}", from=1-1, to=2-1]
\arrow["\Omega_{F}", sqarw=2em, from=1-1, to=3-1]
\arrow["\Omega_{F}\otimes(\Omega_{F}\otimes\Omega_{F})"', from=1-2, to=2-2]
\arrow["\Omega_{F}"', sqare=2em, from=1-2, to=3-2]
\arrow["\alpha^{\otimes}", from=2-1, to=2-2]
\arrow["\Omega_{F}", from=2-1, to=3-1]
\arrow["\Omega_{F}"', from=2-2, to=3-2]
\arrow["F\alpha^{\otimes}", from=3-1, to=3-2]
\end{tikzcd}

%% file: sec2.tex
\section{Functor categories}

For categories $\mathcal{C}$ and $\mathcal{C}'$, denote by $\mathcal{C}'^{\mathcal{C}}$
the category with functors from $\mathcal{C}$ to $\mathcal{C}'$
as objects, and natural transformations between them as morphisms.

If some property holds for every component of the natural transformation
$\alpha\colon F\Rightarrow G$, then we say that $\alpha$ enjoys
this property \emph{componentwise}. For example, the obvious natural
transformation $C_{0}([0,\infty),\mmph)\Rightarrow C_{b}([0,\infty),\mmph)$
between the functors from the introduction is a componentwise ideal
inclusion, and $C_{b}([0,\infty),\mmph)\Rightarrow\mathfrak{A}$ is
a componentwise quotient projection onto the asymptotic
algebra functor.
\begin{lem}\label{lem:cw-epic-epic}
If a natural transformation is componentwise
monic (resp. componentwise epic), then it is monic (resp. epic).
\end{lem}

\begin{proof}
Straightforward.
\end{proof}

\subsection{Limits in functor categories}

Here
we recall a few well-known facts concerning
functor categories, some of which can be found in e.g.~\cite{Leinster2014basicCT}
or~\cite{Riehl2016CategoryInContext}.

Let $\mathcal{J}$, $\mathcal{K}$ and $\mathcal{C}$ be categories,
with $\mathcal{K}$ and $\mathcal{J}$ small, and $\mathcal{C}$ complete.
 Denote by $e_{j}\colon\mathcal{C}^{\mathcal{J}}\to\mathcal{C}$
the functor of evaluation at $j\in\in\mathcal{J}$. We also have the
functor of $\lim\colon\mathcal{C}^{\mathcal{J}}\to\mathcal{C}$ and
the natural transformation
\begin{equation}
\pr_{j}\colon\lim\Rightarrow e_{j},\qquad\pr_{j}D\colon\lim_{\mathcal{J}}D\to D_{j}\label{eq:lim-eval}
\end{equation}
where $\pr_{j}D\colon\lim_{\mathcal{J}}D\to D_{j}$ stands for the
$j$th leg of the limit cone over the diagram $D\colon\mathcal{J}\to\mathcal{C}$.
By abuse of notation, we shall write $\pr_{j}$ instead of $\pr_{j}D$.
\begin{lem}\label{lem:iter-lims}
Let $H\colon\mathcal{J}\times\mathcal{K}\to\mathcal{C}$
be a small diagram. Then the outer diagram below commutes,
\[
\input{diagrams/d266.tikz}
\]
and the dashed arrow induced by the universal property of $\lim_{(j,k)}H(j,k)$
is an isomorphism.
\end{lem}

\begin{defn}\label{def:ab}
Let $\mathcal{C}$, $\mathcal{C}'$, $\mathcal{C}''$
and $\mathcal{J}$ be categories. Define the \emph{canonical arrows}
(provided the limits below exist)
\begin{alignat}{3}
\rc\colon(\lim_{j}F_{j})G & \Rightarrow & \lim_{j}(F_{j}G), &\ 
& \textrm{ where } & F_{\bullet}\colon\mathcal{J}\to\mathcal{C}''^{\mathcal{C}'},~G\colon\mathcal{C}\to\mathcal{C}',\label{eq:canar1}\\
\lc\colon G(\lim_{j}F_{j}) & \Rightarrow & \lim_{j}(GF_{j}), &\ 
& \textrm{ where } & F_{\bullet}\colon\mathcal{J}\to\mathcal{C}'^{\mathcal{C}},~G\colon\mathcal{C}'\to\mathcal{C}''\label{eq:canar11}
\end{alignat}
as the unique arrows making the following diagrams commutative:
\begin{alignat}{2}
\input{diagrams/d226.tikz} & \qquad & \input{diagrams/d310.tikz}\label{eq:canar2}
\end{alignat}
\end{defn}

\begin{rem}
If we choose $\mathcal{C}$ in (\ref{eq:canar1}) and (\ref{eq:canar2})
to be the terminal category\footnote{\emph{The terminal category} is a category with one object and its
identity morphism.}, we automatically get the following arrows, which we also call canonical
arrows:
\begin{alignat*}{3}
\rc\colon(\lim_{j}F_{j})B & \to & \lim_{j}(F_{j}B), &\ 
& \textrm{ where } & F_{\bullet}\colon\mathcal{J}\to\mathcal{C}'^{\mathcal{C}},~B\in\in\mathcal{C};\\
\lc\colon G(\lim_{j}B_{j}) & \to & \lim_{j}(GB_{j}), &\ 
& \textrm{ where } & G\colon\mathcal{C}\to\mathcal{C}',~B_{\bullet}\colon\mathcal{J}\to\mathcal{C}.
\end{alignat*}
\end{rem}

\begin{lem}\label{lem:almost-bimon}
Canonical arrows satisfy the following associativity
laws:
\begin{alignat}{2}
\input{diagrams/d315.tikz} & \qquad & \input{diagrams/d316.tikz}\label{eq:d315}
\end{alignat}
\begin{equation}
\input{diagrams/d317.tikz}\label{eq:d317}
\end{equation}
\end{lem}

\begin{proof}
The commutativity of (\ref{eq:d317}) is proved by the following commuting
diagram:
\begin{equation}
\input{diagrams/d341.tikz}\label{eq:limitsoffunctors-1}
\end{equation}
in which the triangles commute by the definition of $\rc$ and $\lc$,
and the outer rectangle commutes by the universal property of $\lim_{i}(G_{1}F_{i}G_{2})$.
The proof for the diagrams in~(\ref{eq:d315}) is analogous.
\end{proof}

The following lemma says that limits in functor categories are precisely
componentwise limits. This fact is well known and we state it without
proof.
\begin{lem}\label{lem:cw-lim}
Let $\mathcal{J}$, $\mathcal{C}$ and $\mathcal{C}'$
be
categories with $\mathcal{J}$ small and $\mathcal{C}'$
complete, and let $F_{\bullet}\colon\mathcal{J}\to\mathcal{C}'^{\mathcal{C}}$
be a diagram. Then there is a functor
\[
\Lambda\colon\mathcal{C}\to\mathcal{C}'\colon(A\xrightarrow{\varphi}B)\mapsto(\lim_{i}(F_{i}A)\xrightarrow{\lim(F_{i}\varphi\circ\pr_{i})}\lim_{i}(F_{i}B))
\]
and natural transformations
\[
\pi_{j}\colon\Lambda\Rightarrow F_{j},\qquad\pi_{j}A\colon\lim_{i}(F_{i}A)\xrightarrow{\pr_{j}}F_{j}A
\]
for all $j\in\in\mathcal{J}$, such that $\{\pi_{j}\colon\Lambda\Rightarrow F_{j}\}_{j\in\in\mathcal{J}}$
is a limit cone over $F_{\bullet}$.
\end{lem}

\begin{lem}\label{lem:rc-isom-id}
Let
$\mathcal{J}$, $\mathcal{C}$,
$\mathcal{C}'$, and $F_{\bullet}$ be as in the previous lemma, let
$\mathcal{D}$ be a category, and let $G\colon\mathcal{D}\to\mathcal{C}$
be a functor. Then the canonical arrow $(\lim F_{\bullet})G\xRightarrow{\rc}\lim(F_{\bullet}G)$
is an isomorphism. Furthermore, one can choose $\lim F_{\bullet}$
such that this canonical arrow is the identity natural transformation.
\end{lem}

\begin{proof}
We claim that $\Lambda$ from the previous lemma is the choice of
$\lim F_{\bullet}$ making $\rc$ the identity natural transformation.
Indeed, it is easy to obtain this when $\mathcal{D}$ is a terminal
category. This fact being combined with the left-hand side diagram
of (\ref{eq:d315}) in Lemma~\ref{lem:almost-bimon} proves the claim.

The first statement follows from the fact that the limit cone is unique
up to isomorphism.
\end{proof}

\begin{lem}\label{lem:lim-prod}
Let $\mathcal{C}$, $\mathcal{C}'$ and $\mathcal{C}''$
be categories with $\mathcal{C}'$ and $\mathcal{C}''$ complete.
Let $G_{\bullet}\colon\mathcal{K}\to\mathcal{C}'^{\mathcal{C}}$ and
$F_{\bullet}\colon\mathcal{J}\to\mathcal{C}''^{\mathcal{C}'}$ be
two small diagrams. Then the diagram
\begin{equation}
\input{diagrams/d318.tikz}\label{eq:d318}
\end{equation}
(in which the unnamed arrows both denote the isomorphisms from Lemma~\ref{lem:iter-lims})
commutes, and the arrows marked with $\cong$ are isomorphisms.
\end{lem}

\begin{proof}
We claim that the following diagram commutes:
\begin{equation}
\input{diagrams/d488.tikz}\label{eq:lims1-1}
\end{equation}
Indeed, triangle~$\tc 1$ commutes by the definition of the canonical
arrow. Triangle~$\tc 2$ commutes since it is obtained by applying
the functor
$\lim_{\mathcal{J}}$ to the following
commutative triangle of $\mathcal{J}$-shaped diagrams:
\begin{equation}
\input{diagrams/d340.tikz}\label{eq:nat-canonar-1}
\end{equation}
Parallelogram~$\tc 3$ commutes by the naturality of $\pr_{j}$.
 Finally,~$\tc 4$ commutes by Lemma~\ref{lem:iter-lims}, and~$\tc 5$
commutes trivially.

By a similar argument, we can prove the commutativity of
\begin{equation}
\input{diagrams/d267.tikz}\label{eq:lims2-1}
\end{equation}
The commutativity of~(\ref{eq:d318}) follows from diagrams~(\ref{eq:lims1-1})
and~(\ref{eq:lims2-1}) together with the universal property of $\lim_{(j,k)}F_{j}G_{k}$.

The statement about isomorphisms follows from Lemma~\ref{lem:rc-isom-id}.
\end{proof}

\begin{defn}\label{def:Delta}
Let $\mathcal{J}$, $\mathcal{C}$ be categories,
$\mathcal{J}$ small, $\mathcal{C}$ complete, and let $A\in\in\mathcal{C}$.
Denote by $\Delta$ the unique arrow making the following diagram
commutative for all $j\in\in\mathcal{J}$:
\[
\input{diagrams/d406.tikz}
\]
\end{defn}

\begin{lem}\label{lem:Delta}
Let $\mathcal{C}$, $\mathcal{C}'$, $\mathcal{C}''$
and $\mathcal{J}$ be categories with $\mathcal{C}'$ and $\mathcal{C}''$
complete and $\mathcal{J}$ small, and let $G\colon\mathcal{C}\to\mathcal{C}'$
and $F\colon\mathcal{C}'\to\mathcal{C}''$ be functors. Then the following
diagrams commute:
\begin{alignat*}{2}
\input{diagrams/d407.tikz} & \qquad & \input{diagrams/d435.tikz}
\end{alignat*}
\end{lem}

\begin{proof}
The proof is similar to that of Lemma~\ref{lem:almost-bimon}.
\end{proof}

\begin{lem}\label{lem:comonoid}
Let $(\mathcal{C},\oplus,\tobj,\alpha^{\oplus},\lambda^{\oplus},\rho^{\oplus},\xi^{\oplus})$
be a cartesian symmetric monoidal category, and let $A\in\in\mathcal{C}$.
Then the following diagram commutes:
\[
\input{diagrams/d424.tikz}
\]
\end{lem}

\begin{proof}
Straightforward.
\end{proof}

\begin{rem}
In fact, every object $A$ in a cartesian symmetric monoidal category
gives rise to a commutative comonoid $(A,\Delta,!)$, where $!$ is
the unique arrow $A\to\tobj$. However, we shall not need this fact.
\end{rem}

\subsection{Cartesian monoidal structure}

In this section, we lift the cartesian monoidal structure to the level
of the corresponding endofunctors. Throughout the section, we assume
that $\mathcal{C}$ is a category with binary products $\oplus$ and
a terminal object $\tobj$.

\begin{defn}\label{def:oplus}
For endofunctors $F,F',G,G'\colon\mathcal{C}\to\mathcal{C}$
and natural transformations $\alpha\colon F\Rightarrow F'$ and $\beta\colon G\Rightarrow G'$,
let us introduce:
\begin{itemize}
\item the endofunctor
$F\oplus G$ defined as
\[
F\oplus G\colon(A\xrightarrow{\varphi}B)\mapsto(FA\oplus GA\xrightarrow{F\varphi\oplus G\varphi}FB\oplus GB);
\]
\item the natural transformation $\alpha\oplus\beta\colon F\oplus G\Rightarrow F'\oplus G'$
defined at component $B\in\in\mathcal{C}$ as
\[
(\alpha\oplus\beta)B\coloneqq\alpha B\oplus\beta B.
\]
\end{itemize}
It is straightforward to check that $\oplus$ is well-defined and
functorial.
\end{defn}

\begin{rem}\label{rem:rc-eq}
It follows
from Lemmas~\ref{lem:cw-lim}
and~\ref{lem:rc-isom-id} that:
\begin{itemize}
\item $F\oplus F'$ is a product of $F$ and $F'$ in $\mathcal{C}^{\mathcal{C}}$;
\item the canonical arrow $(F\oplus F')G\xRightarrow{\rc}FG\oplus F'G$
is the identity natural transformation for any $G\in\in\mathcal{C}^{\mathcal{C}}$.
\end{itemize}
\end{rem}

\begin{defn}
For $A\in\in\mathcal{C}$, we introduce the constant endofunctor $\Const_{A}\colon\mathcal{C}\to\mathcal{C}$
defined as
\[
\Const_{A}\colon(B\to B')\mapsto(A\xrightarrow{=}A).
\]
For every morphism $\varphi\colon A\to A'$, we also define the obvious
natural transformation $\Const_{\varphi}\colon\Const_{A}\Rightarrow\Const_{A'}$.
\end{defn}

\begin{defn}\label{def:tbpp}
 We
denote by $\bpp{\mathcal{C}}$
the class of endofunctors in $\mathcal{C}$ which preserve binary
products and terminal objects, i.e., the class of all $F\colon\mathcal{C}\to\mathcal{C}$
such that
\begin{alignat*}{2}
\lc\colon F(A\oplus B)\to FA\oplus FB, & \qquad & \qquad!\colon F\tobj\to\tobj
\end{alignat*}
are isomorphisms for all $A,B\in\in\mathcal{C}$, where $\terminal$
denotes the unique arrow into the terminal object.
\end{defn}

Since $\Const_{\tobj}$ is a terminal object in $\mathcal{C}^{\mathcal{C}}$,
let us denote by $\Terminal^{l}$ and $\Terminal^{r}$ the unique
arrows
\begin{alignat*}{2}
\Const_{\tobj}F\xRightarrow{\Terminal^{l}}\Const_{\tobj} & ,\qquad\qquad & F\Const_{\tobj}\xRightarrow{\Terminal^{r}}\Const_{\tobj} & .
\end{alignat*}
Note also that $\Terminal^{l}$ is an identity natural transformation,
and that
\begin{equation}
\input{diagrams/d363.tikz}\label{eq:TerminalFunctor}
\end{equation}
commutes.

It is easy to check that $(\bpp{\mathcal{C}},\singledot,\Id)$ is
a monoidal subcategory of $(\mathcal{C},\singledot,\Id)$ from Example~\ref{exa:str-moncat}.
\begin{lem}\label{lem:bpp-codist}
$(\bpp{\mathcal{C}},\singledot,\Id)$ is codistributive.
\end{lem}

\begin{proof}
The terminal object $\Const_{\tobj}$ is contained in $\bpp{\mathcal{C}}$
since $\tobj\oplus\tobj\cong\tobj$. For $F_{1},F_{2},G_{1},G_{2}\in\in\bpp{\mathcal{C}}$,
the diagram
\[
\input{diagrams/d357.tikz}
\]
commutes by Lemma~\ref{lem:lim-prod}, and the solid arrows are isomorphisms
by Definition~\ref{def:tbpp} and Lemma~\ref{lem:rc-isom-id}. Thus,
the dashed arrow is also an isomorphism, and so $F_{1}\oplus F_{2}\in\in\bpp{\mathcal{C}}$.
We have proved that $\bpp{\mathcal{C}}$ has binary products.

We need to show that $\lc$, $\rc$, $\Terminal^{l}$, $\Terminal^{r}$
are isomorphisms. Indeed, $\rc$ is an identity by Lemma~\ref{lem:rc-isom-id};
$\lc$~is an isomorphism by Definition~\ref{def:tbpp}; $\Terminal^{l}$~is
an identity; and $\Terminal^{r}$ is an isomorphism by the commutativity
of~(\ref{eq:TerminalFunctor}).
\end{proof}

\begin{thm}\label{thm:funcat-bimonoidal}
Every cartesian
category
$\mathcal{C}$ gives rise to the tight bimonoidal category
\begin{equation}
(\bpp{\mathcal{C}},(\oplus,\Const_{\tobj},\au,\lu,\ru,\xi^{\oplus}),(\singledot,\Id),(\Terminal^{l},\Terminal^{r}),(\lc,\rc))\label{eq:bimon-tbpp}
\end{equation}
with the cartesian additive structure.
\end{thm}

\begin{proof}
This follows from Lemma~\ref{lem:bpp-codist} and Proposition~\ref{prop:cart-bicat}.
\end{proof}

\begin{lem}\label{lem:coher-component}
The coherence isomorphisms of the additive
structure of (\ref{eq:bimon-tbpp}) are related to those of $\mathcal{C}$
as follows:
\begin{alignat}{3}
\alpha^{\oplus}_{F,G,H}B= & \alpha^{\oplus}_{FB,GB,HB}, & \qquad\qquad &\ 
& \xi^{\oplus}_{F,G}B= & \xi^{\oplus}_{FB,GB},\label{eq:1414}\\
\lambda^{\oplus}_{F}B= & \lambda^{\oplus}_{FB}, &\ 
&\ 
& \rho^{\oplus}_{F}B= & \rho^{\oplus}_{FB}.\label{eq:1515}
\end{alignat}
\end{lem}

\begin{proof}
Combining (\ref{eq:canar2}) from Definition~\ref{def:ab} with Definition~\ref{def:oplus}
we get the commutative diagram
\[
\input{diagrams/d437.tikz}
\]
from which one can deduce~(\ref{eq:1414}). The proof of~(\ref{eq:1515})
is left as an exercise.
\end{proof}
\subsection{Tensoring}
\begin{defn}\label{def:tensoring}
Let $\mathcal{C}$ be a bimonoidal category.
For $A\in\in\mathcal{C}$, define the \emph{left tensoring endofunctor}
\[
\mathfrak{O}_{A}\coloneqq A\otimes\mmph\colon\mathcal{C}\to\mathcal{C}\colon(B\xrightarrow{\varphi}B')\mapsto(A\otimes B\xrightarrow{A\otimes\varphi}A\otimes B'),
\]
and for $\varphi\in\mathcal{C}(A,A')$, define the natural transformation
\[
\mathfrak{O}_{\varphi}\colon\mathfrak{O}_{A}\Rightarrow\mathfrak{O}_{A'},\qquad\mathfrak{O}_{\varphi}B\colon A\otimes B\xrightarrow{\varphi\otimes B}A'\otimes B.
\]
\end{defn}
\begin{prop}\label{prop:bimon-O}
Let $\mathcal{C}$ be a codistributive monoidal
category. There is a bimonoidal functor
\[
(\mathfrak{O},\mathfrak{O}^{2},\mathfrak{O}^{0},\mathfrak{O}^{2}_{\oplus},\mathfrak{O}^{0}_{\oplus})\colon\mathcal{C}\to\bpp{\mathcal{C}}
\]
where $\mathcal{C}$ and $\bpp{\mathcal{C}}$ stand respectively for
the tight bimonoidal categories
\begin{alignat*}{1}
 & (\mathcal{C},(\oplus,\tobj,\alpha^{\oplus},\lambda^{\oplus},\rho^{\oplus},\xi^{\oplus}),(\otimes,\tensorunit,\alpha^{\otimes},\lambda^{\otimes},\rho^{\otimes}),(\terminal^{l},\terminal^{r}),(\delta^{l},\delta^{r}))),\\
 & (\bpp{\mathcal{C}},(\oplus,\Const_{\tobj},\au,\lu,\ru,\xi^{\oplus}),(\Id,\singledot),(\Terminal^{l},\Terminal^{r}),(\lc,\rc)),
\end{alignat*}
and where:
\begin{itemize}
\item $\mathfrak{O}\colon\mathcal{C}\to\bpp{\mathcal{C}}\colon(A\xrightarrow{\varphi}A')\mapsto(\mathfrak{O}_{A}\xRightarrow{\mathfrak{O}_{\varphi}}\mathfrak{O}_{A'})$;
\item $(\mathfrak{O}^{2})_{A,A'}\colon\mathfrak{O}_{A}\mathfrak{O}_{A'}\Rightarrow\mathfrak{O}_{A\otimes A'},~(\mathfrak{O}^{2})_{A,A'}B\colon A\otimes(A'\otimes B)\xrightarrow{(\alpha^{\otimes})^{-1}}(A\otimes A')\otimes B$;
\item $\mathfrak{O}^{0}\colon\Id\Rightarrow\mathfrak{O}_{\tensorunit},~\mathfrak{O}^{0}B\colon B\xrightarrow{(\lambda^{\otimes})^{-1}}\tensorunit\otimes B$;
\item $(\mathfrak{O}^{2}_{\oplus})_{A,A'}\colon\mathfrak{O}_{A}\oplus\mathfrak{O}_{A'}\Rightarrow\mathfrak{O}_{A\oplus A'},~(\mathfrak{O}^{2}_{\oplus})_{A,A'}B\colon(A\otimes B)\oplus(A'\otimes B)\xrightarrow{(\delta^{r})^{-1}}(A\oplus A')\otimes B$;
\item $\mathfrak{O}^{0}_{\oplus}\colon\Const_{\tobj}\Rightarrow\mathfrak{O}_{\tobj},~\mathfrak{O}^{0}_{\oplus}B\colon\tobj\xrightarrow{(\terminal^{l})^{-1}}\tobj\otimes B$.
\end{itemize}
\end{prop}

\begin{proof}
It follows from Definition~\ref{def:codistributive} that $\mathfrak{O}_{A}\in\in\bpp{\mathcal{C}}$
for every $A\in\in\mathcal{C}$. One can deduce from the universal
properties of products and terminal objects that $(\mathfrak{O},\mathfrak{O}^{2}_{\oplus},\mathfrak{O}^{0}_{\oplus})$
is a symmetric monoidal functor. It follows from the coherence theorem
for the multiplicative structure of $\mathcal{C}$ that $(\mathfrak{O},\mathfrak{O}^{2},\mathfrak{O}^{0})$
is a monoidal functor. The commutativity of the four diagrams in Definition~\ref{def:bimon-fun}
follows from Laplaza axioms\footnote{To be exact, we need the right-hand side diagrams of~(\ref{eq:d557})
and (\ref{eq:d558}), as well as the central and the right-hand side
diagrams of~(\ref{eq:d568}) in the Appendix.} for~$\mathcal{C}$.
\end{proof}

\begin{rem}\label{rem:bimon-ext}
The
previous proposition is
a special case of a more general fact: any monoidal functor between
the multiplicative structures of two codistributive bimonoidal categories
with cartesian additive structure gives rise to a bimonoidal functor
between these categories, cf.~Proposition~5.1.11 from~\cite{JY_bimon24}.
\end{rem}

\begin{lem}\label{lem:218}
In the setting of Proposition~\ref{prop:bimon-O},
the following diagram commutes for all $A_{1},A_{2}\in\in\mathcal{C}$
and~$j=1,2$:
\[
\input{diagrams/d369.tikz}
\]
\end{lem}

\begin{proof}
Straightforward.
\end{proof}

%% file: diagrams/d266.tikz
\begin{tikzcd}[ampersand replacement=\&]
 \& {\lim_{k}H(j,k)} \\
{\lim_{j}\lim_{k}H(j,k)} \& {\lim_{(j,k)}H(j,k)} \& {H(j,k)} \\
 \& {\lim_{j}H(j,k)}
\arrow["\pr_{k}", curve={height=-12pt}, from=1-2, to=2-3]
\arrow["\pr_{j}", curve={height=-12pt}, from=2-1, to=1-2]
\arrow["\cong", dashed, from=2-1, to=2-2]
\arrow["\lim_{j}\pr_{k}"', curve={height=12pt}, from=2-1, to=3-2]
\arrow["\pr_{(j,k)}", from=2-2, to=2-3]
\arrow["\pr_{j}"', curve={height=12pt}, from=3-2, to=2-3]
\end{tikzcd}

%% file: diagrams/d226.tikz
\begin{tikzcd}[ampersand replacement=\&]
{(\lim_{j}F_{j})G} \& {\lim_{j}(F_{j}G)} \\
 \& {F_{j}G\mmdc}
\arrow["\rc", Rightarrow, dashed, from=1-1, to=1-2]
\arrow["\pr_{j}G"', Rightarrow, from=1-1, to=2-2]
\arrow["\pr_{j}", Rightarrow, from=1-2, to=2-2]
\end{tikzcd}

%% file: diagrams/d310.tikz
\begin{tikzcd}[ampersand replacement=\&]
{G(\lim_{j}F_{j})} \& {\lim_{j}(GF_{j})} \\
 \& {GF_{j}\mmdd}
\arrow["\lc", Rightarrow, dashed, from=1-1, to=1-2]
\arrow["G\pr_{j}"', Rightarrow, from=1-1, to=2-2]
\arrow["\pr_{j}", Rightarrow, from=1-2, to=2-2]
\end{tikzcd}

%% file: diagrams/d315.tikz
\begin{tikzcd}[ampersand replacement=\&]
 \& {(\lim_{i}F_{i}G_{1})G_{2}} \\
{(\lim_{i}F_{i})G_{1}G_{2}} \& {\lim_{i}(F_{i}G_{1}G_{2})\mmdc}
\arrow["\rc", Rightarrow, from=1-2, to=2-2]
\arrow["\rc G_{2}", Rightarrow, from=2-1, to=1-2]
\arrow["\rc"', Rightarrow, from=2-1, to=2-2]
\end{tikzcd}

%% file: diagrams/d316.tikz
\begin{tikzcd}[ampersand replacement=\&]
 \& {G_{1}(\lim_{i}G_{2}F_{i})} \\
{G_{1}G_{2}\lim_{i}F_{i}} \& {\lim_{i}(G_{1}G_{2}F_{i})\mmdc}
\arrow["\lc", Rightarrow, from=1-2, to=2-2]
\arrow["G_{1}\lc", Rightarrow, from=2-1, to=1-2]
\arrow["\lc"', Rightarrow, from=2-1, to=2-2]
\end{tikzcd}

%% file: diagrams/d317.tikz
\begin{tikzcd}[ampersand replacement=\&]
{G_{1}(\lim_{i}F_{i})G_{2}} \&  \& {(\lim_{i}G_{1}F_{i})G_{2}} \\
{G_{1}\lim_{i}(F_{i}G_{2})} \&  \& {\lim_{i}(G_{1}F_{i}G_{2})\mmdd}
\arrow["\lc G_{2}", Rightarrow, from=1-1, to=1-3]
\arrow[""', draw=none, from=1-1, to=1-3]
\arrow["G_{1}\rc"', Rightarrow, from=1-1, to=2-1]
\arrow["\rc", Rightarrow, from=1-3, to=2-3]
\arrow["\lc"', Rightarrow, from=2-1, to=2-3]
\arrow[""', draw=none, from=2-1, to=2-3]
\end{tikzcd}

%% file: diagrams/d341.tikz
\begin{tikzcd}[column sep=small, ampersand replacement=\&]  %d341
{G_{1}(\lim_{i}F_{i})G_{2}} \&  \&  \&  \&  \&  \& {G_{1}\lim_{i}(F_{i}G_{2})} \\
 \&  \&  \& {G_{1}F_{j}G_{2}} \\
{\lim_{i}(G_{1}F_{i})G_{2}} \&  \&  \&  \&  \&  \& {\lim_{i}(G_{1}F_{i}G_{2})}
\arrow["G_{1}\rc", Rightarrow, from=1-1, to=1-7]
\arrow[draw=none, from=1-1, to=1-7]
\arrow["G_{1}\pr_{j}G_{2}"'{pos=0.4}, Rightarrow, from=1-1, to=2-4]
\arrow["\lc G_{2}"', Rightarrow, from=1-1, to=3-1]
\arrow["G_{1}\pr_{j}"{pos=0.4}, Rightarrow, from=1-7, to=2-4]
\arrow["\lc", Rightarrow, from=1-7, to=3-7]
\arrow["\pr_{j}G_{2}"{pos=0.4}, Rightarrow, from=3-1, to=2-4]
\arrow["\rc"', Rightarrow, from=3-1, to=3-7]
\arrow[draw=none, from=3-1, to=3-7]
\arrow["G_{1}\pr_{j}G_{2}"'{pos=0.4}, Rightarrow, from=3-7, to=2-4]
\end{tikzcd}

%% file: diagrams/d318.tikz
\begin{tikzcd}[ampersand replacement=\&]
{{(\lim_{j}F_{j})(\lim_{k}G_{k})}} \&  \& {{\lim_{j}(F_{j}\lim_{k}G_{k})}} \&  \& {{\lim_{j}\lim_{k}(F_{j}G_{k})}} \\
{{\lim_{k}((\lim_{j}F_{j})G_{k})}} \&  \& {{\lim_{k}(\lim_{j}F_{j}G_{k})}} \&  \& {{\lim_{(j,k)}F_{j}G_{k}}}
\arrow["\rc", Rightarrow, from=1-1, to=1-3]
\arrow["\cong"', draw=none, from=1-1, to=1-3]
\arrow["\lc"', Rightarrow, from=1-1, to=2-1]
\arrow["{\lim_{j}\lc}", Rightarrow, from=1-3, to=1-5]
\arrow["\cong"', Rightarrow, from=1-5, to=2-5]
\arrow["\cong", Rightarrow, from=2-1, to=2-3]
\arrow["{\lim_{k}\rc}"', draw=none, from=2-1, to=2-3]
\arrow["\cong", Rightarrow, from=2-3, to=2-5]
\end{tikzcd}

%% file: diagrams/d488.tikz
\begin{tikzcd}[ampersand replacement=\&]
{(\lim_{j}F_{j})(\lim_{k}G_{k})} \& {\lim_{j}(F_{j}\lim_{k}G_{k})} \&  \& {\lim_{j}\lim_{k}(F_{j}G_{k})} \\
 \& {F_{j}\lim_{k}G_{k}} \& {\lim_{j}(F_{j}G_{k})} \\
{\,} \&  \& {F_{j}G_{k}} \& {\lim_{(j,k)}F_{j}G_{k}\mmdd}
\arrow["\rc", Rightarrow, from=1-1, to=1-2]
\arrow[""{name=0, anchor=center, inner sep=0}, "\pr_{j}\lim_{k}G_{k}"', Rightarrow, from=1-1, to=2-2]
\arrow["\pr_{j}\pr_{k}"{pos=0.7}, manv, Rightarrow, from=1-1, to=3-3]
\arrow[""{name=1, anchor=center, inner sep=0}, "\lim_{j}\lc", Rightarrow, from=1-2, to=1-4]
\arrow["\pr_{j}", Rightarrow, from=1-2, to=2-2]
\arrow[""{name=2, anchor=center, inner sep=0}, "\lim_{j}(F_{j}\pr_{k})", Rightarrow, from=1-2, to=2-3]
\arrow[""{name=3, anchor=center, inner sep=0}, "\lim_{j}\pr_{k}", Rightarrow, from=1-4, to=2-3]
\arrow["\cong", Rightarrow, from=1-4, to=3-4]
\arrow[""{name=4, anchor=center, inner sep=0}, "F_{j}\pr_{k}"', Rightarrow, from=2-2, to=3-3]
\arrow["\pr_{j}", Rightarrow, from=2-3, to=3-3]
\arrow["\tc 5"{pos=0.8,description}, draw=none, from=3-1, to=2-2]
\arrow["\pr_{(j,k)}", Rightarrow, from=3-4, to=3-3]
\arrow["\tc 1"{description, pos=0.4}, draw=none, from=1-2, to=0]
\arrow["\tc 2"{description}, draw=none, from=1, to=3]
\arrow["\tc 3"{description}, draw=none, from=2, to=4]
\arrow["\tc 4"{description}, shift right=2, draw=none, from=3, to=3-4]
\end{tikzcd}

%% file: diagrams/d340.tikz
\begin{tikzcd}[ampersand replacement=\&]
{F_{\bullet}\lim_{k}G_{k}} \&  \& {\lim_{k}(F_{\bullet}G_{k})} \\
 \& {F_{\bullet}G_{k}\mmdd}
\arrow["\lc", Rightarrow, from=1-1, to=1-3]
\arrow["F_{\bullet}\pr_{k}"', Rightarrow, from=1-1, to=2-2]
\arrow["\pr_{k}", Rightarrow, from=1-3, to=2-2]
\end{tikzcd}

%% file: diagrams/d267.tikz
\begin{tikzcd}[ampersand replacement=\&]
{(\lim_{j}F_{j})(\lim_{k}G_{k})} \& {\lim_{k}((\lim_{j}F_{j})G_{k})} \&  \& {\lim_{k}\lim_{j}(F_{j}G_{k})} \\
 \& {(\lim_{j}F_{j})G_{k}} \& {\lim_{k}(F_{j}G_{k})} \\
 \&  \& {F_{j}G_{k}} \& {\lim_{(j,k)}F_{j}G_{k}\mmdd}
\arrow["\lc", Rightarrow, from=1-1, to=1-2]
\arrow["(\lim_{j}F_{j})\pr_{k}"', Rightarrow, from=1-1, to=2-2]
\arrow["\pr_{j}\pr_{k}"{pos=0.7}, Rightarrow, manv, from=1-1, to=3-3]
\arrow["\lim_{k}\rc", Rightarrow, from=1-2, to=1-4]
\arrow["\pr_{k}", Rightarrow, from=1-2, to=2-2]
\arrow["\lim_{k}\pr_{j}G_{k}", Rightarrow, from=1-2, to=2-3]
\arrow["\lim_{k}\pr_{j}", Rightarrow, from=1-4, to=2-3]
\arrow["\cong", Rightarrow, from=1-4, to=3-4]
\arrow["\pr_{j}G_{k}"', Rightarrow, from=2-2, to=3-3]
\arrow["\pr_{k}", Rightarrow, from=2-3, to=3-3]
\arrow["\pr_{(j,k)}", Rightarrow, from=3-4, to=3-3]
\end{tikzcd}

%% file: diagrams/d406.tikz
\begin{tikzcd}[ampersand replacement=\&]
{A} \& {\lim_{\mathcal{J}}A} \\
 \& {A.}
\arrow["\Delta", dashed, from=1-1, to=1-2]
\arrow[" ="', from=1-1, to=2-2]
\arrow["\pr_{j}", from=1-2, to=2-2]
\end{tikzcd}

%% file: diagrams/d407.tikz
\begin{tikzcd}[column sep=small, ampersand replacement=\&]  %d407
{FG} \&  \& {\lim_{\mathcal{J}}(FG)} \\
 \& {F\lim_{\mathcal{J}}G\mmdc}
\arrow["\Delta", Rightarrow, from=1-1, to=1-3]
\arrow["F\Delta"', Rightarrow, from=1-1, to=2-2]
\arrow["\lc"', Rightarrow, from=2-2, to=1-3]
\end{tikzcd}

%% file: diagrams/d435.tikz
\begin{tikzcd}[column sep=small, ampersand replacement=\&]  %d435
{FG} \&  \& {\lim_{\mathcal{J}}(F)G} \\
 \& {(\lim_{\mathcal{J}}F)G\mmdd}
\arrow["\Delta", Rightarrow, from=1-1, to=1-3]
\arrow["\Delta G"', Rightarrow, from=1-1, to=2-2]
\arrow["\rc"', Rightarrow, from=2-2, to=1-3]
\end{tikzcd}

%% file: diagrams/d424.tikz
\begin{tikzcd}[ampersand replacement=\&]
{A\oplus A} \& {A} \& {A\oplus A} \\
{(A\oplus A)\oplus A} \&  \& {A\oplus(A\oplus A)\mmdd}
\arrow["\Delta\oplus A"', from=1-1, to=2-1]
\arrow["\Delta"', from=1-2, to=1-1]
\arrow["\Delta", from=1-2, to=1-3]
\arrow["A\oplus\Delta", from=1-3, to=2-3]
\arrow["\alpha^{\oplus}", from=2-1, to=2-3]
\end{tikzcd}

%% file: diagrams/d363.tikz
\begin{tikzcd}[ampersand replacement=\&]
{F\Const_{\tobj}} \& {\Const_{\tobj}} \\
{\Const_{F\tobj}}
\arrow["\Terminal^{r}", Rightarrow, from=1-1, to=1-2]
\arrow[equals, from=1-1, to=2-1]
\arrow["\Const_{\terminal}"', curve={height=12pt}, Rightarrow, from=2-1, to=1-2]
\arrow["\cong", curve={height=12pt}, draw=none, from=2-1, to=1-2]
\end{tikzcd}

%% file: diagrams/d357.tikz
\begin{tikzcd}[ampersand replacement=\&]
{(F_{1}\oplus F_{2})(G_{1}\oplus G_{2})} \& {(F_{1}\oplus F_{2})G_{1}\oplus(F_{1}\oplus F_{2})G_{2}} \\
{F_{1}(G_{1}\oplus G_{2})\oplus F_{2}(G_{1}\oplus G_{2})} \& {(F_{1}G_{1}\oplus F_{2}G_{1})\oplus(F_{1}G_{2}\oplus F_{2}G_{2})}
\arrow["\lc", Rightarrow, dashed, from=1-1, to=1-2]
\arrow["\rc"', Rightarrow, from=1-1, to=2-1]
\arrow["\cong", draw=none, from=1-1, to=2-1]
\arrow["\rc\oplus\rc", Rightarrow, from=1-2, to=2-2]
\arrow["\cong"', draw=none, from=1-2, to=2-2]
\arrow["\lc\oplus\lc"', Rightarrow, from=2-1, to=2-2]
\arrow["\cong", draw=none, from=2-1, to=2-2]
\end{tikzcd}

%% file: diagrams/d437.tikz
\begin{tikzcd}[ampersand replacement=\&]
{(F_{1}\oplus F_{2})B} \&  \& {F_{1}B\oplus F_{2}B} \\
 \& {F_{j}G}
\arrow["\rc"', Rightarrow, from=1-1, to=1-3]
\arrow[" =", draw=none, from=1-1, to=1-3]
\arrow["\pr_{j}B"', Rightarrow, from=1-1, to=2-2]
\arrow["\pr_{j}", Rightarrow, from=1-3, to=2-2]
\end{tikzcd}

%% file: diagrams/d369.tikz
\begin{tikzcd}[ampersand replacement=\&]
{\mathfrak{O}_{A_{1}}\oplus\mathfrak{O}_{A_{2}}} \& {\mathfrak{O}_{A_{j}}} \\
{\mathfrak{O}_{A_{1}\oplus A_{2}}\mmdd}
\arrow["\pr_{j}", Rightarrow, from=1-1, to=1-2]
\arrow["\mathfrak{O}_{\oplus}^{2}"', Rightarrow, from=1-1, to=2-1]
\arrow["\mathfrak{O}_{\pr_{j}}"', Rightarrow, from=2-1, to=1-2]
\end{tikzcd}

%% file: sec3.tex
\section{$C^{*}$-algebras and decent endofunctors}

In this section we start studying the
category of
$C^{*}$-algebras\footnote{Not necessarily unital.} and $*$-homomorphisms~\cite{murphy1990},
its homotopy version, and the corresponding functor categories.

\subsection{\label{subsec:cal}The category of $C^{*}$-algebras}

Denote by $\Cstar$ the category of $C^{*}$-algebras and $*$-homomorphisms~\cite{murphy1990}.
The algebraic tensor product over $\mathbb{C}$ of two $C^{*}$-algebras
is not complete (and hence is not a $C^{*}$-algebra) in general,
so in order to define tensoring in the category of $C^{*}$-algebras,
given an algebraic tensor product we need to take its completion with
respect to some possibly non-unique $C^{*}$-norm. Thus, in $C^{*}$-algebra
theory there exist various different tensor products.\ 
The most well-known
examples are the \emph{maximal} and \emph{minimal} tensor products,
denoted respectively by $\otimes_{\max}$ and $\otimes_{\min}$. In
this paper we mostly work with $\otimes_{\max}$ which will be abbreviated
to $\otimes$.

The maximal tensor product distributes over direct sums, in particular
there exist isomorphisms
\begin{alignat*}{2}
\delta^{r}_{A,B_{1},B_{2}}\colon & A\otimes(B_{1}\oplus B_{2}) & \xrightarrow{\cong}(A\otimes B_{1})\oplus(A\otimes B_{2}),\\
\delta^{l}_{B_{1},B_{2},A}\colon & (B_{1}\oplus B_{2})\otimes A & \xrightarrow{\cong}(B_{1}\otimes A)\oplus(B_{2}\otimes A)
\end{alignat*}
which are natural in $A,B_{1},B_{2}\in\in\Cstar$\footnote{In fact, $\otimes$ distributes over infinite direct sums as well.}.

It is also well-known that $\Cstar$ is complete, and that finite
direct sums of $C^{*}$-algebras are products in the categorical sense.
If $\varphi_{1}\colon A_{1}\to B$ and $\varphi_{2}\colon A_{2}\to B$
are two $*$-homomorphisms, then the explicit formula for the pullback
along $\varphi_{1}$ and $\varphi_{2}$ is
\begin{equation}
A_{1}\oplus_{B}A_{2}=\left\{ (a_{1},a_{2})\in A_{1}\oplus A_{2}\mid\varphi_{1}(a_{1})=\varphi_{2}(a_{2})\right\} .\label{eq:explicit-pullback}
\end{equation}

\begin{lem}
[Functional calculus]\label{lem:spectrum} Let
$A$ be a $C^{*}$-algebra, and $a\in A$ a normal element\footnote{An element $a$ of a $C^{*}$-algebra is called \emph{normal} if $aa^{*}=a^{*}a$.}.
Denote by $C^{*}(a)$ the $C^{*}$-algebra generated by $\{a\}$,
and denote
by $\spec(a)\subset\mathbb{C}$ the spectrum
of $a$. Then there exists a unique $*$-isomorphism $C_{0}(\spec(a)\setminus\{0\})\cong C^{*}(a)$
which maps the identity function on $\spec(a)\setminus\{0\}$ to~$a$.
\end{lem}

Denote by $\bbzero$ the zero $C^{*}$-algebra, which is clearly a
zero
object in $\Cstar$.
\begin{prop}\label{prop:cstar-mult-str}
The
category $\Cstar$
carries the structure of a codistributive symmetric monoidal category
\[
(\Cstar,\otimes,\mathbb{C},\alpha^{\otimes},\lambda^{\otimes},\rho^{\otimes},\xi^{\otimes}),
\]
 where:
\begin{itemize}
\item $\otimes\colon\Cstar\times\Cstar\to\Cstar$\textup{ is the maximal
tensor product;}
\item $\alpha^{\otimes}_{A,B,C}\colon(A\otimes B)\otimes C\to A\otimes(B\otimes C)\colon(a\otimes b)\otimes c\mapsto a\otimes(b\otimes c)$;
\item $\lambda^{\otimes}_{A}\colon\mathbb{C}\otimes A\to A\colon z\otimes a\mapsto za$;
\item $\rho^{\otimes}_{A}\colon A\otimes\mathbb{C}\to A\colon a\otimes z\mapsto za$;
\item $\xi^{\otimes}_{A,B}\colon A\otimes B\to B\otimes A\colon a\otimes b\mapsto b\otimes a$.
\end{itemize}
\end{prop}

\begin{proof}
Straightforward.
\end{proof}

The codistributivity of $\Cstar$ allows us to turn it into a bimonoidal
category.
\begin{defn}
Introduce the tight symmetric bimonoidal category
\begin{alignat}{1}
(\Cstar,(\bbzero,\oplus,\au,\lu,\ru,\xi^{\oplus}),(\mathbb{C},\otimes,\alpha^{\otimes},\lambda^{\otimes},\rho^{\otimes},\xi^{\otimes}),(\terminal^{l},\terminal^{r}),(\delta^{l},\delta^{r})),\label{eq:bimon-Cstar}
\end{alignat}
 as in Proposition~\ref{prop:cart-bicat}.
\end{defn}

\begin{defn}
Let $A$ and $B$ be $C^{*}$-algebras. Two $*$-homomorphisms $\varphi_{0},\varphi_{1}\colon A\to B$
are called \emph{homotopic} if there is a $*$-homomorphism $\Phi\colon A\to IB$
such that the following diagram commutes for $j=0,1$:
\[
\input{diagrams/d472.tikz}
\]
where $IB\coloneqq C([0,1],B)$ denotes
the $C^{*}$-algebra
of continuous functions on $[0,1]$ taking values in~$B$, and $e_{t}\colon IB\to B$
stands for the $*$-homomorphism of evaluation at $t\in[0,1]$.

It is straightforward to verify that $\simeq$ is a bimonoidal congruence.
\end{defn}

\begin{defn}
We call the symmetric bimonoidal quotient category $\hCstar\coloneqq\Cstar\modrel{\simeq}$
the \emph{homotopy category} \emph{of $C^{*}$-algebras}.

\end{defn}

\subsection{\label{subsec:rig-comp}The rig of compact operators}

Denote by $\oM_{n}$ the $C^{*}$-algebra of $n$-by-$n$ matrices
with entries in $\mathbb{C}$, and by $\oK(\mathcal{H})$ the $C^{*}$-algebra
of compact operators on a separable infinite-dimensional Hilbert space
$\mathcal{H}$. In most cases, the choice of $\mathcal{H}$ is not
important and we often use the abbreviation $\oK\coloneqq\oK(\mathcal{H})$.
We write $\epsilon_{i,j}$ for the $(i,j)$th matrix unit in $\oK$
or in~$\oM_{n}$.

In this section,
we prove that $\oK$ carries the
structure of a rig in $\hCstar$. In Section~\ref{subsec:starig},
we promote this construction to the level of endofunctors and use
it in Section~\ref{sec:genhom} to define bilinear composition on
generalized morphisms.
\begin{lem}\label{lem:Ads-homot}
Let $\mathcal{H}$ and $\mathcal{H}'$ be separable
infinite dimensional Hilbert spaces, and let $V_{j}\colon\mathcal{H}\to\mathcal{H}'$
be isometries for~$j=0,1$. Then the following two $*$-homomorphisms
are homotopic:
\[
\Ad V_{j}\colon\oK(\mathcal{H})\to\oK(\mathcal{H}')\colon k\mapsto V_{j}kV^{*}_{j},\qquad j=0,1.
\]
\end{lem}

\begin{proof}
This follows from~\cite[Lemma~1.3.7]{Jensen-Thomsen}.
\end{proof}

\begin{cor}\label{cor:ae}
The symmetry isomorphism $\xi^{\otimes}_{\oK,\oK}\colon\oK\otimes\oK\to\oK\otimes\oK$
is homotopic to the identity $\id_{\oK\otimes\oK}$.
\end{cor}

\begin{proof}
This follows from Lemma~\ref{lem:Ads-homot}.
\end{proof}

We now introduce the following $*$-homomorphisms:
\begin{itemize}
\item $\widehat{\iota}\colon\mathbb{C}\oplus\mathbb{C}\to\oM_{2}\colon(z_{1},z_{2})\mapsto\diag(z_{1},z_{2})$
is the diagonal embedding;
\item $\widehat{\iota}_{00}\colon\mathbb{C}\to\oK\colon z\mapsto\diag(z,0,0,\dots)$
is the corner embedding;
\item $\widehat{\theta}\colon\oK\otimes\oK\to\oK$ is the isomorphism\footnote{$\widehat{\theta}$ and $\widehat{\theta}_{2}$ are often called \emph{stability
isomorphisms}.} induced by some unitary isomorphism $l_{2}(\mathbb{N})\otimes l_{2}(\mathbb{N})\xrightarrow{\cong}l_{2}(\mathbb{N})$;
\item $\widehat{\theta}_{2}\colon\oM_{2}\otimes\oK\to\oK$ is the isomorphism
induced by some unitary isomorphism $l_{2}(\{0,1\})\otimes l_{2}(\mathbb{N})\xrightarrow{\cong}l_{2}(\mathbb{N})$;
\item $\widehat{\mu}\colon\oK\oplus\oK\xrightarrow{(\lambda^{\otimes})^{-1}\oplus(\lambda^{\otimes})^{-1}}(\mathbb{C}\otimes\oK)\oplus(\mathbb{C}\otimes\oK)\xrightarrow{(\delta^{r})^{-1}}(\mathbb{C}\oplus\mathbb{C})\otimes\oK\xrightarrow{\widehat{\iota}\otimes\oK}\oM_{2}\otimes\oK\xrightarrow{\widehat{\theta}_{2}}\oK$.
\end{itemize}
We also introduce the $*$-homomorphism
\[
\widehat{d}\colon\oK\oplus\oK\to\oM_{2}\otimes\oK\colon(k_{1},k_{2})\mapsto\epsilon_{11}\otimes k_{1}+\epsilon_{22}\otimes k_{2}
\]
which, as one can directly check, satisfies the equalities $\widehat{d}=(\widehat{\iota}\otimes\oK)\circ(\delta^{r})^{-1}\circ((\lambda^{\otimes})^{-1}\oplus(\lambda^{\otimes})^{-1})$
and $\widehat{\mu}=\widehat{\theta}_{2}\circ\widehat{d}$.
\begin{prop}\label{prop:rig-K-1}
 $(\oK,[\widehat{\theta}],[\widehat{\iota}_{00}],[\widehat{\mu}],[0])$
is a commutative
rig in $\hCstar$.
\end{prop}

\begin{proof}
One can easily show using Lemma~\ref{lem:Ads-homot} that $(\oK,[\widehat{\theta}],[\widehat{\iota}_{00}])$
is a commutative monoid in $\hCstar$. That $(\oK,[\widehat{\mu}],[0])$
is a commutative monoid follows from~\cite[Lemma~1.3.12]{Jensen-Thomsen}.

Now, we need to check the commutativity of the four diagrams from
Definition~\ref{def:rig}. The right-hand side diagrams in (\ref{eq:rig-1})
and (\ref{eq:rig-2}) commute since $\bbzero$ is an initial object.
The commutativity of the left-hand side diagram in (\ref{eq:rig-1})
follows from that of
\begin{equation}
\input{diagrams/d358.tikz}\label{eq:d358}
\end{equation}
in which the central vertical arrow is an obvious composite involving
associativity and symmetry isomorphisms. The upper and lower triangles
in (\ref{eq:d358}) commute by the definition of $\widehat{\mu}$
and $\widehat{d}$; the left and the bottom-left subdiagrams commute
by a direct computation. To see the commutativity of the right subdiagram,
note that both legs are induced by unitary isomorphisms, and hence
are homotopic by Lemma~\ref{lem:Ads-homot}. The proof for the left-hand
side diagram in (\ref{eq:rig-2}) is similar.
\end{proof}

\subsection{\label{subsec:decent-end}Decent endofunctors}

We begin our study of endofunctors in the category of $C^{*}$-algebras
and $*$-homomorphisms. We restrict our attention only to the class
of sufficiently nice endofunctors which is just a slight elaboration
of binary product-preserving endofunctors considered in the previous
section. Henceforth, the term ``endofunctor'' will always refer
to an endofunctor of $\Cstar$.

To
clarify the motivation for the definition below,
we refer the reader to property~\ref{enu:o1} and the subsequent
paragraph in the introduction.

\begin{defn}\label{def:decentFunctor}
 We call $F\in\in\Cstar^{\Cstar}$ a \emph{decent
endofunctor} if:
\begin{itemize}
\item for every pullback
\[
\input{diagrams/d249.tikz}
\]
where $\varphi_{2}$ is a split epimorphism\footnote{An arrow is called a \emph{split epimorphism} if it has a pre-inverse.},
the canonical arrow $F(B_{1}\underset{B}{\oplus}B_{2})\xrightarrow{\lc}FB_{1}\underset{FB}{\oplus}FB_{2}$
is an isomorphism;
\item The unique arrow $F\bbzero\xrightarrow{\terminal}\bbzero$ is an isomorphism.
\end{itemize}
We denote by $\DEFC$ the full subcategory in $\EFC$ generated by
decent endofunctors. It is clear that $\DEFC$ is in fact a full subcategory
in $\bpp{\Cstar}$.
\end{defn}

We refer the reader to the beginning of Section~\ref{subsec:stacont}
for basic examples of decent endofunctors.
\begin{prop}
$(\DEFC,\singledot,\Id)$ is a well-defined strict monoidal category.
\end{prop}

\begin{proof}
Using the fact that any functor being applied to a split epimorphism
yields again a split epimorphism, we can deduce that the class of
decent endofunctors is closed under composition, from which the statement
follows.
\end{proof}

\begin{lem}\label{lem:pullback-functors-and-c-algebras}
Let $G\in\in\DEFC$,
$\alpha_{1}\in\EFC(F_{1},F_{0})$, $\alpha_{2}\in\EFC(F_{2},F_{0})$
such that $\alpha_{2}$ is a componentwise split epimorphism. Then
the canonical arrow $G(F_{1}\underset{F_{0}}{\oplus}F_{2})\xRightarrow{\lc}GF_{1}\underset{GF_{0}}{\oplus}GF_{2}$
is an isomorphism in $\EFC$.
\end{lem}

\begin{proof}
Let $F_{\bullet}$ denote the diagram $F_{1}\xRightarrow{\alpha_{1}}F_{0}\xLeftarrow{\alpha_{2}}F_{2}$.
We need to show that for every $B\in\in\Cstar$ the canonical arrow
$G(\lim F_{\bullet})B\xrightarrow{\lc B}\lim(GF_{\bullet})B$ is an
isomorphism. Diagram~(\ref{eq:d317}) from Lemma~\ref{lem:almost-bimon}
gives rise to the commuting diagram\footnote{The object $B$ here is regarded as a functor from the terminal category.}
\[
\input{diagrams/d343.tikz}
\]
in which the right vertical arrow is an isomorphism since $G$ is
decent, and the horizontal arrows are isomorphisms by Lemma~\ref{lem:rc-isom-id}.
Therefore, the left vertical arrow is also an isomorphism.
\end{proof}

\begin{prop}\label{prop:decent-complete}
The category $\DEFC$ is complete.
\end{prop}

\begin{proof}
Let $F_{\bullet}\colon\mathcal{J}\to\EFC$ be a small diagram. Denote
by $\mathcal{K}$ the category $\bullet\to\bullet\leftarrow\bullet$,
and fix a diagram of $C^{*}$-algebras $B_{1}\xrightarrow{\varphi_{1}}B_{0}\xleftarrow{\varphi_{2}}B_{2}$
of shape $\mathcal{K}$ with $\varphi_{2}$ a split epimorphism. The
limit $\lim_{j}F_{j}$ exists since $\EFC$ is complete (see~Lemma~\ref{lem:cw-lim}),
and we only need to show that the canonical arrow $(\lim_{j}F_{j})(\lim_{k}B_{k})\to\lim_{k}((\lim_{j}F_{j})B_{k})$
is an isomorphism. Lemma~\ref{lem:lim-prod} yields the commuting
diagram
\[
\input{diagrams/d342.tikz}
\]
in which $\lim_{j}\lc$ is an isomorphism since $F_{j}$ is decent
for all $j\in\in\mathcal{J}$, hence, the leftmost vertical arrow
is also an isomorphism, and the proof is complete.
\end{proof}

\begin{lem}\label{lem:bimon-defc}
The strict monoidal category $(\DEFC,\singledot,\Id)$
is codistributive.
\end{lem}

\begin{proof}
The statement follows from Lemma~\ref{lem:bpp-codist} and the fact
that $\DEFC$ is a full subcategory of $\bpp{\Cstar}$.
\end{proof}

Lemma~\ref{lem:bimon-defc} and Theorem~\ref{thm:funcat-bimonoidal}
allow us to endow $\DEFC$ with a bimonoidal structure.
\begin{defn}\label{def:bimon-defc}
Introduce the bimonoidal category
\[
(\DEFC,(\bzero,\oplus,\au,\lu,\ru,\xi^{\oplus}),(\Id,\singledot),(\Terminal^{l},\Terminal^{l}),(\lc,\rc)),
\]
with $\bzero\coloneqq\Const_{\bbzero}$ the zero endofunctor, and
$(\bzero,\oplus,\au,\lu,\ru,\xi^{\oplus})$ the cartesian additive
structure.
\end{defn}

\subsection{Tensor embedding}

For $A\in\in\Cstar$, denote by $\mathfrak{O}_{A}\in\in\EFC$ the
endofunctor of the left (maximal) tensor multiplication as in Definition~\ref{def:tensoring}:
\[
\mathfrak{O}_{A}\coloneqq A\otimes\mmph\colon(B\xrightarrow{\varphi}B')\mapsto(A\otimes B\xrightarrow{A\otimes\varphi}A\otimes B').
\]

\begin{prop}\label{prop:O-decent}
For every $C^{*}$-algebra $A$, the endofunctor
$\mathfrak{O}_{A}$ is decent and
exact\footnote{A functor is called \emph{exact} if it preserves short exact sequences.}.
\end{prop}

\begin{proof}
It is clear that $\mathfrak{O}_{A}$ preserves terminal objects. Furthermore,
it is known that $\mathfrak{O}_{A}$ preserves all pullbacks (see
the remark at the top of page~264 of~\cite{Ped_pullback99}). Thus,
$\mathfrak{O}_{A}$ is decent. For the proof of exactness, see e.g.~\cite[Lemma~4.1]{GHT}.
\end{proof}

\begin{defn}\label{def:bimon-C-DEFC}
Introduce the bimonoidal functor
\[
(\mathfrak{O},\mathfrak{O}^{2},\mathfrak{O}^{0},\mathfrak{O}^{2}_{\oplus},\mathfrak{O}^{0}_{\oplus})\colon\Cstar\to\DEFC
\]
defined as in Proposition~\ref{prop:bimon-O}. We call it the \emph{tensor
embedding.}
\end{defn}

In order to justify the word ``embedding'', we are going to show
that $\mathfrak{O}$ is fully faithful. To do that we shall need a
lemma which uses the distinguishing feature of $C^{*}$-algebras---the
Gelfand isomorphism.

\begin{lem}\label{lem:well-ppointed}
Let $\varphi,\psi\colon A\otimes B\to C$
be $*$-homomorphisms, such that for any $\mathsf{X}\subset\mathbb{R}$
and any $*$-homomorphism $\zeta\colon C_{0}(\mathsf{X})\to A$ we
have the equality
\[
\varphi\circ(\zeta\otimes\id_{B})=\psi\circ(\zeta\otimes\id_{B})\colon C_{0}(\mathsf{X})\otimes B\to C.
\]
Then $\varphi=\psi$.
\end{lem}

\begin{proof}
Take an arbitrary $b\in B$ and a self-adjoint $a\in A$. Denote by
$C^{*}(a)$ the commutative $C^{*}$-algebra generated by $\{a\}$,
and by $z$ the identity function on the spectrum $\spec(a)\subset\mathbb{R}$.
By Lemma~\ref{lem:spectrum} there is a $*$-isomorphism $\zeta\colon C_{0}(\spec(a)\setminus\{0\})\xrightarrow{\cong}C^{*}(a)$
which maps $z$ to $a$. The statement of the lemma follows from the
equality
\[
\varphi(a\otimes b)=\varphi((\zeta\otimes\id)(z\otimes b))=\psi((\zeta\otimes\id)(z\otimes b))=\psi(a\otimes b).
\]
\end{proof}

\begin{thm}\label{thm:O-ffaith}
The functor $\mathfrak{O}\colon\Cstar\to\EFC$
is fully faithful.
\end{thm}

\begin{proof}
Faithfulness is clear, and we only need show
fullness.
Let $\Phi\colon\mathfrak{O}_{A}\Rightarrow\mathfrak{O}_{A'}$ be a
natural transformation. Set $\varphi\coloneqq\rho^{\otimes}_{A'}\circ\Phi\mathbb{C}\circ(\rho^{\otimes}_{A})^{-1}$.
We need to prove that $\Phi B=\mathfrak{O}_{\varphi}B$. To this end,
take some subspace $\mathsf{X}\subset\mathbb{R}$, denote by $e_{x}\colon C_{0}(\mathsf{X})\to\mathbb{C}$
the $*$-homomorphism of evaluation at $x\in\mathsf{X}$, and consider
the diagram
\[
\input{diagrams/d275.tikz}
\]
in which the outermost rectangle commutes by the naturality of $\Phi$;
the right subdiagram commutes by the definition of $\varphi$; the
middle subdiagram commutes by a straightforward computation. Thus,
we have the equality
\[
(\id_{A'}\otimes e_{x})\circ(\varphi\otimes\id_{C_{0}(\mathsf{X})})=(\id_{A'}\otimes e_{x})\circ(\Phi C_{0}(\mathsf{X}))
\]
for all $x\in\mathsf{X}$, and hence $\varphi\otimes\id_{C_{0}(\mathsf{X})}=\Phi C_{0}(\mathsf{X})$
by Corollary~\ref{cor:teval} from the next subsection.

Let now $\zeta\colon C_{0}(\mathsf{X})\to B$ be an arbitrary $*$-homomorphism.
The outermost diagram in
\[
\input{diagrams/d276.tikz}
\]
commutes by the naturality of $\Phi$, the left subdiagram commutes
by what we have just proved; the middle subdiagram commutes by a direct
computation. Thus, we have the equality
\[
\Phi B\circ(\id_{A}\otimes\zeta)=(\varphi\otimes\id_{B})\circ(\id_{A}\otimes\zeta)
\]
which by Lemma~\ref{lem:well-ppointed} (and using that $\otimes$
is symmetric) implies that $\Phi B=\varphi\otimes\id_{B}=\mathfrak{O}_{\varphi}B$
for all $B\in\in\Cstar$, from which the statement follows.
\end{proof}

\subsection{\label{subsec:stacont}Stabilization and continuous functions}

\begin{defn}\label{def:Cb}
Let $C_{b}(\mathsf{X},B)$ be the $C^{*}$-algebras
of continuous bounded functions on a locally compact Hausdorff space
$\mathsf{X}$ taking values in a $C^{*}$-algebra
$B$, and let $C_{0}(\mathsf{X},B)$ be the ideal in $C_{b}(\mathsf{X},B)$
comprised of those functions that vanish at infinity. These $C^{*}$-algebras
give rise to the endofunctors $\mathfrak{C}^{b}_{\mathsf{X}}\coloneqq C_{b}(\mathsf{X},\mmph)$
and $\C_{\mathsf{X}}\coloneqq C_{0}(\mathsf{X},\mmph)$.

Every continuous (resp. proper\footnote{A map is called \emph{proper} if the inverse image of a precompact
set is precompact.} continuous) map $f\colon\mathsf{X}\to\mathsf{Y}$ induces the natural
transformation $\mathfrak{C}^{b}_{f}\colon\mathfrak{C}^{b}_{\mathsf{Y}}\Rightarrow\mathfrak{C}^{b}_{\mathsf{X}}$
(resp. $\C_{f}\colon\C_{\mathsf{Y}}\Rightarrow\C_{\mathsf{X}}$) given
by the formula $g\mapsto g\circ f$.
\end{defn}

\begin{defn}\label{def:bbK}
Let $\mathcal{L}_{B}(\mathcal{H}_{B})$ be the $C^{*}$-algebra
of adjointable operators in the standard Hilbert $B$-module $\mathcal{H}_{B}$~\cite{lance1995}.
Denote by $\fK B$ the closure of the $*$-subalgebra in $\mathcal{L}_{B}(\mathcal{H}_{B})$
consisting of $\mathbb{N}$-by-$\mathbb{N}$ matrices with entries
in $B$, in which all but finitely many entries are zero. In other
words, $\fK B$ is just the $C^{*}$-algebra $\mathcal{K}_{B}(\mathcal{H}_{B})$
of compact operators in the sense of Hilbert $C^{*}$-modules. It
is easy to show that $\fK$ is functorial, and we call $\fK$ the
\emph{stabilization} functor. Denote by $\iota_{00}\colon\Id\Rightarrow\fK$
the obvious corner embedding.

Similarly, we define the endofunctor $\fM_{n}$ which sends a $C^{*}$-algebra
$B$ to the $C^{*}$-algebra of $n$-by-$n$ matrices with entries
in $B$, and the endofunctor $\mathfrak{M}^{u}_{\mathbb{Z}}$ which
sends a $C^{*}$-algebra $B$ to the $C^{*}$-algebra generated by
all $\mathbb{Z}$-by-$\mathbb{Z}$ matrices $(b_{i,j})$ with $b_{i,j}\in B$
such that $\sup_{b_{i,j}\neq0}|i-j|<\infty$. We call $\mathfrak{M}^{u}_{\mathbb{Z}}$
the \emph{uniform Roe functor} for $\mathbb{Z}$; see~\cite{makeev_asadj}.
\end{defn}

\begin{lem}\label{lem:decent-ex}
The endofunctors $\fM_{n}$, $\fK$, $\mathfrak{M}^{u}_{\mathbb{Z}}$,
$\mathfrak{C}^{b}_{\mathsf{X}}$, and $\C_{\mathsf{X}}$ are decent.
\end{lem}

\begin{proof}
Let us prove the statement for $\C_{\mathsf{X}}$.
Bearing in mind the explicit formula~(\ref{eq:explicit-pullback})
for pullbacks in~$\Cstar$, we conclude that the canonical morphism
$\lc$ is given by the formula
\[
\lc\colon\C_{\mathsf{X}}(B_{1}\underset{B}{\oplus}B_{2})\to\C_{\mathsf{X}}B_{1}\underset{\C_{\mathsf{X}}B}{\oplus}\C_{\mathsf{X}}B_{2}:(x\mapsto(f_{1}(x),f_{2}(x)))\mapsto(x\mapsto f_{1}(x),x\mapsto f_{2}(x))
\]
and is clearly injective and surjective. Furthermore, $\C_{\mathsf{X}}$
obviously preserves zero objects, and hence, it is decent.

The proof carries over verbatim for $\mathfrak{C}^{b}_{\mathsf{X}}$,
$\fM_{n}$, $\mathfrak{M}^{u}_{\mathbb{Z}}$, and $\fK$.
\end{proof}

Note that $\oK=\fK\mathbb{C}$ and $\oM_{n}=\fM_{n}\mathbb{C}$ where
$\oK$ and $\oM_{n}$ are the $C^{*}$-algebras discussed in Section~\ref{subsec:rig-comp}.
The following fact is well known, and we state it without proof.

\begin{lem}\label{lem:39}
Let $\mathsf{X}$ be a locally compact Hausdorff topological
space. There exist natural isomorphisms $\mathfrak{O}_{C_{0}(\mathsf{X})}\xRightarrow{\cong}\C_{\mathsf{X}}$,
$\mathfrak{O}_{\oM_{n}}\xRightarrow{\cong}\fM_{n}$, and $\mathfrak{O}_{\oK}\xRightarrow{\cong}\fK$
given by the formulas
\begin{alignat}{1}
C_{0}(\mathsf{X})\otimes B\to C_{0}(\mathsf{X},B)\colon & f\otimes b\mapsto(x\mapsto f(x)b),\label{eq:al}\\
\oM_{n}\otimes B\to\fM_{n}B\colon & \epsilon_{i,j}\otimes b\mapsto\epsilon_{i,j}b,\nonumber \\
\oK\otimes B\to\fK B\colon & \epsilon_{i,j}\otimes b\mapsto\epsilon_{i,j}b.\nonumber
\end{alignat}
\end{lem}

\begin{lem}\label{lem:alt-homot-Cstar}
Two $*$-homomorphisms $\varphi_{0},\varphi_{1}\colon A\to B$
are homotopic if and only if there is a $*$-homomorphism $\Phi\colon A\to B\otimes C[0,1]$
such that the following diagram commutes for $j=0,1$:
\[
\input{diagrams/d366.tikz}
\]
\end{lem}

\begin{proof}
The statement is straightforward to prove using the isomorphism (\ref{eq:al}).
\end{proof}

\subsection{\label{subsec:tt}Tensor-type endofunctors}

In this section we introduce a subcategory of decent endofunctors
that behave like tensoring, and which we call \emph{tensor-type} endofunctors.
This category yields a strictification theorem
for
the category of $C^{*}$-algebras and $*$-homomorphisms. The construction
we use is analogous to that of~\cite{street-joyal-93}, also known
as the category of \emph{left-module endofunctors}~\cite{EGNO15,becerra2024strictification},
but adapted to leverage some special properties of $C^{*}$-algebras.
\begin{defn}\label{def:tt}
We say that $T\in\in\DEFC$ is \emph{tensor-type} if
there is an isomorphism $\omega_{T}\colon T\Rightarrow\mathfrak{O}_{T\mathbb{C}}$
such that the following diagram commutes:
\begin{equation}
\input{diagrams/d581.tikz}\label{eq:tt}
\end{equation}
Denote by $\DEFC_{\ttt}$ the full subcategory in $\DEFC$ generated
by tensor-type endofunctors.
\end{defn}

\begin{lem}
Let $\mathsf{X}$ be a locally compact Hausdorff topological space,
and $A$ a $C^{*}$-algebra. The endofunctors $\bzero$, $\Id$, $\fM_{n}$,
$\fK$, $\C_{\mathsf{X}}$, $\mathfrak{O}_{A}$ are tensor-type.
\end{lem}

\begin{proof}
Define $\omega_{\Id}$, $\omega_{\mathfrak{O}_{A}}$, $\omega_{\bzero}$
as
\begin{alignat*}{3}
\omega_{\Id} & \colon\Id\Rightarrow\mathfrak{O}_{\mathbb{C}}, & \qquad &\ 
& \omega_{\Id}B\colon & B\xrightarrow{(\lambda^{\otimes})^{-1}}\mathbb{C}\otimes B;\\
\omega_{\mathfrak{O}_{A}} & \colon\mathfrak{O}_{A}\Rightarrow\mathfrak{O}_{A\otimes\mathbb{C}}, & \qquad &\ 
& \omega_{\mathfrak{O}_{A}}B\colon & A\otimes B\xrightarrow{(\rho^{\otimes})^{-1}\otimes B}(A\otimes\mathbb{C})\otimes B;\\
\omega_{\bzero} & \colon\bzero\Rightarrow\mathfrak{O}_{\bbzero}, & \qquad &\ 
& \omega_{\Id}B\colon & \bbzero\xrightarrow{(\lambda^{\bullet})^{-1}}\bbzero\otimes B,
\end{alignat*}
and define $\omega_{\C_{\mathsf{X}}}$, $\mathfrak{O}_{\fM_{n}}$,
$\mathfrak{O}_{\fK}$ as the inverse of the natural isomorphisms from
Lemma~\ref{lem:39}. We leave it to the reader to check that (\ref{eq:tt})
holds for these natural transformations.
\end{proof}

\begin{rem}
One
can show that $\mathfrak{C}^{b}_{\mathsf{X}}$
is not tensor-type for any non-compact locally compact metric space~$\mathsf{X}$.
To do this, assume for contradiction that $\mathfrak{C}^{b}_{\mathsf{X}}$
is tensor-type. We leave it to the reader to show that in this case
the isomorphism $\omega^{-1}_{\mathfrak{C}^{b}_{\mathsf{X}}}$ is
given by the formula
\begin{equation}
\omega^{-1}_{\mathfrak{C}^{b}_{\mathsf{X}}}B\colon\mathfrak{C}^{b}_{\mathsf{X}}\mathbb{C}\otimes B\to\mathfrak{C}^{b}_{\mathsf{X}}B\colon f\otimes b\mapsto(x\mapsto f(x)b).\label{eq:rem-nontt}
\end{equation}
(Hint: Using Definition~\ref{def:tt}, naturality of $\omega^{-1}_{\mathfrak{C}^{b}_{\mathsf{X}}}$,
and Lemma~\ref{lem:spectrum}, successively prove the statement for
$B=\mathbb{C}$; $B=C_{0}(\mathsf{Y})$ where $\mathsf{Y}\subset\mathbb{R}$;
$B$ generated by a self-adjoint element; and arbitrary $B$.) Now,
using the properties of $\mathsf{X}$, we can find a sequence of balls\footnote{$B_{r}(x)$ denotes the ball of radius $r$ centered at $x$.}
$\left\{ B_{r_{j}}(x_{j})\right\} ^{\infty}_{j=1}$ such that every
$y\in\mathsf{X}$ has a neighborhood intersecting at most one ball
from this family, which allows us to define the function $g\in\mathfrak{C}^{b}_{\mathsf{X}}\oK$
by the formula
\[
g\colon\mathsf{X}\to\oK\colon x\mapsto\begin{cases}
\left(1-\frac{\dist(x,x_{j})}{r_{j}}\right)\epsilon_{j,j}, & \textrm{if }x\in B_{r_{j}}(x_{j});\\
0, & \textrm{otherwise},
\end{cases}
\]
where $\epsilon_{i,j}\in\oK$ stands for the $(i,j)$th matrix unit.
Using formula~(\ref{eq:rem-nontt}), we conclude that $\norm{g-h}=1$
for every $h$ in the image of $\omega^{-1}_{\mathfrak{C}^{b}_{\mathsf{X}}}\oK$,
and hence $\omega^{-1}_{\mathfrak{C}^{b}_{\mathsf{X}}}\oK$ is not
an isomorphism.
\end{rem}

\begin{lem}\label{lem:26}
The following diagram commutes:
\[
\input{diagrams/d438.tikz}
\]
\end{lem}

\begin{proof}
Straightforward.
\end{proof}

\begin{cor}\label{cor:teval}
Let $A$ be a $C^{*}$-algebra, and $\mathsf{X}$
a locally compact Hausdorff space. If $m\in C_{0}(\mathsf{X})\otimes A$
is such that $(e_{x}\otimes\id_{A})(m)=0$ for all $x\in\mathsf{X}$,
then $m=0$.
\end{cor}

\begin{proof}
Evaluating the diagram from Lemma~\ref{lem:26}
at $A$, we conclude that $\ev_{x}A\left(\omega^{-1}_{\C_{\mathsf{X}}}A(m)\right)=0$
 for all $x\in\mathsf{X}$, which implies that $\omega^{-1}_{\C_{\mathsf{X}}}A(m)=0$,
and hence, $m=0$.
\end{proof}

\begin{thm}\label{thm:aoa}
Let $\alpha\in\DEFC_{\ttt}(T,S)$. Then the following
diagram commutes:
\[
\input{diagrams/d379.tikz}
\]

\end{thm}

\begin{proof}
Let $B$ be a $C^{*}$-algebra, $\mathsf{X}\subset\mathbb{R}$, and
$x\in\mathsf{X}$. The top, bottom, left, right, and back faces in
\[
\input{diagrams/d380.tikz}
\]
commute by the naturality of $\omega_{T}$, $\omega_{S}$, $\alpha$,
$\mathfrak{O}_{\alpha\mathbb{C}}$, and $\rho^{\otimes}$, respectively.
Hence, the front face commutes by Corollary~\ref{cor:teval}.

Let $\zeta\colon C_{0}(\mathsf{X})\to B$ be an arbitrary $*$-homomorphism.
The top, bottom, left and right faces in
\[
\input{diagrams/d381.tikz}
\]
commute by naturality, and the front face commutes by what has been
proved above. Hence, the back face commutes by Lemma~\ref{lem:well-ppointed},
and the proof is complete.
\end{proof}

\begin{cor}
There exists at most one natural transformation $\omega_{T}$ which
turns $T\in\in\DEFC$ into a tensor-type endofunctor.
\end{cor}

\begin{proof}
Set $S\coloneqq T$, and $\alpha\coloneqq\id_{T}$ in Theorem~\ref{thm:aoa},
and deduce that $\omega_{T}=\omega_{S}$.
\end{proof}

\begin{lem}\label{lem:def-omegaF}
 $\DEFC_{\ttt}$ is closed under taking compositions
and binary products. Specifically, for $T,S\in\in\DEFC_{\ttt}$ the
natural transformations $\omega_{TS}$ and $\omega_{T\oplus S}$ are
defined as the following composites:
\begin{alignat}{1}
\omega_{TS}\colon & TS\xRightarrow{\omega_{T}\omega_{S}}\mathfrak{O}_{T\mathbb{C}}\mathfrak{O}_{S\mathbb{C}}\xRightarrow{\mathfrak{O}^{2}}\mathfrak{O}_{T\mathbb{C}\otimes S\mathbb{C}}\xRightarrow{\mathfrak{O}_{\omega^{-1}_{T}S\mathbb{C}}}\mathfrak{O}_{TS\mathbb{C}},\label{eq:omega1}\\
\omega_{T\oplus S}\colon & T\oplus S\xRightarrow{\omega_{T}\oplus\omega_{S}}\mathfrak{O}_{T\mathbb{C}}\oplus\mathfrak{O}_{S\mathbb{C}}\xRightarrow{\mathfrak{O}^{2}_{\oplus}}\mathfrak{O}_{T\mathbb{C}\oplus S\mathbb{C}}=\mathfrak{O}_{(T\oplus S)\mathbb{C}}.\label{eq:omega2}
\end{alignat}
\end{lem}

\begin{proof}
The equality $\omega_{TS}\mathbb{C}=(\rho^{\otimes})^{-1}$ follows
from the diagram
\[
\input{diagrams/d382.tikz}
\]
in which the left subdiagram commutes by the naturality of $\rho^{\otimes}$,
and the right triangle commutes by coherence. The equality $\omega_{T\oplus S}\mathbb{C}=(\rho^{\otimes})^{-1}$
follows from the diagram
\[
\input{diagrams/d384.tikz}
\]
in which the rectangles commute by the definitions of tensor-type
endofunctors and $\mathfrak{O}^{2}$, and the lower triangle commutes
by one of the Laplaza axioms (see the right-hand side diagram of~(\ref{eq:d560})
in the Appendix).
\end{proof}

\begin{cor}
$\DEFC_{\ttt}$ is a tight bimonoidal category.
\end{cor}

\begin{proof}
This follows from the fact that $\DEFC$ is a tight bimonoidal category
together with Lemma~\ref{lem:def-omegaF}.
\end{proof}

\begin{lem}\label{lem:aa}
For every $T\in\in\DEFC_{\ttt}$ and $A_{1},A_{2}\in\in\Cstar$,
the following diagram commutes:
\[
\input{diagrams/d415.tikz}
\]
\end{lem}

\begin{proof}
All the inner subdiagrams in
\[
\input{diagrams/d416.tikz}
\]
commute for $j=1,2$ either by definition or by naturality. Hence,
by the universal property of $\oplus$, the outermost rectangle also
commutes.
\end{proof}

\begin{prop}\label{prop:EV}
There
is a fully faithful unitary
bimonoidal functor $(\EV,\EV^{2}_{\otimes},\EV^{2}_{\oplus})\colon\DEFC_{\ttt}\to\Cstar$,
with
\begin{alignat*}{1}
\EV\colon & (T\xRightarrow{\alpha}S)\mapsto(T\mathbb{C}\xRightarrow{\alpha\mathbb{C}}S\mathbb{C}),\\
\EV^{2}_{\otimes}\colon & T\mathbb{C}\otimes S\mathbb{C}\xrightarrow{\omega^{-1}_{T}S\mathbb{C}}TS\mathbb{C},\\
\EV^{2}_{\oplus}\colon & T\mathbb{C}\oplus S\mathbb{C}\xrightarrow{=}(T\oplus S)\mathbb{C}.
\end{alignat*}
\end{prop}

\begin{proof}
Let us first show that $(\EV,\EV^{2}_{\otimes})$ and $(\EV,\EV^{2}_{\oplus})$
are respectively monoidal and symmetric monoidal functors. Indeed,
for $T,S,R\in\in\DEFC_{\ttt}$ the diagrams
\begin{alignat*}{2}
\input{diagrams/d392.tikz} & \qquad\qquad & \input{diagrams/d411.tikz}
\end{alignat*}
commute by the definitions of $\omega_{TS}$ (see the proof of Lemma~\ref{lem:def-omegaF})
and by Lemma~\ref{lem:coher-component}. The diagrams
\begin{alignat*}{2}
\input{diagrams/d394.tikz} & \qquad & \input{diagrams/d395.tikz}\\
\input{diagrams/d413.tikz} &\ 
& \input{diagrams/d412.tikz}
\end{alignat*}
commute since $\omega_{\Id}=(\lambda^{\otimes})^{-1}$, $\omega_{T}\mathbb{C}=(\rho^{\otimes})^{-1}$,
and by Lemma~\ref{lem:coher-component}. It is straightforward to
show that $(\EV,\EV^{2}_{\oplus})$ is symmetric, and we are done.

To finish the proof, we need to check the commutativity of four diagrams
as in~(\ref{eq:bmf1}) and~(\ref{eq:bmf2}). The left-hand side
diagram in~(\ref{eq:bmf2}) takes the form
\[
\input{diagrams/d486.tikz}
\]
where $\tc 1$ commutes by the definition of $\omega_{T\oplus S}$
(see Lemma~\ref{lem:def-omegaF}); $\tc 2$ by the definition of
$\mathfrak{O}^{2}_{\oplus}$ (see Proposition~\ref{prop:bimon-O});
and $\tc 3$ commutes by the definition of $\oplus$ (see Definition~\ref{def:oplus}).
The left-hand side diagram in~(\ref{eq:bmf1}) takes the form
\[
\input{diagrams/d414.tikz}
\]
where $\tc 1$ commutes by the definition of $\oplus$ (see Definition~\ref{def:oplus});
$\tc 2$ by Lemma~\ref{lem:aa} and Remark~\ref{rem:rc-eq}; and
$\tc 3$ by Lemma~\ref{lem:almost-bimon}. Finally, the right-hand
side diagrams as in (\ref{eq:bmf1}) and (\ref{eq:bmf2}) commute
since $\bbzero$ is a terminal object.

It is easy to show using Theorems~\ref{thm:aoa} and~\ref{thm:O-ffaith}
that $\EV$ is fully faithful.
\end{proof}

\begin{thm}\label{thm:strictif-defc}
There
is a bimonoidal
equivalence of the categories $\Cstar$ and $\DEFC_{\ttt}$ witnessed
by the bimonoidal functors
\[
\mathfrak{O}\colon\Cstar\leftrightarrows\DEFC_{\ttt}\noloc\EV.
\]
\end{thm}

\begin{proof}
We claim that the morphisms $\left\{ \omega_{T}\colon T\Rightarrow\mathfrak{O}_{T\mathbb{C}}\right\} _{T\in\in\DEFC_{\ttt}}$
from Definition~\ref{def:tt} are in fact the components of the bimonoidal
isomorphism
\[
\omega\colon1_{\DEFC_{\ttt}}\to\mathfrak{O}\circ\EV.
\]
 Indeed, naturality of $\omega$ follows from Theorem~\ref{thm:aoa}.
We only need show that $\omega$ is a monoidal natural transformation
with respect to both the multiplicative and additive structures of
$\DEFC_{\ttt}$. The first statement follows from the outermost parts
of the following two diagrams
\begin{alignat*}{2}
 & \input{diagrams/d576.tikz} & \qquad\qquad & \input{diagrams/d577.tikz}
\end{alignat*}
where the left-hand side one commutes by Definition~\ref{def:monfun-compos},
Lemma~\ref{lem:def-omegaF}, and Proposition~\ref{prop:EV}, and
the right-hand side one commutes trivially. The second statement is
proved similarly.

Now, we claim that the right unity isomorphisms $\left\{ \rho^{\otimes}_{A}\colon A\otimes\mathbb{C}\to A\right\} _{A\in\in\Cstar}$
are precisely the components of a bimonoidal natural isomorphism
\[
\rho\colon\EV\circ\mathfrak{O}\to1_{\Cstar}.
\]
Indeed, naturality in $A$ is part of the definition of $\rho^{\otimes}_{A}$.
Bimonoidality is left as an exercise.
\end{proof}

\begin{rem}\label{rem:strictif-defc}
We can regard Theorem~\ref{thm:strictif-defc}
as a strictification result for the multiplicative structure of~$\Cstar$.
\end{rem}

\begin{cor}\label{cor:EV-DEFC}
There is a fully faithful strict bimonoidal functor
\[
\overline{\EV}\colon\Br(\DEFC_{\ttt})\to\Br(\Cstar)
\]
defined as in Theorem~\ref{thm:defFbar}.
\end{cor}

\begin{proof}
This follows from Proposition~\ref{prop:EV} and Theorem~\ref{thm:defFbar}.
\end{proof}

Now, we are going to find an explicit formula for the $\overline{\EV}$-preimages
of morphisms in $\Br(\Cstar)$. For every $T\in\in\Br(\DEFC_{\ttt})$,
let us define the natural isomorphism $\Lambda_{T}\colon\u(T)\Rightarrow\mathfrak{O}_{\u(\overline{\EV}T)}$
by the following inductive relations:
\begin{alignat*}{1}
\Lambda_{T}\colon & \u(T)\xRightarrow{\omega_{\u(T)}}\mathfrak{O}_{\u(T)\mathbb{C}}=\mathfrak{O}_{\u\overline{\EV}(T)},\qquad\textrm{if }|T|=1;\\
\Lambda_{TS}\colon & \u(TS)=\u(T)\u(S)\xRightarrow{\Lambda_{T}\Lambda_{S}}\mathfrak{O}_{\u\overline{\EV}(T)}\mathfrak{O}_{\u\overline{\EV}(S)}\xRightarrow{\mathfrak{O}^{2}}\mathfrak{O}_{\u\overline{\EV}(T)\otimes\u\overline{\EV}(T)}=\mathfrak{O}_{\u\overline{\EV}(TS)};\\
\Lambda_{T\oplus S}\colon & \u(T\oplus S)=\u(T)\oplus\u(S)\xRightarrow{\Lambda_{T}\oplus\Lambda_{S}}\mathfrak{O}_{\u\overline{\EV}(T)}\oplus\mathfrak{O}_{\u\overline{\EV}(S)}\xRightarrow{\mathfrak{O}^{2}_{\oplus}}\mathfrak{O}_{\u\overline{\EV}(T)\oplus\u\overline{\EV}(T)}=\mathfrak{O}_{\u\overline{\EV}(T\oplus S)}.
\end{alignat*}

\begin{lem}\label{lem:ba}
The following diagram commutes for all $T\in\in\Br(\DEFC_{\ttt})$:
\[
\input{diagrams/d439.tikz}
\]
where $\Omega_{\EV}$ is defined as in Proposition~\ref{prop:natOmega}.
\end{lem}

\begin{proof}
When $|T|=1$, the statement of the lemma follows from the definition
of $\Lambda$, $\omega$, and $\Omega$. So, it suffices to show that
if the statement holds for both $T$ and $S$, then it also holds
for $T\cdot S$ and $T\oplus S$. The former is proved by the diagram
\[
\input{diagrams/d473.tikz}
\]
in which the upper and lower subdiagrams commute by the definitions
of $\Lambda_{TS}$ and $\omega_{TS}$, respectively; the left rectangle
commutes by the induction hypothesis, the middle one by naturality
of $\mathfrak{O}^{2}$, and the right one by the definition of $\EV$
and $\Omega_{\EV}$. The proof for $T\oplus S$ is similar.
\end{proof}

\begin{prop}\label{prop:bb}
The family $\{\Lambda_{T}\}_{T\in\in\Br(\DEFC_{\ttt})}$
gives rise to the natural isomorphism
\[
\Lambda\colon1_{\Br(\DEFC_{\ttt})}\xrightarrow{\cong}\i\mathfrak{O}_{\u\overline{\EV}(\mmph)}
\]
where $1_{\Br(\DEFC_{\ttt})}\colon\Br(\DEFC_{\ttt})\to\Br(\DEFC_{\ttt})$
stands for the identity functor.
\end{prop}

\begin{proof}
Let $\alpha\in\Br(\DEFC_{\ttt})(T,S)$. Recall that the arrows in
$\Br(\DEFC_{\ttt})$ are the natural transformations between the underlying
endofunctors, and so it suffices to show the commutativity of the
outermost rectangle in
\begin{equation}
\input{diagrams/d441.tikz}\label{eq:d441}
\end{equation}
in which the left and right triangles commute by Lemma~\ref{lem:ba};
the upper subdiagram by Theorem~\ref{thm:aoa}; and the lower subdiagram
by the naturality of $\Omega_{\EV}$ (see Proposition~\ref{prop:natOmega}).
The outermost rectangle commutes since $\Omega_{\EV}$ is an isomorphism.
\end{proof}

\begin{rem}
By abuse of notation, later on we shall regard $\Lambda_{T}$ as the
arrows in $\Br(\DEFC_{\ttt})$, i.e., $\Lambda_{T}\colon T\Rightarrow\i\mathfrak{O}_{\u(T)}$.
\end{rem}

Proposition~\ref{prop:bb} gives us a very convenient way to describe
$\overline{\EV}$-preimages of morphisms in $\Br(\Cstar)$.
\begin{lem}\label{lem:bc}
If $\varphi\in\Br(\Cstar)(\overline{\EV}(T),\overline{\EV}(S))$,
then $\Phi\coloneqq\Lambda^{-1}_{S}\circ\i\mathfrak{O}_{\varphi}\circ\Lambda_{T}\in\Br(\DEFC_{\ttt})(T,S)$
is the $\overline{\EV}$-preimage of $\varphi$.
\end{lem}

\begin{proof}
It follows from diagram~(\ref{eq:d441}) that $\i\mathfrak{O}_{\varphi}=\i\mathfrak{O}_{\overline{\EV}(\Phi)}$,
thus, $\varphi=\overline{\EV}(\Phi)$ since $\mathfrak{O}$ and $\i$
are faithful.
\end{proof}

\begin{defn}\label{def:xitt}
For every $T,S\in\in\DEFC_{\ttt}$, define
$\overline{\xi}_{T,S}\colon T\cdot S\Rightarrow S\cdot T$ as the
$\overline{\EV}$-preimage of the symmetry isomorphism
\[
\xi^{\otimes}_{T\mathbb{C},S\mathbb{C}}\colon T\mathbb{C}\otimes S\mathbb{C}\to S\mathbb{C}\otimes T\mathbb{C}
\]
regarded as a morphism in $\Br(\Cstar)$; and denote by $\xi$ the
underlying natural transformation $\xi\coloneqq\u(\overline{\xi})$.
It follows from Lemma~\ref{lem:bc} that $\xi_{T,S}$ is given by
the composite:
\[
\xi_{T,S}\colon TS\xRightarrow{\omega_{T}\omega_{S}}\mathfrak{O}_{T\mathbb{C}}\mathfrak{O}_{S\mathbb{C}}\xRightarrow{\mathfrak{O}^{2}}\mathfrak{O}_{T\mathbb{C}\otimes S\mathbb{C}}\xRightarrow{\mathfrak{O}_{\xi^{\otimes}}}\mathfrak{O}_{S\mathbb{C}\otimes T\mathbb{C}}\xRightarrow{(\mathfrak{O}^{2})^{-1}}\mathfrak{O}_{S\mathbb{C}}\mathfrak{O}_{T\mathbb{C}}\xRightarrow{\omega^{-1}_{S}\omega^{-1}_{T}}ST.
\]
\end{defn}

\begin{prop}\label{prop:xitt}
The
natural isomorphism $\xi$
defined above turns $\DEFC_{\ttt}$ into a symmetric bimonoidal category.
\end{prop}

\begin{proof}
To show that the diagrams as in Definition~\ref{def:symmoncat} commute
in $\DEFC_{\ttt}$, just note that applying the strict monoidal functor
$\overline{\EV}$ to them gives us exactly the diagrams in $\Cstar$
as in Definition~\ref{def:symmoncat} which commute since $\Cstar$
is a symmetric bimonoidal category. Hence, the initial diagrams also
commute since $\overline{\EV}$ is faithful.
\end{proof}

\begin{rem}\label{rem:xiid}
For
every $T\in\in\DEFC_{\ttt}$,
the equality $\xi_{\Id,T}=\id_{T}$ holds.
\end{rem}

%% file: diagrams/d472.tikz
\begin{tikzcd}[ampersand replacement=\&]
{A} \& {IB} \\
 \& {B,}
\arrow["\Phi", from=1-1, to=1-2]
\arrow["\varphi_{j}"', from=1-1, to=2-2]
\arrow["e_{j}", from=1-2, to=2-2]
\end{tikzcd}

%% file: diagrams/d358.tikz
\begin{tikzcd}[column sep=small, ampersand replacement=\&]  %d358
{\oK\otimes(\oK\oplus\oK)} \&  \&  \&  \& {\oK\otimes\oK} \\
 \& {\oK\otimes(\oM_{2}\otimes\oK)} \\
 \& {\oM_{2}\otimes(\oK\otimes\oK)} \&  \& {\oM_{2}\otimes\oK} \\
{(\oK\otimes\oK)\oplus(\oK\otimes\oK)} \&  \& {\oK\oplus\oK} \&  \& {\oK}
\arrow["\oK\otimes\widehat{\mu}", from=1-1, to=1-5]
\arrow["\oK\otimes\widehat{d}", from=1-1, to=2-2]
\arrow["\delta^{l}"', from=1-1, to=4-1]
\arrow["\widehat{\theta}", from=1-5, to=4-5]
\arrow["\oK\otimes\widehat{\theta}_{2}"', from=2-2, to=1-5]
\arrow["\cong", from=2-2, to=3-2]
\arrow["\oM_{2}\otimes\widehat{\theta}", from=3-2, to=3-4]
\arrow["\widehat{\theta}_{2}", from=3-4, to=4-5]
\arrow["\widehat{d}", from=4-1, to=3-2]
\arrow["\widehat{\theta}\oplus\widehat{\theta}"', from=4-1, to=4-3]
\arrow["\widehat{d}", from=4-3, to=3-4]
\arrow["\widehat{\mu}"', from=4-3, to=4-5]
\end{tikzcd}

%% file: diagrams/d249.tikz
\begin{tikzcd}[ampersand replacement=\&]
{B_{1}\oplus_{B}B_{2}} \& {B_{2}} \\
{B_{1}} \& {B}
\arrow[from=1-1, to=1-2]
\arrow[""', from=1-1, to=2-1]
\arrow["\varphi_{2}", two heads, from=1-2, to=2-2]
\arrow["\varphi_{1}", from=2-1, to=2-2]
\end{tikzcd}

%% file: diagrams/d343.tikz
\begin{tikzcd}[ampersand replacement=\&]
{G(\lim_{i}F_{i})B} \&  \& {G\lim_{i}(F_{i}B)} \\
{\lim_{i}(GF_{i})B} \&  \& {\lim_{i}(GF_{i}B)}
\arrow["G\rc", from=1-1, to=1-3]
\arrow["\cong"', draw=none, from=1-1, to=1-3]
\arrow["\lc B"', from=1-1, to=2-1]
\arrow["\lc", from=1-3, to=2-3]
\arrow["\cong"', draw=none, from=1-3, to=2-3]
\arrow["\rc"', from=2-1, to=2-3]
\arrow["\cong", draw=none, from=2-1, to=2-3]
\end{tikzcd}

%% file: diagrams/d342.tikz
\begin{tikzcd}[ampersand replacement=\&]
{(\lim_{j}F_{j})(\lim_{k}B_{k})} \&  \& {\lim_{j}(F_{j}\lim_{k}B_{k})} \&  \& {\lim_{j}\lim_{k}(F_{j}B_{k})} \\
{\lim_{k}((\lim_{j}F_{j})B_{k})} \&  \& {\lim_{k}((\lim_{j}F_{j}B_{k}))} \&  \& {\lim_{(j,k)}F_{j}B_{k}}
\arrow["\cong"', from=1-1, to=1-3]
\arrow["\rc", draw=none, from=1-1, to=1-3]
\arrow["\lc"', from=1-1, to=2-1]
\arrow["\lim_{j}\lc", from=1-3, to=1-5]
\arrow["\cong"', draw=none, from=1-3, to=1-5]
\arrow["\cong"', from=1-5, to=2-5]
\arrow["\lim_{k}\rc"', from=2-1, to=2-3]
\arrow["\cong", no body, from=2-1, to=2-3]
\arrow["\cong", from=2-3, to=2-5]
\end{tikzcd}

%% file: diagrams/d275.tikz
\begin{tikzcd}[row sep=small, ampersand replacement=\&]  %d275
{A\otimes C_{0}(\mathsf{X})} \&  \&  \&  \& {A\otimes\mathbb{C}} \\
 \&  \&  \& {A} \\
 \&  \&  \& {A'} \\
{A'\otimes C_{0}(\mathsf{X})} \&  \&  \&  \& {A'\otimes\mathbb{C}}
\arrow["\id_{A}\otimes e_{x}=\mathfrak{O}_{A}e_{x}", from=1-1, to=1-5]
\arrow[""{name=0, anchor=center, inner sep=0}, "\Phi C_{0}(\mathsf{X})"', from=1-1, to=4-1]
\arrow[""{name=1, anchor=center, inner sep=0}, "\varphi\otimes\id_{C_{0}(\mathsf{X})}", shift left=3, curve={height=-30pt}, from=1-1, to=4-1]
\arrow["\rho_{A}", from=1-5, to=2-4]
\arrow[""{name=2, anchor=center, inner sep=0}, "\Phi\mathbb{C}", from=1-5, to=4-5]
\arrow[""{name=3, anchor=center, inner sep=0}, "\varphi", from=2-4, to=3-4]
\arrow["\id_{A'}\otimes e_{x}=\mathfrak{O}_{A'}e_{x}"', from=4-1, to=4-5]
\arrow["\rho_{A'}"', from=4-5, to=3-4]
\arrow[""{description, pos=0.7}, draw=none, from=1, to=3]
\arrow[""{description}, draw=none, from=0, to=1]
\arrow[""{description}, draw=none, from=3, to=2]
\end{tikzcd}

%% file: diagrams/d276.tikz
\begin{tikzcd}[ampersand replacement=\&]
{A\otimes C_{0}(\mathsf{X})} \&  \&  \&  \&  \& {A\otimes B} \\
 \\
{A'\otimes C_{0}(\mathsf{X})} \&  \&  \&  \&  \& {A'\otimes B}
\arrow[""{name=0, anchor=center, inner sep=0}, "\id_{A}\otimes\zeta=\mathfrak{O}_{A}\zeta", from=1-1, to=1-6]
\arrow["\Phi C_{0}(\mathsf{X})"', curve={height=24pt}, from=1-1, to=3-1]
\arrow["\varphi\otimes\id_{C_{0}(\mathsf{X})}", curve={height=-24pt}, from=1-1, to=3-1]
\arrow["\Phi B", curve={height=-24pt}, from=1-6, to=3-6]
\arrow["\varphi\otimes\id_{B}"', curve={height=24pt}, from=1-6, to=3-6]
\arrow[""{description}, draw=none, from=1-6, to=3-6]
\arrow[""{description}, draw=none, from=3-1, to=1-1]
\arrow[""{name=1, anchor=center, inner sep=0}, "\id_{A'}\otimes\zeta=\mathfrak{O}_{A'}\zeta"', from=3-1, to=3-6]
\arrow[""{description}, draw=none, from=0, to=1]
\end{tikzcd}

%% file: diagrams/d366.tikz
\begin{tikzcd}[ampersand replacement=\&]
{A} \& {B\otimes C[0,1]} \\
{B} \& {B\otimes\mathbb{C}\mmdd}
\arrow["\Phi", from=1-1, to=1-2]
\arrow["\varphi_{j}"', from=1-1, to=2-1]
\arrow["B\otimes e_{j}", from=1-2, to=2-2]
\arrow["(\rho^{\otimes})^{-1}"', from=2-1, to=2-2]
\end{tikzcd}

%% file: diagrams/d581.tikz
\begin{tikzcd}[ampersand replacement=\&]
{T\mathbb{C}} \&  \& {T\mathbb{C}\otimes\mathbb{C}.}
\arrow["\omega_{T}\mathbb{C}"', shift right, curve={height=6pt}, from=1-1, to=1-3]
\arrow["(\rho_{T\mathbb{C}}^{\otimes})^{-1}"{pos=0.6}, shift left, curve={height=-6pt}, from=1-1, to=1-3]
\end{tikzcd}

%% file: diagrams/d438.tikz
\begin{tikzcd}[ampersand replacement=\&]
{\C_{\mathsf{X}}} \& {\mathfrak{O}_{C_{0}(\mathsf{X})}} \\
{\Id} \& {\mathfrak{O}_{\mathbb{C}}\mmdd}
\arrow["\omega_{\C_{\mathsf{X}}}", Rightarrow, from=1-1, to=1-2]
\arrow["\ev_{x}"', Rightarrow, from=1-1, to=2-1]
\arrow["\mathfrak{O}_{\ev_{x}\mathbb{C}}", Rightarrow, from=1-2, to=2-2]
\arrow["\omega_{\Id}"', Rightarrow, from=2-1, to=2-2]
\end{tikzcd}

%% file: diagrams/d379.tikz
\begin{tikzcd}[ampersand replacement=\&]
{T} \& {\mathfrak{O}_{T\mathbb{C}}} \\
{S} \& {\mathfrak{O}_{S\mathbb{C}}\mmdd}
\arrow["\omega_{T}", Rightarrow, from=1-1, to=1-2]
\arrow["\alpha"', Rightarrow, from=1-1, to=2-1]
\arrow["\mathfrak{O}_{\alpha\mathbb{C}}", Rightarrow, from=1-2, to=2-2]
\arrow["\omega_{S}"', Rightarrow, from=2-1, to=2-2]
\end{tikzcd}

%% file: diagrams/d380.tikz
\begin{tikzcd}[ampersand replacement=\&]
 \& {T\mathbb{C}} \&  \&  \& {\mathfrak{O}_{T\mathbb{C}}\mathbb{C}} \\
{TC_{0}(\mathsf{X})} \&  \&  \& {\mathfrak{O}_{T\mathbb{C}}C_{0}(\mathsf{X})} \\
 \& {S\mathbb{C}} \&  \&  \& {\mathfrak{O}_{S\mathbb{C}}\mathbb{C}} \\
{SC_{0}(\mathsf{X})} \&  \&  \& {\mathfrak{O}_{S\mathbb{C}}C_{0}(\mathsf{X})}
\arrow["\omega_{T}\mathbb{C}", from=1-2, to=1-5]
\arrow["\alpha\mathbb{C}"'{pos=0.7}, from=1-2, to=3-2]
\arrow["\mathfrak{O}_{\alpha\mathbb{C}}\mathbb{C}"{pos=0.4}, from=1-5, to=3-5]
\arrow["Te_{x}"{pos=0.4}, from=2-1, to=1-2]
\arrow["\alpha C_{0}(\mathsf{X})"', from=2-1, to=4-1]
\arrow["\mathfrak{O}_{T\mathbb{C}}e_{x}"{pos=0.4}, from=2-4, to=1-5]
\arrow["\omega_{S}\mathbb{C}"'{pos=0.3}, from=3-2, to=3-5]
\arrow["Se_{x}", from=4-1, to=3-2]
\arrow["\omega_{S}C_{0}(\mathsf{X})"', from=4-1, to=4-4]
\arrow["\mathfrak{O}_{S\mathbb{C}}e_{x}"', from=4-4, to=3-5]
\arrow["\omega_{T}C_{0}(\mathsf{X})"{pos=0.7}, crossing over, from=2-1, to=2-4]
\arrow["\mathfrak{O}_{\alpha\mathbb{C}}C_{0}(\mathsf{X})"{pos=0.3}, crossing over, from=2-4, to=4-4]
\end{tikzcd}

%% file: diagrams/d381.tikz
\begin{tikzcd}[ampersand replacement=\&]
 \& {TB} \&  \&  \& {\mathfrak{O}_{T\mathbb{C}}B} \\
{TC_{0}(\mathsf{X})} \&  \&  \& {\mathfrak{O}_{T\mathbb{C}}C_{0}(\mathsf{X})} \\
 \& {SB} \&  \&  \& {\mathfrak{O}_{S\mathbb{C}}B} \\
{SC_{0}(\mathsf{X})} \&  \&  \& {\mathfrak{O}_{S\mathbb{C}}C_{0}(\mathsf{X})}
\arrow["\alpha B"'{pos=0.7}, from=1-2, to=3-2]
\arrow["\omega_{T}^{-1}B"', from=1-5, to=1-2]
\arrow["\mathfrak{O}_{\alpha\mathbb{C}}B", from=1-5, to=3-5]
\arrow["T\zeta", from=2-1, to=1-2]
\arrow["\alpha C_{0}(\mathsf{X})"', from=2-1, to=4-1]
\arrow["\mathfrak{O}_{T\mathbb{C}}\zeta", from=2-4, to=1-5]
\arrow["\omega_{S}^{-1}B"{pos=0.6}, from=3-5, to=3-2]
\arrow["S\zeta", from=4-1, to=3-2]
\arrow["\mathfrak{O}_{S\mathbb{C}}\zeta"', from=4-4, to=3-5]
\arrow["\omega_{S}^{-1}C_{0}(\mathsf{X})", from=4-4, to=4-1]
\arrow["\omega_{T}^{-1}C_{0}(\mathsf{X})"'{pos=0.4}, crossing over, from=2-4, to=2-1]
\arrow["\mathfrak{O}_{\alpha\mathbb{C}}C_{0}(\mathsf{X})"{pos=0.3}, crossing over, from=2-4, to=4-4]
\end{tikzcd}

%% file: diagrams/d382.tikz
\begin{tikzcd}[ampersand replacement=\&]
{TS\mathbb{C}} \&  \& {T\mathbb{C}\otimes S\mathbb{C}} \&  \& {T\mathbb{C}\otimes(S\mathbb{C}\otimes\mathbb{C})} \\
{TS\mathbb{C}\otimes\mathbb{C}} \&  \&  \& {(T\mathbb{C}\otimes S\mathbb{C})\otimes\mathbb{C}}
\arrow["\omega_{T}S\mathbb{C}", from=1-1, to=1-3]
\arrow["\omega_{T}\omega_{S}\mathbb{C}", sqarn=1.5em, from=1-1, to=1-5]
\arrow["\rho_{\otimes}^{-1}"', from=1-1, to=2-1]
\arrow["T\mathbb{C}\otimes\omega_{S}\mathbb{C}=T\mathbb{C}\otimes\rho_{\otimes}^{-1}", from=1-3, to=1-5]
\arrow["\rho_{\otimes}^{-1}"'{pos=0.3}, from=1-3, to=2-4]
\arrow["\alpha_{\otimes}^{-1}=\mathfrak{O}^{2}\mathbb{C}"{pos=0.3}, from=1-5, to=2-4]
\arrow["\omega_{T}^{-1}S\mathbb{C}\otimes\mathbb{C}", from=2-4, to=2-1]
\end{tikzcd}

%% file: diagrams/d384.tikz
\begin{tikzcd}[ampersand replacement=\&]
{(T\oplus S)\mathbb{C}} \&  \&  \& {(\mathfrak{O}_{T\mathbb{C}}\oplus\mathfrak{O}_{S\mathbb{C}})\mathbb{C}} \& {\mathfrak{O}_{T\mathbb{C}\oplus S\mathbb{C}}\mathbb{C}} \\
{T\mathbb{C}\oplus S\mathbb{C}} \&  \&  \& {(T\mathbb{C}\otimes\mathbb{C})\oplus(S\mathbb{C}\otimes\mathbb{C})} \& {(T\mathbb{C}\oplus S\mathbb{C})\otimes\mathbb{C}}
\arrow["(\omega_{T}\oplus\omega_{S})\mathbb{C}", from=1-1, to=1-4]
\arrow[equals, from=1-1, to=2-1]
\arrow["\mathfrak{O}_{\oplus}^{2}\mathbb{C}", from=1-4, to=1-5]
\arrow[equals, from=1-4, to=2-4]
\arrow[equals, from=1-5, to=2-5]
\arrow["(\rho^{\otimes})^{-1}\oplus(\rho^{\otimes})^{-1}", from=2-1, to=2-4]
\arrow["(\rho^{\otimes})^{-1}"', sqars=1.5em, from=2-1, to=2-5]
\arrow["(\delta^{r})^{-1}", from=2-4, to=2-5]
\end{tikzcd}

%% file: diagrams/d415.tikz
\begin{tikzcd}[ampersand replacement=\&]
{{T\left(A_{1}\oplus A_{2}\right)}} \&  \& {{TA_{1}\oplus TA_{2}}} \\
{{T\mathbb{C}\otimes\left(A_{1}\oplus A_{2}\right)}} \&  \& {{\left(T\mathbb{C}\otimes A_{1}\right)\oplus\left(T\mathbb{C}\otimes A_{2}\right)\mmdd}}
\arrow["\lc", from=1-1, to=1-3]
\arrow["{\omega_{T}(A_{1}\oplus A_{2})}"', from=1-1, to=2-1]
\arrow["{\omega_{T}A_{1}\oplus\omega_{T}A_{2}}", from=1-3, to=2-3]
\arrow["{\delta^{l}}"', from=2-1, to=2-3]
\end{tikzcd}

%% file: diagrams/d416.tikz
\begin{tikzcd}[ampersand replacement=\&]
{{T\left(A_{1}\oplus A_{2}\right)}} \&  \& {{TA_{1}\oplus TA_{2}}} \\
 \& {{TA_{j}}} \\
 \& {{T\mathbb{C}\otimes A_{j}}} \\
{{T\mathbb{C}\otimes\left(A_{1}\oplus A_{2}\right)}} \&  \& {{\left(T\mathbb{C}\otimes A_{1}\right)\oplus\left(T\mathbb{C}\otimes A_{2}\right)}}
\arrow["\lc", from=1-1, to=1-3]
\arrow["{T\pr_{j}}"', from=1-1, to=2-2]
\arrow["{\omega_{T}(A_{1}\oplus A_{2})}"', from=1-1, to=4-1]
\arrow["{\pr_{j}}", from=1-3, to=2-2]
\arrow["{\omega_{T}A_{1}\oplus\omega_{T}A_{2}}", from=1-3, to=4-3]
\arrow["{\omega_{T}A_{j}}", from=2-2, to=3-2]
\arrow["{T\mathbb{C}\otimes\pr_{j}}", from=4-1, to=3-2]
\arrow["{\delta^{l}}"', from=4-1, to=4-3]
\arrow["{\pr_{j}}"', from=4-3, to=3-2]
\end{tikzcd}

%% file: diagrams/d392.tikz
\begin{tikzcd}[ampersand replacement=\&]
{\left(T\mathbb{C}\otimes S\mathbb{C}\right)\otimes R\mathbb{C}} \& {T\mathbb{C}\otimes\left(S\mathbb{C}\otimes R\mathbb{C}\right)} \\
{\left(TS\mathbb{C}\right)\otimes R\mathbb{C}} \& {T\mathbb{C}\otimes SR\mathbb{C}} \\
{TSR\mathbb{C}} \& {TSR\mathbb{C}\mmdc}
\arrow["\alpha^{\otimes}", from=1-1, to=1-2]
\arrow["\omega_{T}^{-1}S\mathbb{C}\otimes R\mathbb{C}"', from=1-1, to=2-1]
\arrow["T\mathbb{C}\otimes\omega_{S}^{-1}R\mathbb{C}", from=1-2, to=2-2]
\arrow["\omega_{TS}^{-1}R\mathbb{C}"', from=2-1, to=3-1]
\arrow["\omega_{T}^{-1}SR\mathbb{C}", from=2-2, to=3-2]
\arrow[equals, from=3-1, to=3-2]
\end{tikzcd}

%% file: diagrams/d411.tikz
\begin{tikzcd}[ampersand replacement=\&]
{(T\mathbb{C}\oplus S\mathbb{C})\oplus R\mathbb{C}} \& {T\mathbb{C}\oplus(S\mathbb{C}\oplus R\mathbb{C})} \\
{(T\oplus S)\mathbb{C}\oplus R\mathbb{C}} \& {T\mathbb{C}\oplus(S\oplus R)\mathbb{C}} \\
{((T\oplus S)\oplus R)\mathbb{C}} \& {(T\oplus(S\oplus R))\mathbb{C}\mmdd}
\arrow["\alpha^{\oplus}", from=1-1, to=1-2]
\arrow[equals, from=1-1, to=2-1]
\arrow[equals, from=1-2, to=2-2]
\arrow[equals, from=2-1, to=3-1]
\arrow[equals, from=2-2, to=3-2]
\arrow["\alpha^{\oplus}\mathbb{C}", from=3-1, to=3-2]
\end{tikzcd}

%% file: diagrams/d394.tikz
\begin{tikzcd}[ampersand replacement=\&]
{\mathbb{C}\otimes T\mathbb{C}} \&  \& {T\mathbb{C}} \\
{\Id\mathbb{C}\otimes T\mathbb{C}} \&  \& {\Id T\mathbb{C}\mmdc}
\arrow["\lambda^{\otimes}", from=1-1, to=1-3]
\arrow[equals, from=1-1, to=2-1]
\arrow[equals, from=1-3, to=2-3]
\arrow["\omega_{\Id}^{-1}T\mathbb{C}", from=2-1, to=2-3]
\end{tikzcd}

%% file: diagrams/d395.tikz
\begin{tikzcd}[ampersand replacement=\&]
{\bbzero\oplus T\mathbb{C}} \&  \& {T\mathbb{C}} \\
{\bzero\mathbb{C}\oplus T\mathbb{C}} \&  \& {(\bzero\oplus T)\mathbb{C}\mmdc}
\arrow["\lambda^{\oplus}", from=1-1, to=1-3]
\arrow[equals, from=1-1, to=2-1]
\arrow[equals, from=2-1, to=2-3]
\arrow["\lambda^{\oplus}\mathbb{C}"', from=2-3, to=1-3]
\end{tikzcd}

%% file: diagrams/d413.tikz
\begin{tikzcd}[ampersand replacement=\&]
{T\mathbb{C}\otimes\mathbb{C}} \&  \& {T\mathbb{C}} \\
{T\mathbb{C}\otimes\Id\mathbb{C}} \&  \& {T\Id\mathbb{C}\mmdc}
\arrow["\rho^{\otimes}", from=1-1, to=1-3]
\arrow[equals, from=1-1, to=2-1]
\arrow[equals, from=1-3, to=2-3]
\arrow["\omega_{T}^{-1}\Id\mathbb{C}", from=2-1, to=2-3]
\end{tikzcd}

%% file: diagrams/d412.tikz
\begin{tikzcd}[ampersand replacement=\&]
{T\mathbb{C}\oplus\bbzero} \&  \& {T\mathbb{C}} \\
{T\mathbb{C}\oplus\bzero\mathbb{C}} \&  \& {(T\oplus\bzero)\mathbb{C}}
\arrow["\rho^{\oplus}", from=1-1, to=1-3]
\arrow[equals, from=1-1, to=2-1]
\arrow[equals, from=2-1, to=2-3]
\arrow["\rho^{\oplus}\mathbb{C}"', from=2-3, to=1-3]
\end{tikzcd}

%% file: diagrams/d486.tikz
\begin{tikzcd}[ampersand replacement=\&]
{\left(T\mathbb{C}\oplus S\mathbb{C}\right)\otimes R\mathbb{C}} \&  \&  \&  \& {\left(T\mathbb{C}\otimes R\mathbb{C}\right)\oplus\left(S\mathbb{C}\otimes R\mathbb{C}\right)} \\
{\left(T\oplus S\right)\mathbb{C}\otimes R\mathbb{C}} \&  \& {\mathfrak{O}_{T\mathbb{C}\oplus S\mathbb{C}}R\mathbb{C}} \&  \& {TR\mathbb{C}\oplus SR\mathbb{C}} \\
{\left(T\oplus S\right)R\mathbb{C}} \&  \& {(\mathfrak{O}_{T\mathbb{C}}\oplus\mathfrak{O}_{S\mathbb{C}})R\mathbb{C}} \&  \& {\left(TR\oplus SR\right)\mathbb{C}}
\arrow["\delta^{r}", from=1-1, to=1-5]
\arrow[equals, from=1-1, to=2-1]
\arrow["\omega_{T}^{-1}R\mathbb{C}\oplus\omega_{S}^{-1}R\mathbb{C}", from=1-5, to=2-5]
\arrow[equals, from=2-1, to=2-3]
\arrow[""{name=0, anchor=center, inner sep=0}, "\omega_{T\oplus S}^{-1}"', from=2-1, to=3-1]
\arrow["\tc 2"{description}, draw=none, from=2-3, to=1-5]
\arrow[""{name=1, anchor=center, inner sep=0}, "\mathfrak{O}_{\oplus}^{2}R\mathbb{C}"', from=2-3, to=3-3]
\arrow[sqars=1.0em, equals, from=3-1, to=3-5]
\arrow[equals, from=3-3, to=1-5]
\arrow["(\omega_{T}^{-1}\oplus\omega_{S}^{-1})R\mathbb{C}"', from=3-3, to=3-1]
\arrow["\tc 3"{description}, shift left=5, draw=none, from=3-3, to=3-5]
\arrow[equals, from=3-5, to=2-5]
\arrow["\tc 1"{description}, shift left=2, draw=none, from=0, to=1]
\end{tikzcd}

%% file: diagrams/d414.tikz
\begin{tikzcd}[ampersand replacement=\&]
{R\mathbb{C}\otimes\left(T\mathbb{C}\oplus S\mathbb{C}\right)} \&  \& {\left(R\mathbb{C}\otimes T\mathbb{C}\right)\oplus\left(R\mathbb{C}\otimes S\mathbb{C}\right)} \\
{R\mathbb{C}\otimes\left(T\oplus S\right)\mathbb{C}} \& {R\left(T\mathbb{C}\oplus S\mathbb{C}\right)} \& {RT\mathbb{C}\oplus RS\mathbb{C}} \\
{R\left(T\oplus S\right)\mathbb{C}} \&  \& {\left(RT\oplus RS\right)\mathbb{C}}
\arrow[""{name=0, anchor=center, inner sep=0}, "\delta^{l}", from=1-1, to=1-3]
\arrow[equals, from=1-1, to=2-1]
\arrow["\omega_{R}^{-1}(T\mathbb{C}\oplus S\mathbb{C})", from=1-1, to=2-2]
\arrow["\omega_{R}^{-1}TC\oplus\omega_{R}^{-1}S\mathbb{C}", from=1-3, to=2-3]
\arrow["\tc 1"{description}, draw=none, from=2-1, to=2-2]
\arrow["\omega_{R}^{-1}(T\oplus S)\mathbb{C}"', from=2-1, to=3-1]
\arrow[""{name=1, anchor=center, inner sep=0}, "\lc"', from=2-2, to=2-3]
\arrow[equals, from=2-2, to=3-1]
\arrow[""{name=2, anchor=center, inner sep=0}, "\lc\mathbb{C}"', from=3-1, to=3-3]
\arrow[equals, from=3-3, to=2-3]
\arrow["\tc 2"{description}, draw=none, from=0, to=1]
\arrow["\tc 3"{description}, draw=none, from=1, to=2]
\end{tikzcd}

%% file: diagrams/d576.tikz
\begin{tikzcd}[ampersand replacement=\&]
{TS} \&  \& {\mathfrak{O}_{S\mathbb{C}}\mathfrak{O}_{S\mathbb{C}}} \\
 \& {\mathfrak{O}_{T\mathbb{C}\otimes S\mathbb{C}}} \\
{TS} \&  \& {\mathfrak{O}_{TS\mathbb{C}},}
\arrow["\omega_{T}\omega_{S}", Rightarrow, from=1-1, to=1-3]
\arrow[equals, from=1-1, to=3-1]
\arrow["\mathfrak{O}^{2}", Rightarrow, from=1-3, to=2-2]
\arrow["(\mathfrak{O}\circ\EV)^{2}", curve={height=-6pt}, Rightarrow, from=1-3, to=3-3]
\arrow["\mathfrak{O}_{\omega_{T}^{-1}S\mathbb{C}}", Rightarrow, from=2-2, to=3-3]
\arrow["\omega_{TS}"', Rightarrow, from=3-1, to=3-3]
\end{tikzcd}

%% file: diagrams/d577.tikz
\begin{tikzcd}[ampersand replacement=\&]
 \& {\Id} \\
 \& {\phantom{\mathfrak{O}_{T}}} \\
{\Id} \&  \& {\mathfrak{O}_{\mathbb{C}},}
\arrow["\omega_{\Id}", Rightarrow, from=1-2, to=3-3]
\arrow[equals, from=3-1, to=1-2]
\arrow["\omega_{\Id}"', Rightarrow, from=3-1, to=3-3]
\end{tikzcd}

%% file: diagrams/d439.tikz
\begin{tikzcd}[ampersand replacement=\&]
{\u(T)} \&  \& {\mathfrak{O}_{\u\overline{\EV}(T)}} \\
{\mathfrak{O}_{\u(T)\mathbb{C}}} \&  \& {\mathfrak{O}_{\EV\u(T)}}
\arrow["\Lambda_{T}", Rightarrow, from=1-1, to=1-3]
\arrow["\omega_{\u(T)}"', Rightarrow, from=1-1, to=2-1]
\arrow["\mathfrak{O}_{\Omega_{\EV}}", Rightarrow, from=1-3, to=2-3]
\arrow[equals, from=2-1, to=2-3]
\end{tikzcd}

%% file: diagrams/d473.tikz
\begin{tikzcd}[column sep=tiny,row sep=scriptsize, ampersand replacement=\&]  %d473
{{\u(TS)}} \& {{\u(T)\u(S)}} \& {{\hspc 4}} \& {{\mathfrak{O}_{\u\overline{\EV}(T)}\mathfrak{O}_{\u\overline{\EV}(S)}}} \& {{\hspc 0}} \& {{\mathfrak{O}_{\u\overline{\EV}(T)\otimes\u\overline{\EV}(S)}}} \& {{\hspc 4}} \& {{\mathfrak{O}_{\u\overline{\EV}(T\cdot S)}}} \\
 \\
 \& {{\u(T)\u(S)}} \&  \& {{\mathfrak{O}_{\u(T)\mathbb{C}}\mathfrak{O}_{\u(S)\mathbb{C}}}} \&  \& {{\mathfrak{O}_{\u(T)\mathbb{C}\otimes\u(S)\mathbb{C}}}} \&  \& {{\mathfrak{O}_{\u(TS)\mathbb{C}}}}
\arrow[equals, from=1-1, to=1-2]
\arrow["{\Lambda_{TS}}", sqarn=1.5em, Rightarrow, from=1-1, to=1-8]
\arrow["{\Lambda_{T}\Lambda_{S}}", Rightarrow, from=1-2, to=1-4]
\arrow["{\mathfrak{O}^{2}}", Rightarrow, from=1-4, to=1-6]
\arrow["{\mathfrak{O}_{\Omega_{\EV}}\mathfrak{O}_{\Omega_{\EV}}}"', Rightarrow, from=1-4, to=3-4]
\arrow[equals, from=1-6, to=1-8]
\arrow["{\mathfrak{O}_{\Omega_{\EV}\otimes\Omega_{\EV}}}", Rightarrow, from=1-6, to=3-6]
\arrow["{\Omega_{\EV}}", Rightarrow, from=1-8, to=3-8]
\arrow[equals, from=3-2, to=1-2]
\arrow["{\omega_{\u(T)}\omega_{\u(S)}}", Rightarrow, from=3-2, to=3-4]
\arrow["{\omega_{\u(T)\u(S)}}", sqars=1.5em, Rightarrow, from=3-2, to=3-8]
\arrow["{\mathfrak{O}^{2}}", Rightarrow, from=3-4, to=3-6]
\arrow["{\mathfrak{O}_{\omega^{-1}_{\u(T)}\u(S)\mathbb{C}}}", Rightarrow, from=3-6, to=3-8]
\end{tikzcd}

%% file: diagrams/d441.tikz
\begin{tikzcd}[ampersand replacement=\&]
{\u(T)} \&  \&  \& {\u(S)} \\
 \& {\mathfrak{O}_{\u(T)\mathbb{C}}} \& {\mathfrak{O}_{\u(S)\mathbb{C}}} \\
{\mathfrak{O}_{\u\overline{\EV}(T)}} \&  \&  \& {\mathfrak{O}_{\u\overline{\EV}(S)}}
\arrow["\alpha", Rightarrow, from=1-1, to=1-4]
\arrow["\omega_{\u(T)}", Rightarrow, from=1-1, to=2-2]
\arrow["\Lambda_{T}"', Rightarrow, from=1-1, to=3-1]
\arrow["\omega_{\u(S)}"', Rightarrow, from=1-4, to=2-3]
\arrow["\Lambda_{S}", Rightarrow, from=1-4, to=3-4]
\arrow["\mathfrak{O}_{\alpha\mathbb{C}}"', Rightarrow, from=2-2, to=2-3]
\arrow["\mathfrak{O}_{\Omega_{\EV}}"', Rightarrow, from=3-1, to=2-2]
\arrow["\mathfrak{O}_{\overline{\EV}(\alpha)}"', Rightarrow, from=3-1, to=3-4]
\arrow["\mathfrak{O}_{\Omega_{\EV}}", Rightarrow, from=3-4, to=2-3]
\end{tikzcd}

%% file: sec4.tex
\section{Good endofunctors}

In this section, we start the discussion of good endofunctors, which
are designed to mimic the behavior of the asymptotic algebra functor.

\subsection{\label{subsec:labeled-end}Labeled endofunctors}

We
have already discussed in the introduction
that there exists a natural transformation $I\mathfrak{A}\Rightarrow\mathfrak{A}I$
which yields a well-defined composition of generalized morphisms (see
property~\ref{enu:o2} and the subsequent paragraph in the introduction).
We now formalize this property as an axiom, allowing us to swap endofunctors
of $C^{*}$-algebras not only with the cylinder $I$ but also with
an arbitrary tensor-type endofunctor. We shall use this to prove that
 the aforementioned composition is bilinear (see Proposition~\ref{prop:bullet-bilin-assoc}).
\begin{defn}\label{def:labeled-endof}
A \emph{labeled endofunctor} is a pair $(F,\{\kappa^{A,F}\}_{A\in\in\Cstar})$
where $F\in\in\text{\ensuremath{\Cstar}}^{\text{\ensuremath{\Cstar}}}$
and $\{\kappa^{A,F}\}$ is a family of natural transformations $\left\{ \kappa^{A,F}\colon\mathfrak{O}_{A}F\Rightarrow F\mathfrak{O}_{A}\right\} _{A\in\in\Cstar}$
referred to as the \emph{labeling}, such that for every $*$-homomorphism
$\varphi\colon A\to A'$ the following diagrams commute:
\begin{equation}
\input{diagrams/d245.tikz}\label{eq:kappa-def-2}
\end{equation}
\begin{alignat}{2}
\input{diagrams/d246.tikz} & \qquad & \input{diagrams/d247.tikz}\label{eq:kappa-def22}
\end{alignat}
The composition of labeled endofunctors is defined as
\[
(F,\{\kappa^{A,F}\}_{A})(G,\{\kappa^{A,G}\}_{A})\coloneqq(FG,\{\mathfrak{O}_{A}FG\xRightarrow{\kappa^{A,F}G}F\mathfrak{O}_{A}G\xRightarrow{F\kappa^{A,G}}FG\mathfrak{O}_{A}\}_{A}).
\]
It is straightforward to check that the composition is associative.

We shall sometimes abuse notation and write just~$F$ instead of
$(F,\{\kappa^{A,F}\}_{A\in\in\Cstar})$.
\end{defn}

\begin{defn}\label{def:admissible}
Let
$(F,\{\kappa^{A,F}\}_{A})$
and $(G,\{\kappa^{A,G}\}_{A})$ be two labeled endofunctors. A \emph{labeled
natural transformation}
\[
\alpha\colon(F,\{\kappa^{A,F}\}_{A})\Rightarrow(G,\{\kappa^{A,G}\}_{A})
\]
is a natural transformation $\alpha\colon F\Rightarrow G$ of the
underlying endofunctors, such that for all $A\in\in\Cstar$ the following
diagram commutes:
\[
\input{diagrams/d248.tikz}
\]
\end{defn}

\begin{rem}\label{rem:bicat-interpr}
There
is a concise bicategorical
interpretation of labeled endofunctors (which is slightly beyond the
scope of this paper, and we refer the reader to~\cite{JY2021} or~\cite{leinster1998basicbicategories}
for the necessary background). Monoidal categories are exactly the
bicategories with a single object. Thus, the monoidal functor
\[
(\mathfrak{O},\mathfrak{O}^{2},\mathfrak{O}^{0})\colon(\Cstar,\otimes,\mathbb{C},\alpha^{\otimes},\lambda^{\otimes},\rho^{\otimes})\to(\EFC,\singledot,\Id).
\]
can be viewed as a lax functor between bicategories. Labeled endofunctors
are exactly the lax transformations from~$(\mathfrak{O},\mathfrak{O}^{2},\mathfrak{O}^{0})$
to itself, and labeled natural transformations between labeled endofunctors
are exactly the modifications between the corresponding lax transformations.
\end{rem}

\begin{example}
The identity functor can be viewed as a labeled endofunctor: $(\Id,\{\mathfrak{O}_{A}\Id\xRightarrow{=}\Id\mathfrak{O}_{A}\}_{A})$.
\end{example}

\begin{example}\label{exa:CX-labeled}
The endofunctors $\mathfrak{C}^{b}_{\mathsf{X}}$
and $\C_{\mathsf{X}}$ from the previous section can be endowed with
the natural labelings $\{\kappa^{A,\mathfrak{C}^{b}_{\mathsf{X}}}\}$
and $\{\kappa^{A,\C_{\mathsf{X}}}\}$ given by the formula $a\otimes f\mapsto(x\mapsto a\otimes f(x))$.
It will be shown in Proposition~\ref{prop:unique-labeling} that
these labelings are unique. For any continuous (resp. proper continuous)
function $f\colon\mathsf{X}\to\mathsf{Y}$, the natural transformation
$\mathfrak{C}^{b}_{f}\colon\mathfrak{C}^{b}_{\mathsf{Y}}\Rightarrow\mathfrak{C}^{b}_{\mathsf{X}}$
(resp. $\C_{f}\colon\C_{\mathsf{Y}}\Rightarrow\C_{\mathsf{X}}$) is
obviously labeled.
\end{example}

\begin{example}\label{exa:mat-labeled}
$\fM_{n}$, $\fK$, and $\mathfrak{M}^{u}_{\mathbb{Z}}$
from Section~\ref{subsec:stacont} carry the structure of labeled
endofunctors with $\kappa^{A,\fM_{n}}$, $\kappa^{A,\fK}$, and $\kappa^{A,\mathfrak{M}^{u}_{\mathbb{Z}}}$
given by the formula
\[
a\otimes(b_{i,j})\mapsto(a\otimes b_{i,j}).
\]
\end{example}

\begin{defn}\label{def:lefc-2cat}
Define the category $\LEFC$ with objects labeled
endofunctors and arrows labeled natural transformations.
\end{defn}

\begin{lem}\label{lem:adm-vh-closed}
Labeled
natural transformations
are closed under vertical and horizontal compositions.
\end{lem}

\begin{proof}
Straightforward.
\end{proof}

\begin{lem}
There is a strict monoidal category $(\LEFC,\singledot,\Id)$ where
$\singledot$ stands for the composition of labeled endofunctors (which
will be denoted again by concatenation).
\end{lem}

\begin{proof}
Straightforward.
\end{proof}

\begin{rem}
Generally speaking, the statement of the previous lemma follows from
Remark~\ref{rem:bicat-interpr}, since lax functors between bicategories
$\mathcal{C}$ and $\mathcal{C}'$, lax transformations, and modifications
form a bicategory, which is a $2$-category whenever $\mathcal{C}'$
is a $2$-category.
\end{rem}

\begin{lem}\label{lem:labeled-isom}
A labeled natural transformation $\alpha$
is an isomorphism in $\LEFC$ if and only if $\alpha$ is a natural
isomorphism of the underlying endofunctors.
\end{lem}

\begin{proof}
Straightforward.
\end{proof}

\begin{prop}\label{prop:lefc-complete}
The category $\LEFC$ is complete. Specifically,
if $\overline{F}_{\bullet}\colon\mathcal{J}\to\LEFC$ is a small diagram,
and $F_{\bullet}\colon\mathcal{J}\to\EFC$ is the diagram of the underlying
endofunctors, then $\lim\overline{F}_{\bullet}$ is equal to
\[
\left(\lim_{j}F_{j},\left\{ \mathfrak{O}_{A}\lim_{j}(F_{j})\xRightarrow{\lc}\lim_{j}(\mathfrak{O}_{A}F_{j})\xRightarrow{\lim_{j}\kappa^{A,F_{j}}}\lim_{j}(F_{j}\mathfrak{O}_{A})\xRightarrow{\rc^{-1}}\lim_{j}(F_{j})\mathfrak{O}_{A}\right\} _{A\in\in\Cstar}\right).
\]
\end{prop}

\begin{proof}
Denote by $F_{\infty}\coloneqq\lim_{j}F_{j}$ the limit of $F_{\bullet}$,
and recall that $\pr_{j}\colon F_{\infty}\Rightarrow F_{j}$ stands
for the legs of the limit cone. It follows from Definition~\ref{def:ab}
and Lemma~\ref{lem:rc-isom-id} that for all $D\in\in\Cstar$,
\begin{equation}
\{F_{\infty}\mathfrak{O}_{D}\xRightarrow{\pr_{j}\mathfrak{O}_{D}}F_{j}\mathfrak{O}_{D}\}\label{eq:uni-cone}
\end{equation}
is the limit cone of the diagram $F_{\bullet}\mathfrak{O}_{D}$.

Let $\varphi\colon A\to A'$ be a $*$-homomorphism. The solid faces
of the cube
\[
\input{diagrams/d268.tikz}
\]
clearly commute. Using the universal property of the cone (\ref{eq:uni-cone}),
we can define $\kappa^{A,F_{\infty}}$ and $\kappa^{A',F_{\infty}}$
to be the unique arrows making respectively the front and the back
faces commutative for all $j\in\in\mathcal{J}$, and conclude that
the left face also commutes. Thus, we have just obtained a labeled
endofunctor $(F_{\infty},\{\kappa^{A,F_{\infty}}\}_{A})$. Furthermore,
$\pr_{j}\colon F_{\infty}\Rightarrow F_{j}$ is labeled for all $j\in\in\mathcal{J}$
since the front face commutes.

Let us check that $(F_{\infty},\{\kappa^{A,F_{\infty}}\}_{A})$ is
a limit. To do this, let $\tau_{j}\colon(G,\{\kappa^{A,G}\})\Rightarrow(F_{j},\{\kappa^{A,F_{j}}\})$
be a cone in $\LEFC$ over $\overline{F}_{\bullet}$. Forgetting about
the labeling for the time being, we can consider the cone $\tau_{j}\colon G\Rightarrow F_{j}$
which factors uniquely through the limit cone over~$F_{\bullet}$:
\[
\input{diagrams/d270.tikz}
\]
Now we only need to verify that $\tau_{\infty}$ is labeled. In other
words, we should show that the top face of
\[
\input{diagrams/d269.tikz}
\]
commutes, which in turn follows from the commutativity of the rest
of the faces and the universal property of the cone~(\ref{eq:uni-cone}).
\end{proof}

\begin{prop}\label{prop:underlying-canon}
Let
$\overline{F}_{\bullet}\colon\mathcal{J}\to\LEFC$
be a small diagram, let $\overline{G}\in\in\LEFC$, let $\{\overline{\tau}_{j}\colon\overline{G}\Rightarrow\overline{F}_{j}\}_{j\in\in\mathcal{J}}$
be a cone over~$\overline{F}_{\bullet}$, and let $\overline{\tau}_{\infty}\colon\overline{G}\Rightarrow\lim\overline{F}_{\bullet}$
be the unique labeled natural transformation such that $\overline{\tau}_{j}=\pr_{j}\circ\overline{\tau}_{\infty}$
for all $j\in\in\mathcal{J}$. Let $G,F_{j}\in\in\Cstar$ and $\tau_{j}\colon G\Rightarrow F_{j}$
be the endofunctors and natural transformations underlying respectively
$\overline{G}$, $\overline{F}_{j}$, and $\overline{\tau}_{j}$,
for all $j\in\in\mathcal{J}$. Denote by $\tau_{\infty}\colon G\Rightarrow\lim F_{\bullet}$
the unique labeled natural transformation such that $\tau_{j}=\pr_{j}\circ\tau_{\infty}$
for all $j\in\in\mathcal{J}$. Then:%--!--
\begin{enumerate}[label=\textup{(\roman*)}]

\item\label{enu:undercanon1}$\tau_{\infty}$ is the natural transformation
underlying $\overline{\tau}_{\infty}$;
\item\label{enu:undercanon2}$\{\tau_{j}\}$ is a limit cone if and only
if $\{\overline{\tau}_{j}\}$ is.
\end{enumerate}
\end{prop}

\begin{proof}
The first statement follows from the proof of Proposition~\ref{prop:lefc-complete}.
To prove the second statement, note that if $\{\tau_{j}\}$ is a limit
cone in $\EFC$, then $\tau_{\infty}$ is an isomorphism in $\EFC$,
and apply Lemma~\ref{lem:labeled-isom}.
\end{proof}

\begin{defn}
Let $(F_{1},\kappa_{1})$ and $(F_{2},\kappa_{2})$ be labeled endofunctors.
Define their binary product
\[
(F,\{\kappa^{A,F}\})\oplus(G,\{\kappa^{A,G}\})\coloneqq(F\oplus G,\{\kappa^{A,F\oplus G}\})
\]
 where $\kappa^{A,F\oplus G}$ is the composite
\[
\kappa^{A,F\oplus G}\colon\mathfrak{O}_{A}(F\oplus G)\xRightarrow{\lc}\mathfrak{O}_{A}F\oplus\mathfrak{O}_{A}G\xRightarrow{\kappa^{A,F}\oplus\kappa^{A,G}}F\mathfrak{O}_{A}\oplus G\mathfrak{O}_{A}\xRightarrow{\rc^{-1}}(F\oplus G)\mathfrak{O}_{A}.
\]
It follows from Proposition~\ref{prop:lefc-complete} that $\oplus$
is indeed a binary product.
\end{defn}

\begin{prop}\label{prop:OisLabaled}
For every $D\in\in\Cstar$,
\[
\{\kappa^{A,\mathfrak{O}_{D}}\colon\mathfrak{O}_{A}\mathfrak{O}_{D}\xRightarrow{\mathfrak{O}^{2}}\mathfrak{O}_{A\otimes D}\xRightarrow{\mathfrak{O}_{\xi^{\otimes}}}\mathfrak{O}_{D\otimes A}\xRightarrow{(\mathfrak{O}^{2})^{-1}}\mathfrak{O}_{D}\mathfrak{O}_{A}\}_{A\in\in\Cstar}
\]
is the unique family which turns $\mathfrak{O}_{D}$ into a labeled
endofunctor. Furthermore, every $*$-homomorphism $\varphi\colon D\to D'$
induces a labeled natural transformation
\[
\mathfrak{O}_{\varphi}\colon(\mathfrak{O}_{D},\{\kappa^{A,\mathfrak{O}_{D}}\}_{A})\Rightarrow(\mathfrak{O}_{D'},\{\kappa^{A,\mathfrak{O}_{D'}}\}_{A}).
\]
\end{prop}

\begin{proof}
That $(\mathfrak{O}_{D},\{\kappa^{A,\mathfrak{O}_{D}}\})$ and $\mathfrak{O}_{\varphi}$
are labeled follows from the coherence theorem for the multiplicative
structure of the symmetric bimonoidal category~$\Cstar$. The uniqueness
of $\{\kappa^{A,\mathfrak{O}_{D}}\}$ follows from Lemma~\ref{lem:OA-wp}
and Proposition~\ref{prop:unique-labeling}, which we shall prove
later.
\end{proof}

\begin{prop}\label{prop:kappa-labeled}
Let $(F,\{\kappa^{A,F}\}_{A})$ be a labeled
endofunctor. Then for all $A\in\in\Cstar$ the natural transformation
$\kappa^{A,F}\colon\mathfrak{O}_{A}F\Rightarrow F\mathfrak{O}_{A}$
is labeled.
\end{prop}

\begin{proof}
We need only check that the central rectangle in the diagram
\[
\input{diagrams/d491.tikz}
\]
commutes, which follows from the fact that all the other inner subdiagrams
and the outermost rectangle commute by the definition of labeled endofunctors
and their composition, and the definition of $\kappa^{B,\mathfrak{O}_{A}}$
(see Proposition~\ref{prop:OisLabaled}).
\end{proof}

We finish this subsection with a couple of easy lemmas which we shall
use in future papers.
\begin{lem}\label{lem:abg}
Let $F,G,H\in\in\LEFC$, and let $\alpha$, $\beta$,
$\gamma$ some natural transformations between the underlying endofunctors
as in the diagram
\[
\input{diagrams/d467.tikz}
\]
such that $\alpha$ is a componentwise monomorphism. Then $\gamma$
is labeled whenever $\alpha$ and $\beta$ are.
\end{lem}

\begin{proof}
The upper face of
\[
\input{diagrams/d468.tikz}
\]
is commutative since the rest of the faces are and $\alpha$ is a
componentwise monomorphism.
\end{proof}

We finish this subsection with a technical lemma which we shall need
in~\cite{makeev_asadj}.
\begin{lem}\label{lem:27}
Let $F_{0},F,F_{1},G\in\in\LEFC$, and let $\iota$,
$\pi$, $\alpha$ be natural transformations between the underlying
endofunctors as in the solid arrow part of the diagram
\[
\input{diagrams/d469.tikz}
\]
such that the upper row is a componentwise short exact sequence, and
the left-hand triangle commutes. Then:
\begin{itemize}
\item there exists a unique $\beta\colon F_{1}\Rightarrow G$ making the
right-hand triangle commutative;
\item if $\iota$, $\pi$ and $\alpha$ are labeled, then $\beta$ is also
labeled.
\end{itemize}
\end{lem}

\begin{proof}
Left as an exercise.
\end{proof}

\subsection{Good endofunctors}
\begin{defn}
A labeled endofunctor $(F,\{\kappa^{A,F}\}_{A\in\in\Cstar})$ where
$F$ is decent will be called a \emph{good endofunctor}.
\end{defn}

\begin{example}\label{exa:good-endfs}
$\fM_{n}$,
$\fK$, $\mathfrak{M}^{u}_{\mathbb{Z}}$,
$\mathfrak{C}^{b}_{\mathsf{X}}$, and $\C_{\mathsf{X}}$ carry the
structure of good endofunctors by Lemma~\ref{lem:decent-ex} and
Examples~\ref{exa:CX-labeled} and~\ref{exa:mat-labeled}.
\end{example}

\begin{defn}
Define the category $\GEFC$ with objects good endofunctors and arrows
labeled natural transformations. As usual, $\GEFC$ naturally carries
the structure of a strict monoidal category $(\GEFC,\singledot,\Id)$,
where $\singledot$ stands for the composition of labeled endofunctors
and $\Id$ is the identity endofunctor with the obvious labeling.
\end{defn}

\begin{prop}\label{prop:ac}
$\GEFC$ is a complete category.
\end{prop}

\begin{proof}
The statement follows from Propositions~\ref{prop:lefc-complete}
and~\ref{prop:decent-complete}.
\end{proof}

\begin{prop}\label{prop:gefc-codist}
The strict monoidal category $(\GEFC,\singledot,\Id)$
is codistributive.
\end{prop}

\begin{proof}
It follows from Proposition~\ref{prop:ac} that $\GEFC$ has binary
products (denoted again by $\oplus$). Let $F,G_{1},G_{2}\in\in\GEFC$.
The canonical arrows $F(G_{1}\oplus G_{2})\xRightarrow{\lc}FG_{1}\oplus FG_{2}$
and $(G_{1}\oplus G_{2})F\xRightarrow{\rc}G_{1}F\oplus G_{2}F$ are
labeled natural transformations by Proposition~\ref{prop:underlying-canon}(\ref{enu:undercanon1}).
Furthermore, they are natural isomorphisms of the underlying endofunctors,
and hence by Lemma~\ref{lem:labeled-isom} are isomorphisms in $\LEFC$,
and therefore in $\GEFC$.

It is straightforward to show that the zero endofunctor $\bzero$
with the labeling
\[
\kappa^{A,\bzero}\colon\mathfrak{O}_{A}\bzero\xRightarrow{\cong}\bzero\xRightarrow{\cong}\bzero\mathfrak{O}_{A}
\]
is a terminal object in $\GEFC$, and that the unique natural isomorphisms
$\Terminal^{l}\colon\bzero F\Rightarrow\bzero$ and $\Terminal^{r}\colon F\bzero\Rightarrow\bzero$
are labeled.
\end{proof}

\begin{cor}
There is a tight bimonoidal category
\[
(\GEFC,(\bzero,\oplus,\au,\lu,\ru,\xi^{\oplus}),(\Id,\singledot),(\Terminal^{l},\Terminal^{r}),(\lc,\rc))
\]
with the cartesian additive structure.
\end{cor}

\begin{proof}
This follows from Proposition~\ref{prop:gefc-codist} and Theorem~\ref{thm:funcat-bimonoidal}.
\end{proof}

\begin{prop}
There is a bimonoidal functor
\begin{alignat}{1}
(\mathfrak{O},\mathfrak{O}^{2},\mathfrak{O}^{0},\mathfrak{O}^{2}_{\oplus},\mathfrak{O}^{0}_{\oplus})\colon & \Cstar\to\GEFC\label{eq:ogefc}
\end{alignat}
as in Definition~\ref{def:bimon-C-DEFC}.
\end{prop}

\begin{proof}
Bearing in mind Proposition~\ref{prop:bimon-O}, we only need to
check that $(\mathfrak{O}^{2})_{A,B}$, $(\mathfrak{O}^{2}_{\oplus})_{A,B}$
and $\mathfrak{O}^{0}$, $\mathfrak{O}^{0}_{\oplus}$ are labeled
natural transformations for all $A,B\in\in\Cstar$, which follows
respectively from the coherence theorem for the (symmetric) multiplicative
structure of~$\Cstar$, and from Laplaza axioms for $\Cstar$.
\end{proof}
We conclude this subsection with a technical lemma
which will be used in the proof of Proposition~\ref{prop:bullet-bilin-assoc}
via Proposition~\ref{prop:kappa-prop}.
\begin{lem}\label{lem:d296}
Let $A_{1},A_{2}\in\in\Cstar$ and $F\in\in\GEFC$.
The outermost rectangle in the diagram below commutes:
\[
\input{diagrams/d475.tikz}
\]
\end{lem}

\begin{proof}
The upper triangles commute by the definition of $\lc$ and $\rc$;
the lower triangles commute by Lemma~\ref{lem:218}; the bottom middle
subdiagram commutes since $F$ is labeled; and the top middle one
commutes obviously. Thus, the whole diagram commutes for $j=1,2$,
and so by the universal property of $\oplus$ the outermost rectangle
also commutes.
\end{proof}

\subsection{\label{subsec:quot}Quotients of good endofunctors}

\begin{defn}\label{def:quot}
Let $F,G\in\in\EFC$ and let $\iota\colon F\Rightarrow G$
be a componentwise ideal inclusion. Define the \emph{quotient endofunctor}
$G/F$ at the object $B$ as $(G/F)B\coloneqq GB/FB$ and at the morphism
$\varphi\colon A\to B$ by a diagram chase in
\begin{equation}
\input{diagrams/d171.tikz}\label{eq:diagram-search-quotient-functor}
\end{equation}
Note that the diagram above also defines the componentwise quotient
projection, which we denote by $q\colon G\Rightarrow G/F$.
\end{defn}

\begin{example}\label{exa:frakA}
For
good endofunctors $\mathfrak{T}\coloneqq\mathfrak{C}^{b}_{[0,\infty)}$
and $\mathfrak{T}_{0}\coloneqq\mathfrak{C}_{[0,\infty)}$ as in Definition~\ref{def:Cb},
consider the obvious componentwise ideal inclusion $\mathfrak{T}_{0}\Rightarrow\mathfrak{T}$.
We call the quotient $\mathfrak{A}\coloneqq\mathfrak{T}/\mathfrak{T}_{0}$
the \emph{asymptotic algebra} functor.
\end{example}

We proceed to introduce quotients of good endofunctors. This construction
requires some preparation.

Denote by $\Ch(\Cstar)$ (resp. $\chex(\Cstar)$) the category with
objects chain complexes (resp. exact chain complexes) of $C^{*}$-algebras
and arrows chain maps. For any chain complex we write
$\delta\colon A^{k}\to A^{k+1}$ for its boundary map; and for a chain
map $\varphi\colon A\to B$, denote by $\varphi^{k}\colon A^{k}\to B^{k}$
its degree $k$ component.
\begin{lem}\label{lem:chain-pullback}
Let $\varphi_{1}\colon A_{1}\to A_{0}$
and $\varphi_{2}\colon A_{2}\to A_{0}$ be a pair of morphisms in
$\chex(\Cstar)$ such that $\varphi^{k}_{1}\colon A^{k}_{1}\to A^{k}_{0}$
is epic for all $k\in\mathbb{Z}$. Then there exists $P\in\in\chex(\Cstar)$
and $p_{j}\colon P\to A_{j}$, $j=1,2$, such that the left-hand side
diagram in
\begin{alignat}{2}
\input{diagrams/d192.tikz} & \qquad & \input{diagrams/d201.tikz}\label{eq:chain-pullback}
\end{alignat}
is a pullback in $\Ch(\Cstar)$, and the right-hand side diagram is
a pullback in $\Cstar$ for all~$k\in\mathbb{Z}$.
\end{lem}

\begin{proof}
Denote by $P^{k}$ the pullbacks as in the right-hand side diagram
in~(\ref{eq:chain-pullback}). Consider the complex $\dots\xrightarrow{\delta}P^{n}\xrightarrow{\delta}P^{n+1}\xrightarrow{\delta}\dots$,
where $\delta\colon P^{n}\to P^{n+1}$ is induced by the pair $\delta\circ p^{n}_{j}$,
$j=1,2$ (in other words $\delta(h_{1},h_{2})=(\delta h_{1},\delta h_{2})$\footnote{Here we use the explicit form of pullbacks in $\Cstar$, see (\ref{eq:explicit-pullback})
in Section~\ref{subsec:cal}.}). It is clear that $\delta\delta=0$, and that the left-hand side
diagram in~(\ref{eq:chain-pullback}) commutes in~$\Ch(\Cstar)$.

To show that $P$ is exact, take $(h_{1},h_{2})\in P^{n+1}$, such
that $\delta(h_{1},h_{2})=0$, and choose $g_{1}\in A^{n}_{1}$ and
$g_{2}\in A^{n}_{2}$, such that $\delta g_{1}=h_{1}$, $\delta g_{2}=h_{2}$.
Note that $\varphi_{2}(g_{2})-\varphi_{1}(g_{1})\in\ker\delta$. Hence,
there is $a_{0}\in A^{n-1}_{0}$, such that $\delta a_{0}=\varphi_{2}(g_{2})-\varphi_{1}(g_{1})$.
Since all the components of $\varphi_{1}$ are epimorphisms, there
is $a_{1}\in A^{n-1}_{1}$, such that $\varphi_{1}(a_{1})=a_{0}$.
It follows from the equality
\[
\varphi_{1}(g_{1}+\delta a_{1})=\varphi_{1}(g_{1})+\delta\varphi_{1}(a_{1})=\varphi_{1}(g_{1})+\delta a_{0}=\varphi_{2}(g_{2}),
\]
that $(g_{1}+\delta a_{1},g_{2})\in P^{n}$. So, $(h_{1},h_{2})=\delta(g_{1}+\delta a_{1},g_{2})$,
and the complex $P$ is exact.

Let $B\in\in\Ch(\Cstar)$ be a chain complex and let $q_{1}\colon B\to A_{1}$,
$q_{2}\colon B\to A_{2}$ be a pair of chain maps, such that $\varphi_{1}\circ q_{1}=\varphi_{2}\circ q_{2}$.
For every $n\in\mathbb{Z}$, denote by $f^{n}\colon B^{n}\to P^{n}$
the unique
$*$-homomorphisms induced by the pair
$(q^{n}_{1},q^{n}_{2})$. The two $*$-homomorphisms
\[
f^{n}\circ\delta,~\delta\circ f^{n-1}\colon B^{n}\to P^{n+1}
\]
are both induced by $(\delta\circ q^{n}_{1},\delta\circ q^{n}_{2})$,
and hence coincide. Thus, $\{f^{n}\}$ is a chain map which is unique
by construction, and the proof is complete.
\end{proof}

\begin{lem}\label{lem:quot-decent}
Let $F,G\in\in\DEFC$ and let $\iota\colon F\Rightarrow G$
be a componentwise ideal inclusion. Then:
\begin{itemize}
\item $G/F\in\in\DEFC$;
\item if $F$ and $G$ are exact, then so is $G/F$.
\end{itemize}
\end{lem}

\begin{proof}
For brevity, we set $H\coloneqq G/F$, and define the functor $\Omega\colon\Cstar\to\chex(\Cstar)$
by the formula
\[
A\mapsto(0\to FA\to GA\to HA\to0).
\]
Consider the diagrams
\begin{alignat*}{2}
\input{diagrams/d172.tikz} & \qquad & \input{diagrams/d168.tikz}
\end{alignat*}
where the left-hand side diagram is a pullback with $\varphi_{2}$
a split epimorphism; the perimeter of the right-hand side diagram
is obtained by applying $\Omega$ to the left-hand side one; the right-bottom
rectangle is a pullback square which by Lemma~\ref{lem:chain-pullback}
exists and is computed degree-wise (i.e., each degree $k$ component
is a pullback in $\Cstar$); the dashed arrow is induced by the universal
property of the pullback. Expanding the definition of $\Omega$, we
rewrite the dashed arrow in the following form:
\[
\input{diagrams/d173.tikz}
\]
where the vertical arrows are easily seen to be the canonical ones
(see Definition~\ref{def:ab}). Furthermore, the left and the middle
ones are isomorphisms since $F$ and $G$ are decent. Thus, the rightmost
one is also an isomorphism, and so $H$ is decent.

Let now $F$ and $G$ be exact, and let
\[
0\to J\to A\to B\to0
\]
 be a short exact sequence. The exactness of $H$ follows from a diagram
chase on
\[
\input{diagrams/d476.tikz}
\]
\end{proof}

\begin{lem}\label{lem:quotient-criterion}
For a componentwise ideal inclusion
$\iota\in\GEFC(F,G)$, there exists a unique family
\[
\{\kappa^{A,G/F}\colon\mathfrak{O}_{A}(G/F)\Rightarrow(G/F)\mathfrak{O}_{A}\}_{A\in\in\Cstar}
\]
 satisfying the following properties:
\begin{itemize}
\item $\{\kappa^{A,G/F}\}_{A}$ turns $G/F$ into a good endofunctor;
\item $q\colon G\Rightarrow G/F$ from Definition~\ref{def:quot} is labeled,
i.e., $q\in\GEFC(G,G/F)$.
\end{itemize}
\end{lem}

\begin{proof}
For every $B\in\in\Cstar$, define $\kappa^{A,H}_{B}$ by a diagram
chase in
\begin{equation}
\input{diagrams/d260.tikz}\label{eq:quot-d0}
\end{equation}
To prove the naturality of $\kappa^{A,H}_{B}$ in $B$, we should
show that for a $*$-homomorphism $\varphi\colon B\to B'$, the right
face in
\begin{equation}
\input{diagrams/d261.tikz}\label{eq:d261}
\end{equation}
commutes. But this follows from the fact that the other faces commute
(by the definition of $\kappa^{A,H}_{B}$, $\kappa^{A,H}_{B'}$ and
$q$) and that $\mathfrak{O}_{A}qB$ is epic (since $\mathfrak{O}_{A}$
is exact by Lemma~\ref{prop:O-decent}, and $q$ is a quotient projection).
Thus, we have defined the natural transformation $\kappa^{A,H}\colon\mathfrak{O}_{A}H\Rightarrow H\mathfrak{O}_{A}$.

Note that $G/F$ is decent by Lemma~\ref{lem:quot-decent}. So, to
finish the proof of the first statement, we need to show the commutativity
of the three diagrams from Definition~\ref{def:labeled-endof}. The
first one is exactly the back face in
\[
\input{diagrams/d262.tikz}
\]
which commutes since the rest of the faces do and $\mathfrak{O}_{A}q$
is epic by Lemma~\ref{lem:cw-epic-epic}. The proof for the other
two diagrams is similar.

The second statement follows from the commutativity of the right square
in~(\ref{eq:quot-d0}).
\end{proof}

\begin{example}\label{exa:frakA-labeling}
Recall
the asymptotic
algebra endofunctor $\mathfrak{A}$ defined in~Example~\ref{exa:frakA}.
Lemma~\ref{lem:quotient-criterion} yields the labeling $\{\kappa^{A,\mathfrak{A}}\}$
which turn $\mathfrak{A}$ into a good endofunctor. The components
of this labeling are given by the formula
\[
a\otimes[t\mapsto f(t)]\mapsto[t\mapsto a\otimes f(t)]
\]
 (cf. Example~\ref{exa:CX-labeled}), where the square brackets denote
the class of a function in the asymptotic algebra.
\end{example}

\begin{rem}
Henceforth, the quotient $G/F$ of two good endofunctors (whenever
it is defined) will always stand for the good endofunctor with the
labeling $\{\kappa^{A,G/F}\}_{A}$ as in Lemma~\ref{lem:quotient-criterion}.
\end{rem}

\subsection{Tensor-type endofunctors revisited}

In this section, we extend the notion of labeling and allow us to
swap labeled endofunctors with an arbitrary tensor-type endofunctor
rather than just $\mathfrak{O}_{A}$.
\begin{defn}\label{def:kappaTT}
For $T\in\in\DEFC_{\ttt}$ and $F\in\in\LEFC$,
define the natural transformation $\kappa^{T,F}$ as the composite
\begin{equation}
\kappa^{T,F}\colon TF\xRightarrow{\omega_{T}F}\mathfrak{O}_{T\mathbb{C}}F\xRightarrow{\kappa^{T\mathbb{C},F}}F\mathfrak{O}_{T\mathbb{C}}\xRightarrow{F\omega^{-1}_{T}}FT.\label{eq:kappa37}
\end{equation}
\end{defn}

\begin{lem}\label{lem:DEFCttProp}
For every $T\in\in\DEFC_{\ttt}$,
\[
\{\kappa^{A,T}:\mathfrak{O}_{A}T\xRightarrow{\mathfrak{O}_{A}\omega_{T}}\mathfrak{O}_{A}\mathfrak{O}_{T\mathbb{C}}\xRightarrow{\kappa^{A,\mathfrak{O}_{T\mathbb{C}}}}\mathfrak{O}_{T\mathbb{C}}\mathfrak{O}_{A}\xRightarrow{(\omega_{T})^{-1}\mathfrak{O}_{A}}T\mathfrak{O}_{A}\}_{A\in\in\Cstar},
\]
is the unique family of natural transformations turning $T$ into
a good endofunctor. Furthermore, every $\alpha\in\DEFC_{\ttt}(T_{1},T_{2})$
is labeled (with respect to the unique labelings on $T_{1}$ and~$T_{2}$).
\end{lem}

\begin{proof}
This follows from Proposition~\ref{prop:OisLabaled}.
\end{proof}

\begin{defn}\label{def:GEFCtt}
Denote by $\GEFC_{\ttt}$ the full subcategory
of $\GEFC$ generated by the objects whose underlying endofunctors
are
tensor-type.
\end{defn}

It follows from Lemma~\ref{lem:DEFCttProp} that we can identify
$\GEFC_{\ttt}$ with $\DEFC_{\ttt}$. Bearing this in mind, let us
introduce the fully faithful unitary bimonoidal functor $\EV\colon\GEFC_{\ttt}\to\Cstar$
as in Proposition~\ref{prop:EV} and rewrite Corollary~\ref{cor:EV-DEFC}
in the following form.
\begin{cor}\label{cor:EV-GEFC}
There is a fully faithful strict bimonoidal functor
$\overline{\EV}\colon\Br(\GEFC_{\ttt})\to\Br(\Cstar)$.
\end{cor}

\begin{lem}
Let $T\in\in\GEFC_{\ttt}$, and $F\in\in\GEFC$. Then $\kappa^{T,F}\colon TF\Rightarrow FT$
is labeled.
\end{lem}

\begin{proof}
Follows from Proposition~\ref{prop:kappa-labeled}.
\end{proof}

\begin{lem}\label{lem:kappa-tensor-type}
Let $\alpha\in\GEFC_{\ttt}(T,S)$, and
$\beta\in\GEFC(F,G)$. Then the following diagram commutes:
\[
\input{diagrams/d385.tikz}
\]
\end{lem}

\begin{proof}
It suffices to prove the statement in two cases: (a) $F=G$, $\beta=\id_{F}$,
and (b) $T=S$, $\alpha=\id_{T}$. We prove only (a) and leave (b)
to the reader. Consider the diagram
\[
\input{diagrams/d417.tikz}
\]
in which the top and bottom faces commute by the definitions of $\kappa^{T,F}$
and $\kappa^{S,F}$; the back face commutes since $F$ is labeled;
and the left and right faces commute by Theorem~\ref{thm:aoa}. Thus,
the front face also commutes, and the proof is complete.
\end{proof}

The previous lemma gives us a useful corollary for the functors $\C_{\mathsf{X}}$
and $\fK$ from Section~\ref{subsec:stacont}.
\begin{cor}\label{cor:332}
Let $F\in\in\GEFC$, $\mathsf{X}$ a locally compact
Hausdorff space, and $x\in\mathsf{X}$. Then the following diagrams
commute:
\begin{alignat}{2}
\input{diagrams/d273.tikz} & \qquad & \input{diagrams/d354.tikz}\label{eq:d273}
\end{alignat}
\end{cor}

\begin{proof}
This follows from Lemma~\ref{lem:kappa-tensor-type} and the fact
that $\kappa^{\Id,F}=\id_{F}$.
\end{proof}

\begin{prop}\label{prop:kappa-prop}
Let $T,S\in\in\GEFC_{\ttt}$ and $F\in\in\GEFC$.
The following diagrams commute in $\GEFC$:
\begin{alignat*}{2}
\input{diagrams/d397.tikz} & \qquad & \input{diagrams/d396.tikz}
\end{alignat*}
\end{prop}

\begin{proof}
The commutativity of the left-hand side diagram follows from that
of the outermost rectangle in
\[
\input{diagrams/d398.tikz}
\]
in which the upper two subdiagrams commute by the definitions of $\kappa^{T,F}$
and $\kappa^{S,F}$; the left and right subdiagrams commute by the
definition of $\omega_{TS}$ (see (\ref{eq:omega1}) in the proof
of Lemma~\ref{lem:def-omegaF}); and the middle and the bottom subdiagrams
commute by the axioms of labeled endofunctors.

The commutativity of the right-hand side diagram follows similarly
from
\[
\input{diagrams/d484.tikz}
\]
where the central rectangle commutes by Lemma~\ref{lem:d296}, and
the rest of the subdiagrams commute by the definition of $\omega_{T\oplus S}$,
the axioms of labeled endofunctors, and naturality
of $\lc$ and $\rc$.
\end{proof}

\begin{lem}\label{lem:d376}
For every $G\in\in\GEFC$,
the
following diagram commutes in $\GEFC$:
\[
\input{diagrams/d376.tikz}
\]
\end{lem}

\begin{proof}
This follows from Lemma~\ref{prop:kappa-prop} and the fact that
$\kappa^{\Id,G}=\id_{G}$.
\end{proof}

The next
lemma is a good example of how bracketing
can simplify the statement by avoiding the use of coherence morphisms.
\begin{lem}\label{lem:af}
For
$T,S\in\in\GEFC_{\ttt}$, let
us consider $\kappa^{T,S}\colon TS\Rightarrow ST$ as an arrow in
$\Br(\GEFC_{\ttt})$:
\[
\kappa^{T,S}\colon((T,S),\mmph\cdot\mmph)\Rightarrow((S,T),\mmph\cdot\mmph).
\]
Then $\u\overline{\EV}(\kappa^{T,S})=\xi^{\otimes}_{T\mathbb{C},S\mathbb{C}}$.
\end{lem}

\begin{proof}
The diagrams
\[
\input{diagrams/d420.tikz}
\]
\[
\input{diagrams/d482.tikz}
\]
commute by the definitions of $\overline{\EV}$, $\kappa^{T,S}$,
$\kappa^{T\mathbb{C},\mathfrak{O}_{S\mathbb{C}}}$, $\kappa^{T,\mathfrak{O}_{S\mathbb{C}}}$,
the definition of tensor-type endofunctors, and the naturality of~$\rho^{\otimes}$.
Hence,
\[
\u\overline{\EV}(\kappa^{T,S})=\omega_{S}\omega^{-1}_{T}\mathbb{C}\circ\kappa^{T\mathbb{C},S}\mathbb{C}=\xi^{\otimes}_{T\mathbb{C},S\mathbb{C}}.
\]
\end{proof}

\subsection{Well-pointed endofunctors}
\begin{defn}
We say that an endofunctor $F\in\in\DEFC$ is \emph{well-pointed}
if for every locally compact Hausdorff space $\mathsf{X}$, the natural
transformation $F\C_{\mathsf{X}}\Rightarrow\prod_{\mathsf{X}}F$ induced
by $\{F\ev_{x}\}_{x\in\mathsf{X}}$ is monic.
\end{defn}

\begin{lem}\label{lem:OA-wp}
The following endofunctors are well-pointed:
\begin{itemize}
\item $\mathfrak{O}_{A}$ for all $A\in\in\Cstar$;
\item all $T\in\in\DEFC_{\ttt}$;
\item $\mathfrak{C}^{b}_{\mathsf{X}}$ for every topological space $\mathsf{X}$.
\end{itemize}
\end{lem}

\begin{proof}
Let us prove the first statement. Consider the diagram
\[
\input{diagrams/d430.tikz}
\]
in which the dashed arrows are induced by the universal property of
the product, and $\xi_{\mathfrak{O}_{A},\C_{\mathsf{X}}}$ is as in
Proposition~\ref{prop:xitt}. The upper triangle commutes by the
naturality of $\xi$ and Remark~\ref{rem:xiid}. The right-hand dashed
arrow is obviously monic, and $\xi_{\mathfrak{O}_{A},\C_{\mathsf{X}}}$
is an isomorphism, thus the left-hand dashed arrow is also monic,
and hence, $\mathfrak{O}_{A}$ is well-pointed.

The second statement follows from the first. The third is obvious.
\end{proof}

\begin{rem}
Note that the quotient of two well-pointed endofunctors is not necessarily
well-pointed. For example, one can show that $\mathfrak{A}\coloneqq\mathfrak{C}^{b}_{\mathbb{R}_{+}}/\C_{\mathbb{R}_{+}}$
is not well-pointed.
\end{rem}

\begin{prop}\label{prop:unique-labeling}
If $F\in\in\DEFC$ is a well-pointed
endofunctor, then there is at most one family $\{\kappa^{A,F}\colon\mathfrak{O}_{A}F\Rightarrow F\mathfrak{O}_{A}\}_{A\in\in\Cstar}$
which turns $F$ into a labeled endofunctor.
\end{prop}

\begin{proof}
Suppose we have two such families, say $\{\kappa^{A,F}_{1}\}$ and
$\{\kappa^{A,F}_{2}\}$. For an arbitrary $\mathsf{X}\subset\mathbb{R}$
and $j=1,2$, denote by $\kappa^{\C_{\mathsf{X}},F}_{j}$ the natural
transformation obtained from $\{\kappa^{A,F}_{j}\}$ as in Definition~\ref{def:kappaTT}.
It follows from the left-hand diagram in Corollary~\ref{cor:332}
and the definition of well-pointed endofunctors that $\kappa^{\C_{\mathsf{X}},F}_{1}=\kappa^{\C_{\mathsf{X}},F}_{2}$,
and therefore $\kappa^{C_{0}(\mathsf{X}),F}_{1}=\kappa^{C_{0}(\mathsf{X}),F}_{2}$.
For an arbitrary $*$-homomorphism $\zeta\colon C_{0}(\mathsf{X})\to A$,
there is a commuting diagram
\[
\input{diagrams/d274.tikz}
\]
for $j=1,2$. Hence, for every $B\in\in\Cstar$, we have $\kappa^{A,F}_{1}B\circ\mathfrak{O}_{\zeta}FB=\kappa^{A,F}_{2}B\circ\mathfrak{O}_{\zeta}FB$,
which by Lemma~\ref{lem:well-ppointed} implies that $\kappa^{A,F}_{1}B=\kappa^{A,F}_{2}B$.
\end{proof}

\begin{prop}\label{prop:wp-lbl}
Let $F,G\in\in\GEFC$, with $G$ well-pointed.
Then every natural transformation $\alpha\colon F\Rightarrow G$ between
the underlying endofunctors is labeled.
\end{prop}

\begin{proof}
For an arbitrary $\mathsf{X}\subset\mathbb{R}$, consider the diagram
\begin{equation}
\input{diagrams/d278.tikz}\label{eq:d278}
\end{equation}
in which the left and right subdiagrams clearly commute, and the upper
and lower triangles commute by Corollary~\ref{cor:332}. Hence, by
the definition of well-pointed endofunctors, the outermost rectangle
also commutes.

Now for an arbitrary $*$-homomorphism $\zeta\colon C_{0}(\mathsf{X})\to A$,
consider the diagram
\[
\input{diagrams/d279.tikz}
\]
in which the front face commutes by what has already been proved (the
outermost rectangle in (\ref{eq:d278})); the top and bottom faces
commute since $F$ and $G$ are labeled; the left and right commute
obviously. It follows from Lemma~\ref{lem:well-ppointed} that the
back face also commutes, and the proof is complete.
\end{proof}

It would be nice to have a way to construct new well-pointed endofunctors
from old ones. To this end, we introduce a slightly stronger property
which is closed under taking compositions.
\begin{defn}
Let $\mathcal{C}$ and $\mathcal{C}'$ be complete categories. We
say that a functor $F\colon\mathcal{C}\to\mathcal{C}'$ is
\begin{itemize}
\item \emph{product-separating} if $F\prod_{j}B_{j}\xrightarrow{\lc}\prod_{j}FB_{j}$
is a monomorphism for any family $\{B_{j}\}$ of objects in $\mathcal{C}$;
\item \emph{monic-preserving} if it preserves monomorphisms.
\end{itemize}
\end{defn}

\begin{example}
$\mathfrak{C}^{b}_{\mathsf{X}}$, $\C_{\mathsf{X}}$, and $\fK$ defined
in the previous subsection are obviously product-separating and monic-preserving,
but $\mathfrak{O}_{A}$ is not monic-preserving in general~\cite[p.~254]{Ped_pullback99}.
\end{example}

\begin{lem}
Product-separating monic-preserving endofunctors are closed under
composition.
\end{lem}

\begin{proof}
Straightforward.
\end{proof}

\begin{lem}\label{lem:eval-at-dense}
Product-separating and monic-preserving
endofunctors in $\DEFC$ are well-pointed.
\end{lem}

\begin{proof}
Let $\mathsf{X}$ be a locally compact Hausdorff space, and let $\alpha_{1},\alpha_{2}\colon G\Rightarrow F$
be two natural transformations, such that $\ev_{x}\circ\alpha_{1}=\ev_{x}\circ\alpha_{2}$
for all $x\in\mathsf{X}$. It is clear that the $*$-homomorphism
$\C_{\mathsf{X}}B\to\prod_{\mathsf{X}}B$ generated by $\{\ev_{x}B\}_{x\in\mathsf{X}}$
is monic for all $B\in\in\Cstar$. For all $x\in\mathsf{X}$, the
dashed arrow in
\[
\input{diagrams/d365.tikz}
\]
is independent of $j\in\{1,2\}$. Hence, the horizontal composite
is also independent of $j$ by the universal property of the product
$\prod_{\mathsf{X}}FB$, and so $\alpha_{1}B=\alpha_{2}B$ since the
middle and right horizontal arrows are monomorphisms.
\end{proof}

We finish this subsection with an example involving \emph{minimal
tensor products} and \emph{nuclear} $C^{*}$-algebras, whose definitions
the reader can find in e.g.~\cite{murphy1990}.
\begin{example}
For $A\in\in\Cstar$, the functor $A\otimes_{\min}\mmph$ is product
separating and monic-preserving. Indeed, let $\{B_{j}\}$ be a family
of $C^{*}$-algebras, and let $A\to\mathcal{H}$, $B_{j}\to\mathcal{H}_{j}$
be some faithful representations. Then by the definition of the minimal
tensor product, the vertical arrows in
\[
\input{diagrams/d348.tikz}
\]
 are monomorphisms, and so the upper horizontal arrow is a monomorphism.
The fact that $A\otimes_{\min}\mmph$ is monic-preserving is well-known~\cite{Ped_pullback99}.

As a consequence, we may conclude that $\mathfrak{O}_{A}$ is product
separating and monic-preserving for every nuclear $C^{*}$-algebra
$A$.
\end{example}

\subsection{\label{subsec:homnattr}Homotopy category of natural transformations}

We
are now ready to introduce homotopies of natural
transformations between good endofunctors. To this end, we shall need
the \emph{cylinder functor} $I\coloneqq\C_{[0,1]}$ (see Definition~\ref{def:Cb})
which carries the structure of a good endofunctor (see Example~\ref{exa:good-endfs}).
\begin{defn}
We call two labeled natural transformations $\gamma_{0},\gamma_{1}\in\GEFC(F,G)$
\emph{homotopic} (written $\gamma_{0}\simeq\gamma_{1}$) if there
is $\Gamma\in\GEFC(F,GI)$ such that the following diagram commutes
for~$j=0,1$:
\[
\input{diagrams/d224.tikz}
\]
It will be shown in the next proposition that homotopy is an equivalence
relation, and we shall denote by $[\gamma]$ the homotopy class of
$\gamma\in\GEFC$.
\end{defn}

\begin{thm}\label{thm:homot-eqrel}
Homotopy
is an equivalence
relation on morphisms of $\GEFC$.
\end{thm}

\begin{proof}
Reflexivity and symmetry are obvious, so we only need to prove transitivity.
To this end, denote by $r_{0}\colon I\Rightarrow I$ and $r_{1}\colon I\Rightarrow I$
the natural transformations induced respectively by the functions
$[0,1]\to[0,1]\colon t\mapsto t/2$ and $[0,1]\to[0,1]\colon t\mapsto1/2+t/2$
(roughly speaking, $r_{0}$ and $r_{1}$ are the restrictions to the
left and right halves of the unit interval). Consider the following
diagrams:
\begin{alignat}{2}
\input{diagrams/d271.tikz} & \qquad & \input{diagrams/d272.tikz}\label{eq:homot-eqrel}
\end{alignat}
The left-hand diagram is clearly a pullback in $\DEFC$. Bearing in
mind that $\ev_{0}$ is a componentwise split epimorphism, we conclude
from Lemma~\ref{lem:pullback-functors-and-c-algebras} that the right-hand
diagram is also a pullback in $\DEFC$. Furthermore, $G\ev_{0}$,
$G\ev_{1}$, $Gr_{0}$, and $Gr_{1}$ are labeled, and hence, by Proposition~\ref{prop:underlying-canon}(\ref{enu:undercanon2}),
the right-hand diagram is a pullback in $\GEFC$.

Let now $\alpha_{0},\alpha_{1},\alpha_{2}\in\GEFC(F,G)$, and let
$\alpha_{01},\alpha_{12}\in\GEFC(F,IG)$ be some homotopies connecting
$\alpha_{0}$ with $\alpha_{1}$ and $\alpha_{1}$ with $\alpha_{2}$.
The left-hand diagram in
\begin{alignat*}{2}
\input{diagrams/d194.tikz} & \qquad & \input{diagrams/d195.tikz}
\end{alignat*}
gives rise to the right-hand one by the universal property of the
pullback $GI$. The dashed arrow is a homotopy connecting $\alpha_{0}$
with $\alpha_{2}$.
\end{proof}

\begin{prop}\label{prop:homot-bimon-congr}
Homotopy is a bimonoidal congruence
relation on $\GEFC$.
\end{prop}

\begin{proof}
The commutative diagrams
\begin{alignat*}{2}
\input{diagrams/d196.tikz} & \qquad\qquad & \input{diagrams/d199.tikz}
\end{alignat*}
for $j=0,1$, show that $\beta_{0}\circ\alpha\simeq\beta_{1}\circ\alpha$
whenever $\beta_{0}\simeq\beta_{1}$, and $\beta\circ\alpha_{0}\simeq\beta\circ\alpha_{1}$
whenever $\alpha_{0}\simeq\alpha_{1}$.
Now let $\beta\in\GEFC(G,G')$, $\alpha_{0},\alpha_{1}\in\GEFC(F,F')$,
and let $\alpha\in\GEFC(F,IF')$ be a homotopy connecting $\alpha_{0}$
with $\alpha_{1}$, and consider the diagrams
\begin{alignat*}{2}
\phantom{j=1,2} & \input{diagrams/d222.tikz} & \qquad\qquad & \input{diagrams/d223.tikz}
\end{alignat*}
for $j=0,1$, which obviously commute (the upper triangle in the right-hand
side diagram commutes by Corollary~\ref{cor:332}). These two diagrams
show that $\beta\alpha_{0}\simeq\beta\alpha_{1}$ whenever $\alpha_{0}\simeq\alpha_{1}$,
and $\beta_{0}\alpha\simeq\beta_{1}\alpha$ whenever $\beta_{0}\simeq\beta_{1}$.

Finally, let $\alpha_{0},\alpha_{1}\in\in\GEFC(F,F')$ and $\beta_{0},\beta_{1}\in\in\GEFC(G,G')$,
and let $\alpha$ and $\beta$ be the homotopies connecting $\alpha_{0}$
with $\alpha_{1}$ and $\beta_{0}$ with $\beta_{1}$. We need to
check that $\alpha_{0}\oplus\beta_{0}\simeq\alpha_{1}\oplus\beta_{1}$.
The left triangle in
\[
\input{diagrams/d288.tikz}
\]
commutes since $\oplus$ is functorial; the right triangle commutes
by naturality
of $\rc$. Thus, $\rc^{-1}\circ(\alpha\oplus\beta)$
is the required homotopy.
\end{proof}

The previous proposition allows us to equip the category of good endofunctors
and homotopy classes of labeled natural transformations with a bimonoidal
structure (see Definition~\ref{def:aj}).

\begin{defn}
Introduce the bimonoidal category $\hGEFC\coloneqq\GEFC\modrel{\simeq}$
which we call the \emph{homotopy category of natural transformations}.
It follows from Remark~\ref{rem:quot-tight-strict} that $\hGEFC$
is tight, and its multiplicative structure is strict.
\end{defn}

\begin{lem}\label{lem:br-homot}
Let
$\alpha_{0},\alpha_{1}\in\Br(\GEFC_{\ttt})(T,S)$.
Then $\u(\alpha_{0})\simeq\u(\alpha_{1})$ if and only if $\u\overline{\EV}(\alpha_{0})\simeq\u\overline{\EV}(\alpha_{1})$.
\end{lem}

\begin{proof}
For brevity, let us write $I$ and $I\mathbb{C}$ respectively for
the bracketed objects $\i(I)\in\in\Br(\GEFC_{\ttt})$ and $\i(I\mathbb{C})\in\in\Br(\Cstar)$.
For some $\alpha\in\Br(\GEFC_{\ttt})(T,SI)$ let us consider the diagrams
\begin{alignat*}{2}
\input{diagrams/d402.tikz} & \qquad\qquad & \input{diagrams/d403.tikz}
\end{alignat*}
respectively in $\Br(\GEFC_{\ttt})$ and $\Br(\Cstar)$. Bearing in
mind that $\overline{\EV}$ is a fully faithful strict bimonoidal
functor (see Corollary~\ref{cor:EV-GEFC}), and that $\overline{\EV}(\ev_{j})=\ev_{j}\mathbb{C}$,
we have the following chain of equivalences:
\begin{alignat*}{3}
\u(\alpha_{0})\simeq\u(\alpha_{1})~ & \Longleftrightarrow~ & \textrm{the left-hand diagram above commutes for some }\alpha\\
 & \Longleftrightarrow & \textrm{the right-hand diagram above commutes for some }\alpha & ~\Longleftrightarrow~ & \u\overline{\EV}(\alpha_{0})\simeq\u\overline{\EV}(\alpha_{1}) & .
\end{alignat*}
\end{proof}

The previous lemma yields a few corollaries, the first of which is
the construction of the homotopy versions of $\mathfrak{O}$ and $\EV$.
\begin{lem}\label{lem:O-pres-homot}
The functors $\mathfrak{O}\colon\Cstar\rightleftarrows\GEFC_{\ttt}\noloc\EV$
preserve homotopies, i.e,
\begin{enumerate}
\item\label{enu:pres-homot-O}if $\varphi_{0},\varphi_{1}\in\Cstar(A,B)$,
then $\varphi_{0}\simeq\varphi_{1}$ if and only if $\mathfrak{O}_{\varphi_{0}}\simeq\mathfrak{O}_{\varphi_{1}}$;
\item\label{enu:pres-homot-EV}if $\alpha_{0},\alpha_{1}\in\DEFC_{\ttt}(T,S)$,
then $\alpha_{0}\simeq\alpha_{1}$ if and only if $\alpha_{0}\mathbb{C}\simeq\alpha_{1}\mathbb{C}$.
\end{enumerate}
\end{lem}

\begin{proof}
It follows from Lemma~\ref{lem:br-homot} and the definition of $\overline{\EV}$
(see Theorem~\ref{thm:defFbar}) that
\begin{alignat*}{2}
 &\ 
&\ 
& \mathfrak{O}_{\varphi_{0}}\simeq\mathfrak{O}_{\varphi_{1}}\Leftrightarrow\u\i\mathfrak{O}_{\varphi_{0}}\simeq\u\i\mathfrak{O}_{\varphi_{1}}\Leftrightarrow\u\overline{\EV}\i\mathfrak{O}_{\varphi_{0}}\simeq\u\overline{\EV}\i\mathfrak{O}_{\varphi_{1}}\Leftrightarrow\varphi_{0}\otimes\mathbb{C}\simeq\varphi_{1}\otimes\mathbb{C}\Leftrightarrow\varphi_{0}\simeq\varphi_{1},\\
\textrm{and} & \qquad &\ 
& \alpha_{0}\simeq\alpha_{1}\Leftrightarrow\u\i(\alpha_{0})\simeq\u\i(\alpha_{1})\Leftrightarrow\u\overline{\EV}\i(\alpha_{0})\simeq\u\overline{\EV}\i(\alpha_{1})\Leftrightarrow\EV(\alpha_{0})\simeq\EV(\alpha_{0}).
\end{alignat*}

\end{proof}

\begin{thm}\label{thm:OEequivs}
There
exist bimonoidal functors
$\mathfrak{O}_{\simeq}\colon\hCstar\to\hGEFC_{\ttt}$ and $\EV_{\simeq}\colon\hGEFC_{\ttt}\to\hCstar$
rendering commutative the diagram
\begin{equation}
\input{diagrams/d578.tikz}\label{eq:d578}
\end{equation}
in which the vertical arrows are the quotients, and both rows yield
bimonoidal equivalences of $\Cstar$ and $\GEFC_{\ttt}$, and $\hCstar$
and $\hGEFC_{\ttt}$.
\end{thm}

\begin{proof}
It follows from Lemmas~\ref{lem:ak} and~\ref{lem:O-pres-homot}
that $\mathfrak{O}_{\simeq}$ and $\EV_{\simeq}$ exist and are defined
by the formulas
\begin{alignat}{4}
 &\ 
& \mathfrak{O}_{\simeq}[\varphi] & \coloneqq[\mathfrak{O}\varphi], & \qquad\qquad &\ 
& \EV_{\simeq}[\alpha] & \coloneqq[\EV\alpha],\label{eq:oeeq1}\\
 &\ 
& \mathfrak{O}^{2}_{\simeq} & \coloneqq[\mathfrak{O}^{2}], &\ 
&\ 
& \EV^{2}_{\simeq} & \coloneqq[\EV^{2}],\label{eq:oeeq2}\\
 &\ 
& \mathfrak{O}^{0}_{\simeq} & \coloneqq[\mathfrak{O}^{0}], &\ 
&\ 
& \EV^{0}_{\simeq} & \coloneqq[\EV^{0}],\label{eq:oeeq3}
\end{alignat}
The upper row in~(\ref{eq:d578}) is a bimonoidal equivalence by
Theorem~\ref{thm:strictif-defc} and the remark after Definition~\ref{def:GEFCtt}.
Using this fact together with~(\ref{eq:oeeq1})--(\ref{eq:oeeq3}),
one can show that the lower row is also a bimonoidal equivalence.
\end{proof}

The second corollary yields a strictified version of the functor $\EV_{\simeq}$.
\begin{lem}\label{lem:EV-hGEFC}
There
is a fully faithful strict
bimonoidal functor
\[
\overline{\EV}_{\simeq}\colon\Br(\hGEFC_{\ttt})\to\Br(\hCstar).
\]
\end{lem}

\begin{proof}
It follows from Proposition~\ref{prop:brquot} that $\Br(\GEFC)$
and $\Br(\Cstar)$ inherit
the homotopy congruences
(both denoted by~$\simeq$) from the underlying categories $\GEFC$
and $\Cstar$, and that the equalities $\Br(\hGEFC)=\Br(\GEFC)\modrel{\simeq}$
and $\Br(\hCstar)=\Br(\Cstar)\modrel{\simeq}$ hold. Using Lemmas~\ref{lem:ak}
and~\ref{lem:br-homot} we obtain the dashed arrow in
\[
\input{diagrams/d436.tikz}
\]
Bearing in mind that the top-right leg of the diagram above is a full
functor, we conclude from Lemma~\ref{lem:ak} that $\overline{\EV}_{\simeq}$
is also full. That $\overline{\EV}_{\simeq}$ is faithful readily
follows from Lemma~\ref{lem:br-homot}.
\end{proof}

The last corollary of Lemma~\ref{lem:br-homot} characterizes the
homotopy class of the natural transformation $\kappa^{\fK,\fK}$.

\begin{cor}\label{cor:ai}
$\kappa^{\fK,\fK}\simeq\id_{\fK^{2}}$.
\end{cor}

\begin{proof}
This follows from Corollary~\ref{cor:ae}, Lemma~\ref{lem:af},
and Lemma~\ref{lem:br-homot}.
\end{proof}

\subsection{\label{subsec:starig}The stabilization rig}

Recall the following $*$-homomorphisms from Section~\ref{subsec:rig-comp}:
\begin{alignat*}{2}
 &\ 
& \widehat{\iota}\colon & \mathbb{C}\oplus\mathbb{C}\to\oM_{2},\hspc 4\phantom{\mathbb{K}\mathbb{K}\mathbb{K}\widehat{\theta}_{n}}\quad\widehat{\iota}_{00}\colon\mathbb{C}\to\oK,\\
 &\ 
& \widehat{\theta}\colon & \oK\otimes\oK\to\oK,\hspc 4\phantom{\mathbb{C}\mathbb{C}\mathbb{M}_{2}\widehat{\iota}_{00}}\quad\widehat{\theta}_{n}\colon\oM_{2}\otimes\oK\to\oK,\\
 &\ 
& \widehat{\mu}\colon & \oK\oplus\oK\xrightarrow{\lambda^{-1}\oplus\lambda^{-1}}(\mathbb{C}\otimes\oK)\oplus(\mathbb{C}\otimes\oK)\xrightarrow{(\delta^{r})^{-1}}(\mathbb{C}\otimes\mathbb{C})\otimes\oK\xrightarrow{\widehat{\iota}\otimes\oK}\oM_{2}\otimes\oK\xrightarrow{\widehat{\theta}_{2}}\oK.
\end{alignat*}
Let us consider them as arrows in $\Br(\hCstar)$, endowing the corresponding
objects with the obvious bracketing. By Corollary~\ref{lem:EV-hGEFC},
we can find their $\overline{\EV}_{\simeq}$-preimages in $\Br(\hGEFC_{\ttt})$:
\begin{alignat*}{2}
 & \iota\colon\Id\oplus\Id\Rightarrow\fM_{2},\hspc 3\phantom{\fK\fK\fK\theta_{n}}\quad\iota_{00}\colon\Id\Rightarrow\fK,\\
 & \theta\colon\fK\fK\Rightarrow\fK,\hspc 3\phantom{\Id\oplus\Id\fM_{2}\iota_{00}}\quad\theta_{n}\colon\fM_{n}\fK\Rightarrow\fK,\\
 & \mu\colon\fK\oplus\fK=\Id\fK\oplus\Id\fK\xRightarrow{\rc^{-1}}(\Id\oplus\Id)\fK\xRightarrow{\iota\fK}\fM_{2}\fK\xRightarrow{\theta_{2}}\fK. & \hspc{12}
\end{alignat*}
It is straightforward to check that the components of $\iota$ and
$\iota_{00}$ at $B\in\in\Cstar$ are given by
\[
\iota B\colon(b_{1},b_{2})\mapsto\left(\begin{array}{cc}
b_{1} & 0\\
0 & b_{2}
\end{array}\right),\qquad\qquad\iota_{00}B\colon b\mapsto\left(\begin{array}{ccc}
b & 0 & \ldots\\
0 & 0 & \ldots\\
\vdots & \vdots & \ddots
\end{array}\right).
\]

\begin{prop}\label{prop:rig-K}
$(\fK,[\theta],[\iota_{00}],[\mu],[0])$ is a rig
in $\hGEFC$.
\end{prop}

\begin{proof}
We need to check the commutativity of nine diagrams in $\hGEFC$ (see
Definition~\ref{def:rig}). Note that it suffices to do it regarding
them as diagrams in $\Br(\hGEFC_{\ttt})$. Applying the strict bimonoidal
functor $\overline{\EV}_{\simeq}$ (see Corollary~\ref{lem:EV-hGEFC})
to these diagrams yields exactly the commutative diagrams witnessing
that $(\oK,[\widehat{\theta}],[\widehat{\iota}_{00}],[\widehat{\mu}],[0])$
is a rig in $\hCstar$. Thus, the original diagrams are also commutative
since $\overline{\EV}_{\simeq}$ is faithful by Corollary~\ref{lem:EV-hGEFC}.
\end{proof}
The rig structure of $\fK$ will be necessary for
the definition of generalized morphisms (see Section~\ref{sec:genhom}).

%% file: diagrams/d245.tikz
\begin{tikzcd}[ampersand replacement=\&]
{\mathfrak{O}_{A}F} \&  \& {F\mathfrak{O}_{A}} \\
{\mathfrak{O}_{A'}F} \&  \& {F\mathfrak{O}_{A'}\mmdc}
\arrow["\kappa^{A,F}", Rightarrow, from=1-1, to=1-3]
\arrow["\mathfrak{O}_{\varphi}F"', Rightarrow, from=1-1, to=2-1]
\arrow["F\mathfrak{O}_{\varphi}", Rightarrow, from=1-3, to=2-3]
\arrow["\kappa^{A',F}", Rightarrow, from=2-1, to=2-3]
\end{tikzcd}

%% file: diagrams/d246.tikz
\begin{tikzcd}[column sep=large, ampersand replacement=\&]  %d246
{\mathfrak{O}_{A}\mathfrak{O}_{A'}F} \& {\mathfrak{O}_{A}F\mathfrak{O}_{A'}} \& {F\mathfrak{O}_{A}\mathfrak{O}_{A'}} \\
{\mathfrak{O}_{A\otimes A'}F} \&  \& {F\mathfrak{O}_{A\otimes A'}\mmdc}
\arrow["\mathfrak{O}_{A}\kappa^{A',F}", Rightarrow, from=1-1, to=1-2]
\arrow["\mathfrak{O}_{A,A'}^{2}F"', Rightarrow, from=1-1, to=2-1]
\arrow["\kappa^{A,F}\mathfrak{O}_{A'}", Rightarrow, from=1-2, to=1-3]
\arrow["F\mathfrak{O}_{A,A'}^{2}", Rightarrow, from=1-3, to=2-3]
\arrow["\kappa^{A\otimes A',F}", Rightarrow, from=2-1, to=2-3]
\end{tikzcd}

%% file: diagrams/d247.tikz
\begin{tikzcd}[column sep=small, ampersand replacement=\&]  %d247
{\Id F} \& {F} \& {F\Id} \\
{\mathfrak{O}_{\mathbb{C}}F} \&  \& {F\mathfrak{O}_{\mathbb{C}}\mmdd}
\arrow["\mathfrak{O}^{0}F"', Rightarrow, from=1-1, to=2-1]
\arrow[" ="', Rightarrow, from=1-2, to=1-1]
\arrow[" =", Rightarrow, from=1-2, to=1-3]
\arrow["F\mathfrak{O}^{0}", Rightarrow, from=1-3, to=2-3]
\arrow["\kappa^{\mathbb{C},F}", Rightarrow, from=2-1, to=2-3]
\end{tikzcd}

%% file: diagrams/d248.tikz
\begin{tikzcd}[ampersand replacement=\&]
{\mathfrak{O}_{A}F} \&  \& {F\mathfrak{O}_{A}} \\
{\mathfrak{O}_{A}G} \&  \& {G\mathfrak{O}_{A}\mmdd}
\arrow["\kappa^{A,F}", Rightarrow, from=1-1, to=1-3]
\arrow["\mathfrak{O}_{A}\alpha"', Rightarrow, from=1-1, to=2-1]
\arrow["\alpha\mathfrak{O}_{A}", Rightarrow, from=1-3, to=2-3]
\arrow["\kappa^{A,G}"', Rightarrow, from=2-1, to=2-3]
\end{tikzcd}

%% file: diagrams/d268.tikz
\begin{tikzcd}[ampersand replacement=\&]
 \&  \& {\mathfrak{O}_{A'}F_{\infty}} \&  \&  \& {\mathfrak{O}_{A'}F_{j}} \\
{\mathfrak{O}_{A}F_{\infty}} \&  \&  \& {\mathfrak{O}_{A}F_{j}} \\
 \&  \& {F_{\infty}\mathfrak{O}_{A'}} \&  \&  \& {F_{j}\mathfrak{O}_{A'}} \\
{F_{\infty}\mathfrak{O}_{A}} \&  \&  \& {F_{j}\mathfrak{O}_{A}}
\arrow["\mathfrak{O}_{A'}\pr_{j}", Rightarrow, from=1-3, to=1-6]
\arrow["\kappa^{A',F_{\infty}}"{pos=0.7}, Rightarrow, dashed, from=1-3, to=3-3]
\arrow["\kappa^{A',F_{j}}", Rightarrow, from=1-6, to=3-6]
\arrow["\mathfrak{O}_{\varphi}F_{\infty}", Rightarrow, from=2-1, to=1-3]
\arrow["\kappa^{A,F_{\infty}}"', Rightarrow, dashed, from=2-1, to=4-1]
\arrow["\mathfrak{O}_{\varphi}F_{j}"{pos=0.3}, Rightarrow, from=2-4, to=1-6]
\arrow["\pr_{j}\mathfrak{O}_{A'}"{pos=0.7}, Rightarrow, from=3-3, to=3-6]
\arrow["F_{\infty}\mathfrak{O}_{\varphi}"{pos=0.4}, Rightarrow, from=4-1, to=3-3]
\arrow["\pr_{j}\mathfrak{O}_{A}"', Rightarrow, from=4-1, to=4-4]
\arrow["F_{j}\mathfrak{O}_{\varphi}"', Rightarrow, from=4-4, to=3-6]
\arrow["\mathfrak{O}_{A}\pr_{j}"{pos=0.4}, crossing over, Rightarrow, from=2-1, to=2-4]
\arrow["\kappa^{A,F_{j}}"{pos=0.2}, crossing over, Rightarrow, from=2-4, to=4-4]
\end{tikzcd}

%% file: diagrams/d270.tikz
\begin{tikzcd}[ampersand replacement=\&]
{G} \&  \& {F_{\infty}} \\
 \&  \& {F_{j}\mmdd}
\arrow["\tau_{\infty}", Rightarrow, dashed, from=1-1, to=1-3]
\arrow["\tau_{j}"', Rightarrow, from=1-1, to=2-3]
\arrow["\pr_{j}"{pos=0.4}, Rightarrow, from=1-3, to=2-3]
\end{tikzcd}

%% file: diagrams/d269.tikz
\begin{tikzcd}[ampersand replacement=\&]
 \&  \&  \& {G\mathfrak{O}_{A}} \&  \&  \& {F_{\infty}\mathfrak{O}_{A}} \\
{\mathfrak{O}_{A}G} \&  \&  \& {\mathfrak{O}_{A}F_{\infty}} \\
 \&  \&  \&  \&  \&  \& {F_{j}\mathfrak{O}_{A}} \\
 \&  \&  \& {\mathfrak{O}_{A}F_{j}}
\arrow["\mathfrak{O}_{A}\tau_{\infty}", Rightarrow, dashed, from=1-4, to=1-7]
\arrow["\tau_{j}\mathfrak{O}_{A}"', curve={height=12pt}, Rightarrow, from=1-4, to=3-7]
\arrow["\pr_{j}\mathfrak{O}_{A}", Rightarrow, from=1-7, to=3-7]
\arrow["\kappa^{A,G}", Rightarrow, from=2-1, to=1-4]
\arrow["\tau_{\infty}\mathfrak{O}_{A}", Rightarrow, dashed, from=2-1, to=2-4]
\arrow["\mathfrak{O}_{A}\tau_{j}"', curve={height=12pt}, Rightarrow, from=2-1, to=4-4]
\arrow["\mathfrak{O}_{A}\pr_{j}"{pos=0.4}, Rightarrow, from=2-4, to=4-4]
\arrow["\kappa^{A,F_{j}}"', Rightarrow, from=4-4, to=3-7]
\arrow["\kappa^{A,F_{\infty}}"'{pos=0.6}, crossing over, Rightarrow, from=2-4, to=1-7]
\end{tikzcd}

%% file: diagrams/d491.tikz
\begin{tikzcd}[ampersand replacement=\&]
{\mathfrak{O}_{B\otimes A}F} \&  \&  \&  \&  \&  \& {F\mathfrak{O}_{B\otimes A}} \\
 \&  \& {\mathfrak{O}_{B}\mathfrak{O}_{A}F} \&  \& {\mathfrak{O}_{B}F\mathfrak{O}_{A}} \&  \& {F\mathfrak{O}_{B}\mathfrak{O}_{A}} \\
 \\
{\mathfrak{O}_{A}\mathfrak{O}_{B}F} \&  \& {\mathfrak{O}_{A}F\mathfrak{O}_{B}} \&  \& {F\mathfrak{O}_{A}\mathfrak{O}_{B}} \\
{\mathfrak{O}_{A\otimes B}F} \&  \&  \&  \&  \&  \& {F\mathfrak{O}_{A\otimes B}}
\arrow["\kappa^{B\otimes A,F}", Rightarrow, from=1-1, to=1-7]
\arrow["\mathfrak{O}_{\xi^{\otimes}}F"{pos=0.15}, sqarw=1.5em, Rightarrow, from=1-1, to=5-1]
\arrow["F\mathfrak{O}_{\xi^{\otimes}}"'{pos=0.35}, sqare=1.5em, Rightarrow, from=1-7, to=5-7]
\arrow["\mathfrak{O}^{2}F"', Rightarrow, from=2-3, to=1-1]
\arrow["\cong", draw=none, from=2-3, to=1-1]
\arrow["\mathfrak{O}_{B}\kappa^{A,F}", Rightarrow, from=2-3, to=2-5]
\arrow["\kappa^{B,\mathfrak{O}_{A}}F"', Rightarrow, from=2-3, to=4-1]
\arrow["\kappa^{B,\mathfrak{O}_{A}F}", Rightarrow, from=2-3, to=4-3]
\arrow["\kappa^{B,F}\mathfrak{O}_{A}", Rightarrow, from=2-5, to=2-7]
\arrow["\kappa^{B,F\mathfrak{O}_{A}}"', Rightarrow, from=2-5, to=4-5]
\arrow["F\mathfrak{O}^{2}", Rightarrow, from=2-7, to=1-7]
\arrow["F\kappa^{B,\mathfrak{O}_{A}}", Rightarrow, from=2-7, to=4-5]
\arrow["\mathfrak{O}_{A}\kappa^{B,F}", Rightarrow, from=4-1, to=4-3]
\arrow["\kappa^{A,F}\mathfrak{O}_{B}", Rightarrow, from=4-3, to=4-5]
\arrow["(\mathfrak{O}^{2})^{-1}F"', Rightarrow, from=5-1, to=4-1]
\arrow["\kappa^{A\otimes B,F}", Rightarrow, from=5-1, to=5-7]
\arrow["F(\mathfrak{O}^{2})^{-1}", Rightarrow, from=5-7, to=4-5]
\arrow["\cong"', draw=none, from=5-7, to=4-5]
\end{tikzcd}

%% file: diagrams/d467.tikz
\begin{tikzcd}[ampersand replacement=\&]
{F} \& {G} \\
 \& {H}
\arrow["\gamma", Rightarrow, from=1-1, to=1-2]
\arrow["\beta"', Rightarrow, from=1-1, to=2-2]
\arrow["\alpha", Rightarrow, from=1-2, to=2-2]
\end{tikzcd}

%% file: diagrams/d468.tikz
\begin{tikzcd}[row sep=scriptsize, ampersand replacement=\&]  %d468
 \&  \&  \& {FA} \&  \&  \& {GA} \\
{AF} \&  \&  \& {AG} \\
 \&  \&  \&  \&  \&  \& {HA} \\
 \&  \&  \& {AH}
\arrow["\gamma\mathfrak{O}_{A}", Rightarrow, from=1-4, to=1-7]
\arrow["\beta\mathfrak{O}_{A}"', curve={height=12pt}, Rightarrow, from=1-4, to=3-7]
\arrow["\alpha\mathfrak{O}_{A}", Rightarrow, from=1-7, to=3-7]
\arrow["\kappa^{A,F}", Rightarrow, from=2-1, to=1-4]
\arrow["\mathfrak{O}_{A}\gamma"{pos=0.6}, Rightarrow, from=2-1, to=2-4]
\arrow["\mathfrak{O}_{A}\beta"', curve={height=12pt}, Rightarrow, from=2-1, to=4-4]
\arrow["\mathfrak{O}_{A}\alpha"'{pos=0.4}, Rightarrow, from=2-4, to=4-4]
\arrow["\kappa^{A,H}"', Rightarrow, from=4-4, to=3-7]
\arrow["\kappa^{A,G}"{pos=0.4}, crossing over, Rightarrow, from=2-4, to=1-7]
\end{tikzcd}

%% file: diagrams/d469.tikz
\begin{tikzcd}[ampersand replacement=\&]
{0} \& {F_{0}} \& {F} \& {F_{1}} \& {0} \\
 \&  \& {G}
\arrow[Rightarrow, from=1-1, to=1-2]
\arrow["\iota", Rightarrow, from=1-2, to=1-3]
\arrow[" 0"', curve={height=6pt}, Rightarrow, from=1-2, to=2-3]
\arrow["\pi", Rightarrow, from=1-3, to=1-4]
\arrow["\alpha", Rightarrow, from=1-3, to=2-3]
\arrow[Rightarrow, from=1-4, to=1-5]
\arrow["\beta", curve={height=-6pt}, Rightarrow, dashed, from=1-4, to=2-3]
\end{tikzcd}

%% file: diagrams/d475.tikz
\begin{tikzcd}[column sep=tiny, ampersand replacement=\&]  %d475
{\mathfrak{O}_{A_{1}}F\oplus\mathfrak{O}_{A_{2}}F} \&  \&  \&  \&  \&  \& {F\mathfrak{O}_{A_{1}}\oplus F\mathfrak{O}_{A_{2}}} \\
{(\mathfrak{O}_{A_{1}}\oplus\mathfrak{O}_{A_{2}})F} \& {\hspc 4} \& {\mathfrak{O}_{A_{j}}F} \& {\hspc{2.5}} \& {F\mathfrak{O}_{A_{j}}} \& {\hspc 4} \& {F(\mathfrak{O}_{A_{1}}\oplus\mathfrak{O}_{A_{2}})} \\
{\mathfrak{O}_{A_{1}\oplus A_{2}}F} \&  \&  \&  \&  \&  \& {F\mathfrak{O}_{A_{1}\oplus A_{2}}\mmdd}
\arrow["\kappa^{A_{1},F}\oplus\kappa^{A_{2},F}", Rightarrow, from=1-1, to=1-7]
\arrow["\pr_{j}", Rightarrow, from=1-1, to=2-3]
\arrow["\pr_{j}"', Rightarrow, from=1-7, to=2-5]
\arrow["\rc", Rightarrow, from=2-1, to=1-1]
\arrow["\pr_{j}F"{pos=0.3}, Rightarrow, from=2-1, to=2-3]
\arrow["\kappa^{A_{j},F}", Rightarrow, from=2-3, to=2-5]
\arrow["\lc"', Rightarrow, from=2-7, to=1-7]
\arrow["F\pr_{j}"'{pos=0.3}, Rightarrow, from=2-7, to=2-5]
\arrow["(\mathfrak{O}_{\oplus}^{2})^{-1}F", Rightarrow, from=3-1, to=2-1]
\arrow["\mathfrak{O}_{\pr_{j}}F"', Rightarrow, from=3-1, to=2-3]
\arrow["\kappa^{A_{1}\oplus A_{2},F}"', Rightarrow, from=3-1, to=3-7]
\arrow["F\mathfrak{O}_{\pr_{j}}", Rightarrow, from=3-7, to=2-5]
\arrow["F(\mathfrak{O}_{\oplus}^{2})^{-1}"', Rightarrow, from=3-7, to=2-7]
\end{tikzcd}

%% file: diagrams/d171.tikz
\begin{tikzcd}[ampersand replacement=\&]
{0} \& {FA} \& {GA} \& {(G/F)A} \& {0} \\
{0} \& {FB} \& {GB} \& {(G/F)B} \& {0.}
\arrow[from=1-1, to=1-2]
\arrow["\iota A", from=1-2, to=1-3]
\arrow["F\varphi", from=1-2, to=2-2]
\arrow[from=1-3, to=1-4]
\arrow["G\varphi", from=1-3, to=2-3]
\arrow[from=1-4, to=1-5]
\arrow[dashed, from=1-4, to=2-4]
\arrow[from=2-1, to=2-2]
\arrow["\iota B", from=2-2, to=2-3]
\arrow[from=2-3, to=2-4]
\arrow[from=2-4, to=2-5]
\end{tikzcd}

%% file: diagrams/d192.tikz
\begin{tikzcd}[ampersand replacement=\&]
{P} \& {A_{1}} \\
{A_{2}} \& {A_{0}\mmdc}
\arrow["p_{1}", , from=1-1, to=1-2]
\arrow["p_{2}"', , from=1-1, to=2-1]
\arrow["\varphi_{1}", , from=1-2, to=2-2]
\arrow["\varphi_{2}"', , from=2-1, to=2-2]
\end{tikzcd}

%% file: diagrams/d201.tikz
\begin{tikzcd}[ampersand replacement=\&]
{P^{k}} \& {A_{1}^{k}} \\
{A_{2}^{k}} \& {A_{0}^{k}}
\arrow["p_{1}^{k}", , from=1-1, to=1-2]
\arrow["p_{2}^{k}"', , from=1-1, to=2-1]
\arrow["\varphi_{1}^{k}", , from=1-2, to=2-2]
\arrow["\varphi_{2}^{k}"', , from=2-1, to=2-2]
\end{tikzcd}

%% file: diagrams/d172.tikz
\begin{tikzcd}[ampersand replacement=\&]
{B_{1}\oplus_{B_{0}}B_{2}} \\
 \&  \& {B_{2}} \\
 \& {B_{1}} \& {B_{0}\mmdc}
\arrow["p_{2}", curve={height=-12pt}, from=1-1, to=2-3]
\arrow["p_{1}"', curve={height=12pt}, from=1-1, to=3-2]
\arrow["\varphi_{2}", two heads, from=2-3, to=3-3]
\arrow["\varphi_{1}", from=3-2, to=3-3]
\end{tikzcd}

%% file: diagrams/d168.tikz
\begin{tikzcd}[ampersand replacement=\&]
{\Omega(B_{1}\oplus_{B_{0}}B_{2})} \\
 \& {\Omega B_{1}\oplus_{\Omega B_{0}}\Omega B_{2}} \& {\Omega B_{2}} \\
 \& {\Omega B_{1}} \& {\Omega B_{0}}
\arrow[dashed, from=1-1, to=2-2]
\arrow["\Omega p_{2}", curve={height=-12pt}, from=1-1, to=2-3]
\arrow["\Omega p_{1}"', curve={height=12pt}, from=1-1, to=3-2]
\arrow[""', from=2-2, to=2-3]
\arrow[from=2-2, to=3-2]
\arrow["\Omega\varphi_{2}", two heads, from=2-3, to=3-3]
\arrow["\Omega\varphi_{1}"', from=3-2, to=3-3]
\end{tikzcd}

%% file: diagrams/d173.tikz
\begin{tikzcd}[ampersand replacement=\&]
{0} \& {F(B_{1}\underset{B_{0}}{\oplus}B_{2})} \& {G(B_{1}\underset{B_{0}}{\oplus}B_{2})} \& {H(B_{1}\underset{B_{0}}{\oplus}B_{2})} \& {0} \\
{0} \& {FB_{1}\underset{FB_{0}}{\oplus}FB_{2}} \& {GB_{1}\underset{GB_{0}}{\oplus}GB_{2}} \& {HB_{1}\underset{HB_{0}}{\oplus}HB_{2}} \& {0,}
\arrow[from=1-1, to=1-2]
\arrow[from=1-2, to=1-3]
\arrow["\lc", from=1-2, to=2-2]
\arrow["\cong"', draw=none, from=1-2, to=2-2]
\arrow[from=1-3, to=1-4]
\arrow["\lc", from=1-3, to=2-3]
\arrow["\cong"', draw=none, from=1-3, to=2-3]
\arrow[from=1-4, to=1-5]
\arrow["\lc", from=1-4, to=2-4]
\arrow[from=2-1, to=2-2]
\arrow[from=2-2, to=2-3]
\arrow[from=2-3, to=2-4]
\arrow[from=2-4, to=2-5]
\end{tikzcd}

%% file: diagrams/d476.tikz
\begin{tikzcd}[row sep=small, ampersand replacement=\&]  %d476
 \& {0} \& {0} \& {0} \\
{0} \& {FJ} \& {GJ} \& {HJ} \& {0} \\
{0} \& {FA} \& {GA} \& {HA} \& {0} \\
{0} \& {FB} \& {GB} \& {HB} \& {0} \\
 \& {0} \& {0} \& {0\mmdd}
\arrow[from=1-2, to=2-2]
\arrow[from=1-3, to=2-3]
\arrow[from=1-4, to=2-4]
\arrow[from=2-1, to=2-2]
\arrow[from=2-2, to=2-3]
\arrow[from=2-2, to=3-2]
\arrow[from=2-3, to=2-4]
\arrow[from=2-3, to=3-3]
\arrow[from=2-4, to=2-5]
\arrow[from=2-4, to=3-4]
\arrow[from=3-1, to=3-2]
\arrow[from=3-2, to=3-3]
\arrow[from=3-2, to=4-2]
\arrow[from=3-3, to=3-4]
\arrow[from=3-3, to=4-3]
\arrow[from=3-4, to=3-5]
\arrow[from=3-4, to=4-4]
\arrow[from=4-1, to=4-2]
\arrow[from=4-2, to=5-2]
\arrow[from=4-2, to=4-3]
\arrow[from=4-3, to=5-3]
\arrow[from=4-3, to=4-4]
\arrow[from=4-4, to=5-4]
\arrow[from=4-4, to=4-5]
\end{tikzcd}

%% file: diagrams/d260.tikz
\begin{tikzcd}[ampersand replacement=\&]
{0} \& {\mathfrak{O}_{A}FB} \& {\mathfrak{O}_{A}GB} \& {\mathfrak{O}_{A}HB} \& {0} \\
{0} \& {F\mathfrak{O}_{A}B} \& {G\mathfrak{O}_{A}B} \& {H\mathfrak{O}_{A}B} \& {0.}
\arrow[from=1-1, to=1-2]
\arrow[from=1-2, to=1-3]
\arrow["\kappa^{A,F}B", from=1-2, to=2-2]
\arrow[from=1-3, to=1-4]
\arrow["\kappa^{A,G}B", from=1-3, to=2-3]
\arrow[from=1-4, to=1-5]
\arrow["\kappa_{B}^{A,H}", dashed, from=1-4, to=2-4]
\arrow[from=2-1, to=2-2]
\arrow[from=2-2, to=2-3]
\arrow[from=2-3, to=2-4]
\arrow[from=2-4, to=2-5]
\end{tikzcd}

%% file: diagrams/d261.tikz
\begin{tikzcd}[ampersand replacement=\&]
 \& {\mathfrak{O}_{A}GB'} \&  \&  \& {\mathfrak{O}_{A}HB'} \\
{\mathfrak{O}_{A}GB} \&  \&  \& {\mathfrak{O}_{A}HB} \\
 \& {G\mathfrak{O}_{A}B'} \&  \&  \& {H\mathfrak{O}_{A}B'} \\
{G\mathfrak{O}_{A}B} \&  \&  \& {H\mathfrak{O}_{A}B}
\arrow["\mathfrak{O}_{A}qB'", from=1-2, to=1-5]
\arrow["\kappa^{A,G}B'"'{pos=0.7}, from=1-2, to=3-2]
\arrow["\kappa_{B'}^{A,H}"{pos=0.4}, from=1-5, to=3-5]
\arrow["\mathfrak{O}_{A}G\varphi"{pos=0.4}, from=2-1, to=1-2]
\arrow["\kappa^{A,G}B"'{pos=0.4}, from=2-1, to=4-1]
\arrow["\mathfrak{O}_{A}H\varphi"{pos=0.3}, from=2-4, to=1-5]
\arrow["q\mathfrak{O}_{A}B'"'{pos=0.2}, from=3-2, to=3-5]
\arrow["G\mathfrak{O}_{A}\varphi", from=4-1, to=3-2]
\arrow["q\mathfrak{O}_{A}B"', from=4-1, to=4-4]
\arrow["H\mathfrak{O}_{A}\varphi"', from=4-4, to=3-5]
\arrow["\mathfrak{O}_{A}qB"{pos=0.7}, crossing over, from=2-1, to=2-4]
\arrow["\kappa_{B}^{A,H}"{pos=0.3}, crossing over, from=2-4, to=4-4]
\end{tikzcd}

%% file: diagrams/d262.tikz
\begin{tikzcd}[ampersand replacement=\&]
 \& {\mathfrak{O}_{A}H} \&  \&  \& {H\mathfrak{O}_{A}} \\
{\mathfrak{O}_{A}G} \&  \&  \& {G\mathfrak{O}_{A}} \\
 \& {\mathfrak{O}_{A'}H} \&  \&  \& {H\mathfrak{O}_{A'}} \\
{\mathfrak{O}_{A'}G} \&  \&  \& {G\mathfrak{O}_{A'}}
\arrow["\kappa^{A,H}", Rightarrow, from=1-2, to=1-5]
\arrow["\mathfrak{O}_{\varphi}H"'{pos=0.7}, Rightarrow, from=1-2, to=3-2]
\arrow["H\mathfrak{O}_{\varphi}", Rightarrow, from=1-5, to=3-5]
\arrow["\mathfrak{O}_{A}q", Rightarrow, from=2-1, to=1-2]
\arrow["\mathfrak{O}_{\varphi}G"', Rightarrow, from=2-1, to=4-1]
\arrow["q\mathfrak{O}_{A}", Rightarrow, from=2-4, to=1-5]
\arrow["\kappa^{A',H}"{pos=0.3}, Rightarrow, from=3-2, to=3-5]
\arrow["\mathfrak{O}_{A'}q", Rightarrow, from=4-1, to=3-2]
\arrow["\kappa^{A',G}"', Rightarrow, from=4-1, to=4-4]
\arrow["q\mathfrak{O}_{A'}"', Rightarrow, from=4-4, to=3-5]
\arrow["\kappa^{A,G}"{pos=0.6}, crossing over, Rightarrow, from=2-1, to=2-4]
\arrow["G\mathfrak{O}_{\varphi}"{pos=0.3}, crossing over, Rightarrow, from=2-4, to=4-4]
\end{tikzcd}

%% file: diagrams/d385.tikz
\begin{tikzcd}[ampersand replacement=\&]
{TF} \&  \& {FT} \\
{SG} \&  \& {GS\mmdd}
\arrow["\kappa^{T,F}", Rightarrow, from=1-1, to=1-3]
\arrow["\alpha\beta"', Rightarrow, from=1-1, to=2-1]
\arrow["\beta\alpha", Rightarrow, from=1-3, to=2-3]
\arrow["\kappa^{S,G}", Rightarrow, from=2-1, to=2-3]
\end{tikzcd}

%% file: diagrams/d417.tikz
\begin{tikzcd}[ampersand replacement=\&]
 \& {\mathfrak{O}_{T\mathbb{C}}F} \&  \&  \& {F\mathfrak{O}_{T\mathbb{C}}} \\
{TF} \&  \&  \& {FT} \\
 \& {\mathfrak{O}_{S\mathbb{C}}F} \&  \&  \& {F\mathfrak{O}_{S\mathbb{C}}} \\
{SF} \&  \&  \& {FS}
\arrow["\kappa^{T\mathbb{C},F}", Rightarrow, from=1-2, to=1-5]
\arrow["\mathfrak{O}_{\alpha\mathbb{C}}F"'{pos=0.7}, Rightarrow, from=1-2, to=3-2]
\arrow["F\mathfrak{O}_{\alpha\mathbb{C}}", Rightarrow, from=1-5, to=3-5]
\arrow["\omega_{T}F", Rightarrow, from=2-1, to=1-2]
\arrow["\alpha F"', Rightarrow, from=2-1, to=4-1]
\arrow["F\omega_{T}", Rightarrow, from=2-4, to=1-5]
\arrow["\kappa^{S\mathbb{C},F}"{pos=0.3}, Rightarrow, from=3-2, to=3-5]
\arrow["\omega_{S}F"{pos=0.7}, Rightarrow, from=4-1, to=3-2]
\arrow["\kappa^{S,F}"', Rightarrow, from=4-1, to=4-4]
\arrow["F\omega_{S}"', Rightarrow, from=4-4, to=3-5]
\arrow["\kappa^{T,F}"{pos=0.6}, crossing over, Rightarrow, from=2-1, to=2-4]
\arrow["F\alpha"{pos=0.3}, crossing over, Rightarrow, from=2-4, to=4-4]
\end{tikzcd}

%% file: diagrams/d273.tikz
\begin{tikzcd}[ampersand replacement=\&]
{\C_{\mathsf{X}}F} \&  \& {F\C_{\mathsf{X}}} \\
 \& {F\mmdc}
\arrow["\kappa^{\C_{\mathsf{X}},F}", Rightarrow, from=1-1, to=1-3]
\arrow["\ev_{x}F"', Rightarrow, from=1-1, to=2-2]
\arrow["F\ev_{x}", Rightarrow, from=1-3, to=2-2]
\end{tikzcd}

%% file: diagrams/d354.tikz
\begin{tikzcd}[ampersand replacement=\&]
 \& {F} \\
{\fK F} \&  \& {F\fK\mmdd}
\arrow["\iota_{00}F"', Rightarrow, from=1-2, to=2-1]
\arrow["F\iota_{00}", Rightarrow, from=1-2, to=2-3]
\arrow["\kappa^{\fK,F}"', Rightarrow, from=2-1, to=2-3]
\end{tikzcd}

%% file: diagrams/d397.tikz
\begin{tikzcd}[column sep=small, ampersand replacement=\&]  %d397
{TSF} \&  \& {FTS} \\
 \& {TFS\mmdc}
\arrow["\kappa^{TS,F}", Rightarrow, from=1-1, to=1-3]
\arrow["T\kappa^{S,F}"', Rightarrow, from=1-1, to=2-2]
\arrow["\kappa^{T,F}S"', Rightarrow, from=2-2, to=1-3]
\end{tikzcd}

%% file: diagrams/d396.tikz
\begin{tikzcd}[ampersand replacement=\&]
{(T\oplus S)F} \&  \& {F(T\oplus S)} \\
{TF\oplus SF} \&  \& {TF\oplus SF\mmdd}
\arrow["\kappa^{T\oplus S,F}", Rightarrow, from=1-1, to=1-3]
\arrow["\rc"', Rightarrow, from=1-1, to=2-1]
\arrow["\lc", Rightarrow, from=1-3, to=2-3]
\arrow["\kappa^{T,F}\oplus\kappa^{S,F}", Rightarrow, from=2-1, to=2-3]
\end{tikzcd}

%% file: diagrams/d398.tikz
\begin{tikzcd}[ampersand replacement=\&]
{TSF} \&  \& {TFS} \&  \& {FTS} \\
 \& {\mathfrak{O}_{T\mathbb{C}}\mathfrak{O}_{S\mathbb{C}}F} \& {\mathfrak{O}_{T\mathbb{C}}F\mathfrak{O}_{S\mathbb{C}}} \& {F\mathfrak{O}_{T\mathbb{C}}\mathfrak{O}_{S\mathbb{C}}} \\
 \& {\mathfrak{O}_{T\mathbb{C}\otimes S\mathbb{C}}F} \&  \& {F\mathfrak{O}_{T\mathbb{C}\otimes S\mathbb{C}}} \\
{\mathfrak{O}_{TS\mathbb{C}}F} \&  \&  \&  \& {F\mathfrak{O}_{TS\mathbb{C}}}
\arrow["T\kappa^{S,F}", Rightarrow, from=1-1, to=1-3]
\arrow["\omega_{T}\omega_{S}F"{pos=0.6}, Rightarrow, from=1-1, to=2-2]
\arrow["\omega_{TS}F"', Rightarrow, from=1-1, to=4-1]
\arrow["\kappa^{T,F}S", Rightarrow, from=1-3, to=1-5]
\arrow["\omega_{T}F\omega_{S}", Rightarrow, from=1-3, to=2-3]
\arrow["F\omega_{T}\omega_{S}"'{pos=0.6}, Rightarrow, from=1-5, to=2-4]
\arrow["F\omega_{TS}", Rightarrow, from=1-5, to=4-5]
\arrow[Rightarrow, from=2-2, to=2-3]
\arrow["\mathfrak{O}^{2}F", Rightarrow, from=2-2, to=3-2]
\arrow[Rightarrow, from=2-3, to=2-4]
\arrow["F\mathfrak{O}^{2}"', Rightarrow, from=2-4, to=3-4]
\arrow["\kappa^{T\mathbb{C}\otimes S\mathbb{C},F}", Rightarrow, from=3-2, to=3-4]
\arrow["\mathfrak{O}_{\omega_{T}^{-1}S\mathbb{C}}F"{pos=0.4}, Rightarrow, from=3-2, to=4-1]
\arrow["F\mathfrak{O}_{\omega_{T}^{-1}S\mathbb{C}}"'{pos=0.4}, Rightarrow, from=3-4, to=4-5]
\arrow["\kappa^{TS\mathbb{C},F}"', Rightarrow, from=4-1, to=4-5]
\end{tikzcd}

%% file: diagrams/d484.tikz
\begin{tikzcd}[ampersand replacement=\&]
{(T\oplus S)F} \& {TF\oplus SF} \&  \& {FT\oplus FS} \& {F(T\oplus S)} \\
 \& {\mathfrak{O}_{T\mathbb{C}}F\oplus\mathfrak{O}_{S\mathbb{C}}F} \&  \& {F\mathfrak{O}_{T\mathbb{C}}\oplus F\mathfrak{O}_{S\mathbb{C}}} \\
 \& {(\mathfrak{O}_{T\mathbb{C}}\oplus\mathfrak{O}_{S\mathbb{C}})F} \&  \& {F(\mathfrak{O}_{T\mathbb{C}}\oplus\mathfrak{O}_{S\mathbb{C}})} \\
{\mathfrak{O}_{(T\oplus S)\mathbb{C}}F} \& {\mathfrak{O}_{T\mathbb{C}\oplus S\mathbb{C}}F} \&  \& {F\mathfrak{O}_{T\mathbb{C}\oplus S\mathbb{C}}} \& {F\mathfrak{O}_{(T\oplus S)\mathbb{C}}}
\arrow["\rc", Rightarrow, from=1-1, to=1-2]
\arrow["(\omega_{T}\oplus\omega_{S})F"', curve={height=12pt}, Rightarrow, from=1-1, to=3-2]
\arrow["\omega_{T\oplus S}F"', curve={height=30pt}, Rightarrow, from=1-1, to=4-1]
\arrow["\kappa^{T,F}\oplus\kappa^{S,F}", Rightarrow, from=1-2, to=1-4]
\arrow["\omega_{T}F\oplus\omega_{S}F"', Rightarrow, from=1-2, to=2-2]
\arrow["\lc^{-1}", Rightarrow, from=1-4, to=1-5]
\arrow["F\omega_{T}\oplus F\omega_{S}", Rightarrow, from=1-4, to=2-4]
\arrow["F(\omega_{T}\oplus\omega_{S})", curve={height=-12pt}, Rightarrow, from=1-5, to=3-4]
\arrow["F\omega_{T\oplus S}", curve={height=-30pt}, Rightarrow, from=1-5, to=4-5]
\arrow["\kappa^{T\mathbb{C},F}\oplus\kappa^{S\mathbb{C},F}", Rightarrow, from=2-2, to=2-4]
\arrow["\rc^{-1}"', Rightarrow, from=2-2, to=3-2]
\arrow["\lc^{-1}", Rightarrow, from=2-4, to=3-4]
\arrow["\mathfrak{O}_{\oplus}^{2}F"', Rightarrow, from=3-2, to=4-2]
\arrow["F\mathfrak{O}_{\oplus}^{2}", Rightarrow, from=3-4, to=4-4]
\arrow["\kappa^{(T\oplus S)\mathbb{C},F}", sqars=1.5em, Rightarrow, from=4-1, to=4-5]
\arrow[equals, from=4-2, to=4-1]
\arrow["\kappa^{T\mathbb{C}\oplus S\mathbb{C},F}", Rightarrow, from=4-2, to=4-4]
\arrow[equals, from=4-4, to=4-5]
\end{tikzcd}

%% file: diagrams/d376.tikz
\begin{tikzcd}[ampersand replacement=\&]
 \& {G\oplus G} \\
{(\Id\oplus\Id)G} \&  \& {G(\Id\oplus\Id)\mmdd}
\arrow["\lc^{-1}", Rightarrow, from=1-2, to=2-3]
\arrow["\rc", Rightarrow, from=2-1, to=1-2]
\arrow["\kappa^{\Id\oplus\Id,G}"', Rightarrow, from=2-1, to=2-3]
\end{tikzcd}

%% file: diagrams/d420.tikz
\begin{tikzcd}[ampersand replacement=\&]
{\mathfrak{O}_{T\mathbb{C}}S\mathbb{C}} \&  \&  \&  \&  \&  \& {S\mathfrak{O}_{T\mathbb{C}}\mathbb{C}} \\
 \&  \& {TS\mathbb{C}} \&  \& {ST\mathbb{C}} \&  \& {\mathfrak{O}_{S\mathbb{C}}T\mathbb{C}} \\
{T\mathbb{C}\otimes S\mathbb{C}} \&  \&  \&  \&  \&  \& {S\mathbb{C}\otimes T\mathbb{C}\mmdc}
\arrow["\kappa^{T\mathbb{C},S}\mathbb{C}", from=1-1, to=1-7]
\arrow[equals, from=1-1, to=3-1]
\arrow["\omega_{S}\omega_{T}^{-1}\mathbb{C}", from=1-7, to=2-7]
\arrow["\omega_{T}S\mathbb{C}"'{pos=0.4}, from=2-3, to=1-1]
\arrow["\kappa^{T,S}\mathbb{C}"', from=2-3, to=2-5]
\arrow["S\omega_{T}\mathbb{C}"{pos=0.4}, from=2-5, to=1-7]
\arrow["\omega_{S}T\mathbb{C}"', from=2-5, to=2-7]
\arrow[equals, from=2-7, to=3-7]
\arrow["\u\overline{\EV}(\kappa^{S,T})\mathbb{C}", from=3-1, to=3-7]
\end{tikzcd}

%% file: diagrams/d482.tikz
\begin{tikzcd}[ampersand replacement=\&]
{\mathfrak{O}_{T\mathbb{C}}S\mathbb{C}} \&  \& {\mathfrak{O}_{T\mathbb{C}}\mathfrak{O}_{S\mathbb{C}}\mathbb{C}} \&  \& {\mathfrak{O}_{S\mathbb{C}}\mathfrak{O}_{T\mathbb{C}}\mathbb{C}} \&  \& {S\mathfrak{O}_{T\mathbb{C}}\mathbb{C}} \\
 \&  \& {T\mathbb{C}\otimes(S\mathbb{C}\otimes\mathbb{C})} \&  \& {S\mathbb{C}\otimes(T\mathbb{C}\otimes\mathbb{C})} \\
{T\mathbb{C}\otimes S\mathbb{C}} \&  \& {(T\mathbb{C}\otimes S\mathbb{C})\otimes\mathbb{C}} \&  \& {(S\mathbb{C}\otimes T\mathbb{C})\otimes\mathbb{C}} \&  \& {S\mathbb{C}\otimes T\mathbb{C}}
\arrow["\mathfrak{O}_{T\mathbb{C}}\omega_{S}\mathbb{C}"{pos=0.6}, from=1-1, to=1-3]
\arrow["\kappa^{T\mathbb{C},S}\mathbb{C}", sqarn=1.5em, from=1-1, to=1-7]
\arrow[equals, from=1-1, to=3-1]
\arrow["\kappa^{T\mathbb{C},\mathfrak{O}_{S\mathbb{C}}}\mathbb{C}", from=1-3, to=1-5]
\arrow[equals, from=1-3, to=2-3]
\arrow[equals, from=1-5, to=2-5]
\arrow["\mathfrak{O}_{S\mathbb{C}}\omega_{T}^{-1}\mathbb{C}"{pos=0.6}, curve={height=-12pt}, from=1-5, to=3-7]
\arrow["\omega_{S}^{-1}\mathfrak{O}_{T\mathbb{C}}\mathbb{C}", from=1-5, to=1-7]
\arrow["\omega_{S}\omega_{T}^{-1}\mathbb{C}", curve={height=-6pt}, from=1-7, to=3-7]
\arrow["(\alpha^{\otimes})^{-1}", from=2-3, to=3-3]
\arrow["(\alpha^{\otimes})^{-1}"', from=2-5, to=3-5]
\arrow["S\mathbb{C}\otimes\rho^{\otimes}"{pos=0.4}, from=2-5, to=3-7]
\arrow["\xi^{\otimes}", sqars=1.5em, from=3-1, to=3-7]
\arrow["\rho^{\otimes}"', from=3-3, to=3-1]
\arrow["\xi^{\otimes}\otimes\mathbb{C}", from=3-3, to=3-5]
\arrow["\rho^{\otimes}", from=3-5, to=3-7]
\end{tikzcd}

%% file: diagrams/d430.tikz
\begin{tikzcd}[ampersand replacement=\&]
{\mathfrak{O}_{A}\C_{\mathsf{X}}} \&  \&  \&  \& {\C_{\mathsf{X}}\mathfrak{O}_{A}} \\
 \&  \& {\mathfrak{O}_{A}} \\
 \&  \& {\prod_{\mathsf{X}}\mathfrak{O}_{A}}
\arrow["\xi_{\mathfrak{O}_{A},\C_{\mathsf{X}}}", Rightarrow, from=1-1, to=1-5]
\arrow["\mathfrak{O}_{A}\ev_{x}", Rightarrow, from=1-1, to=2-3]
\arrow[curve={height=18pt}, dashed, Rightarrow, from=1-1, to=3-3]
\arrow["\ev_{x}\mathfrak{O}_{A}"', Rightarrow, from=1-5, to=2-3]
\arrow[curve={height=-18pt}, dashed, Rightarrow, from=1-5, to=3-3]
\arrow["\pr_{x}", Rightarrow, from=3-3, to=2-3]
\end{tikzcd}

%% file: diagrams/d274.tikz
\begin{tikzcd}[ampersand replacement=\&]
{\mathfrak{O}_{C_{0}(\mathsf{X})}F} \&  \&  \&  \& {F\mathfrak{O}_{C_{0}(\mathsf{X})}} \\
{\mathfrak{O}_{A}F} \&  \&  \&  \& {F\mathfrak{O}_{A}}
\arrow["\kappa_{1}^{C_{0}(\mathsf{X}),F}=\kappa_{2}^{C_{0}(\mathsf{X}),F}", Rightarrow, from=1-1, to=1-5]
\arrow["\mathfrak{O}_{\zeta}F"', Rightarrow, from=1-1, to=2-1]
\arrow["F\mathfrak{O}_{\zeta}", Rightarrow, from=1-5, to=2-5]
\arrow["\kappa_{j}^{A,F}", Rightarrow, from=2-1, to=2-5]
\end{tikzcd}

%% file: diagrams/d278.tikz
\begin{tikzcd}[ampersand replacement=\&]
{\C_{\mathsf{X}}F} \&  \&  \&  \& {F\C_{\mathsf{X}}} \\
 \&  \& {F} \\
 \&  \& {G} \\
{\C_{\mathsf{X}}G} \&  \&  \&  \& {G\C_{\mathsf{X}}}
\arrow["\kappa^{\C_{\mathsf{X}},F}", Rightarrow, from=1-1, to=1-5]
\arrow["\ev_{x}F"', Rightarrow, from=1-1, to=2-3]
\arrow["\C_{\mathsf{X}}\alpha"', Rightarrow, from=1-1, to=4-1]
\arrow["F\ev_{x}", Rightarrow, from=1-5, to=2-3]
\arrow["\alpha\C_{\mathsf{X}}", Rightarrow, from=1-5, to=4-5]
\arrow["\alpha", Rightarrow, from=2-3, to=3-3]
\arrow["\ev_{x}G", Rightarrow, from=4-1, to=3-3]
\arrow["\kappa^{\C_{\mathsf{X}},G}", Rightarrow, from=4-1, to=4-5]
\arrow["G\ev_{x}"', Rightarrow, from=4-5, to=3-3]
\end{tikzcd}

%% file: diagrams/d279.tikz
\begin{tikzcd}[ampersand replacement=\&]
 \&  \& {\mathfrak{O}_{A}F} \&  \&  \& {F\mathfrak{O}_{A}} \\
{\mathfrak{O}_{C_{0}(\mathsf{X})}F} \&  \&  \& {F\mathfrak{O}_{C_{0}(\mathsf{X})}} \\
 \&  \& {\mathfrak{O}_{A}G} \&  \&  \& {G\mathfrak{O}_{A}} \\
{\mathfrak{O}_{C_{0}(\mathsf{X})}G} \&  \&  \& {G\mathfrak{O}_{C_{0}(\mathsf{X})}}
\arrow["\kappa^{A,F}", Rightarrow, from=1-3, to=1-6]
\arrow["\mathfrak{O}_{A}\alpha"{pos=0.7}, Rightarrow, from=1-3, to=3-3]
\arrow["\alpha\mathfrak{O}_{A}", Rightarrow, from=1-6, to=3-6]
\arrow["\mathfrak{O}_{\zeta}F", Rightarrow, from=2-1, to=1-3]
\arrow["\mathfrak{O}_{C_{0}(\mathsf{X})}\alpha"', Rightarrow, from=2-1, to=4-1]
\arrow["F\mathfrak{O}_{\zeta}"{pos=0.3}, Rightarrow, from=2-4, to=1-6]
\arrow["\kappa^{A,G}"{pos=0.6}, Rightarrow, from=3-3, to=3-6]
\arrow["\mathfrak{O}_{\zeta}G"{pos=0.4}, Rightarrow, from=4-1, to=3-3]
\arrow["\kappa^{C_{0}(\mathsf{X}),G}"', Rightarrow, from=4-1, to=4-4]
\arrow["G\mathfrak{O}_{\zeta}"', Rightarrow, from=4-4, to=3-6]
\arrow["\kappa^{{C_{0}(\mathsf{X})},F}"{pos=0.4}, crossing over, Rightarrow, from=2-1, to=2-4]
\arrow["\alpha\mathfrak{O}_{C_{0}(\mathsf{X})}"'{pos=0.7}, crossing over, Rightarrow, from=2-4, to=4-4]
\end{tikzcd}

%% file: diagrams/d365.tikz
\begin{tikzcd}[column sep=large, ampersand replacement=\&]  %d365
{GB} \& {F\C_{\mathsf{X}}B} \& {F\prod_{\mathsf{X}}B} \& {\prod_{\mathsf{X}}FB} \\
 \&  \& {FB}
\arrow["\alpha_{j}B", from=1-1, to=1-2]
\arrow[curve={height=18pt}, dashed, from=1-1, to=2-3]
\arrow[hook, from=1-2, to=1-3]
\arrow["F\ev_{x}B"'{pos=0.3}, from=1-2, to=2-3]
\arrow[hook, from=1-3, to=1-4]
\arrow["F\pr_{x}", from=1-3, to=2-3]
\arrow["\pr_{x}"{pos=0.3}, curve={height=-6pt}, from=1-4, to=2-3]
\end{tikzcd}

%% file: diagrams/d348.tikz
\begin{tikzcd}[ampersand replacement=\&]
{A\otimes_{\min}\prod_{j}B_{j}} \& {\prod_{j}(A\otimes_{\min}B_{j})} \\
{\mathcal{H}\otimes(\bigoplus_{j}\mathcal{H}_{j})} \& {\bigoplus_{j}(\mathcal{H}\otimes\mathcal{H}_{j})}
\arrow[from=1-1, to=1-2]
\arrow[hook, from=1-1, to=2-1]
\arrow[hook, from=1-2, to=2-2]
\arrow["\cong", from=2-1, to=2-2]
\end{tikzcd}

%% file: diagrams/d224.tikz
\begin{tikzcd}[ampersand replacement=\&]
{F} \& {GI} \\
 \& {G\mmdd}
\arrow["\gamma_{j}"', Rightarrow, from=1-1, to=2-2]
\arrow["\Gamma", Rightarrow, from=1-1, to=1-2]
\arrow["G\ev_{j}", Rightarrow, from=1-2, to=2-2]
\end{tikzcd}

%% file: diagrams/d271.tikz
\begin{tikzcd}[ampersand replacement=\&]
{I} \& {I} \\
{I} \& {\Id\mmdc}
\arrow["r_{1}", Rightarrow, from=1-1, to=1-2]
\arrow["r_{0}"', Rightarrow, from=1-1, to=2-1]
\arrow["\ev_{0}", Rightarrow, from=1-2, to=2-2]
\arrow["\ev_{1}"', Rightarrow, from=2-1, to=2-2]
\end{tikzcd}

%% file: diagrams/d272.tikz
\begin{tikzcd}[ampersand replacement=\&]
{GI} \& {GI} \\
{GI} \& {G\mmdd}
\arrow["Gr_{1}", Rightarrow, from=1-1, to=1-2]
\arrow["Gr_{0}"', Rightarrow, from=1-1, to=2-1]
\arrow["G\ev_{0}", Rightarrow, from=1-2, to=2-2]
\arrow["G\ev_{1}"', Rightarrow, from=2-1, to=2-2]
\end{tikzcd}

%% file: diagrams/d194.tikz
\begin{tikzcd}[ampersand replacement=\&]
{F} \&  \& {G} \\
 \&  \& {GI} \\
{G} \& {GI} \& {G\mmdc}
\arrow["\alpha_{2}", Rightarrow, from=1-1, to=1-3]
\arrow["\alpha_{12}", Rightarrow, from=1-1, to=2-3]
\arrow["\alpha_{0}"', Rightarrow, from=1-1, to=3-1]
\arrow["\alpha_{01}"', Rightarrow, from=1-1, to=3-2]
\arrow["\alpha_{1}"{pos=0.7}, Rightarrow, from=1-1, to=3-3]
\arrow["G\ev_{1}"', Rightarrow, from=2-3, to=1-3]
\arrow["G\ev_{0}", Rightarrow, from=2-3, to=3-3]
\arrow["G\ev_{0}", Rightarrow, from=3-2, to=3-1]
\arrow["G\ev_{1}"', Rightarrow, from=3-2, to=3-3]
\end{tikzcd}

%% file: diagrams/d195.tikz
\begin{tikzcd}[ampersand replacement=\&]
{F} \&  \& {G} \\
 \& {GI} \& {GI} \\
{G} \& {GI} \& {G}
\arrow["\alpha_{2}", Rightarrow, from=1-1, to=1-3]
\arrow[Rightarrow, dashed, from=1-1, to=2-2]
\arrow["\alpha_{0}"', Rightarrow, from=1-1, to=3-1]
\arrow["G\ev_{1}"{pos=0.4}, Rightarrow, from=2-2, to=1-3]
\arrow["Gr_{1}", Rightarrow, from=2-2, to=2-3]
\arrow["G\ev_{0}"'{pos=0.4}, Rightarrow, from=2-2, to=3-1]
\arrow["Gr_{0}", Rightarrow, from=2-2, to=3-2]
\arrow["G\ev_{1}"', Rightarrow, from=2-3, to=1-3]
\arrow["G\ev_{0}", Rightarrow, from=2-3, to=3-3]
\arrow["G\ev_{0}", Rightarrow, from=3-2, to=3-1]
\arrow["G\ev_{1}"', Rightarrow, from=3-2, to=3-3]
\end{tikzcd}

%% file: diagrams/d196.tikz
\begin{tikzcd}[ampersand replacement=\&]
{F} \& {G} \& {HI} \\
 \&  \& {H\mmdc}
\arrow["\alpha", Rightarrow, from=1-1, to=1-2]
\arrow["\beta_{j}"', Rightarrow, from=1-2, to=2-3]
\arrow["\beta", Rightarrow, from=1-2, to=1-3]
\arrow["H\ev_{j}", Rightarrow, from=1-3, to=2-3]
\end{tikzcd}

%% file: diagrams/d199.tikz
\begin{tikzcd}[ampersand replacement=\&]
{F} \& {GI} \& {HI} \\
 \& {G} \& {H\mmdc}
\arrow["\alpha_{j}"', Rightarrow, from=1-1, to=2-2]
\arrow["\alpha", Rightarrow, from=1-1, to=1-2]
\arrow["G\ev_{j}", Rightarrow, from=1-2, to=2-2]
\arrow["\beta I", Rightarrow, from=1-2, to=1-3]
\arrow["H\ev_{j}", Rightarrow, from=1-3, to=2-3]
\arrow["\beta"', Rightarrow, from=2-2, to=2-3]
\end{tikzcd}

%% file: diagrams/d222.tikz
\begin{tikzcd}[ampersand replacement=\&]
{GF} \& {GF'I} \& {G'F'I} \\
 \& {GF'} \& {G'F'\mmdc}
\arrow["G\alpha", Rightarrow, from=1-1, to=1-2]
\arrow["G\alpha_{j}"', Rightarrow, from=1-1, to=2-2]
\arrow["\beta F'I", Rightarrow, from=1-2, to=1-3]
\arrow["GF'\ev_{j}", Rightarrow, from=1-2, to=2-2]
\arrow["G'F'\ev_{j}", Rightarrow, from=1-3, to=2-3]
\arrow["\beta F'"', Rightarrow, from=2-2, to=2-3]
\end{tikzcd}

%% file: diagrams/d223.tikz
\begin{tikzcd}[ampersand replacement=\&]
{F'IG} \&  \& {F'GI} \&  \& {F'G'I} \\
{FG} \&  \& {F'G} \&  \& {F'G'\mmdc}
\arrow["F'\kappa^{I,G}", Rightarrow, from=1-1, to=1-3]
\arrow["F'\ev_{j}G", curve={height=6pt}, Rightarrow, from=1-1, to=2-3]
\arrow["F'\beta I", Rightarrow, from=1-3, to=1-5]
\arrow["F'G\ev_{j}", Rightarrow, from=1-3, to=2-3]
\arrow["F'G'\ev_{j}", Rightarrow, from=1-5, to=2-5]
\arrow["\alpha G", Rightarrow, from=2-1, to=1-1]
\arrow["\alpha_{j}G"', Rightarrow, from=2-1, to=2-3]
\arrow["F'\beta"', Rightarrow, from=2-3, to=2-5]
\end{tikzcd}

%% file: diagrams/d288.tikz
\begin{tikzcd}[ampersand replacement=\&]
{F\oplus G} \&  \& {F'I\oplus G'I} \&  \& {(F'\oplus G')I} \\
 \&  \& {F'\oplus G'}
\arrow["\alpha\oplus\beta", Rightarrow, from=1-1, to=1-3]
\arrow["\alpha_{j}\oplus\beta_{j}"', curve={height=12pt}, Rightarrow, from=1-1, to=2-3]
\arrow["\rc^{-1}", Rightarrow, from=1-3, to=1-5]
\arrow["F'\ev_{j}\oplus G'\ev_{j}"', Rightarrow, from=1-3, to=2-3]
\arrow["(F'\oplus G')\ev_{j}", curve={height=-12pt}, Rightarrow, from=1-5, to=2-3]
\end{tikzcd}

%% file: diagrams/d402.tikz
\begin{tikzcd}[ampersand replacement=\&]
{T} \&  \& {S} \\
{SI} \&  \& {S\Id\mmdc}
\arrow["\alpha_{j}", Rightarrow, from=1-1, to=1-3]
\arrow["\alpha"', Rightarrow, from=1-1, to=2-1]
\arrow["S\ev_{j}", Rightarrow, from=2-1, to=2-3]
\arrow[equals, from=2-3, to=1-3]
\end{tikzcd}

%% file: diagrams/d403.tikz
\begin{tikzcd}[column sep=large, ampersand replacement=\&]  %d403
{\overline{\EV}(T)} \&  \& {\overline{\EV}(S)} \\
{\overline{\EV}(S)\otimes I\mathbb{C}} \&  \& {\overline{\EV}(S)\otimes\mathbb{C}}
\arrow["\overline{\EV}(\alpha_{j})", from=1-1, to=1-3]
\arrow["\overline{\EV}(\alpha)"', from=1-1, to=2-1]
\arrow["\overline{\EV}(S)\otimes\ev_{j}\mathbb{C}", from=2-1, to=2-3]
\arrow["\lambda^{\otimes}"', from=2-3, to=1-3]
\end{tikzcd}

%% file: diagrams/d578.tikz
\begin{tikzcd}[ampersand replacement=\&]
{\Cstar} \&  \& {\GEFC_{\ttt}} \&  \& {\Cstar} \\
{\hCstar} \&  \& {\hGEFC_{\ttt}} \&  \& {\hCstar,}
\arrow["\mathfrak{O}", from=1-1, to=1-3]
\arrow[from=1-1, to=2-1]
\arrow["\EV", from=1-3, to=1-5]
\arrow[from=1-3, to=2-3]
\arrow[from=1-5, to=2-5]
\arrow["\mathfrak{O}_{\simeq}", dashed, from=2-1, to=2-3]
\arrow["\EV_{\simeq}", dashed, from=2-3, to=2-5]
\end{tikzcd}

%% file: diagrams/d436.tikz
\begin{tikzcd}[ampersand replacement=\&]
 \& {\Br(\GEFC)} \&  \& {\Br(\Cstar)} \\
{\Br(\hGEFC)} \& {\Br(\GEFC)\modrel{\simeq}} \&  \& {\Br(\Cstar)\modrel{\simeq}} \& {\Br(\hCstar)}
\arrow["\overline{\EV}", from=1-2, to=1-4]
\arrow[from=1-2, to=2-2]
\arrow[from=1-4, to=2-4]
\arrow[equals, from=2-1, to=2-2]
\arrow["\overline{\EV}_{\simeq}", dashed, from=2-2, to=2-4]
\arrow[equals, from=2-4, to=2-5]
\end{tikzcd}

%% file: sec5.tex
\section{\label{sec:genhom}Generalized morphisms and asymptotic adjunction}

In this section, we introduce generalized morphisms between $C^{*}$-algebras
associated with a good endofunctor and endow them with the structure
of a commutative monoid. It will be shown that such monoids admit
a bilinear pairing which generalizes the composition in the $E$-theory
category. We also consider the asymptotic version of this construction,
and introduce asymptotically adjoint good endofunctors.

\subsection{Generalized morphisms}

\begin{defn}
Let $A$, $B$ be $C^{*}$-algebras, let $F\in\in\hGEFC$, and let
$\varphi_{0},\varphi_{1}\colon A\to FB$ be two $*$-homomorphisms.
We say that they are $F$-\emph{homotopic} (written $\varphi_{0}\simeq_{F}\varphi_{1}$)
if there is a $*$-homomorphism $\Phi\colon A\to FIB$, such that
the following diagram commutes for $j=0,1$:
\[
\input{diagrams/d434.tikz}
\]
By an argument similar to that in the proof of Theorem~\ref{thm:homot-eqrel},
one can show that $F$-homotopy is an equivalence relation. We denote
by $\left[A,F,B\right]\coloneqq\hom(A,FB)\modrel{\simeq_{F}}$ the
set of $F$-homotopy classes of $*$-homomorphisms from $A$ to $FB$,
and write $[\varphi]$ for the class of $\varphi\colon A\to FB$ in
$[A,F,B]$.
\end{defn}

\begin{rem}
Note that the $n$-homotopy classes from~\cite{GHT} discussed in
the introduction are precisely the classes of $\mathfrak{A}^{n}$-homotopy:
$[A,\mathfrak{A}^{n},B]=\bbrackets{A,B}_{n}$.
\end{rem}

\begin{defn}\label{def:53}
Define the binary operation $+$ on $[A,F\fK,B]$ by
$[\varphi]+[\psi]\coloneqq[\varphi\uu\psi]$, where $\varphi\uu\psi$
is the following composite:
\[
\varphi\uu\psi\colon A\xrightarrow{\Delta}A\oplus A\xrightarrow{\varphi\oplus\psi}F\fK B\oplus F\fK B=(F\fK\oplus F\fK)B\xrightarrow{\lc^{-1}B}F(\fK\oplus\fK)B\xrightarrow{F\mu B}F\fK B,
\]
in which $\Delta$ is as in Definition~\ref{def:Delta} (and which
is explicitly defined by the formula $\Delta(a)=(a,a)$), and $\mu$
is as in Section~\ref{subsec:starig}. It is straightforward to check
that $+$ is well-defined: if $\varphi_{0}$, $\varphi_{1}$ and $\psi_{0}$,
$\psi_{1}$ are connected respectively by $F$-homotopies $\Phi$
and $\Psi$, then $\varphi_{0}\uu\psi_{0}$ and $\varphi_{1}\uu\psi_{1}$
are connected by the $F$-homotopy $\Phi\uu\Psi$.
\end{defn}

\begin{prop}
$([A,F\fK,B],+,[0])$ is a commutative monoid.
\end{prop}

\begin{proof}
We claim that the following diagrams commute in $\hGEFC$:
\begin{equation}
\input{diagrams/d492.tikz}\label{eq:d280}
\end{equation}
\begin{equation}
\input{diagrams/d326.tikz}\label{eq:326}
\end{equation}
\begin{equation}
\input{diagrams/d323.tikz}\label{eq:323}
\end{equation}
Indeed, the subdiagrams marked by $\tc 1$ commute by the bimonoidality
of $\hGEFC$; those marked with $\tc 2$ commute up to homotopy since
$(\fK,[\mu],\bzero)$ is a rig in $\hGEFC$; and $\tc 3$ commutes
since $\bzero$ is an initial object in $\GEFC$.

We leave it to the reader to deduce from the commutativity of the
perimeters of these three diagrams that the binary operation $+$
is associative (diagram~(\ref{eq:d280}) and Lemma~\ref{lem:comonoid}),
symmetric (diagram~(\ref{eq:326})), and that $[0]$ is the neutral
element (diagram~(\ref{eq:323})).
\end{proof}

\begin{prop}\label{prop:bullet-bilin-assoc}
There is a well-defined bilinear associative
operation
\[
\bullet\colon[B,G\fK,C]\times[A,F\fK,B]\to[A,FG\fK,C]\colon([\psi],[\varphi])\mapsto[\psi]\bullet[\varphi]\coloneqq[\psi\bullet\varphi]
\]
where $\psi\bullet\varphi$ is the following composite:
\begin{equation}
\psi\bullet\varphi\colon A\xrightarrow{\varphi}F\fK B\xrightarrow{F\fK\psi}F\fK G\fK C\xrightarrow{F\kappa^{\fK,G}\fK C}FG\fK^{2}C\xrightarrow{FG\theta C}FG\fK C.\label{bullet}
\end{equation}
Furthermore, the corner embedding $\iota_{00}\colon\Id\Rightarrow\fK$
gives rise to the neutral element, i.e.,
\[
[\varphi]\bullet[\iota_{00}A]=[\varphi],\qquad\textrm{and}\qquad[\iota_{00}B]\bullet[\varphi]=[\varphi].
\]
\end{prop}

\begin{rem}\label{rem:57}
It follows from Corollary~\ref{cor:ai} that $\kappa^{\fK,G\fK}\simeq\kappa^{\fK,G}\fK$.
So, we can use $F\kappa^{\fK,G\fK}C$ instead of $F\kappa^{\fK,G}\fK C$
in (\ref{bullet}).
\end{rem}

\begin{proof}
[Proof of Proposition~\ref{prop:bullet-bilin-assoc}]To prove that
$\bullet$ is well-defined we need to show that $\psi\bullet\varphi_{0}\simeq\psi\bullet\varphi_{1}$
whenever $\varphi_{0}\simeq\varphi_{1}$, and $\psi_{0}\bullet\varphi\simeq\psi_{1}\bullet\varphi$
whenever $\psi_{0}\simeq\psi_{1}$. The first implication is proved
by the diagram
\[
\input{diagrams/d490.tikz}
\]
where $\varphi$ is an $F\fK$-homotopy connecting $\varphi_{0}$
and $\varphi_{1}$, and where $\tc 1$ and $\tc 4$ commute by naturality;
$\tc 2$~by Lemma~\ref{lem:kappa-tensor-type}; $\tc 3$~by Corollary~\ref{cor:332}.
The proof of the second implication is similar.

Associativity essentially follows from the outer part of the diagram
\[
\input{diagrams/d304.tikz}
\]
in which the left-bottom rectangle commutes by the definition of labeled
endofunctors, the right-bottom rectangle commutes up to homotopy since
$(\fK,[\theta],[\iota_{00}])$ is a monoid in $(\hGEFC,\singledot,\Id)$,
and the top triangle commutes trivially.

We now proceed to prove bilinearity. The equality $([\varphi_{0}]+[\varphi_{1}])\bullet[\psi]=[\psi\bullet\varphi_{0}]+[\psi\bullet\varphi_{1}]$
follows from the outer part of the diagram
\[
\input{diagrams/d400.tikz}
\]
in which the subdiagrams marked with $\tc 1$ commute by the bimonoidality
of $\GEFC$; $\tc 2$ commutes by Proposition~\ref{prop:kappa-prop};
$\tc 3$ by Lemma~\ref{lem:kappa-tensor-type}; the commutativity
up to homotopy of $\tc 4$ follows from the outer diagram in
\[
\input{diagrams/d401.tikz}
\]
in which $\tc 5$ and $\tc 6$ commute by the bimonoidality of $\GEFC$,
and $\tc 7$ commutes up to homotopy by Proposition~\ref{prop:rig-K}.
The proof of the equality $[\varphi]\bullet([\psi_{0}]+[\psi_{1}])=[\varphi]\bullet[\psi_{0}]+[\varphi]\bullet[\psi_{1}]$
is similar.

We leave it to the reader to prove that $[\varphi]\bullet[\iota_{00}A]=[\varphi]$
and $[\iota_{00}B]\bullet[\varphi]=[\varphi]$.
\end{proof}

\begin{defn}
Let $F,G\in\in\GEFC$, and let $\alpha\colon F\Rightarrow G$ be a
natural transformation (not necessarily labeled) between the underlying
endofunctors. Denote by $\alpha_{*}$ the following map of sets:
\begin{alignat*}{4}
\alpha_{*}\colon &\ 
& [A,F,B] & ~ & \to & ~ &\ 
& [A,G,B]\\
 &\ 
& [A\xrightarrow{\varphi}FB] &\ 
& \mapsto &\ 
&\ 
& [A\xrightarrow{\varphi}FB\xrightarrow{\alpha B}GB].
\end{alignat*}
\end{defn}

\begin{lem}\label{lem:weak-rangle}
The following statements hold:
\begin{enumerate}
\item The map $\alpha_{*}$ is well-defined;
\item $(\beta\circ\alpha)_{*}=\beta_{*}\circ\alpha{}_{*}$ for $\alpha\colon F\Rightarrow G$
and $\beta\colon G\Rightarrow H$.
\end{enumerate}
\end{lem}

\begin{proof}
Straightforward.
\end{proof}

\begin{defn}
For $A,B\in\in\Cstar$ and $\alpha\in\GEFC(F,G)$ define the map
\[
\left\langle \alpha\right\rangle \coloneqq(\alpha\fK)_{*}\colon[A,F\fK,B]\to[A,G\fK,B]\colon[\varphi]\mapsto[\alpha\fK B\circ\varphi].
\]
\end{defn}

\begin{thm}\label{thm:58}
The following statements hold:%--!--
\begin{enumerate}[label=\textup{(\roman*)}]

\item\label{enu:mm}The map $\left\langle \alpha\right\rangle $ is a morphism
of monoids;
\item\label{enu:hmt}The map $\alpha\mapsto\left\langle \alpha\right\rangle $
is homotopy invariant, i.e. $\left\langle \alpha\right\rangle =\left\langle \beta\right\rangle $
whenever $\alpha\simeq\beta$;
\item\label{enu:v}If $\alpha\in\GEFC(F,G)$, $\beta\in\GEFC(G,H)$, then
$\left\langle \beta\right\rangle \circ\left\langle \alpha\right\rangle =\left\langle \beta\circ\alpha\right\rangle $;
\item\label{enu:nat}If $\alpha\in\GEFC(F,F')$, $\beta\in\GEFC(G,G')$,
then the following diagram commutes:
\[
\input{diagrams/d287.tikz}
\]
\end{enumerate}
\end{thm}

\begin{proof}
Statements~\ref{enu:hmt},~\ref{enu:v} are straightforward to verify.

Let $\varphi_{1},\varphi_{2}\colon A\to F\fK B$, and let $\varphi'_{j}\coloneqq\alpha\fK B\circ\varphi_{j}$,
for $j=1,2$. To prove statement~\ref{enu:mm},
it suffices to show that $\alpha\fK B\circ(\varphi_{1}\uu\varphi_{2})=\varphi'_{1}\uu\varphi'_{2}$,
which in turn follows from the diagram
\[
\input{diagrams/d281.tikz}
\]
commuting by the definition of $\varphi'$ and $\psi'$, and the bimonoidality
of $\GEFC$.

Now, let $\varphi\colon A\to F\fK B$, $\psi\colon B\to F\fK C$,
and let $\varphi'\coloneqq\alpha\fK B\circ\varphi$, $\psi'\coloneqq\alpha\fK B\circ\psi$.
Statement~\ref{enu:nat} follows from the commutativity of the outermost
diagram in
\[
\input{diagrams/d302.tikz}
\]
which can be proved similarly.
\end{proof}

\begin{lem}\label{lem:78}
Let $F,G\in\in\GEFC$, $A,A',B,B'\in\in\Cstar$, and
let $\zeta\colon A'\to A$, $\xi\colon B\to B'$ be $*$-homomorphisms.
The following maps are well-defined morphisms of monoids:
\begin{alignat*}{1}
 & [A,F\fK,B]\to[A',F\fK,B]\colon[\varphi]\mapsto[\varphi\circ\zeta],\\
 & [A,F\fK,B]\to[A,F\fK,B']\colon[\varphi]\mapsto[F\fK\xi\circ\varphi],\\
 & [A,F\fK,B]\to[GA,GF\fK,B]\colon[\varphi]\mapsto[G\varphi],\\
 & [A,F\fK^{2},B]\to[A,F\fK,B]\colon[\varphi]\mapsto[F\theta B\circ\varphi].
\end{alignat*}
\end{lem}

\begin{proof}
It is straightforward to check that these maps are well-defined, and
that the first two maps are morphisms of monoids.

The fact that the third map is a morphism of monoids follows from
\[
\input{diagrams/d425.tikz}
\]
where the left triangle commutes by Lemma~\ref{lem:Delta}, the right
triangle commutes by the bimonoidality of $\GEFC$, and the commutativity
of the left rectangle can be deduced from the universal property of
$\oplus$.

To prove that the last map is a morphism of monoids, we need to check
that if $\varphi,\psi\colon A\to F\fK^{2}B$, $\varphi'\coloneqq F\theta B\circ\varphi$,
and $\psi'\coloneqq F\theta B\circ\psi$, then $F\theta B\circ(\varphi\uu\psi)\simeq\varphi'\uu\psi'$.
This follows from the diagram
\[
\input{diagrams/d305.tikz}
\]
where $\tc 1$ and $\tc 2$ commute by the bimonoidality of $\GEFC$;
$\tc 3$ commutes up to homotopy by Proposition~\ref{prop:rig-K};
and the unlabeled subdiagrams clearly commute.
\end{proof}

\subsection{Asymptotic morphisms and asymptotic adjunction}

We now proceed to define an asymptotic version of generalized morphisms.
Recall
the asymptotic algebra functor $\mathfrak{A}\coloneqq\mathfrak{T}/\mathfrak{T}_{0}\in\in\GEFC$
from Examples~\ref{exa:frakA} and~\ref{exa:frakA-labeling}, and
introduce the following labeled natural transformations:
\begin{itemize}
\item $\const\colon\Id\Rightarrow\mathfrak{T},~\const B\colon B\to\mathfrak{T}B\colon b\mapsto[t\mapsto b]$;
\item $\as\colon\mathfrak{T}\Rightarrow\mathfrak{A}$ the componentwise
quotient projection;
\item $\alpha\coloneqq\as\circ\const\colon\Id\Rightarrow\mathfrak{A}$;
\item $\iota_{00}\colon\Id\Rightarrow\fK$ the corner embedding (see Section~\ref{subsec:starig}).
\end{itemize}
\begin{lem}\label{lem:alpha-A}
 There exists a homotopy $\alpha\mathfrak{A}\simeq\mathfrak{A}\alpha$.
\end{lem}

\begin{proof}
See the proof of~\cite[Proposition~2.8]{GHT}.
\end{proof}

\begin{defn}
For
$A,B\in\in\Cstar$ and $F\in\in\hGEFC$, introduce
the following monoid:
\begin{alignat*}{1}
\bbrackets{A,F,B} & \coloneqq\colim\left([A,F\mathfrak{A}\fK,B]\xrightarrow{\left\langle F\alpha\right\rangle }[A,F\mathfrak{A}^{2}\fK,B]\xrightarrow{\left\langle F\alpha\mathfrak{A}\right\rangle }[A,F\mathfrak{A}^{3}\fK,B]\to\cdots\right).
\end{alignat*}
For a $*$-homomorphism $\varphi\colon A\to F\mathfrak{A}^{n}\fK B$
denote by $\bbrackets{\varphi}$ its class in $\bbrackets{A,F,B}$.
\end{defn}

\begin{rem}
Note that the hom-sets of the category $\hAsy$ from the introduction
is a special case of the monoids defined above: $\bbrackets{A,B}=\bbrackets{A,\Id,B}$.
\end{rem}

\begin{defn}\label{def:asadj}
Let
$S,N\in\in\hGEFC$. We say
that there is an \emph{asymptotic adjunction} $S\asadj N$ if there
exist two labeled natural transformations $\eta\colon\Id\Rightarrow NS$
and $\varepsilon\colon SN\Rightarrow\mathfrak{A}\fK$ (which we call
the \emph{unit} and \emph{counit}) such that the following diagrams
commute in~$\hGEFC$:
\begin{alignat*}{2}
\input{diagrams/d202.tikz} & \qquad\qquad & \input{diagrams/d203.tikz}
\end{alignat*}
\end{defn}

\begin{example}\label{exa:asadj}
A basic example of asymptotically adjoint good endofunctors
is $\C_{\mathbb{R}}\asadj\mathfrak{M}^{u}_{\mathbb{Z}}$ (see Section~\ref{subsec:stacont}
for the definitions). We refer the reader to~\cite{makeev_asadj}
for further details.
\end{example}

\begin{thm}\label{thm:adj2brackets}
Let $A$ and $B$ be $C^{*}$-algebras, let
$S,N\in\in\hGEFC$, and let $S\asadj N$ be an asymptotic adjunction
witnessed by a unit $\eta$ and a counit~$\varepsilon$. Then the
formulas
\begin{align}
\Phi\colon\bbrackets{SA,\Id,B}\to\bbrackets{A,N,B} & \colon\bbrackets{SA\xrightarrow{\varphi}\mathfrak{A}^{n}\fK B}\mapsto\bbrackets{A\xrightarrow{\eta A}NSA\xrightarrow{N\varphi}N\mathfrak{A}^{n}\fK B},\label{eq:mut-izom}\\
\Psi\colon\bbrackets{A,N,B}\to\bbrackets{SA,\Id,B} & \colon\bbrackets{A\xrightarrow{\psi}N\mathfrak{A}^{n}\fK B}\nonumber \\
\mapsto\myllbracket SA\xrightarrow{S\psi}SN\mathfrak{A}^{n}\fK B & \xrightarrow{\varepsilon\mathfrak{A}^{n}\fK B}\mathfrak{A}\fK\mathfrak{A}^{n}\fK B\xrightarrow{\mathfrak{A}\kappa^{\fK,\mathfrak{A}^{n}}\fK B}\mathfrak{A}^{n+1}\fK^{2}B\xrightarrow{\mathfrak{A}^{n+1}\theta B}\mathfrak{A}^{n+1}\fK B\myrrbracket\nonumber
\end{align}
define mutually inverse isomorphisms of monoids natural in~$A$ and~$B$.
\end{thm}

\begin{proof}
Note first that $\Phi$ and $\Psi$ admit the following factorizations
into morphisms of monoids described in Lemma~\ref{lem:78} and Theorem~\ref{thm:58}:
\begin{alignat*}{1}
\Phi\colon & [SA,\mathfrak{A}^{n}\fK,B]\to[NSA,N\mathfrak{A}^{n}\fK,B]\to[A,N\mathfrak{A}^{n}\fK,B],\\
\Psi\colon & [A,N\mathfrak{A}^{n}\fK,B]\to[SA,SN\mathfrak{A}^{n}\fK,B]\to[SA,\mathfrak{A}\fK\mathfrak{A}^{n}\fK,B]\to[SA,\mathfrak{A}^{n+1}\fK^{2},B]\to[SA,\mathfrak{A}^{n+1}\fK,B].
\end{alignat*}
Hence, $\Phi$ and $\Psi$ are also morphisms of monoids.

The equalities $\Psi\circ\Phi=\id$ and $\Phi\circ\Psi=\id$ follow
respectively from the commutativity of the following
diagrams:
\[
\input{diagrams/d308.tikz}
\]
\[
\input{diagrams/d309.tikz}
\]
whose lower rows represent the classes of $\varphi$ and $\psi$ respectively
in $\bbrackets{SA,\Id,B}$ and $\bbrackets{A,N,B}$ by definition.
The triangles marked by $\tc 1$ commute since $S\asadj N$; $\tc 2$
by Lemma~\ref{lem:kappa-tensor-type}; $\tc 3$ since $(\fK,[\theta],[\iota_{00}])$
is a monoid in $(\hGEFC,\singledot,\Id)$; and the unlabeled subdiagrams
by naturality.
\end{proof}

\subsection{Compositions of asymptotic morphisms}

Suppose that $T\in\in\GEFC_{\ttt}$ has a right asymptotic adjoint
$T^{\star}$
(cf.~Example~\ref{exa:asadj}), and
denote by $\eta$ and $\varepsilon$ a unit and a counit of the asymptotic
adjunction. For every $F\in\in\GEFC$ let us introduce the labeled
natural transformation
\[
\kapparev^{F,T^{\star}}\colon FT^{\star}\xRightarrow{\eta FT^{\star}}T^{\star}TFT^{\star}\xRightarrow{T^{\star}\kappa^{T,F}T^{\star}}T^{\star}FTT^{\star}\xRightarrow{T^{\star}F\varepsilon}T^{\star}F\mathfrak{A}\fK.
\]

\begin{lem}\label{lem:kapparev}
Let $\alpha\in\GEFC(F,G)$, $T\in\in\GEFC_{\ttt}$,
and let $T\asadj T^{\star}$. Then the following diagram commutes:
\[
\input{diagrams/d426.tikz}
\]
\end{lem}

\begin{proof}
The statement follows from
\[
\input{diagrams/d427.tikz}
\]
in which the left and right rectangles commute by naturality, and
the middle one commutes since $\alpha$ is labeled.
\end{proof}

For $T$ and $T^{\star}$ as above and $n,k\in\mathbb{N}$, let us
define the natural transformation $\xi_{n,k}$ as the composite
\[
\xi_{n,k}\colon\mathfrak{A}^{n}T^{\star}\mathfrak{A}^{k}\xRightarrow{\kapparev^{\mathfrak{A}^{n},T^{\star}}\mathfrak{A}^{k}}T^{\star}\mathfrak{A}^{n}\mathfrak{A}\fK\mathfrak{A}^{k}\xRightarrow{T^{\star}\mathfrak{A}^{n+1}\kappa^{\fK,\mathfrak{A}^{k}}}T^{\star}\mathfrak{A}^{n+k+1}\fK,
\]
and introduce the following
map:
\begin{alignat}{1}
[B,T^{\star}\mathfrak{A}^{k}\fK,C]\times[A,F\mathfrak{A}^{n}\fK,B] & \xrightarrow{\bullet}[A,F\mathfrak{A}^{n}T^{\star}\mathfrak{A}^{k}\fK,C]\label{eq:ashom-pairing}\\
 & \xrightarrow{\left\langle F\xi_{n,k}\right\rangle }[A,FT^{\star}\mathfrak{A}^{n+k+1}\fK^{2},C]\xrightarrow{(FT^{\star}\mathfrak{A}^{n+k+1}\theta)_{*}}[A,FT^{\star}\mathfrak{A}^{n+k+1}\fK,C].\nonumber
\end{alignat}

\begin{prop}\label{prop:518}
Let
$F\in\in\GEFC$, $T\in\in\GEFC_{\ttt}$,
and let $T\asadj T^{\star}$. Then the formula (\ref{eq:ashom-pairing})
gives rise to the following well-defined bilinear map:
\[
\bullet_{\as}\colon\bbrackets{B,T^{\star},C}\times\bbrackets{A,F,B}\to\bbrackets{A,FT^{\star},C}.
\]
\end{prop}

\begin{proof}
To show that $\bullet_{\as}$ is well-defined it is sufficient to
show that for $n=m+k$ and $n'=m'+k'$ the following diagram commutes:
\begin{equation}
\input{diagrams/d429.tikz}\label{eq:429}
\end{equation}
where $\alpha^{k}$ stands for $\alpha\alpha\dots\alpha$ ($k$ times).
The upper rectangle commutes by Theorem~\ref{thm:58}; and the lower
rectangle commutes by Lemma~\ref{lem:weak-rangle}. In order to show
the commutativity of the middle rectangle, consider the diagram
\begin{equation}
\input{diagrams/d481.tikz}\label{eq:307}
\end{equation}
in which the upper rectangle of (\ref{eq:307}) commutes by Lemma~\ref{lem:kapparev};
and the bottom one commutes by Lemma~\ref{lem:kappa-tensor-type}.
It follows from Theorem~\ref{thm:58}, Lemma~\ref{lem:alpha-A}
and the outer rectangle in (\ref{eq:307}) that the middle rectangle
in (\ref{eq:429}) commutes.

The bilinearity of $\bullet_{\as}$ follows from that of $\bullet$.
\end{proof}

\subsection{Invertibility}

Here we discuss a class of good endofunctors which give rise to abelian
groups (rather than monoids) at the level of generalized morphisms.
Recall that a morphism is called a \emph{zero morphism} if it factors
through a zero object (or equivalently, the chosen zero object).
If $0$ is a zero morphism in $\GEFC$, then $0G$ and $G0$ are also
zero morphisms in $\GEFC$ because $\bzero G$ and $G\bzero$ are
zero objects.
\begin{defn}\label{def:admitinv}
We say that $F\in\in\GEFC$ is
\emph
{homotopy symmetric} (cf.~\cite[Section~5]{dadarlat1994})
if there is a natural transformation $\inv_{F}\colon F\Rightarrow F$
such that
\[
F\xRightarrow{\Delta}F\oplus F\xRightarrow{F\oplus\inv_{F}}F\oplus F\xRightarrow{\lc^{-1}}F(\Id\oplus\Id)\xRightarrow{F\iota}F\fM_{2}
\]
is null-homotopic\footnote{That is to say, homotopic to a zero morphism in $\GEFC$.},
where $\Delta$ is as in Definition~\ref{def:Delta}, and $\iota\colon\Id\oplus\Id\Rightarrow\fM_{2}$
is as Section~\ref{subsec:starig}.
\end{defn}

\begin{prop}
Let $F,G\in\in\GEFC$, and assume that $F$ is homotopy symmetric.
Then $FG$ and $GF$ are also homotopy symmetric. Specifically, one
can choose $\inv_{FG}\coloneqq\inv_{F}G$ and $\inv_{GF}\coloneqq G\inv_{F}$.
\end{prop}

\begin{proof}
That $GF$ is homotopy symmetric follows from the outer diagram in
\[
\input{diagrams/d374.tikz}
\]
which commutes up to homotopy by Lemma~\ref{lem:Delta}, the bimonoidality
of $\GEFC$, and Definition~\ref{def:admitinv}.

That $FG$ is homotopy symmetric is proved by
\[
\input{diagrams/d375.tikz}
\]
which commutes up to homotopy by the bimonoidality of $\GEFC$, Lemmas~\ref{lem:Delta},~\ref{lem:kappa-tensor-type},
Corollary~\ref{lem:d376}, and Definition~\ref{def:admitinv}.
\end{proof}

\begin{thm}
If $F\in\in\GEFC$ is homotopy symmetric, then $[A,F,B]$ is an abelian
group. The inverse of $[\varphi]\in[A,F,B]$ is given by the formula
\[
-[\varphi]\coloneqq[-\varphi\colon A\xrightarrow{\varphi}FB\xrightarrow{\inv_{F}B}FB].
\]
\end{thm}

\begin{proof}
The statement follows from the diagram
\[
\input{diagrams/d480.tikz}
\]
which commutes up to homotopy by the bimonoidality of $\GEFC$, the
definition of $\mu$, Lemma~\ref{lem:Delta}, and Definition~\ref{def:admitinv}.
\end{proof}

\begin{prop}
The suspension functor $S\coloneqq\C_{(0,1)}$ is homotopy symmetric.
Specifically, one can define $\inv_{S}$ by the formula
\[
\inv_{S}B\colon C_{0}((0,1),B)\to C_{0}((0,1),B)\colon f\mapsto(t\mapsto f(1-t)).
\]
\end{prop}

\begin{proof}
For $0\leq a\leq b\leq1$ and $f\in C_{0}((0,1),B)$,
let us define
\[
T_{a,b}(f)(t)\coloneqq\begin{cases}
f\left(\frac{t-a}{b-a}\right), & t\in[a,b]\\
0, & \textrm{otherwise}.
\end{cases}
\]
For all $f,g\in C([0,1],B)$, we define the functions $\overline{f},f*g\in C([0,1],B)$
by the formulas
\begin{alignat*}{3}
 & \overline{f}(t)\coloneqq f(1-t), & \qquad\qquad\qquad\qquad & f*g(t) & \coloneqq\begin{cases}
f(2t), & t\in\left[0,\frac{1}{2}\right];\\
g(2t-1), & t\in\left[\frac{1}{2},1\right].
\end{cases}
\end{alignat*}
The statement of the proposition follows from the chain of homotopies
\[
\left[f\mapsto\left(\begin{array}{cc}
f & 0\\
0 & \overline{f}
\end{array}\right)\right]\simeq\left[f\mapsto\left(\begin{array}{cc}
T_{0,\frac{1}{3}}(f) & 0\\
0 & T_{\frac{2}{3},1}(\overline{f})
\end{array}\right)\right]\simeq\left[f\mapsto\left(\begin{array}{cc}
T_{0,\frac{1}{3}}(f)*T_{\frac{2}{3},1}(\overline{f}) & 0\\
0 & 0
\end{array}\right)\right]\simeq0
\]
given by the following formulas:
\begin{alignat*}{2}
 & H_{1}(f)(s)\coloneqq &\ 
& \left(\begin{array}{cc}
T_{0,\frac{1}{3}s+(1-s)}(f) & 0\\
0 & T_{\frac{2}{3}s,1}(\overline{f})
\end{array}\right),\\
 & H_{2}(f)(s)\coloneqq &\ 
& V_{s}\left(\begin{array}{cc}
T_{0,\frac{1}{3}}(f) & 0\\
0 & T_{\frac{2}{3},1}(\overline{f})
\end{array}\right)V^{*}_{s},\\
 & H_{3}(f)(s)\coloneqq &\ 
& \left(\begin{array}{cc}
T_{0,\frac{1}{3}+\frac{s}{1-s}}(f)*T_{\frac{2}{3}-\frac{s}{1-s},1}(\overline{f}) & 0\\
0 & 0
\end{array}\right),
\end{alignat*}
where $s\in[0,1]$ is the homotopy variable, and where $s\mapsto V_{s}$
is a continuous path of unitaries such that $V_{s}=\left(\begin{array}{cc}
1 & 0\\
0 & 1
\end{array}\right)$ for $s\in\left[0,\frac{1}{3}\right]$ and $V_{s}=\left(\begin{array}{cc}
0 & 1\\
-1 & 0
\end{array}\right)$ for $s\in\left[\frac{2}{3},1\right]$.

\end{proof}

\appendix

%% file: diagrams/d434.tikz
\begin{tikzcd}[ampersand replacement=\&]
{A} \& {FIB} \\
 \& {FB\mmdd}
\arrow["\Phi", from=1-1, to=1-2]
\arrow["\varphi_{j}"', from=1-1, to=2-2]
\arrow["F\ev_{j}B", from=1-2, to=2-2]
\end{tikzcd}

%% file: diagrams/d492.tikz
\begin{tikzcd}[column sep=1.0em,row sep=0.1em, ampersand replacement=\&]  %d492
 \& {F(\fK\oplus\fK)\oplus F\fK} \&  \& {F\fK\oplus F\fK} \\
 \& {\phantom{F}} \\
{(F\fK\oplus F\fK)\oplus F\fK} \& {F((\fK\oplus\fK)\oplus\fK)} \& {\hspc 3} \& {F(\fK\oplus\fK)} \\
 \&  \&  \&  \&  \& {F\fK} \\
{F\fK\oplus(F\fK\oplus F\fK)} \& {F(\fK\oplus(\fK\oplus\fK))} \&  \& {F(\fK\oplus\fK)} \\
 \& {\phantom{F}} \\
 \& {F\fK\oplus F(\fK\oplus\fK)} \&  \& {F\fK\oplus F\fK\mmdc}
\arrow[""{name=0, anchor=center, inner sep=0}, "F\mu\oplus F\fK", Rightarrow, from=1-2, to=1-4]
\arrow[Rightarrow, from=1-2, to=3-2]
\arrow["\lc^{-1}", Rightarrow, from=1-4, to=3-4]
\arrow["\lc^{-1}\oplus F\fK"{pos=0.7}, manv, Rightarrow, from=3-1, to=1-2]
\arrow[""{name=1, anchor=center, inner sep=0}, "\alpha^{\oplus}"', Rightarrow, from=3-1, to=5-1]
\arrow[""{name=2, anchor=center, inner sep=0}, Rightarrow, from=3-2, to=3-4]
\arrow[""{name=3, anchor=center, inner sep=0}, "F\alpha^{\oplus}"', Rightarrow, from=3-2, to=5-2]
\arrow["F\mu"{pos=0.3}, manh, Rightarrow, from=3-4, to=4-6]
\arrow["F\fK\oplus\lc^{-1}"'{pos=0.7}, manv, Rightarrow, from=5-1, to=7-2]
\arrow[""{name=4, anchor=center, inner sep=0}, Rightarrow, from=5-2, to=5-4]
\arrow["F\mu"'{pos=0.3}, manh, Rightarrow, from=5-4, to=4-6]
\arrow[Rightarrow, from=7-2, to=5-2]
\arrow[""{name=5, anchor=center, inner sep=0}, "F\fK\oplus F\mu"', Rightarrow, from=7-2, to=7-4]
\arrow["\lc^{-1}"', Rightarrow, from=7-4, to=5-4]
\arrow["\tc 1"{description}, draw=none, from=0, to=2]
\arrow["\tc 1"{description}, draw=none, from=1, to=3]
\arrow["\tc 2"{description}, draw=none, from=2, to=4]
\arrow["\tc 1"{description}, draw=none, from=4, to=5]
\end{tikzcd}

%% file: diagrams/d326.tikz
\begin{tikzcd}[ampersand replacement=\&]
{F\fK\oplus F\fK} \&  \& {F(\fK\oplus\fK)} \\
{F\fK\oplus F\fK} \&  \& {F(\fK\oplus\fK)} \&  \& {F\fK\mmdc}
\arrow["\lc^{-1}", Rightarrow, from=1-1, to=1-3]
\arrow[""{name=0, anchor=center, inner sep=0}, "\xi^{\oplus}"', Rightarrow, from=1-1, to=2-1]
\arrow[""{name=1, anchor=center, inner sep=0}, "F\xi^{\oplus}"', Rightarrow, from=1-3, to=2-3]
\arrow["F\mu", curve={height=-12pt}, Rightarrow, from=1-3, to=2-5]
\arrow["\lc^{-1}"', Rightarrow, from=2-1, to=2-3]
\arrow[""{name=2, anchor=center, inner sep=0}, "F\mu"', Rightarrow, from=2-3, to=2-5]
\arrow["\tc 1"{description}, draw=none, from=0, to=1]
\arrow["\tc 2"{description, pos=0.6}, shift right=2, draw=none, from=2, to=1-3]
\end{tikzcd}

%% file: diagrams/d323.tikz
\begin{tikzcd}[ampersand replacement=\&]
{F\fK} \&  \&  \& {F(\fK\oplus\bzero)} \& {F\fK} \\
{F\fK\oplus\bzero} \&  \&  \& {F\fK\oplus F\bzero} \\
{F\fK\oplus F\fK} \&  \&  \&  \& {F(\fK\oplus\fK)\mmdd}
\arrow[sqarn=1.5em, equals, from=1-1, to=1-5]
\arrow["(\rho^{\oplus})^{-1}"', Rightarrow, from=1-1, to=2-1]
\arrow[""{name=0, anchor=center, inner sep=0}, "F\rho^{\oplus}"', Rightarrow, from=1-4, to=1-1]
\arrow["\tc 2"{description}, draw=none, from=1-4, to=1-5]
\arrow[""{name=1, anchor=center, inner sep=0}, "F(\fK\oplus\bzero)", Rightarrow, from=1-4, to=3-5]
\arrow["F\fK\oplus(\Terminal^{l})^{-1}", Rightarrow, from=2-1, to=2-4]
\arrow["F\fK\oplus\bzero"', Rightarrow, from=2-1, to=3-1]
\arrow["\lc^{-1}", Rightarrow, from=2-4, to=1-4]
\arrow[""{name=2, anchor=center, inner sep=0}, "F\fK\oplus F\bzero", Rightarrow, from=2-4, to=3-1]
\arrow[""{name=3, anchor=center, inner sep=0}, "\lc^{-1}"', Rightarrow, from=3-1, to=3-5]
\arrow["F\mu"', Rightarrow, from=3-5, to=1-5]
\arrow["\tc 1"{description}, draw=none, from=1, to=3]
\arrow["\tc 1"{description}, draw=none, from=2-1, to=0]
\arrow["\tc 3"{description}, draw=none, from=2, to=2-1]
\end{tikzcd}

%% file: diagrams/d490.tikz
\begin{tikzcd}[column sep=tiny, ampersand replacement=\&]  %d490
{A} \& {\hspc 2} \& {F\fK IB} \& {\hspc{2.5}} \& {F\fK IG\fK C} \& {\hspc 6} \& {FG\fK I\fK C} \& {\hspc 6} \& {FG\fK^{2}IC} \& {\hspc 3} \& {FG\fK IC} \\
 \&  \& {F\fK B} \&  \& {F\fK G\fK C} \&  \& {FG\fK^{2}C} \&  \& {FG\fK^{2}C} \&  \& {FG\fK C}
\arrow["\varphi", from=1-1, to=1-3]
\arrow["\varphi_{j}"', curve={height=12pt}, from=1-1, to=2-3]
\arrow["F\fK I\psi", from=1-3, to=1-5]
\arrow[""{name=0, anchor=center, inner sep=0}, "F\fK\ev_{j}B"', from=1-3, to=2-3]
\arrow["F\kappa^{\fK I,G}\fK C", from=1-5, to=1-7]
\arrow[""{name=1, anchor=center, inner sep=0}, "F\fK\ev_{j}G\fK C"', from=1-5, to=2-5]
\arrow["FG\fK\kappa^{I,\fK}C", from=1-7, to=1-9]
\arrow[""{name=2, anchor=center, inner sep=0}, "FG\fK\ev_{j}\fK C"', from=1-7, to=2-7]
\arrow["FG\theta IC", from=1-9, to=1-11]
\arrow[""{name=3, anchor=center, inner sep=0}, "FG\fK^{2}\ev_{j}C"', from=1-9, to=2-9]
\arrow[""{name=4, anchor=center, inner sep=0}, "FG\fK\ev_{j}C"', from=1-11, to=2-11]
\arrow["F\fK\psi"', from=2-3, to=2-5]
\arrow["F\kappa^{\fK,G}\fK C"', from=2-5, to=2-7]
\arrow[equals, from=2-7, to=2-9]
\arrow["F\fK\theta C"', from=2-9, to=2-11]
\arrow["\tc 1"{description, pos=0.2}, draw=none, from=0, to=1]
\arrow["\tc 2"{description, pos=0.3}, draw=none, from=1, to=2]
\arrow["\tc 3"{description, pos=0.7}, draw=none, from=3, to=2]
\arrow["\tc 4"{description, pos=0.7}, draw=none, from=4, to=3]
\end{tikzcd}

%% file: diagrams/d304.tikz
\begin{tikzcd}[ampersand replacement=\&]
 \&  \&  \&  \& {\fK F\fK} \\
{\fK\fK F\fK} \&  \& {\fK F\fK\fK} \&  \& {F\fK\fK\fK} \&  \& {F\fK\fK} \\
{\fK F\fK} \&  \&  \&  \& {F\fK\fK} \&  \& {F\fK}
\arrow["\kappa^{\fK,F}\fK", Rightarrow, from=1-5, to=2-7]
\arrow["\fK\kappa^{\fK,F}\fK", Rightarrow, from=2-1, to=2-3]
\arrow["\theta F\fK"', Rightarrow, from=2-1, to=3-1]
\arrow["\fK F\theta", Rightarrow, from=2-3, to=1-5]
\arrow["\kappa^{\fK,F}\fK\fK", Rightarrow, from=2-3, to=2-5]
\arrow["F\fK\theta", Rightarrow, from=2-5, to=2-7]
\arrow["F\theta\fK", Rightarrow, from=2-5, to=3-5]
\arrow["F\theta", Rightarrow, from=2-7, to=3-7]
\arrow["\kappa^{\fK,F}\fK"', Rightarrow, from=3-1, to=3-5]
\arrow["F\theta"', Rightarrow, from=3-5, to=3-7]
\end{tikzcd}

%% file: diagrams/d400.tikz
\begin{tikzcd}[ampersand replacement=\&]
{(F\fK\oplus F\fK)G\fK} \\
{F\fK G\fK\oplus F\fK G\fK} \& {F(\fK G\fK\oplus\fK G\fK)} \& {F(\fK\oplus\fK)G\fK} \& {F\fK G\fK} \\
 \\
{FG\fK^{2}\oplus FG\fK^{2}} \& {F(G\fK^{2}\oplus G\fK^{2})} \& {FG\fK(\fK\oplus\fK)} \& {FG\fK^{2}} \\
{FG\fK\oplus FG\fK} \& {F(G\fK\oplus G\fK)} \& {FG(\fK\oplus\fK)} \& {FG\fK}
\arrow["\rc"', Rightarrow, from=1-1, to=2-1]
\arrow[""{name=0, anchor=center, inner sep=0}, "\lc^{-1}G\fK", curve={height=-18pt}, Rightarrow, from=1-1, to=2-3]
\arrow["\lc^{-1}"', Rightarrow, from=2-1, to=2-2]
\arrow["F\kappa^{\fK,G\fK}\oplus F\kappa^{\fK,G\fK}"', Rightarrow, from=2-1, to=4-1]
\arrow["F\rc^{-1}"', Rightarrow, from=2-2, to=2-3]
\arrow["F(\kappa^{\fK,G\fK}\oplus\kappa^{\fK,G\fK})"', Rightarrow, from=2-2, to=4-2]
\arrow[""{name=1, anchor=center, inner sep=0}, "F\mu G\fK"', Rightarrow, from=2-3, to=2-4]
\arrow["F\kappa^{\fK\oplus\fK,G\fK}"', Rightarrow, from=2-3, to=4-3]
\arrow["F\kappa^{\fK,G\fK}", Rightarrow, from=2-4, to=4-4]
\arrow[""{name=2, anchor=center, inner sep=0}, "\lc^{-1}", Rightarrow, from=4-1, to=4-2]
\arrow[""{name=3, anchor=center, inner sep=0}, "FG\theta\oplus FG\theta"', Rightarrow, from=4-1, to=5-1]
\arrow[""{name=4, anchor=center, inner sep=0}, "F\lc^{-1}", Rightarrow, from=4-2, to=4-3]
\arrow[""{name=5, anchor=center, inner sep=0}, "F(G\theta\oplus G\theta)", Rightarrow, from=4-2, to=5-2]
\arrow[""{name=6, anchor=center, inner sep=0}, "FG\fK\mu", Rightarrow, from=4-3, to=4-4]
\arrow["\tc 4"{description}, draw=none, from=4-3, to=5-3]
\arrow["FG\theta", Rightarrow, from=4-4, to=5-4]
\arrow["\lc^{-1}"', Rightarrow, from=5-1, to=5-2]
\arrow[""{name=7, anchor=center, inner sep=0}, "\lc^{-1}"', sqars=1.5em, Rightarrow, from=5-1, to=5-3]
\arrow["F\lc^{-1}"', Rightarrow, from=5-2, to=5-3]
\arrow["FG\mu"', Rightarrow, from=5-3, to=5-4]
\arrow["\tc 1"{description}, draw=none, from=0, to=2-1]
\arrow["\tc 1"{description}, draw=none, from=2-1, to=2]
\arrow["\tc 2"{description}, draw=none, from=2-2, to=4]
\arrow["\tc 3"{description}, draw=none, from=1, to=6]
\arrow["\tc 1"{description}, draw=none, from=3, to=5]
\arrow["\tc 1"{description, pos=0.7}, draw=none, from=7, to=5-2]
\end{tikzcd}

%% file: diagrams/d401.tikz
\begin{tikzcd}[ampersand replacement=\&]
{G\fK\oplus G\fK} \&  \& {G\fK^{2}\oplus G\fK^{2}} \\
{G(\fK\oplus\fK)} \&  \& {G(\fK^{2}\oplus\fK^{2})} \&  \& {G\fK(\fK\oplus\fK)} \\
{G\fK} \&  \&  \&  \& {G\fK^{2}}
\arrow[""{name=0, anchor=center, inner sep=0}, "\lc^{-1}", Rightarrow, from=1-1, to=2-1]
\arrow["G\mu\oplus G\mu"', Rightarrow, from=1-3, to=1-1]
\arrow[""{name=1, anchor=center, inner sep=0}, "\lc^{-1}"', Rightarrow, from=1-3, to=2-3]
\arrow["\lc^{-1}", curve={height=-12pt}, Rightarrow, from=1-3, to=2-5]
\arrow["G\mu"', Rightarrow, from=2-1, to=3-1]
\arrow["G(\mu\oplus\mu)"', Rightarrow, from=2-3, to=2-1]
\arrow["G\lc^{-1}"', Rightarrow, from=2-3, to=2-5]
\arrow["\tc 6"{description}, draw=none, from=2-5, to=1-3]
\arrow["G\fK\mu", Rightarrow, from=2-5, to=3-5]
\arrow[""{name=2, anchor=center, inner sep=0}, "G\theta", Rightarrow, from=3-5, to=3-1]
\arrow["\tc 5"{description}, draw=none, from=1, to=0]
\arrow["\tc 7"{description}, draw=none, from=2-3, to=2]
\end{tikzcd}

%% file: diagrams/d287.tikz
\begin{tikzcd}[ampersand replacement=\&]
{[B,G\fK,C]\times[A,F\fK,B]} \&  \& {[A,FG\fK,C]} \\
{[B,G'\fK,C]\times[A,F'\fK,B]} \&  \& {[A,F'G'\fK,C]\mmdd}
\arrow["\bullet", no head, from=1-1, to=1-3]
\arrow["\left\langle \beta\right\rangle \times\left\langle \alpha\right\rangle "', from=1-1, to=2-1]
\arrow["\left\langle \alpha\beta\right\rangle ", from=1-3, to=2-3]
\arrow["\bullet", no head, from=2-1, to=2-3]
\end{tikzcd}

%% file: diagrams/d281.tikz
\begin{tikzcd}[column sep=tiny, ampersand replacement=\&]  %d281
{{A\oplus A}} \& {{\hspc 2}} \& {{F\fK B\oplus F\fK B}} \&  \& {{(F\fK\oplus F\fK)B}} \& {{\hspc 2}} \& {{F(\fK\oplus\fK)B}} \& {{\hspc 2}} \& {{F\fK B}} \\
{{A}} \&  \& {{G\fK B\oplus G\fK B}} \&  \& {{(GB\oplus G\fK)B}} \&  \& {{G(\fK\oplus\fK)B}} \&  \& {{G\fK B}}
\arrow["\varphi_{1}\oplus\varphi_{2}", from=1-1, to=1-3]
\arrow["\varphi'_{1}\oplus\varphi'_{2}"', from=1-1, to=2-3]
\arrow[equals, from=1-3, to=1-5]
\arrow["\alpha\fK B\oplus\alpha\fK B", from=1-3, to=2-3]
\arrow["\lc^{-1}B", from=1-5, to=1-7]
\arrow["(\alpha\fK\oplus\alpha\fK)B", from=1-5, to=2-5]
\arrow["F\mu B", from=1-7, to=1-9]
\arrow["\alpha(\fK\oplus\fK)B", from=1-7, to=2-7]
\arrow["\alpha\fK B", from=1-9, to=2-9]
\arrow["\Delta", from=2-1, to=1-1]
\arrow[equals, from=2-3, to=2-5]
\arrow["\lc^{-1}B"', from=2-5, to=2-7]
\arrow["G\mu B"', from=2-7, to=2-9]
\end{tikzcd}

%% file: diagrams/d302.tikz
\begin{tikzcd}[ampersand replacement=\&]
{A} \& {F\fK B} \& {F\fK G\fK C} \& {FG\fK^{2}C} \& {FG\fK C} \\
 \& {F'\fK B} \& {F'\fK G\fK C} \& {F'G\fK^{2}C} \& {F'G\fK C} \\
{A} \&  \& {F'\fK G'\fK C} \& {F'G'\fK^{2}C} \& {F'G'\fK C,}
\arrow["\varphi", from=1-1, to=1-2]
\arrow["\psi\bullet\varphi", sqarn=1em, from=1-1, to=1-5]
\arrow["\varphi'"', from=1-1, to=2-2]
\arrow[equals, from=1-1, to=3-1]
\arrow["F\fK\psi", from=1-2, to=1-3]
\arrow["\alpha\fK B", from=1-2, to=2-2]
\arrow[from=1-3, to=1-4]
\arrow[from=1-3, to=2-3]
\arrow[from=1-4, to=1-5]
\arrow[from=1-4, to=2-4]
\arrow[from=1-5, to=2-5]
\arrow["\alpha\beta\fK C", sqare=2em, from=1-5, to=3-5]
\arrow[from=2-2, to=2-3]
\arrow["F'\fK\psi'"', from=2-2, to=3-3]
\arrow[from=2-3, to=2-4]
\arrow["F'\fK\beta\fK C", from=2-3, to=3-3]
\arrow[from=2-4, to=2-5]
\arrow[from=2-4, to=3-4]
\arrow[from=2-5, to=3-5]
\arrow["\psi'\bullet\varphi'"', sqars=1em, from=3-1, to=3-5]
\arrow[from=3-3, to=3-4]
\arrow[from=3-4, to=3-5]
\end{tikzcd}

%% file: diagrams/d425.tikz
\begin{tikzcd}[column sep=small, ampersand replacement=\&]  %d425
 \&  \& {GA\oplus GA} \& {\quad} \& {GF\fK B\oplus GF\fK B} \& {(GF\fK\oplus GF\fK)B} \&  \& {GF(\fK\oplus\fK)B} \\
{GA} \&  \& {G(A\oplus A)} \&  \& {G(F\fK B\oplus F\fK B)} \& {G(F\fK\oplus F\fK)B} \&  \& {GF\fK B}
\arrow["G\varphi\oplus G\psi", from=1-3, to=1-5]
\arrow[equals, from=1-5, to=1-6]
\arrow["\lc^{-1}B", from=1-6, to=1-8]
\arrow["GF\mu", from=1-8, to=2-8]
\arrow["\Delta", from=2-1, to=1-3]
\arrow["G\Delta"', from=2-1, to=2-3]
\arrow["\lc"', from=2-3, to=1-3]
\arrow["G(\varphi\oplus\psi)"', from=2-3, to=2-5]
\arrow["\lc", from=2-5, to=1-5]
\arrow[equals, from=2-5, to=2-6]
\arrow["\lc B", from=2-6, to=1-6]
\arrow["G\lc^{-1}B"', from=2-6, to=1-8]
\end{tikzcd}

%% file: diagrams/d305.tikz
\begin{tikzcd}[ampersand replacement=\&]
{A\oplus A} \& {F\fK^{2}B\oplus F\fK^{2}B} \& {(F\fK^{2}\oplus F\fK^{2})B} \& {F\fK(\fK\oplus\fK)B} \& {F\fK^{2}B} \\
 \&  \&  \& {F(\fK^{2}\oplus\fK^{2})B} \& {F\fK B} \\
{A} \& {F\fK B\oplus F\fK B} \& {(F\fK\oplus F\fK)B} \& {F(\fK\oplus\fK)B}
\arrow["\varphi\oplus\psi", from=1-1, to=1-2]
\arrow["\varphi'\oplus\psi'"', from=1-1, to=3-2]
\arrow[equals, from=1-2, to=1-3]
\arrow["F\theta B\oplus F\theta B", from=1-2, to=3-2]
\arrow[""{name=0, anchor=center, inner sep=0}, "\lc^{-1}B", from=1-3, to=1-4]
\arrow["\lc^{-1}B"', from=1-3, to=2-4]
\arrow[""{name=1, anchor=center, inner sep=0}, "(F\theta\oplus F\theta)B"', from=1-3, to=3-3]
\arrow["F\fK\mu B", from=1-4, to=1-5]
\arrow["F\theta B", from=1-5, to=2-5]
\arrow["F\lc^{-1}B"', from=2-4, to=1-4]
\arrow["\tc 3"{description}, draw=none, from=2-4, to=2-5]
\arrow[""{name=2, anchor=center, inner sep=0}, "F(\theta\oplus\theta)B", from=2-4, to=3-4]
\arrow[from=3-1, to=1-1]
\arrow[equals, from=3-2, to=3-3]
\arrow["\lc^{-1}B", from=3-3, to=3-4]
\arrow["F\mu B"', curve={height=12pt}, from=3-4, to=2-5]
\arrow["\tc 2"{description}, draw=none, from=1, to=2]
\arrow["\tc 1"{description}, draw=none, from=2-4, to=0]
\end{tikzcd}

%% file: diagrams/d202.tikz
\begin{tikzcd}[ampersand replacement=\&]
{S} \& {SNS} \\
 \& {\mathfrak{A}\fK S,}
\arrow["S\eta", Rightarrow, from=1-1, to=1-2]
\arrow["\alpha\iota_{00}S"', Rightarrow, from=1-1, to=2-2]
\arrow["\varepsilon S", Rightarrow, from=1-2, to=2-2]
\end{tikzcd}

%% file: diagrams/d203.tikz
\begin{tikzcd}[ampersand replacement=\&]
{N} \& {NSN} \\
 \& {N\mathfrak{A}\fK.}
\arrow["\eta N", Rightarrow, from=1-1, to=1-2]
\arrow["N\alpha\iota_{00}"', Rightarrow, from=1-1, to=2-2]
\arrow["N\varepsilon", Rightarrow, from=1-2, to=2-2]
\end{tikzcd}

%% file: diagrams/d308.tikz
\begin{tikzcd}[ampersand replacement=\&]
 \&  \& {SNSA} \& {SN\mathfrak{A}^{n}\fK B} \\
 \&  \& {\mathfrak{A}\fK SA} \& {\mathfrak{A}\fK\mathfrak{A}^{n}\fK B} \&  \& {\mathfrak{A}\mathfrak{A}^{n}\fK^{2}B} \&  \& {\mathfrak{A}^{n+1}\fK B\mmdc} \\
{SA} \&  \&  \&  \&  \& {\mathfrak{A}^{n}\fK B}
\arrow["SN\varphi", from=1-3, to=1-4]
\arrow["\varepsilon SA"', from=1-3, to=2-3]
\arrow["\varepsilon\mathfrak{A}^{n}\fK B"', from=1-4, to=2-4]
\arrow["\mathfrak{A}\fK\varphi"', from=2-3, to=2-4]
\arrow["\mathfrak{A}\kappa^{\fK,\mathfrak{A}^{n}}\fK B", from=2-4, to=2-6]
\arrow["\mathfrak{A}^{n+1}\theta B", from=2-6, to=2-8]
\arrow[""{name=0, anchor=center, inner sep=0}, "S\eta A", curve={height=-24pt}, from=3-1, to=1-3]
\arrow["\alpha\iota_{00}SA"'{pos=0.8}, from=3-1, to=2-3]
\arrow["\varphi"', from=3-1, to=3-6]
\arrow[""{name=1, anchor=center, inner sep=0}, "\alpha\iota_{00}\mathfrak{A}^{n}\fK B"{pos=0.8}, from=3-6, to=2-4]
\arrow["\alpha\mathfrak{A}^{n}\iota_{00}\fK B"', from=3-6, to=2-6]
\arrow[""{name=2, anchor=center, inner sep=0}, "\alpha\mathfrak{A}^{n}\fK B"', curve={height=12pt}, from=3-6, to=2-8]
\arrow["\tc 1"{description}, draw=none, from=2-3, to=0]
\arrow["\tc 2"{description}, draw=none, from=2-6, to=1]
\arrow["\tc 3"{pos=0.7}, shift left=2, draw=none, from=2-6, to=2]
\end{tikzcd}

%% file: diagrams/d309.tikz
\begin{tikzcd}[ampersand replacement=\&]
{NSA} \&  \& {NSN\mathfrak{A}^{n}\fK B} \&  \&  \&  \& {N\mathfrak{A}\fK\mathfrak{A}^{n}\fK B} \\
 \&  \&  \&  \&  \&  \& {N\mathfrak{A}^{n+1}\fK^{2}B} \\
{A} \&  \& {N\mathfrak{A}^{n}\fK B} \&  \&  \&  \& {N\mathfrak{A}^{n+1}\fK B,}
\arrow["NS\psi", from=1-1, to=1-3]
\arrow["N\varepsilon\mathfrak{A}^{n}\fK B"{pos=0.6}, from=1-3, to=1-7]
\arrow["N\mathfrak{A}\kappa^{\fK,\mathfrak{A}^{n}}\fK B", from=1-7, to=2-7]
\arrow["N\mathfrak{A}^{n+1}\theta B", from=2-7, to=3-7]
\arrow["\eta A", from=3-1, to=1-1]
\arrow["\psi"', from=3-1, to=3-3]
\arrow["\eta N\mathfrak{A}^{n}\fK B", from=3-3, to=1-3]
\arrow[""{name=0, anchor=center, inner sep=0}, "N\alpha\iota_{00}\mathfrak{A}^{n}\fK B"'{pos=0.6}, shift left, curve={height=-12pt}, from=3-3, to=1-7]
\arrow[""{name=1, anchor=center, inner sep=0}, "N\alpha\mathfrak{A}^{n}\iota_{00}\fK B"'{pos=0.6}, curve={height=-6pt}, from=3-3, to=2-7]
\arrow[""{name=2, anchor=center, inner sep=0}, "N\alpha\mathfrak{A}^{n}\fK B"'{pos=0.6}, from=3-3, to=3-7]
\arrow["\tc 1"{description, pos=0.3}, draw=none, from=1-3, to=0]
\arrow["\tc 2"{description}, shift left=3, draw=none, from=1, to=0]
\arrow["\tc 3"{description, pos=0.6}, shift right=3, draw=none, from=1, to=2]
\end{tikzcd}

%% file: diagrams/d426.tikz
\begin{tikzcd}[ampersand replacement=\&]
{FT^{\star}} \&  \& {T^{\star}F\mathfrak{A}\fK} \\
{GT^{\star}} \&  \& {T^{\star}G\mathfrak{A}\fK\mmdd}
\arrow["\kapparev^{F,T^{\star}}", Rightarrow, from=1-1, to=1-3]
\arrow["\alpha T^{\star}"', Rightarrow, from=1-1, to=2-1]
\arrow["T^{\star}\alpha\mathfrak{A}\fK", Rightarrow, from=1-3, to=2-3]
\arrow["\kapparev^{G,T^{\star}}", Rightarrow, from=2-1, to=2-3]
\end{tikzcd}

%% file: diagrams/d427.tikz
\begin{tikzcd}[ampersand replacement=\&]
{FT^{\star}} \&  \& {T^{\star}TFT^{\star}} \&  \& {T^{\star}FTT^{\star}} \&  \& {T^{\star}F\mathfrak{A}\fK} \\
{GT^{\star}} \&  \& {T^{\star}TGT^{\star}} \&  \& {T^{\star}GTT^{\star}} \&  \& {T^{\star}G\mathfrak{A}\fK}
\arrow["\eta FT^{\star}", Rightarrow, from=1-1, to=1-3]
\arrow["\alpha T^{\star}"', Rightarrow, from=1-1, to=2-1]
\arrow["T^{\star}\kappa^{F,T^{\star}}T^{\star}", Rightarrow, from=1-3, to=1-5]
\arrow["T^{\star}T\alpha T^{\star}"', Rightarrow, from=1-3, to=2-3]
\arrow["T^{\star}F\varepsilon", Rightarrow, from=1-5, to=1-7]
\arrow["T^{\star}\alpha TT^{\star}"', Rightarrow, from=1-5, to=2-5]
\arrow["T^{\star}\alpha\mathfrak{A}\fK", Rightarrow, from=1-7, to=2-7]
\arrow["\eta GT^{\star}"', Rightarrow, from=2-1, to=2-3]
\arrow["T^{\star}\kappa^{G,T^{\star}}T^{\star}"', Rightarrow, from=2-3, to=2-5]
\arrow["T^{\star}G\varepsilon"', Rightarrow, from=2-5, to=2-7]
\end{tikzcd}

%% file: diagrams/d429.tikz
\begin{tikzcd}[ampersand replacement=\&]
{[B,F\mathfrak{A}^{m}\fK,C]\times[A,T^{\star}\mathfrak{A}^{m'}\fK,B]} \& {\qquad} \&  \& {[B,F\mathfrak{A}^{n}\fK,C]\times[A,\mathfrak{A}^{n'}\fK,B]} \\
{[A,F\mathfrak{A}^{m}T^{\star}\mathfrak{A}^{m'}\fK,C]} \&  \&  \& {[A,F\mathfrak{A}^{n}T^{\star}\mathfrak{A}^{n'}\fK,C]} \\
{[A,FT^{\star}\mathfrak{A}^{m+m'+1}\fK^{2},C]} \&  \&  \& {[A,FT^{\star}\mathfrak{A}^{n+n'+1}\fK^{2},C]} \\
{[A,FT^{\star}\mathfrak{A}^{m+m'+1}\fK,C]} \&  \&  \& {[A,FT^{\star}\mathfrak{A}^{n+n'+1}\fK,C]}
\arrow["\left\langle F\alpha^{k}\mathfrak{A}^{m}\right\rangle \times\left\langle T^{\star}\alpha^{k'}\mathfrak{A}^{m}\right\rangle ", from=1-1, to=1-4]
\arrow["\bullet"', from=1-1, to=2-1]
\arrow["\bullet", from=1-4, to=2-4]
\arrow["\left\langle F\alpha^{k}\mathfrak{A}^{m}T^{\star}\alpha^{k'}\mathfrak{A}^{m'}\right\rangle ", from=2-1, to=2-4]
\arrow["\left\langle F\xi_{m,m'}\right\rangle "', from=2-1, to=3-1]
\arrow["\left\langle F\xi_{n,n'}\right\rangle ", from=2-4, to=3-4]
\arrow["\left\langle FT^{\star}\alpha^{k+k'}\mathfrak{A}^{m+m'+1}\fK\right\rangle ", from=3-1, to=3-4]
\arrow["(FT^{\star}\mathfrak{A}^{m+m'+1}\theta)_{*}"', from=3-1, to=4-1]
\arrow["(FT^{\star}\mathfrak{A}^{n+n'+1}\theta)_{*}", from=3-4, to=4-4]
\arrow["\left\langle FT^{\star}\alpha^{k+k'}\mathfrak{A}^{m+m'+1}\right\rangle ", from=4-1, to=4-4]
\end{tikzcd}

%% file: diagrams/d481.tikz
\begin{tikzcd}[column sep=huge, ampersand replacement=\&]  %d481
{\mathfrak{A}^{m}T^{\star}\mathfrak{A}^{m'}\fK} \&  \& {\mathfrak{A}^{n}T^{\star}\mathfrak{A}^{n'}\fK} \\
{T^{\star}\mathfrak{A}^{m+1}\fK\mathfrak{A}^{m'}\fK} \&  \& {T^{\star}\mathfrak{A}^{n+1}\fK\mathfrak{A}^{n'}\fK} \\
{T^{\star}\mathfrak{A}^{m+m'+1}\fK^{2}} \&  \& {T^{\star}\mathfrak{A}^{n+n'+1}\fK^{2}}
\arrow["\alpha^{k}\mathfrak{A}^{m}T^{\star}\alpha^{k'}\mathfrak{A}^{m'}\fK", Rightarrow, from=1-1, to=1-3]
\arrow["\kapparev^{\mathfrak{A}^{m},T^{\star}}\mathfrak{A}^{m'}\fK"', Rightarrow, from=1-1, to=2-1]
\arrow["\kapparev^{\mathfrak{A}^{n},T^{\star}}\mathfrak{A}^{n'}\fK", Rightarrow, from=1-3, to=2-3]
\arrow["T^{\star}\alpha^{k}\mathfrak{A}^{m+1}\fK\alpha^{k'}\mathfrak{A}^{m'}\fK", Rightarrow, from=2-1, to=2-3]
\arrow["T^{\star}\mathfrak{A}^{m+1}\kappa^{\fK,\mathfrak{A}^{m'}}\fK"', Rightarrow, from=2-1, to=3-1]
\arrow["T^{\star}\mathfrak{A}^{n+1}\kappa^{\fK,\mathfrak{A}^{n'}}\fK", Rightarrow, from=2-3, to=3-3]
\arrow["T^{\star}\alpha^{k}\mathfrak{A}^{m+1}\alpha^{k'}\mathfrak{A}^{m'}\fK^{2}", Rightarrow, from=3-1, to=3-3]
\end{tikzcd}

%% file: diagrams/d374.tikz
\begin{tikzcd}[ampersand replacement=\&]
 \& {GF\oplus GF} \&  \&  \& {GF(\Id\oplus\Id)} \\
 \& {GF\oplus GF} \&  \& {G(F\oplus F)} \\
{GF} \&  \&  \& {G(F\oplus F)} \& {GF\fM_{2}}
\arrow["\lc^{-1}", Rightarrow, from=1-2, to=1-5]
\arrow["\lc^{-1}", Rightarrow, from=1-2, to=2-4]
\arrow["GF\iota", Rightarrow, from=1-5, to=3-5]
\arrow["GF\oplus G\inv_{F}", Rightarrow, from=2-2, to=1-2]
\arrow["\lc^{-1}", Rightarrow, from=2-2, to=3-4]
\arrow["G\lc^{-1}", Rightarrow, from=2-4, to=1-5]
\arrow["\Delta", Rightarrow, from=3-1, to=2-2]
\arrow["G\Delta"', Rightarrow, from=3-1, to=3-4]
\arrow[" 0"', sqars=1em, Rightarrow, from=3-1, to=3-5]
\arrow["G(F\oplus\inv_{F})"', Rightarrow, from=3-4, to=2-4]
\end{tikzcd}

%% file: diagrams/d375.tikz
\begin{tikzcd}[ampersand replacement=\&]
{FG\oplus FG} \&  \& {FG\oplus FG} \& {F(G\oplus G)} \&  \& {FG(\Id\oplus\Id)} \\
{(F\oplus F)G} \&  \& {(F\oplus F)G} \& {F(\Id\oplus\Id)G} \&  \& {FG\fM_{2}} \\
{FG} \&  \&  \&  \&  \& {F\fM_{2}G}
\arrow["FG\oplus\inv_{F}G", Rightarrow, from=1-1, to=1-3]
\arrow["\rc^{-1}"', Rightarrow, from=1-1, to=2-1]
\arrow["\lc^{-1}"', Rightarrow, from=1-3, to=1-4]
\arrow["\lc^{-1}", sqarn=1em, Rightarrow, from=1-3, to=1-6]
\arrow["\rc^{-1}"', Rightarrow, from=1-3, to=2-3]
\arrow["F\lc^{-1}"', Rightarrow, from=1-4, to=1-6]
\arrow["FG\iota", Rightarrow, from=1-6, to=2-6]
\arrow["(F\oplus\inv_{F})G"', Rightarrow, from=2-1, to=2-3]
\arrow["\lc^{-1}G"', Rightarrow, from=2-3, to=2-4]
\arrow["F\rc"', Rightarrow, from=2-4, to=1-4]
\arrow["F\kappa^{\Id\oplus\Id,G}"', Rightarrow, from=2-4, to=1-6]
\arrow["F\iota G"', Rightarrow, from=2-4, to=3-6]
\arrow["\Delta", sqarw=3em, Rightarrow, from=3-1, to=1-1]
\arrow["\Delta G"', Rightarrow, from=3-1, to=2-1]
\arrow[" 0"', Rightarrow, from=3-1, to=3-6]
\arrow["F\kappa^{\fM_{2},G}"', Rightarrow, from=3-6, to=2-6]
\end{tikzcd}

%% file: diagrams/d480.tikz
\begin{tikzcd}[column sep=small, ampersand replacement=\&]  %d480
{A} \& {\hspc 0} \& {A\oplus A} \\
{F\fK B} \&  \& {F\fK B\oplus F\fK B} \&  \& {F\fK B\oplus F\fK B} \& {\hspc 1} \\
 \&  \& {(F\fK\oplus F\fK)B} \&  \& {(F\fK\oplus F\fK)B} \&  \& {F(\fK\oplus\fK)B} \\
{F\fK B} \&  \& {(F\oplus F)\fK B} \&  \& {(F\oplus F)\fK B} \&  \& {F(\Id\oplus\Id)\fK B} \\
 \&  \&  \& {\hspc 7} \& {F\fM_{2}\fK B} \&  \& {F\fK B}
\arrow["\Delta", from=1-1, to=1-3]
\arrow["\varphi"', from=1-1, to=2-1]
\arrow["\varphi\oplus\varphi"', from=1-3, to=2-3]
\arrow["\varphi\oplus-\varphi", curve={height=-12pt}, from=1-3, to=2-5]
\arrow["\Delta", from=2-1, to=2-3]
\arrow[equals, from=2-1, to=4-1]
\arrow["F\fK B\oplus\inv_{F}\fK B", from=2-3, to=2-5]
\arrow[equals, from=2-3, to=3-3]
\arrow[equals, from=2-5, to=3-5]
\arrow["(F\fK\oplus\inv_{F}\fK)B", from=3-3, to=3-5]
\arrow["\rc^{-1}B", from=3-3, to=4-3]
\arrow["\lc^{-1}B", from=3-5, to=3-7]
\arrow["\rc^{-1}B", from=3-5, to=4-5]
\arrow["F\rc^{-1}B", from=3-7, to=4-7]
\arrow["F\mu B", sqare=2em, from=3-7, to=5-7]
\arrow["\Delta B", from=4-1, to=3-3]
\arrow["\Delta\fK B", from=4-1, to=4-3]
\arrow[" 0"', curve={height=12pt}, from=4-1, to=5-5]
\arrow["(F\oplus\inv_{F})\fK B", from=4-3, to=4-5]
\arrow["\lc^{-1}\fK B", from=4-5, to=4-7]
\arrow["F\iota\fK B", from=4-7, to=5-5]
\arrow["F\theta_{2}B"', from=5-5, to=5-7]
\end{tikzcd}

%% file: sec6.tex
\section*{Appendix. Laplaza axioms}

\addcontentsline{toc}{section}{Appendix. Laplaza axioms}
\renewcommand{\theequation}{A.\arabic{equation}}
\setcounter{equation}{0}

In this appendix we list all the 24 Laplaza axioms of a symmetric
bimonoidal category~\cite[Definition~2.1.2]{JY_bimon24}, expressed
as the commutativity of the diagrams below. To obtain the 22 axioms
for a non-symmetric bimonoidal category, one omits~(\ref{eq:d551})
and (\ref{eq:d565}).

For readability, the operation $\otimes$ will be written here as
concatenation.
\begin{itemize}
\item Distributivity and multiplicative symmetry:
\begin{equation}
\input{diagrams/d551.tikz}\label{eq:d551}
\end{equation}
\item Distributivity and additive symmetry:
\begin{alignat*}{2}
 & \input{diagrams/d552.tikz} & \qquad\qquad & \input{diagrams/d575.tikz}
\end{alignat*}
\item Distributivity and additive associativity:
\begin{alignat*}{2}
 & \input{diagrams/d554.tikz} & \qquad\qquad & \input{diagrams/d555.tikz}
\end{alignat*}
\item Distributivity and multiplicative associativity:
\begin{alignat}{2}
 & \input{diagrams/d556.tikz} & \quad & \input{diagrams/d557.tikz}\label{eq:d557}
\end{alignat}
\begin{equation}
\input{diagrams/d558.tikz}\label{eq:d558}
\end{equation}
\item 2-by-2 distributivity:
\[
\input{diagrams/d559.tikz}
\]
\item Multiplicative zero of 0:
\begin{equation}
\input{diagrams/d560.tikz}\label{eq:d560}
\end{equation}
\item Multiplicative zero of a sum:
\begin{alignat*}{2}
 & \input{diagrams/d561.tikz} & \qquad\qquad & \input{diagrams/d562.tikz}
\end{alignat*}
\item Multiplicative zero and multiplicative unit:
\begin{alignat*}{2}
 & \input{diagrams/d563.tikz} & \qquad\qquad & \input{diagrams/d564.tikz}
\end{alignat*}
\item Symmetry of multiplicative zero:
\begin{equation}
\input{diagrams/d565.tikz}\label{eq:d565}
\end{equation}
\item Multiplicative zero and multiplicative associativity:
\begin{alignat}{3}
 & \input{diagrams/d566.tikz} & \qquad & \input{diagrams/d567.tikz} & \qquad & \input{diagrams/d568.tikz}\label{eq:d568}
\end{alignat}
\item Additive and multiplicative zero:
\begin{alignat*}{2}
 & \input{diagrams/d569.tikz} & \qquad\qquad & \input{diagrams/d570.tikz}\\
 & \input{diagrams/d571.tikz} &\ 
& \input{diagrams/d572.tikz}
\end{alignat*}
\item Distributivity and multiplicative unit:
\begin{alignat*}{2}
 & \input{diagrams/d573.tikz} & \qquad\qquad & \input{diagrams/d574.tikz}
\end{alignat*}
\end{itemize}

%% file: diagrams/d551.tikz
\begin{tikzcd}[ampersand replacement=\&]
{(A\oplus B)C} \& {AC\oplus BC} \\
{C(A\oplus B)} \& {CA\oplus CB}
\arrow["\delta^{r}", from=1-1, to=1-2]
\arrow["\xi^{\otimes}"', from=1-1, to=2-1]
\arrow["\xi^{\otimes}\oplus\xi^{\otimes}", from=1-2, to=2-2]
\arrow["\delta^{l}", from=2-1, to=2-2]
\end{tikzcd}

%% file: diagrams/d552.tikz
\begin{tikzcd}[ampersand replacement=\&]
{A(B\oplus C)} \& {AB\oplus AC} \\
{A(C\oplus B)} \& {AC\oplus AB}
\arrow["\delta^{l}", from=1-1, to=1-2]
\arrow["1_{A}\xi^{\oplus}"', from=1-1, to=2-1]
\arrow["\xi^{\oplus}", from=1-2, to=2-2]
\arrow["\delta^{l}", from=2-1, to=2-2]
\end{tikzcd}

%% file: diagrams/d575.tikz
\begin{tikzcd}[ampersand replacement=\&]
{(A\oplus B)C} \& {AC\oplus BC} \\
{(B\oplus A)C} \& {BC\oplus AC}
\arrow["\delta^{r}", from=1-1, to=1-2]
\arrow["\xi^{\oplus}1_{C}"', from=1-1, to=2-1]
\arrow["\xi^{\oplus}", from=1-2, to=2-2]
\arrow["\delta^{r}", from=2-1, to=2-2]
\end{tikzcd}

%% file: diagrams/d554.tikz
\begin{tikzcd}[column sep=small, ampersand replacement=\&]  %d554
{[(A\oplus B)\oplus C]D} \&  \& {(A\oplus B)D\oplus CD} \\
{[A\oplus(B\oplus C)]D} \&  \& {(AD\oplus BD)\oplus CD} \\
{AD\oplus(B\oplus C)D} \&  \& {AD\oplus(BD\oplus CD)}
\arrow["\delta^{r}", from=1-1, to=1-3]
\arrow["\alpha^{\oplus}1_{D}"', from=1-1, to=2-1]
\arrow["\delta^{r}\oplus1_{CD}", from=1-3, to=2-3]
\arrow["\delta^{r}", from=2-1, to=3-1]
\arrow["\alpha^{\oplus}", from=2-3, to=3-3]
\arrow["1_{AD}\oplus\delta^{r}", from=3-1, to=3-3]
\end{tikzcd}

%% file: diagrams/d555.tikz
\begin{tikzcd}[column sep=small, ampersand replacement=\&]  %d555
{A[(B\oplus C)\oplus D]} \&  \& {A(B\oplus C)\oplus AD} \\
{A[B\oplus(C\oplus D)]} \&  \& {(AB\oplus AC)\oplus AD} \\
{AB\oplus A(C\oplus D)} \&  \& {AB\oplus(AC\oplus AD)}
\arrow["\delta^{l}", from=1-1, to=1-3]
\arrow["1_{A}\alpha^{\oplus}"', from=1-1, to=2-1]
\arrow["\delta^{l}\oplus1_{AD}", from=1-3, to=2-3]
\arrow["\delta^{l}"', from=2-1, to=3-1]
\arrow["\alpha^{\oplus}", from=2-3, to=3-3]
\arrow["1_{AB}\oplus\delta^{l}", from=3-1, to=3-3]
\end{tikzcd}

%% file: diagrams/d556.tikz
\begin{tikzcd}[column sep=0em, ampersand replacement=\&]  %d556
{(AB)(C\oplus D)} \&  \& {(AB)C\oplus(AB)D} \\
{A[B(C\oplus D)]} \&  \& {A(BC)\oplus A(BD)} \\
 \& {\hspc{-4}A(BC\oplus BD)\hspc{-4}}
\arrow["\delta^{l}", from=1-1, to=1-3]
\arrow["\alpha^{\otimes}"', from=1-1, to=2-1]
\arrow["\alpha^{\otimes}\oplus\alpha^{\otimes}", from=1-3, to=2-3]
\arrow["1_{A}\delta^{l}"', from=2-1, to=3-2]
\arrow["\delta^{l}"', from=3-2, to=2-3]
\end{tikzcd}

%% file: diagrams/d557.tikz
\begin{tikzcd}[column sep=0em, ampersand replacement=\&]  %d557
{(A\oplus B)(CD)} \&  \& {A(CD)\oplus B(CD)} \\
{[(A\oplus B)C]D} \&  \& {(AC)D\oplus(BC)D} \\
 \& {\hspc{-4}(AC\oplus BC)D\hspc{-4}}
\arrow["\delta^{r}", from=1-1, to=1-3]
\arrow["\alpha^{\otimes}", from=2-1, to=1-1]
\arrow["\delta^{r}1_{D}"', from=2-1, to=3-2]
\arrow["\alpha^{\otimes}\oplus\alpha^{\otimes}"', from=2-3, to=1-3]
\arrow["\delta^{r}"', from=3-2, to=2-3]
\end{tikzcd}

%% file: diagrams/d558.tikz
\begin{tikzcd}[ampersand replacement=\&]
{[A(B\oplus C)]D} \& {(AB\oplus AC)D} \& {(AB)D\oplus(AC)D} \\
{A[(B\oplus C)D]} \& {A(BD\oplus CD)} \& {A(BD)\oplus A(CD)}
\arrow["\delta^{l}1_{D}", from=1-1, to=1-2]
\arrow["\alpha^{\otimes}"', from=1-1, to=2-1]
\arrow["\delta^{r}", from=1-2, to=1-3]
\arrow["\alpha^{\otimes}\oplus\alpha^{\otimes}", from=1-3, to=2-3]
\arrow["1_{A}\delta^{r}", from=2-1, to=2-2]
\arrow["\delta^{l}", from=2-2, to=2-3]
\end{tikzcd}

%% file: diagrams/d559.tikz
\begin{tikzcd}[ampersand replacement=\&]
{(A\oplus B)(C\oplus D)} \& {A(C\oplus D)\oplus B(C\oplus D)} \\
{(A\oplus B)C\oplus(A\oplus B)D} \& {(AC\oplus AD)\oplus(BC\oplus BD)} \\
{(AC\oplus BC)\oplus(AD\oplus BD)} \& {AC\oplus[AD\oplus(BC\oplus BD)]} \\
{AC\oplus[BC\oplus(AD\oplus BD)]} \& {AC\oplus[(AD\oplus BC)\oplus BD]} \\
{AC\oplus[(BC\oplus AD)\oplus BD]} \& {AC\oplus[(BC\oplus AD)\oplus BD]}
\arrow["\delta^{r}", from=1-1, to=1-2]
\arrow["\delta^{l}"', from=1-1, to=2-1]
\arrow["\delta^{l}\oplus\delta^{l}", from=1-2, to=2-2]
\arrow["\delta^{r}\oplus\delta^{r}"', from=2-1, to=3-1]
\arrow["\alpha^{\oplus}", from=2-2, to=3-2]
\arrow["\alpha^{\oplus}"', from=3-1, to=4-1]
\arrow["1_{AC}\oplus(\alpha^{\oplus})^{-1}", from=3-2, to=4-2]
\arrow["1_{AC}\oplus(\alpha^{\oplus})^{-1}"', from=4-1, to=5-1]
\arrow["1_{AC}\oplus(\xi^{\oplus}\oplus1_{BD})", from=4-2, to=5-2]
\arrow[equals, from=5-1, to=5-2]
\end{tikzcd}

%% file: diagrams/d560.tikz
\begin{tikzcd}[ampersand replacement=\&]
{\bbzero\otimes\bbzero} \& {\bbzero}
\arrow["\lambda^{\bullet}", shift left, from=1-1, to=1-2]
\arrow["\rho^{\bullet}"', shift right, from=1-1, to=1-2]
\end{tikzcd}

%% file: diagrams/d561.tikz
\begin{tikzcd}[ampersand replacement=\&]
{\bbzero(A\oplus B)} \& {\bbzero A\oplus\bbzero B} \\
{\bbzero} \& {\bbzero\oplus\bbzero}
\arrow["\delta^{l}", from=1-1, to=1-2]
\arrow["\lambda^{\bullet}"', from=1-1, to=2-1]
\arrow["\lambda^{\bullet}\oplus\lambda^{\bullet}", from=1-2, to=2-2]
\arrow["\lambda^{\oplus}", from=2-2, to=2-1]
\end{tikzcd}

%% file: diagrams/d562.tikz
\begin{tikzcd}[ampersand replacement=\&]
{(A\oplus B)\bbzero} \& {A\bbzero\oplus B\bbzero} \\
{\bbzero} \& {\bbzero\oplus\bbzero}
\arrow["\delta^{r}", from=1-1, to=1-2]
\arrow["\rho^{\bullet}"', from=1-1, to=2-1]
\arrow["\rho^{\bullet}\oplus\rho^{\bullet}", from=1-2, to=2-2]
\arrow["\lambda^{\oplus}", from=2-2, to=2-1]
\end{tikzcd}

%% file: diagrams/d563.tikz
\begin{tikzcd}[ampersand replacement=\&]
{\bbzero\otimes\tensorunit} \& {\bbzero}
\arrow["\lambda^{\bullet}", shift left, from=1-1, to=1-2]
\arrow["\rho^{\otimes}"', shift right, from=1-1, to=1-2]
\end{tikzcd}

%% file: diagrams/d564.tikz
\begin{tikzcd}[ampersand replacement=\&]
{\tensorunit\otimes\bbzero} \& {\bbzero}
\arrow["\rho^{\bullet}", shift left, from=1-1, to=1-2]
\arrow["\lambda^{\otimes}"', shift right, from=1-1, to=1-2]
\end{tikzcd}

%% file: diagrams/d565.tikz
\begin{tikzcd}[ampersand replacement=\&]
{A\otimes\bbzero} \&  \& {\bbzero\otimes A} \\
 \& {\bbzero}
\arrow["\xi^{\otimes}", from=1-1, to=1-3]
\arrow["\rho^{\bullet}"', from=1-1, to=2-2]
\arrow["\lambda^{\bullet}", from=1-3, to=2-2]
\end{tikzcd}

%% file: diagrams/d566.tikz
\begin{tikzcd}[ampersand replacement=\&]
{(AB)\bbzero} \& {A(B\bbzero)} \\
{\bbzero} \& {A\bbzero}
\arrow["\alpha^{\otimes}", from=1-1, to=1-2]
\arrow["\rho^{\bullet}"', from=1-1, to=2-1]
\arrow["1_{A}\rho^{\bullet}", from=1-2, to=2-2]
\arrow["\rho^{\bullet}", from=2-2, to=2-1]
\end{tikzcd}

%% file: diagrams/d567.tikz
\begin{tikzcd}[ampersand replacement=\&]
{(A\bbzero)B} \&  \& {A(\bbzero B)} \\
{\bbzero B} \& {\bbzero} \& {A\bbzero}
\arrow["\alpha^{\otimes}", from=1-1, to=1-3]
\arrow["\rho^{\bullet}1_{B}"', from=1-1, to=2-1]
\arrow["1_{A}\lambda^{\bullet}", from=1-3, to=2-3]
\arrow["\lambda^{\bullet}"', from=2-1, to=2-2]
\arrow["\rho^{\bullet}", from=2-3, to=2-2]
\end{tikzcd}

%% file: diagrams/d568.tikz
\begin{tikzcd}[ampersand replacement=\&]
{(\bbzero A)B} \& {\bbzero B} \\
{\bbzero(AB)} \& {\bbzero}
\arrow["\lambda^{\bullet}1_{B}", from=1-1, to=1-2]
\arrow["\alpha^{\otimes}"', from=1-1, to=2-1]
\arrow["\lambda^{\bullet}", from=1-2, to=2-2]
\arrow["\lambda^{\bullet}"', from=2-1, to=2-2]
\end{tikzcd}

%% file: diagrams/d569.tikz
\begin{tikzcd}[ampersand replacement=\&]
{A(\bbzero\oplus B)} \& {A\bbzero\oplus AB} \\
{AB} \& {\bbzero\oplus AB}
\arrow["\delta^{l}", from=1-1, to=1-2]
\arrow["1_{A}\lambda^{\oplus}"', from=1-1, to=2-1]
\arrow["\rho^{\bullet}\oplus1_{AB}", from=1-2, to=2-2]
\arrow["\lambda^{\oplus}", from=2-2, to=2-1]
\end{tikzcd}

%% file: diagrams/d570.tikz
\begin{tikzcd}[ampersand replacement=\&]
{(\bbzero\oplus B)A} \& {\bbzero A\oplus BA} \\
{BA} \& {\bbzero\oplus BA}
\arrow["\delta^{r}", from=1-1, to=1-2]
\arrow["\lambda^{\oplus}1_{A}"', from=1-1, to=2-1]
\arrow["\lambda^{\bullet}\oplus1_{BA}", from=1-2, to=2-2]
\arrow["\lambda^{\oplus}", from=2-2, to=2-1]
\end{tikzcd}

%% file: diagrams/d571.tikz
\begin{tikzcd}[ampersand replacement=\&]
{A(B\oplus\bbzero)} \& {AB\oplus A\bbzero} \\
{AB} \& {AB\oplus\bbzero}
\arrow["\delta^{l}", from=1-1, to=1-2]
\arrow["1_{A}\rho^{\oplus}"', from=1-1, to=2-1]
\arrow["1_{AB}\oplus\rho^{\bullet}", from=1-2, to=2-2]
\arrow["\rho^{\oplus}", from=2-2, to=2-1]
\end{tikzcd}

%% file: diagrams/d572.tikz
\begin{tikzcd}[ampersand replacement=\&]
{(B\oplus\bbzero)A} \& {BA\oplus\bbzero A} \\
{BA} \& {BA\oplus\bbzero}
\arrow["\delta^{r}", from=1-1, to=1-2]
\arrow["\rho^{\oplus}1_{A}"', from=1-1, to=2-1]
\arrow["1_{BA}\oplus\lambda^{\bullet}", from=1-2, to=2-2]
\arrow["\rho^{\oplus}", from=2-2, to=2-1]
\end{tikzcd}

%% file: diagrams/d573.tikz
\begin{tikzcd}[column sep=0em, ampersand replacement=\&]  %d573
{\tensorunit(A\oplus B)} \&  \& {\tensorunit A\oplus\tensorunit B} \\
 \& {A\oplus B}
\arrow["\delta^{l}", from=1-1, to=1-3]
\arrow["\lambda^{\otimes}"', from=1-1, to=2-2]
\arrow["\lambda^{\otimes}\oplus\lambda^{\otimes}", from=1-3, to=2-2]
\end{tikzcd}

%% file: diagrams/d574.tikz
\begin{tikzcd}[column sep=0em, ampersand replacement=\&]  %d574
{(A\oplus B)\tensorunit} \&  \& {A\tensorunit\oplus B\tensorunit} \\
 \& {A\oplus B}
\arrow["\delta^{r}", from=1-1, to=1-3]
\arrow["\rho^{\otimes}"', from=1-1, to=2-2]
\arrow["\rho^{\otimes}\oplus\rho^{\otimes}", from=1-3, to=2-2]
\end{tikzcd}